\documentclass[a4paper,11pt]{amsart}
\usepackage{amsmath, amssymb, amsthm}

\newcommand\beqn{\begin{equation}}
\newcommand\eeqn{\end{equation}}
\newcommand\beqny{\begin{eqnarray}}
\newcommand\eeqny{\end{eqnarray}}
\newcommand\beqnyn{\begin{eqnarray*}}
\newcommand\eeqnyn{\end{eqnarray*}}

\usepackage{amsmath,amsfonts,amssymb,amsthm}
\newtheorem{theorem}{Theorem}[section]
\newtheorem{lemma}[theorem]{Lemma}
\newtheorem{proposition}[theorem]{Proposition}

\def\a{\alpha}           
\def\b{\beta}          
\def\g{\gamma}           
\def\d{\delta}         
\def\e{\epsilon}         
\def\z{\zeta}            
\def\th{\theta}          
\def\k{\kappa}

\def\n{\nu}
\def\t{\tau}

\def\r{\rho} 

\def\s{\sigma}

\def\E{\mathcal E}
\def\G{\Gamma}

\def\Th{\Theta}

\begin{document}
\allowdisplaybreaks[1]
\title[A regularity and compactness theory]{A regularity and compactness theory for 
immersed stable minimal hypersurfaces of multiplicity at most 2}
\author{Neshan Wickramasekera}
\thanks{Research partly supported by NSF grant DMS-0406447} 
\date{}

\maketitle 

\tableofcontents

\section{Introduction}\label{introduction}
\setcounter{equation}{0}

Our goal in this paper is to study the local structure of
immersed, possibly branched, stable minimal hypersurfaces of the 
$(n+1)$-dimensional Euclidean space 
for arbitrary $n \geq 2.$ Assuming the singular set 
of such a hypersurface has locally finite $(n-2)$-dimensional Hausdorff measure, we here develop a regularity theory
that is applicable near those points of the hypersurface where the
volume density is less than 3. Our definition of the singular set is such that it consists only of 
``genuine'' singularities, which include possible branch points. Thus, 
the points of self-intersection, where the hypersurface is immersed, are considered regular. (See 
Section~\ref{notation} for the precise 
definition of the singular set.)\\

In particular, we obtain a 
description of the asymptotic behavior of the hypersurface near any of its multiplicity 2 branch points. 
i.e., the asymptotic behavior of the hypersurface near 
any point at which the hypersurface has a multiplicity 2 hyperplane as one of its tangent cones 
while it fails to decompose as the union of 
two regular minimal submanifolds in any neighborhood of the point. Our main regularity result is the following:\\

\begin{theorem} \label{maintheorem}
For each $\delta \in (0, 1),$ there exists a number $\e \in (0, 1)$, depending only 
on $n$ and $\d,$ such that 
the following is true. If $M$ is an orientable immersed stable minimal hypersurface of $B_{2}^{n+1}(0),$ 
with ${\mathcal H}^{n-2} \, ({\rm sing} \, M) < \infty$, $0 \in \overline{M},$
$\frac{{\mathcal H}^{n}(M)}{\omega_{n}2^{n}} \leq 3 - \delta$
and $\int_{M \cap (B_{1}(0) \times {\mathbf R})}|x^{n+1}|^{2} \leq \e,$
then ${M}_{1} \cap (B_{1/2}(0) \times {\mathbf R}) = \mbox{\rm graph} \, u$ 
where ${M}_{1}$ is the connected component of 
${\overline M} \cap (B_{1}(0) \times {\mathbf R})$ containing the origin, 
$u$ is either a single valued or a 2-valued $C^{1, \a}$ function on $B_{1/2}(0)$ satisfying 
$$\|u\|_{C^{1, \, \a}(B_{1/2}(0))} \leq
C\left(\int_{M \cap (B_{1}(0) \times {\mathbf R})} |x^{n+1}|^{2}\right)^{1/2}.$$ 
Here the constants $C$ and $\a \in (0,1)$ depend only on $n$ and $\d$.\\

\end{theorem}

See Section~\ref{notation}  for the definition of the $C^{1, \a}$ ``norm'' of $u$ when $u$ is a 
2-valued function.\\
 
This theorem in particular says that if an $n$-dimensional stable minimal hypersurface with a 
singular set of locally finite $(n-2)$-dimensional Hausdorff measure has a multiplicity 2 
plane as one of its tangent cones at some point, then it is the unique tangent cone to the hypersurface 
at that point. The theorem rules out for example the possibility of having a sequence of 
``necks'' connecting two sheets and accumulating at a branch point.\\
 
A direct consequence of the above theorem is the following decomposition theorem in case 
${\mathcal H}^{n-2} \, (\mbox{sing} \, M) = 0.$

\medskip
 
\begin{theorem} \label{maincorollary}
For each $\delta \in (0, 1),$ there exists a number $\e \in (0, 1)$, depending only 
on $n$ and $\d,$ such that 
the following is true. If $M$ is an orientable immersed stable minimal hypersurface of $B_{2}^{n+1}(0)$ 
with ${\mathcal H}^{n-2} \, ({\rm sing} \, M) = 0$, $0 \in \overline{M},$ 
$\frac{{\mathcal H}^{n}\,(M)}{\omega_{n}2^{n}} \leq 3 - \delta$
and 
$\int_{M \cap (B_{1}(0) \times {\mathbf R})}|x^{n+1}|^{2} \leq \e,$
then either
$M_{1} \cap (B_{1/2}(0) \times {\mathbf R}) = \mbox{\rm graph} \, u^{0}$ or
$M_{1} \cap (B_{1/2}(0) \times {\mathbf R}) = \mbox{\rm graph} \, u^{1} \cup \mbox{\rm graph} \, u^{2}$
where $M_{1}$ is the connected component of ${\overline M} \cap (B_{1}(0) \times {\mathbf R})$ 
containing the origin, $u^{i} \, : \, B_{1/2}(0) \to {\mathbf R}$ are $C^{1, \, \a}$ functions satisfying 
$$\|u^{i}\|_{C^{1, \, \a}(B_{1/2}(0))} \leq
C\left(\int_{M \cap (B_{1}(0) \times {\mathbf R})} |x^{n+1}|^{2}\right)^{1/2}$$ 
for $i = 0, 1, 2$. Here the constants $C$ and $\a \in (0,1)$ depend only on $n$ and $\d$.\\
\end{theorem}

Theorem~\ref{maincorollary}  implies that if $V$ is a varifold arising as the weak limit
of a sequence of stable minimal hypersurfaces having singular sets 
of $(n-2)$-dimensional Hausdorff measure zero, then near every point where $V$ has a tangent cone
equal to the multiplicity 2 varifold associated with a hyperplane, the support of $V$ decomposes as the union of two minimal graphs. 
In particular, classical branching (of multiplicity 2) cannot occur in the weak limit of a sequence of 
smooth, stable minimal hypersurfaces.\\
  
Based on Theorem~\ref{maincorollary}  and the standard dimension reducing principle of Federer, 
we obtained the following compactness result:

\medskip

\begin{theorem}\label{compactnessthm}
Let $\d \in (0, 1).$ There exists $\s = \s(n, \d) \in (0, 1/2)$ such that the following is true. 
Suppose $M_{k}$ is a sequence of orientable stable minimal 
hypersurface immersed in 
$B_{2}^{n+1}(0)$ with $0 \in {\overline M}_{k}$, 
${\mathcal H}^{n-2}(\mbox{\rm sing} \, M_{k} \cap B_{\s}^{n+1}(0)) = 0$ for each $k$ and 
$\limsup_{k \to \infty} \, \frac{{\mathcal H}^{n}(M_{k})}{\omega_{n}2^{n}} \leq 3 - \d$. 
Then there exists a stationary varifold $V$ of 
$B_{2}^{n+1}(0)$ and a closed subset 
$S$ of $\mbox{\rm spt} \, \|V\| \cap B_{\s}^{n+1}(0)$ with $S = \emptyset$ if 
$2\leq n \leq 6$, $S$ finite if $n=7$ and ${\mathcal H}^{n-7+\g}(S) = 0$ for every $\g>0$ if $n\geq 8$ 
such that after passing to 
a subsequence, which we again denote $\{k\}$, $M_{k} \to V$ as 
varifolds and
$(\mbox{\rm spt} \,\|V\| \setminus S) \cap B_{\s}^{n+1}(0)$ is an orientable immersed, smooth, 
stable minimal hypersurface 
of $B_{\s}^{n+1}(0)$.\\ 
\end{theorem}

In low dimensions, the ``smallness of excess'' hypothesis of Theorem~\ref{maintheorem} can 
be dropped provided we assume that the mass ratio is sufficiently close to $2$.  Precisely, we have the following:\\

\begin{theorem}\label{lowdimreg}
There exist fixed constants $\e \in (0, 1),$ $\s \in (0, 1)$, $C \in (0, \infty)$ and 
$\alpha \in (0, 1)$ such that the following holds. If $2 \leq n \leq 6$, $M$ is 
an orientable immersed stable minimal hypersurface of $B_{2}^{n+1}(0)$ with ${\mathcal H}^{n-2} \, 
({\rm sing} \, M) < \infty,$ $0 \in {\overline M}$
and $\frac{{\mathcal H}^{n} \, (M)}{\omega_{n}2^{n}}
\leq 2 + \e$, then for some orthogonal rotation $q$ of 
${\mathbf R}^{n+1}$, either $q \, M_{1} \cap (B_{\s}(0) \times {\mathbf R}) = {\rm graph} \, u$ where
$u$ is either a single valued or a 2-valued $C^{1, \alpha}$ function on $B_{\s}(0)$ with 
$$\|u\|_{C^{1, \, \alpha}(B_{\s}(0))} 
\leq C\left( \int_{q \, M \cap (B_{1}(0) \times {\mathbf R})} |x^{n+1}|^{2}\right)^{1/2} \hspace{.3in}$$
\noindent
 or there exists a pair of transversely intersecting affine hyperplanes $P^{(1)}, P^{(2)}$ of 
${\mathbf R}^{n+1}$ such that  
$q \, M_{1} \cap (B_{\s}(0) \times {\mathbf R}) = {\rm graph} \, (p^{(1)} + u^{(1)}) \cup {\rm graph} \, 
(p^{(2)} + u^{(2)}),$ where $P^{(1)} = {\rm graph} \, p^{(1)}$, $P^{(2)} = {\rm graph} \, p^{(2)}$, 
$p^{(i)} \, : \, {\mathbf R}^{n} \times \{0\} \to {\mathbf R}$, 
$u^{(i)} \in C^{1, \alpha} \, (P^{(i)} \cap (B_{\s}(0) \times {\mathbf R}) ; {\mathbf R})$ with 
$$\|u^{(i)}\|_{C^{1, \, \alpha}(P^{(i)} \cap (B_{\s}(0) \times {\mathbf R}))} 
\leq C\left(\int_{M \cap (B_{1}(0) \times {\mathbf R})}{\rm dist}^{2} \, (x, P)\right)^{1/2}$$
for $i =1, 2.$ Here $M_{1}$ denotes the connected component of ${\overline M} 
\cap (B_{1}(0) \times {\mathbf R})$ containing the origin.\\ 
\end{theorem}

Finally, we mention the following decomposition theorem for the singular set of 
a branched stable minimal hypersurface of the type considered in this paper.\\

\begin{theorem}\label{decomposition}
There exist $\e  = \e(n) \in (0, 1)$ and $\s = \s(n) \in (0, 1)$ such that 
the following holds. If $V$ belongs to the varifold closure of orientable immersed stable minimal 
hypersurfaces $M$ of $B_{2}^{n+1}(0)$ with 
$0 \in {\overline M}$, ${\mathcal H}^{n-2} \, ({\rm sing} \, M) < \infty$ and
$\frac{{\mathcal H}^{n} \, (M)}{\omega_{n}2^{n}} \leq 2 + \e$ then 
$${\rm sing} \, V \cap B_{\s}^{n+1}(0) = B \cup S$$
\noindent
where\\
\begin{itemize}
\item[(a)] $B$ is the set of branch points of $V$ in $B_{\s}^{n+1}(0);$ thus $B$ consists of those points of 
${\rm sing} \, V \cap B_{\s}^{n+1}(0)$ where $V$ has a (unique) multiplicity 2 tangent plane. 
Either $B = \emptyset$ or ${\mathcal H}^{n-2} \, (B) > 0.$\\ 
\item[(b)]$S$ is a relatively closed subset of ${\rm spt} \, \|V\| \cap B_{\s}^{n+1}(0)$ with 
$S \cap B = \emptyset,$ $S = \emptyset$ if $2 \leq n \leq 6$, $S$ a finite set if $n = 7$ and 
${\mathcal H}^{n-7+\g} \, (S) = 0$ for each $\g >0$ if $n \geq 8.$\\ 
\end{itemize}
\end{theorem}
 
The proofs of the above theorems will appear in Sections~\ref{mainthms} and 
\ref{compact} of the paper. Other consequences 
of Theorems~\ref{maincorollary} and \ref{compactnessthm}, which include a pointwise 
curvature estimate and a Bernstein type theorem in dimensions $\leq 6,$  
will appear in Section~\ref{corollaries}.\\

In case the mass bound is $2 - \d$ (instead of $3 - \d$) for some $\d \in (0, 1)$, 
Theorem~\ref{maintheorem} (with the conclusion 
that $M \cap (B_{1/2}(0) \times {\mathbf R})$ is the graph of a single valued function) 
follows from (otherwise much more general) 
interior regularity theorem of W.K. Allard \cite{AW}, \cite{S1}. In case the stable hypersurface
is {\em embedded}, the theorem 
(under the weaker hypothesis of arbitrary mass bound and with the stronger conclusion as in 
Theorem~\ref{maincorollary} with
a finite number of functions $u_{1} < u_{2} < \ldots < u_{k}$ in place of 
$u_{1}, u_{2},$ with $k$ bounded in terms of the mass bound) is due to R. Schoen and L. Simon \cite{SS}.
Schoen-Simon theorem in fact plays an essential role in the present work.\\

The main ingredient in the proof of Theorem~\ref{maintheorem} is a height 
excess decay lemma (Lemma~\ref{mainlemma1}), 
where we show that under the hypotheses of the theorem, the height excess of the 
hypersurface $M$ at a smaller scale, measured relative to a suitable new pair of 
hyperplanes (a transverse pair of hyperplanes or a multiplicity 2 hyperplane) improves by 
a fixed factor. The theorem
follows by iteratively applying  the lemma. At a key stage of the proof of the excess decay lemma, 
we use a type of harmonic approximation, 
where we show that whenever the $L^{2}$-height excess of the hypersurface relative to 
a hyperplane is small in a cylinder, the hypersurface in a smaller cylinder is well approximated  by 
the graph of a certain type of ``2-valued harmonic'' function. F. J. Almgren Jr. \cite{A} used a somewhat  different class of multi-valued harmonic functions in his work 
on area minimizing currents of arbitrary dimension and codimension, where 
harmonic meant Dirichlet energy minimizing. We are working with the weaker assumption of 
stability, so our two-valued harmonic functions do not satisfy this minimizing property. However, 
the codimension 
1 setting we are working in gives them a lot more structure, and we are able to obtain (in
Theorem~\ref{regblowupthm}) sufficiently 
detailed, geometric information about them.\\ 

A feature of our excess decay lemma perhaps worth pointing out here is that 
it gives, at every scale, decay of the excess of the stable minimal hypersurface at one of {\em three} possible, fixed smaller scales.
The reason why excess improvement is exhibited at one of several possible scales in contrast to 
the more familiar scenario where the improvement 
is always seen at a single fixed, smaller scale is partly geometric and partly technical. The  
geometric part of the explanation
is that the way an immersed hypersurface satisfying the hypotheses of the theorem (in particular, 
the mass bound $3 - \d$ which guarantees that it is ``two sheeted'') looks as one goes down in scale 
(fixing a base point) may vary between different possibilities; namely, at any 
given scale, it 
may either look like a pair of distinct, more or less parallel planes (i.e. the hypersurface is embedded) or 
it may look like a pair of transversely 
intersecting planes (i.e. the hypersurface is embedded away from a small tubular 
neighborhood around the axis of a transverse pair of hyperplanes) or it may have many self-intersections 
distributed more or less evenly. 
Different techniques for these different cases are employed in obtaining excess 
improvement. The technical 
part of the reason for the three scales is not having at our disposal, a priori, 
a single decay estimate, valid uniformly at all points of the domain away from the boundary and 
for all scales less than a fixed scale, for the aforementioned approximating 2-valued harmonic 
functions (which arise as blow-ups of sequences of hypersurfaces 
satisfying the hypotheses of the theorem
and converging to multiplicity 2 hyperplanes). Rather, what we obtain (in Theorem~\ref{regblowupthm}) is 
an asymptotic description which gives two alternatives depending on whether the blow-up itself 
has a non-empty interior branch set or not. The presence 
of two such alternatives for the asymptotics of this "linear problem" means that, at the stage 
where knowledge of the asymptotics of 
the linear problem becomes necessary (which is precisely when we are confronted with the 
picture where the minimal hypersurface has many 
self-intersections distributed approximately evenly), the excess improvement we get for the hypersurface 
is,  correspondingly, at one of two different smaller scales.\\

We use methods and results due to L. Simon \cite{S}; R. Hardt and L. Simon \cite{HS}; R. Schoen and 
L. Simon \cite{SS}; F. J. Almgren Jr. \cite{A} 
and the author \cite{WN} at a number of crucial points in the present work. The present work in fact 
should be viewed as a generalization of the results of \cite{WN}. To prove that a
stable minimal hypersurface, when it is weakly close to a multiplicity 2 hyperplane, is well 
approximated by the 
graph of a 2-valued harmonic function of the type aforementioned, we utilize a blow-up argument 
where we blow up sequences of hypersurfaces off affine hyperplanes. This 
blow up procedure is based on the approximate 
graphical decomposition of the hypersurfaces as in \cite{SS}, and is carried out as described in 
\cite{WN}, after making modifications to and replacement of some of the arguments 
of \cite{WN}. The main difference in the present context, as far as this blowing up step is concerned, is that 
we here allow 
the hypersurfaces to be singular unlike in \cite{WN} where they were assumed to be smooth. Consequently, 
in particular, we here need a different argument to establish continuity of the blow-ups. (See 
Proposition~\ref{lipschitz}.)\\ 

A major part  of this paper
is devoted to analyzing the nature of these 2-valued approximating functions. Theorem~\ref{regblowupthm} is the key result in this respect, where we establish 
crucial decay estimates for the two-valued harmonic functions.
Our approach in analyzing these functions has been to use geometric arguments, aimed at proving 
excess decay estimates for the graphs of the functions. To investigate the local regularity properties of 
these functions 
at points where their graphs blow up to transversely intersecting pairs of hyperplanes, and also 
to prove global decay estimates when the base point is a branch point of the function, we use variants 
of powerful techniques developed by Simon~\cite{S} and Hardt and Simon~\cite{HS}. 
In particular, a crucial ingredient is an estimate for the radial derivatives of the blow-up  
(Lemma~\ref{radialexcessbd}) due to Hardt and Simon \cite{HS}.\\    

An important technical tool used in the analysis of the 
2-valued harmonic functions is the monotonicity of a frequency function, an idea used first 
in a geometric setting by F. J. Almgren Jr. \cite{A}.
We here make use of the frequency function directly associated with  
the two-valued function as well as the one associated with 
the single valued function obtained by taking the difference between 
the two values of the two-valued function. Either frequency function, for any given 
center point, is monotonically non-decreasing as a function 
of the radius. Thus, in particular,  we may classify the
points of the domain of the two-valued function according to the values assumed by the limit of the frequency 
function     
associated with the difference function. In a classical setting, e.g. if the function were single 
valued and harmonic, this limit is 
equal to the vanishing order of the function at the point in question. In our setting, it conveys 
analogous information, which 
may be regarded as the order of contact between the ``two sheets'' of the graph of 
the 2-valued function, (although admittedly at a branch point 
one does not have a useful notion of two sheets)   
and it reveals the local geometric picture of the graph; i.e.
whether the graph locally consists of two disjoint harmonic disks, or of two self intersecting 
harmonic disks or 
whether it is branched. Furthermore,  the rate of decay 
of the graph of the two valued function to its (unique) multiplicity 2 tangent plane at a branch point
has a fixed lower bound  independent of the function. 
Said differently, there exists a fixed frequency gap, depending only on $n$ and $\d$ ( 
$\d$ as in Theorem~\ref{maintheorem}), implying that the 
order of contact at a branch point cannot be arbitrarily close to 1.\\

The existence of a rather rich class of 
stable branched minimal immersions of the type studied in this paper has recently been established 
in \cite{SW}.\\

I am very grateful to Leon Simon for several helpful discussions related to this work. I 
also thank David Jerison, Fanghua Lin and Gang Tian for conversations from which I have benefited.\\

\bigskip 
     
\section{Notation and preliminaries}\label{notation}
We shall adopt the following notation, conventions and definitions throughout the paper.\\

${\mathbf R}^{n+1}$ denotes the $(n+1)$-dimensional Euclidean space and 
$(x^{1}, \ldots, x^{n+1})$ denotes a general 
point in ${\mathbf R}^{n+1}$.\\

$B_{\rho}^{n+1}(X)$ denotes the open ball in ${\mathbf R}^{n+1}$ with radius $\rho$ and 
center $X$. For $X \in {\mathbf R}^{n} \times \{0\},$ we let $B_{\rho}(X) = 
B_{\rho}^{n+1}(X) \cap ({\mathbf R}^{n} \times \{0\})$.\\

$\omega_{n}$ denotes the volume of a ball in ${\mathbf R}^{n}$ with radius 1.\\

For compact sets $S, \, T \subseteq {\mathbf R}^{n+1}$, 
${\rm d}_{\mathcal H}(S,T)$ denotes the Hausdorff distance between 
$S$ and $T$.\\

${\mathcal H}^{n} \, (S)$ denotes the $n$-dimensional Hausdorff measure of the set $S.$\\

For $Y \in {\mathbf R}^{n+1}$ and $\r > 0$, $\eta_{Y, \, \r}  \, : \, {\mathbf R}^{n+1} \to {\mathbf R}^{n+1}$ 
is the map defined by $\eta_{Y, \, \r} \, (X) = \frac{X - Y}{\r}.$\\
 
The letter $M$ will always denote an immersed, smooth hypersurface of $B_{2}^{n+1}(0).$ Thus $M$ is a  
subset of $B_{2}^{n+1}(0)$ such that for each $X \in M$, there exists a number $\s > 0$ such that 
$M \cap B_{\s}^{n+1}(X)$ is the union of a finite number of, possibly intersecting, smooth, 
connected, embedded $n$-dimensional submanifolds with no boundary in $B_{\s}^{n+1}(x).$ \\

Let $M$ be a smooth hypersurface of $B_{2}^{n+1}(0)$ with ${\mathcal H}^{n} \, (M) < \infty.$ 
$M$ said to be {\em minimal} (or {\em stationary}) if it has zero first 
variation of volume with respect to deformations by arbitrary $C^{1}$ vector fields of the ambient space 
${\mathbf R}^{n+1}$ having compact support
in $B_{2}^{n+1}(0).$ Minimality of $M$ is equivalent to the condition that 

\begin{equation}\label{firstvar}
\int_{M} {\rm div}_{M} \, \Phi \, d{\mathcal H}^{n} = 0
\end{equation}

\noindent
for every $C^{1}$ vector field $\Phi  = (\Phi^{1}, \Phi^{2}, \ldots, \Phi^{n+1}) 
\, : \, B_{2}^{n+1}(0) \to {\mathbf R}^{n+1}$ with compact support 
in $B_{2}^{n+1}(0).$ (See \cite{S1}, Chapter 2.) Here ${\rm div}_{M} \, \Phi$ is the tangential 
divergence of $\Phi$ with respect to $M$. Thus, ${\rm div} \, \Phi = \sum_{j=1}^{n+1} e_{j} \cdot 
\nabla^{M} \, \Phi^{j}$ where $\nabla^{M}$ denotes the gradient operator on $M$ and  $\{e_{j}\}_{j=1}^{n+1}$ 
is the standard basis of ${\mathbf R}^{n+1}.$\\

A minimal hypersurface $M$ of $B_{2}^{n+1}(0)$ is {\em stable} if it has non-negative second variation 
of volume with respect to deformations as above. Stability of $M$ is equivalent to the statement 
that (\cite{S1}, Chapter 2) 

\begin{equation}\label{secondvar}
\int_{M} |A|^{2} \, \z^{2} \leq \int_{M} |\nabla \z|^{2}
\end{equation}

\noindent
for every $C^{1}$ function $\z$ with compact support in $M.$ Here $A$ denotes the second fundamental 
form of $M$ and $|A|$ its length.\\

For a smooth hypersurface $M$ of $B_{2}^{n+1}(0)$, we say a point $X \in {\overline M} \cap B_{2}^{n+1}(0)$ 
is (an interior) {\em regular point} of $M$ if there exists a number $\s > 0$ such that 
${\overline M} \cap \overline B_{\s}^{n+1}(X)$ is the union of finitely many smooth, 
compact, connected, embedded submanifolds with boundary contained in ${\partial \, B}_{\s}^{n+1}(X).$
We shall redefine $M$ so that if $X \in {\overline M}$ is a regular point of $M$, then 
$X \in M.$ The (interior) {\em singular set} of $M$ is then defined by 

$${\rm sing} \, M = \left({\overline M} \setminus M \right) \cap B_{2}^{n+1}(0).\\$$
  
${\mathcal I}_{b}$ denotes the family of stable minimal hypersurfaces $M$ of $B_{2}^{n+1}(0)$ 
satisfying  ${\mathcal H}^{n-2}({\rm sing} \, M) < \infty.$ (The subscript $b$ in ${\mathcal I}_{b}$ 
indicates that the members $M$ of ${\mathcal I}_{b}$ are allowed to carry {\em branch point} 
singularities; i.e. points $Z \in {\rm sing} \, M$ such that a hyperplane (with multiplicity $>1$) 
occurs as a tangent cone to (the varifold associated with) $M$ at $Z$.)\\

For a stationary, rectifiable $n$-varifold $V$ of some open subset $U$ of 
${\mathbf R}^{n+1}$ and a point $X \in U$, $\Theta \, (\|V\|, X)$ denotes the $n$-dimensional {\em density} 
at $X$ of the weight measure $\|V\|$ on $U$
associated with $V$. We refer the reader to \cite{S1}, Chapters 4 and 8  for an exposition of the 
theory of rectifiable varifolds.\\

For a Radon measure $\mu$ on $U$, ${\rm spt} \, \mu$ denotes the support of $\mu.$\\

If $L$ is an affine hyperplane of ${\mathbf R}^{n+1},$ 
$\pi_{L} \, : \, {\mathbf R}^{n+1} \to L$ denotes 
the orthogonal projection of ${\mathbf R}^{n+1}$ onto $L$. We shall abbreviate $\pi_{{\mathbf R}^{n} \times 
\{0\}}$ as $\pi.$\\ 

Unless stated otherwise, all constants $c,$ $C$ depend only on $n$ and $\d$, where $\d$ is as in 
Theorems~\ref{maintheorem}-\ref{compactnessthm}.\\

A {\em pair of affine hyperplanes} means the union of two not necessarily distinct  
affine hyperplanes of ${\mathbf R}^{n+1}$ neither of which is perpendicular to ${\mathbf R}^{n} \times  \{0\}.$
If $P = P_{1} \cup P_{2}$ is a pair of affine hyperplanes, with 
$P_{1}, P_{2}$ affine hyperplanes, we use the notation 
$p^{+} = \mbox{max} \, \{l_{1}, l_{2}\}$ and $p^{-} = \mbox{min} \, \{l_{1}, l_{2}\}$ where, for $i=1, 2$, 
$l_{i} \, : \, {\mathbf R}^{n} \times \{0\} \to {\mathbf R}$ is the affine function with 
$\mbox{graph} \, l_{i} = P_{i}$, and we set $P^{+} = \mbox{graph} \, p^{+}$ 
and $P^{-} = \mbox{graph} \,p^{-}.$ For such a pair of 
affine hyperplanes $P$, $\angle \, P$ denotes the angle $\th \in [0, \pi)$ between  $P_{1}$ and $P_{2}$ 
given by $\cos \th = \nu_{1} \cdot \nu_{2}$ where, for $i = 1, 2$, $\nu_{i} =  \frac{(-Dl_{i}, 1)}{
\sqrt{1 + |Dl_{i}|^{2}}}.$\\

By a {\em pair of hyperplanes} we mean a pair of affine hyperplanes $P = P_{1} \cup P_{2}$
where $P_{1}$ and $P_{2}$ are hyperplanes (so that $0 \in  P_{1} \cap P_{2}$).\\

 We now briefly explain the basic facts about 2-valued functions needed in this paper. For a detailed 
treatment of multi-valued functions, we refer the reader to \cite{A}.\\

Let $k$ be an integer $\geq 1.$ ($k = 1$ and $k = n$ are the only cases needed in this paper.) 
Denote by ${\bf T} ({\mathbf R}^{k})$ the set of unordered pairs of elements
of ${\mathbf R}^{k}.$ Define a metric ${\mathcal G}$ on ${\bf T}({\mathbf R}^{k})$ by 

$${\mathcal G}(\{v_{1},v_{2}\},\{w_{1}, w_{2}\}) = \min \, \{ \sqrt{|v_{1} - w_{1}|^{2}  + |v_{2} - w_{2}|^{2}},
\sqrt{|v_{1} - w_{2}|^{2}  + |v_{2} - w_{1}|^{2}} \}.$$

If $u \, : \, B_{1}(0) \to {\bf T}({\mathbf R}^{k})$, we say $u$ is 
a 2-valued function on $B_{1}(0)$ with values in ${\bf T}({\mathbf R}^{k}).$ A 2-valued function 
$u \, : \, B_{1}(0) \to {\bf T}({\mathbf R}^{k})$ is continuous if it is continuous with respect to 
the ${\mathcal G}$ metric.\\
 
We say that a 2-valued function $u \, : \, B_{1}(0) \to {\bf T}({\mathbf R}^{k})$ is differentiable 
(or affinely approximable) at 
a point $a \in B_{1}(0)$ if there exist two affine functions $l_{1}^{a}, \, l_{2}^{a} \, 
: \, {\mathbf R}^{n} \to {\mathbf R}^{k}$ 
such that 

$$u(a) = Au(a)(a) \hspace{.2in} {\rm and}$$

$$\lim_{x \to a} \, \frac{{\mathcal G}(u(x), Au(a)(x))}{|x - a|} = 0$$

\noindent
where $Au(a)$ is the 2-valued function defined by $Au(a)(x) = \{l_{1}^{a}(x), l_{2}^{a}(x)\}$ for all 
$x \in {\mathbf R}^{n}.$ 
It follows that if $Au(a)$ exists, it is unique, and that if $u$ is differentiable at $a \in \Omega$ then it is 
continuous at $a$.\\ 

We say that $u$ is differentiable in $B_{1}(0)$ if $u$ is differentiable at $a$ for every $a \in B_{1}(0).$\\

Suppose $u$ is differentiable in $B_{1}(0)$ and $\a \in (0, 1).$ 
We say that $u$ is $C^{1, \, \a}$ in $B_{1}(0),$ and define  
$$\|u\|_{C^{1, \a}(B_{1}(0))} \equiv \sup_{x_{1},x_{2} \in B_{1}(0), \, x_{1} \neq x_{2}} \, \frac{{\mathcal G} \, (\{Dl_{1}^{x_{1}}, Dl_{2}^{x_{1}}\}, 
\{Dl_{1}^{x_{2}}, Dl_{2}^{x_{2}}\})}{|x_{1} - x_{2}|^{\a}},$$

\noindent
provided the right hand side of the above is finite.\\

\section{Blow-ups off affine hyperplanes} \label{blowupprop}
\setcounter{equation}{0}
For $M \in {\mathcal I}_{b}$, $\r \in (0, 3/2]$ and $P$ a pair of affine hyperplanes, define the height 
excess $E_{M}(\r, P)$ of $M$ relative to $P$ at scale $\r$ by\\
 
\begin{equation}\label{fineexcess}
E_{M}^{2}(\r, P) = \r^{-n-2}\int_{M \cap (B_{\r}(0) \times {\mathbf R})} \mbox{dist}^{2} \, (X, P).
\end{equation}

In case $L$ is a single affine hyperplane, we write 

\begin{equation}\label{coarseexcess}
{\hat E}_{M}(\rho, L) = E_{M}(\rho, L).
\end{equation}

Let $\d \in (0, 1)$ be a fixed number, $\{M_{k}\} \subset {\mathcal I}_{b}$  a sequence of 
hypersurfaces such that 

\begin{equation} \label{massbound}
\frac{{\mathcal H}^{n}(M_{k} \cap (B_{2}^{n+1}(0)))}{\omega_{n}2^{n}} \leq 3 - \delta \hspace{.2in} \mbox{and}
\end{equation}

\begin{equation} \label{supconvergence}
{\hat E}_{k} = {\hat E}_{M_{k}}(3/2, L_{k}) \searrow 0
\end{equation}

\noindent
for some sequence $\{L_{k}\}$ of affine hyperplanes of ${\mathbf R}^{n+1}$ converging to
${\mathbf R}^{n} \times \{0\}.$ Note that by a standard argument using the first variation formula~(\ref{firstvar})
(see e.g. proof of inequality (4.18) of \cite{WN}), we then have that for each $\s \in (0, 3/2)$ the 
estimate 

\begin{equation} \label{tiltconvergence}
(E_{M_{k}}^{T}(\s, L_{k}))^{2} \leq \frac{C \, \s^{-n}}{(3/2 - \s)^{2}} \,{\hat E}_{k}^{2} 
\end{equation}

\noindent
where $C = C(n)$ and, for a hypersurface $M \in {\mathcal I}_{b}$ and an affine 
hyperplane $L$, $E_{M}^{T}(\s, L) \equiv \sqrt{\s^{-n}\int_{M \cap (B_{\s}(0) \times {\mathbf R})} 
1 - (\nu \cdot \nu^{L})^{2}}$ is the {\em tilt excess} of $M$ relative to $L$ at scale $\s.$ Here 
$\nu$ and $\nu^{L}$ are the unit normals to $M$ and $L$ respectively.\\ 

We need to blow up the sequence of hypersurfaces $\{M_{k}\}$ off the sequence 
of affine hyperplanes $\{L_{k}\}.$ This is carried out essentially as in \cite{WN}. For convenience, 
we choose here to blow up by the height excess ${\hat E}_{k}$ rather than by the tilt excess 
$E_{M_{k}}^{T}$ 
which was used in \cite{WN}. This is possible in view of (\ref{tiltconvergence}).
Note also that in \cite{WN}, it is assumed that for each $k$, (i) ${\rm sing} \, M_{k}  = \emptyset$ and 
(ii) $M_{k}$ approximates a cone having a singularity at the origin. Here we 
weaken hypothesis (i) to ${\mathcal H}^{n-2}\, ({\rm sing} \, M_{k}) < \infty$ and drop the assumption (ii) altogether. 
The blow up argument of Sections 3 and 4 of \cite{WN}
can however be repeated with some changes to accommodate the weaker hypotheses. We justify 
this assertion as follows:\\

\noindent
\begin{itemize}
\item[(1)] The conclusion of Lemma 3.2 of \cite{WN} holds without change under the present hypotheses. That is to say, 
for each $M \in {\mathcal I}_{b}$ and each bounded, locally Lipschitz function $\varphi$ with 
$\varphi \equiv 0$ in a neighborhood of $M \cap (\partial \, B_{3/2}(0) \times {\mathbf R})$, we have that
for any constant unit vector $\nu_{0}$,

\begin{equation}\label{schoenest}
\int_{M \cap (B_{3/2}(0) \times {\mathbf R})} |A|^{2} \varphi^{2} \leq 
C \int_{M \cap (B_{3/2}(0) \times {\mathbf R})} \left(1 - (\nu \cdot \nu_{0})^{2}\right) |\nabla ^{M} \varphi|^{2}
\end{equation}
 
\noindent
where $A$ denotes the second fundamental  form of $M$, $|A|$ the length of $A$ and $C = C(n).$ 
This estimate was first proved by R. Schoen in \cite{SR}, and 
later used by R. Schoen and L. Simon in \cite{SS} (Lemma 1 of \cite{SS}) under 
the hypothesis that ${\mathcal H}^{n-2}({\rm sing} \, M) = 0.$ We here use an argument of H. Federer and 
W. Ziemer \cite{FZ} (see also \cite{EG}) to justify 
our claim that the estimate in fact continues to hold under the weaker hypothesis ${\mathcal H}^{n-2} \, 
({\rm sing} \, M) < \infty$.\\

First note that by exactly the argument of the proof of Lemma 1 of \cite{SS}, the estimate 
(\ref{schoenest}) holds if $\varphi$ is locally Lipschitz with compact support in 
$M \cap (B_{3/2}(0) \times {\mathbf R}).$ The issue is to argue that it holds for 
bounded, locally Lipschitz $\varphi$ vanishing near $M \cap (\partial \, B_{3/2}(0) \times {\mathbf R})$ 
under the assumption ${\mathcal H}^{n-2} \, ({\rm sing} \, M) < \infty.$ Let $\t \in (0, 1/4)$ be arbitrary. Since 
${\rm sing} \, M \cap ({\overline B_{3/2}(0)} \times {\mathbf R})$ is compact, for 
each $i = 1, \, 2, \, \ldots$ there 
exists a finite number $N_{i}$ and balls $B_{r^{(i)}_{j}}^{n+1}(Z^{(i)}_{j})$, $j = 1, \ldots, N^{(i)}$ with 
$Z^{(i)}_{j} \in {\rm sing} \, M \cap ({\overline B_{3/2}(0)} \times {\mathbf R})$  such that 
${\rm sing} \, M \cap ({\overline B_{3/2}(0)} \times {\mathbf R}) \subset \cup_{j = 1}^{N^{(i)}} \, B_{r^{(i)}_{j}}^{n+1}(Z^{(i)}_{j}),$ 
$\sum_{j=1}^{N^{(i)}} \, \omega_{n-2}(r^{(i)}_{j})^{n-2} \leq K \equiv 1 + 
2^{n-2}{\mathcal H}^{n-2} \, ({\rm sing} \, M \cap (B_{3/2}(0) \times {\mathbf R}))$ and $r^{(i)}_{j} \leq \t^{(i)}.$ 
Here $\t^{(1)} = \t$ and $\t^{(i)} = \frac{1}{4}{\rm dist} \, ({\rm sing} \, M \cap ({\overline B_{3/2}(0)} \times {\mathbf R}), 
{\mathbf R}^{n+1} \setminus U^{(i-1)})$ for $i = 2, \, 3, \, \ldots$, 
where $U^{(i)} = \cup_{j = 1}^{N^{(i)}} B_{r^{(i)}_{j}}^{n+1}(Z^{(i)}_{j}).$
 For each $i = 1, \, 2, \, \ldots$ and each
$j \in \{1, \ldots, N^{(i)}\}$, let $\psi^{(i)}_{j}$ be a $C^{1}$ function on $M$ such that 
$\psi^{(i)}_{j} \equiv 0$ on $B_{r^{(i)}_{j}}^{n+1}(Z^{(i)}_{j}) \cap M$, $\psi^{(i)}_{j} \equiv 1$ on $M \setminus 
B_{2r^{(i)}_{j}}^{n+1}(Z^{(i)}_{j})$, $0 \leq \psi^{(i)}_{j} \leq 1$ everywhere and 
$|\nabla \, \psi^{(i)}_{j}| \leq 2(r^{(i)}_{j})^{-1}.$
Let $\z^{(i)} = {\rm min} \, \{\psi^{(i)}_{1} , \ldots, \psi^{(i)}_{N^{(i)}}\}.$ Then 
${\rm spt} \, |\nabla \, \z^{(i)}| \subset M  \cap (V^{(i)} \setminus V^{(i+1)})$ 
where $V^{(i)} = \cup_{j=1}^{N^{(i)}} B_{2r^{(i)}_{j}}^{n+1}(Z^{(i)}_{j})$ 
and $\int_{M} |\nabla \, \z^{(i)}|^{2} \leq  cK,$  $c = c(n).$ Finally, for $\ell = 1, \, 2, \, \ldots$, let 

\begin{equation*}\label{schoenest0}
\b_{\ell} = \frac{1}{S_{\ell}} \sum_{i=1}^{\ell} \frac{\z^{(i)}}{i}
\end{equation*}

\noindent
where $S_{\ell} = \sum_{i=1}^{\ell} i^{-1}.$ Then, since ${\rm spt} \, \nabla \, \z^{(i)},$ $i = 1, \, 2, \, \ldots$ 
are disjoint, we have that 

\begin{equation}\label{schoenest1}
\int_{M \cap (B_{3/2}(0) \times {\mathbf R})} |\nabla \, \b_{k}|^{2} = 
\frac{1}{S_{\ell}^{2}}\sum_{i=1}^{\ell} \int_{M \cap (B_{3/2}(0) \times {\mathbf R})} i^{-2}|\nabla \, \z^{(i)}|^{2}
\leq \frac{cK}{S^{2}_{\ell}} \sum_{i=1}^{\ell} i^{-2}.
\end{equation}

\noindent
Now, if $\varphi$ is a bounded, locally Lipschitz function vanishing in a neighborhood of 
$M \cap (\partial \, B_{3/2}(0) \times {\mathbf R})$, then for each $\ell$, 
$\b_{\ell}\varphi$ is a locally Lipschitz function with compact support in $M \cap (B_{3/2}(0) \times 
{\mathbf R})$ and hence (\ref{schoenest}) holds with $\b_{\ell}\varphi$ in place of $\varphi.$ Thus

\begin{eqnarray}\label{schoenest2}
\int_{M \cap (B_{3/2}(0) \times {\mathbf R})} |A|^{2} \b^{2}_{\ell} \varphi^{2} &\leq& 
C\int_{M \cap (B_{3/2}(0) \times {\mathbf R})} \left( 1 - (\nu \cdot \nu_{0})^{2}\right)
\b_{\ell}^{2} |\nabla \, \varphi|^{2} \nonumber\\
&& + C \, {\rm sup} \, \varphi^{2} \int_{M \cap (B_{3/2}(0) \times {\mathbf R})} |\nabla \, \b_{\ell}|^{2}.
\end{eqnarray}

Since $\b_{\ell} \leq 1$ and $\b_{\ell} \equiv 1$ on 
$M_{\t} \cap (B_{3/2}(0) \times {\mathbf R})$ where $M_{\t} = 
M \setminus \{X  : \, {\rm dist} \, (X, {\rm sing} \, M) \leq 2\t \}$, 
we conclude from (\ref{schoenest1}) and 
(\ref{schoenest2}) that 

\begin{equation*}\label{schoenest3}
\int_{M_{\t} \cap (B_{3/2}(0) \times {\mathbf R})} |A|^{2}\varphi^{2} \leq 
C \int_{M \cap (B_{3/2}(0) \times {\mathbf R})} \left(1 - (\nu \cdot \nu_{0})^{2}\right)|\nabla \, \varphi|^{2} 
+ \frac{C \, K}{S^{2}_{\ell}} \sum_{i=1}^{\ell} i^{-2}.
\end{equation*}

\noindent
Letting first $\ell \to \infty$ and then $\t \to 0$ in this, we conclude (\ref{schoenest}).\\

\noindent
{\bf Remark}: Note that the validity of (\ref{schoenest}) under the hypothesis 
${\mathcal H}^{n-2} \, ({\rm sing} \, M) < \infty$, as justified above, shows that Schoen-Simon regularity theory \cite{SS}
for embedded stable minimal hypersurfaces $M$ holds under the hypothesis 
${\mathcal H}^{n-2} \, ({\rm sing}  \, M) < \infty$.\\ 
  
\noindent
\item[(2)] Lemma 3.3 of \cite{WN} (which is essentially the same 
as Lemma 2 of \cite{SS}) holds and gives a good approximate graphical decomposition 
of $M_{k}$ {\em relative to the affine hyperplane} $L_{k}$, provided we make 
the minor modification noted in item (3) below, which is necessary due to the presence of 
a singular set. Note that since $L_{k} \to {\mathbf R}^{n} \times \{0\}$, there exists
a sequence of rigid motions $q_{k}$ of ${\mathbf R}^{n+1}$ with $q_{k} 
\to identity$ such that $q_{k}(a_{k}) = \{0\}$ and 
$q_{k} \, L_{k} \equiv {\mathbf R}^{n} \times \{0\},$ where $a_{k}$ is the nearest point of $L_{k}$ 
to $0 \in {\mathbf R}^{n+1}.$
Then, by essentially the same arguments as in \cite{SS}, Section 3 (as detailed in \cite{WN}, Section 3), 
for each given $\s \in (0, 3/2)$ and each sufficiently large $k$
(depending on $\s$), there exists 
a ``good set'' $\widetilde\Omega_{k} = \widetilde\Omega_{k}(\s)  
\subset L_{k} \cap q_{k}^{-1} \, (B_{\s}(0) \times \{0\})$ 
(which corresponds to $\Omega_{k}$ of \cite{WN}, Lemma 3.3),  and two Lipschitz functions 
${\widetilde u}_{k}^{\pm} \, : \, \widetilde\Omega_{k} \to {\mathbf R}$ 
with Lipschitz constants $\leq 1$ (analogous to $u_{k}^{\pm}$ of \cite{WN}, Section 3), 
such that ${\rm graph} \, {\widetilde u}_{k}^{+} \nu^{L_{k}} \cup 
{\rm graph} \, {\widetilde u}_{k}^{-} \nu^{L_{k}} \subseteq M_{k} \cap q_{k}^{-1} (B_{\s}(0) \times {\mathbf R})$
and 

\begin{equation}\label{badsetest}
{\mathcal H}^{n} \, ((M_{k} \setminus ({\rm graph} \, {\widetilde u}_{k}^{+} \nu^{L_{k}} \cup 
{\rm graph} \, {\widetilde u}_{k}^{-} \nu^{L_{k}})) \cap q_{k}^{-1}(B_{\s}(0) \times {\mathbf R})) \leq 
C_{\s}({\hat E}_{k})^{2 + \mu}
\end{equation}

\noindent
where $\nu^{L_{k}}$ denotes the upward pointing unit normal to $L_{k}$, $\mu$ is a fixed constant depending only on 
$n$ and $C_{\s}$ is a constant depending only on $n$ and $\s$. Here 
${\rm graph} \, {\widetilde u}_{k}^{\pm}\nu^{L_{k}} \equiv \{ x + {\widetilde u}_{k}^{\pm}(x)\nu^{L_{k}} \, : 
\, x \in {\widetilde \Omega}_{k}\}.$\\  

In the present paper, we shall use the notation $G_{k}^{\pm} = G_{k}^{\pm}(\s) =
 {\rm graph} \, {\widetilde u}_{k} ^{\pm}\nu^{L_{k}},$ $\Omega_{k}  = q_{k} \, {\widetilde\Omega}_{k}$ and
$u_{k}^{\pm}(x) = {\widetilde u}_{k}^{\pm} \circ q_{k}^{-1}(x)$ for $x \in \Omega_{k}.$\\    

\noindent
\item[(3)] In Lemma 3.3 of \cite{WN}, the definition of $\G_{k}$ needs to be modified to 
$\G_{k} = {\pi}_{L_{k}} \, \{ X \in M_{k} \cap q_{k}^{-1} \, (B_{\s}(0) \times {\mathbf R}) \, : \, g_{k}(X) = \th_{k} \} 
\cup {\pi}_{L_{k}} \, ({\rm sing} \, M_{k}).$ (cf. \cite{SS}.) Here $g_{k}$ and 
$\th_{k}$ are as in \cite{WN}, Section 3, and ${\pi}_{L_{k}}$ is the orthogonal projection of 
${\mathbf R}^{n+1}$ onto $L_{k}$. The conclusions of Lemma 3.3 (with notational 
changes as indicated in item (2) above) hold with this modification and 
with ${\hat E}_{k}$ in place of $\e_{k}$ (where by definition $\e_{k} = $ {\it tilt} excess in \cite{WN}).\\

\noindent
\item[(4)] We may construct cut-off functions ${\widetilde\varphi}_{k}^{0}$, ${\widetilde \psi}_{k}^{(\eta)}$ and 
${\widetilde{\overline \psi}}_{k}$ analogous, respectively, to the cut-off functions 
$\varphi_{k}^{0}$, $\psi_{k}^{(\eta)}$ and ${\overline \psi}_{k}$ of \cite{WN}, Section 3.
The domains of these cut-off functions are 
$q_{k}^{-1} \, (B_{\s}(0) \times \{0\}) \setminus {\pi}_{L_{k}} \, ({\rm sing} \, M_{k}),$
$M_{k} \setminus {\pi}_{L_{k}}^{-1}({\pi}_{L_{k}} \, ({\rm sing} \, M_{k}))$ and 
$q_{k}^{-1}(B_{\s}(0) \times \{0\}) \setminus {\pi}_{L_{k}} \, ({\rm sing} \, M_{k})$ respectively, 
and they take values in ${\mathbf R}.$ We then define ${\overline \psi}_{k} \, : \, B_{\s}(0) \times \{0\}
\setminus q_{k} ({\pi}_{L_{k}} \, {\rm sing} \, M_{k}) 
 \to {\mathbf R}$ by setting ${\overline \psi}_{k}(x) = \widetilde{\overline \psi}_{k} \circ q_{k}^{-1}(x).$
Note that 

\begin{equation}\label{cutoffmeasure}
{\mathcal H}^{n} \, (B_{\s}(0) \setminus \{x \, : \, {\overline \psi}_{k}(x) =1\}) \leq C_{\s}({\hat E}_{k})^{2 + \mu}.
\end{equation}\\ 

\noindent
(See the estimate (3.26) of \cite{WN}, Section 3.)\\
  
\noindent
\item[(5)]  We cannot assume Lemma 3.4 of \cite{WN} in the present context because it 
depends on $M_{k}$ being free of singularities. (Specifically, 
the inequality (3.7) of \cite{WN} assumes that ${\rm sing} \, M_{k} = \emptyset$.) Notice that in \cite{WN}, Lemma 3.4 was used 
precisely at two places; namely,\\

\noindent
\begin{itemize}
\item[$(a)$] to establish the estimate (3.28) of \cite{WN} which bounds the square of the $L^{2}$ norm 
of $|D{\overline \psi}_{k}|$ from 
above by a constant times $(E_{M_{k}}^{T})^{2 + \mu},$ where ${\overline \psi}_{k}$ is the cut-off function 
described in item (4) above
and $\mu = \mu(n) >0$ is as in Lemma 3.5 of \cite{WN}, and\\

\noindent
\item[$(b)$] in the proof of the pointwise gradient estimate for the blow-up (i.e. Lemma 4.9 of \cite{WN}).\\
\end{itemize}

The modifications necessary for $(a)$ above are minor. In fact, it suffices to have the estimate 

\begin{equation}\label{cutoffestimate}
\int_{B_{\s}(0)}|D{\overline \psi}_{k}|^{2} \leq c \, {\hat E}_{k}^{2}, 
\end{equation}

\noindent
$c =c(\s)$, and this weaker estimate follows easily from (\ref{schoenest}) and (\ref{tiltconvergence}) 
in view of the fact that $|D{\overline \psi}_{k}|$
is pointwise bounded  from above by a constant times 
the length of the second fundamental form of $M_{k}.$ (See \cite{WN}, Section 3.) That this weaker estimate 
suffices follows from the fact that $|Du_{k}^{\pm}|$ are bounded, that 
$u_{k}^{\pm} \to 0$ pointwise a.e. and that $Du_{k}^{\pm} \to 0$ in $L^{2}.$\\  

As for $(b)$ above, we shall give an argument in Lemma~\ref{lipschitz} below 
which, under our present (weaker) hypotheses, in fact shows only that the blow-ups are 
continuous and satisfy a Lipschitz condition at points 
where they are single valued. This suffices for proving asymptotic decay estimates for 
the blow-ups later in Section \ref{regblowup}.\\

\noindent    
\item[(6)] Parts $(a), \, (b), \, (f)$ and $(g)$ of Lemma 4.6 of \cite{WN} hold (of course with the 
functions now having domain $B_{\s}(0)$). Thus, letting 
$v_{k}^{\pm} = \frac{\overline{\psi}_{k}u_{k}^{\pm}}{{\hat E}_{k}}$, 
there exist functions $v^{\pm} \in W^{1,2}_{\rm loc}(B_{3/2}(0))$ ---the blow-up of 
$\{M_{k}\}$ off $\{L_{k}\}$---with $v^{+} \geq v^{-}$ such that, 
after passing to a subsequence of $\{k\}$ which we continue to label $\{k\}$, we have

\begin{equation} \label{blowups}
 {\overline \psi}_{k}v_{k}^{\pm} \to v^{\pm}
\end{equation}

\noindent
in $W^{1,2}(B_{\s}(0))$ for each $\s <3/2.$ Note that unlike in 
\cite{WN} (where each $M_{k}$ was assumed to approximate 
a cone arbitrarily closely), $v^{\pm}$ here need not be homogeneous of 
degree 1. Note also that it is easy to see that 
${\overline \psi}_{k}v_{k}^{\pm} \to v^{\pm}$ in $L^{2}(B_{\s}(0))$ 
and weakly in $W^{1,2}(B_{\s}(0))$ 
for each $\s <3/2$ since it follows directly from the definition of ${\hat E}_{k}$ that 
${\overline\psi}_{k}v_{k}^{\pm}$ are uniformly 
bounded in $L^{2}(B_{\s}(0))$, and from the estimates 
(\ref{tiltconvergence}) and (\ref{cutoffestimate}) 
that $D({\overline \psi}_{k}v_{k}^{\pm})$ are uniformly bounded in $L^{2}(B_{\s}(0)).$ 
The proof that the convergence is strong in $W^{1,2}(B_{\s}(0))$ requires only 
some minor modification
of the argument of \cite{WN} used to prove the same assertion (i.e. parts  $(f)$ and $(g)$ of Lemma 4.6,
\cite{WN}.) See item $(8)$ below.\\

\noindent
\item[(7)] $h \equiv \frac{1}{2}(v^{+} + v^{-})$ is harmonic in $B_{3/2}(0).$ The proof of this is as in 
part $(e)$ of Lemma 4.6, \cite{WN}.\\

\noindent
\item[(8)] The necessary modifications to the argument of parts $(f)$ and $(g)$ of Lemma 4.6, \cite{WN} to show 
that the convergence in (\ref{blowups})  is strong in $W^{1,2}(B_{\s}(0))$ 
for each $\s < 3/2$ are as follows: The energy estimate (4.6) of \cite{WN} must be replaced by 

$$\int_{B_{\s}(0) \cap \{|{\overline \psi}_{k}(v_{k}^{+} - h)| \leq \e\}} |D({\overline \psi}_{k} 
(v_{k}^{+} - h))|^{2} 
+ \int_{B_{\s}(0) \cap \{|{\overline \psi}_{k}(v_{k}^{-} - h)| \leq \e\}} |D({\overline \psi}_{k} 
(v_{k}^{-} - h))|^{2}
\leq c\e,$$

\noindent
$c = c(\s),$ and, consequently, the estimate (4.41) of \cite{WN} becomes

$$\int_{B_{\s}(0) \cap \{|v^{+} - h| \leq \e\}} |D(v^{+} - h)|^{2} 
+ \int_{B_{\s}(0) \cap \{|v^{-} - h| \leq \e\}} |D(v^{-} - h)|^{2} \leq c\e,$$

\noindent
$c = c(\s).$ To prove the former estimate, define ${\widetilde h}(x^{\prime}, x^{n+1}) = h(x^{\prime})$ and 
repeat the argument of the proof of the estimate (4.6) of \cite{WN} after replacing 
$\widetilde \z$ in the first variation identity (4.1) of \cite{WN} simply with
$F_{\delta}(x^{n+1} - {\hat E}_{k}{\widetilde h}){\widetilde \z}^{2}$ (rather than with 
$F_{\d}(x^{n+1}){\widetilde \z}^{2}$ which was used in \cite{WN}; here notation is as in \cite{WN}) 
and use the estimate (\ref{tiltconvergence}) above.\\ 

The only other change necessary in the proof of strong convergence is 
that the function $V_{k}^{\e}$ (see paragraph preceding estimate (4.34) of 
\cite{WN})  must now be defined to be
$V_{k}^{\e} = {\overline \psi}_{k}(\g_{\e}(v_{k}^{+}  - h)D(v_{k}^{+} - h) + 
\g_{\e}(v_{k}^{-} - h)D(v_{k}^{-} - h)).$ Of course then subsequent estimates involving 
$V_{k}^{\e}$ need to be modified accordingly in an obvious way.\\ 
\end{itemize}

\noindent
{\bf Remark}: Note that the hypothesis (\ref{massbound}) allows the possibility that 
$M_{k}$ are ``single sheeted,'' in which case, the blow up would be a single valued, harmonic function $v$.
The asymptotic decay estimate as in Lemma~\ref{regblowupthm}, part $(b)$ is much easier to prove and is 
standard in this case. Our analysis throughout the paper however contains this as a special case.\\ 

\noindent
{\bf Definition}: Let ${\mathcal F}_{\delta}$ denote the family of ordered pairs of functions $v = (v^{+}, v^{-})$ 
on the ball $B_{3/2}(0)$ arising as blow-ups of sequences of stable minimal hypersurfaces in the manner 
described above. Precisely, 
each $(v^{+}, v^{-}) \in {\mathcal F}_{\delta}$ 
is the blow-up, as in (\ref{blowups}), 
of a sequence $\{M_{k}\} \subset {\mathcal I}_{b}$ satisfying
(\ref{massbound}) and (\ref{supconvergence}) for some sequence of hyperplanes $L_{k}$ converging to 
${\mathbf R}^{n} \times \{0\}.$ \\

\begin{lemma}\label{compactness}
For each $\s \in (0, 3/2)$, ${\mathcal F}_{\delta}$ is a compact subset of 
$W^{1,2}(B_{\s}(0) \, ; \, {\mathbf R}^{2}).$\\
\end{lemma}

\begin{proof}
The lemma follows directly from the ``diagonal process''. Specifically, let 
$\{(v_{i}^{+}, v_{i}^{-})\}$ be a sequence of functions in ${\mathcal F}_{\delta}.$ Then for each $i$, 
there exists a sequence of hypersurfaces $\{M_{k}^{i}\} \subset {\mathcal I}_{b}$ with 
$\frac{{\mathcal H}^{n} \, (M_{k}^{i} \cap B_{2}^{n+1}(0))}{\omega_{n}2^{n}} \leq 3-\d$ and 
a sequence of affine hyperplanes $L_{k}^{i}$ of ${\mathbf R}^{n+1}$ converging to 
${\mathbf R}^{n} \times \{0\}$ as $k \to \infty$ such that ${\hat E}_{k}^{i} \equiv 
{\hat E}_{M_{k}^{i}}(3/2, L_{k}^{i}) \to 0$ and $(v_{i}^{+}, v_{i}^{-})$ is the blow-up of 
$\{M_{k}^{i}\}$ by ${\hat E}_{k}^{i}.$ Thus, for each $i$, 

\begin{equation} \label{compact1}
\frac{\overline{\psi}_{i, \, k}u_{i, \, k}^{\pm}}{{\hat E}_{k}^{i}} \to v_{i}^{\pm}
\end{equation} 
 
\noindent
as $k \to \infty,$ in $W^{1,2}(B_{\s}(0))$ for each $\s \in (0, 3/2).$ (The notation here is as 
in items (2) and (6) of the discussion at the beginning of this section.) 
Now choose a diagonal sequence $\{M^{i}_{k(i)}\}$, 
$k(1) < k(2) < k(3) < \ldots$ and positive integers $N_{i},$ $i=1, 2, 3, \ldots$ such that 
$L_{k(i)}^{i} \to {\mathbf R}^{n} \times \{0\}$, ${\hat E}^{i}_{k(i)} \to 0$ as $i \to \infty$
and, for each $i$, $\|\frac{\overline{\psi}_{i, \, k(i)}u_{i, \, k(i)}^{\pm}}{E_{k(i)}^{i}} - 
\frac{\overline{\psi}_{i, \, k(j)}u_{i, \, k(j)}^{\pm}}{E_{k(j)}^{i}}\|_{W^{1,2}(B_{\s}(0))} \leq 2^{-i}$ for all 
$j \geq N_{i}.$  (Such $N_{i}$ exist by the convergence (\ref{compact1}).) 
Let $(v^{+}, v^{-}) \in W^{1,2}_{\rm loc}(B_{3/2}(0); {\mathbf R}^{2})$ be the blow-up 
of $\{M^{i}_{k(i)}\}$ by  ${\hat E}^{i}_{k(i)}.$  i.e. for a subsequence $\{i^{\prime}\}$ of $\{i\}$, 
$(v^{+}, v^{-})$ is, for each $\s < 3/2$, the $W^{1,2}(B_{\s}(0); {\mathbf R}^{2})$ limit of the blow-up sequence 
$\Big\{ \, \left(\frac{\overline{\psi}_{i^{\prime}, \, k(i^{\prime})}u_{i^{\prime}, k(i^{\prime})}^{+}}
{{\hat E}^{i^{\prime}}_{k(i^{\prime})}}, \frac{\overline{\psi}_{i^{\prime}, \, 
k(i^{\prime})}u_{i^{\prime}, k(i^{\prime})}^{-}}
{{\hat E}^{i^{\prime}}_{k(i^{\prime})}}\right) \, \Big\}.$  Then, 
by the definition of ${\mathcal F}_{\delta}$, $(v^{+}, v^{-}) \in {\mathcal F}_{\delta},$
and it is easily seen using the triangle inequality that 
$v_{i^{\prime}}^{\pm} \to v^{\pm}$ in $W^{1,2}(B_{\s}(0)).$
\end{proof}

\medskip

\begin{lemma}\label{embedded}
Let $z \in B_{3/2}(0)$ and $\s \in (0, 3/2 - |z|).$ If for all sufficiently 
large $k$, $M_{k} \cap (B_{\s}(z) \times {\mathbf R})$ are embedded, then 
$\left. v^{+}\right|_{B_{\s}(z)}$ and $\left. v^{-} \right|_{B_{\s}(z)}$ are individually (a. e. equal to) harmonic 
functions on $B_{\s}(z).$  
\end{lemma}

\begin{proof}
Under the hypotheses of the lemma, we have that for all sufficiently large $k$, 
$u_{k}^{+} > u_{k}^{-}$ in $\Omega_{k} \cap B_{\s}(z)$ and that $u_{k}^{\pm}$ are (smooth)
solutions of the minimal surface equation:

\begin{equation}\label{mse}
\sum_{j=1}^{n} D_{j} \left(\frac{D_{j}u_{k}^{\pm}}{\sqrt{1 + |Du_{k}^{\pm}|^{2}}}\right) = 0
\end{equation}

\noindent
in $\Omega_{k} \cap B_{\s}(z).$ Let $\z$ be an arbitrary $C^{1}$ function with compact 
support in $B_{\s}(z).$ Multiplying  (\ref{mse}) by ${\overline \psi}_{k} \z$ and integrating 
over $B_{\s}(z)$, we have 

\begin{equation*}
\int_{B_{\s}(z)} \frac{Du_{k}^{\pm} \cdot D({\overline \psi}_{k} \z)}{\sqrt{1 + |Du_{k}^{\pm}|^{2}}} = 0
\end{equation*}

\noindent
which can be written as 

\begin{equation*}
\int_{B_{\s}(z)} \frac{D({\overline \psi}_{k}u_{k}^{\pm})
 \cdot D\z}{\sqrt{1 + |Du_{k}^{\pm}|^{2}}} 
= -\int_{B_{\s}(z)}\frac{\z Du_{k}^{\pm} \cdot D{\overline \psi}_{k}}{\sqrt{1 + |Du_{k}^{\pm}|^{2}}} 
+\int_{B_{\s}(z)}\frac{u_{k}^{\pm} D{\overline \psi}_{k} \cdot D\z}{\sqrt{1 + |Du_{k}^{\pm}|^{2}}}. 
\end{equation*}

\noindent
Dividing this by ${\hat E}_{k}$ and passing to the limit as $k \to \infty$, we conclude using 
the Cauchy-Schwarz inequality and (\ref{cutoffestimate}) that 

$$\int_{B_{\s}(z)} Dv^{\pm} \cdot D\z = 0$$

\noindent
as required.
\end{proof}

\medskip

Any $v = (v^{+}, v^{-}) \in {\mathcal F}_{\delta}$ satisfies the 
properties listed in Propositions~\ref{properties}---\ref{wproperties} below. 
Given $v \in {\mathcal F}_{\d}$, here and subsequently we use the following notation:\\

$$h = \frac{v^{+} + v^{-}}{2}, \hspace{.5in} w = \frac{v^{+} - v^{-}}{2}.$$\\

\medskip

\begin{proposition} \label{properties}
\begin{itemize}
\item[(1)] $h$ is harmonic in $B_{3/2}(0).$\\

\item[(2)] $\int_{B_{3/2}(0)} (v^{+})^{2} + (v^{-})^{2} \leq \left(\frac{3}{2}\right)^{n+2}.$\\

\item[(3)] $\int (|Dv^{+}|^{2} + |Dv^{-}|^{2}) \z = -\int (v^{+}Dv^{+} + v^{-}Dv^{-}) \cdot D\z$ for every 
$\z \in C^{1}_{c}(B_{3/2}(0)).$\\ 
More generally, $\int (|Dv^{+}|^{2} + |Dv^{-}|^{2}) \z = -\int ((v^{+}-y) Dv^{+} + 
(v^{-} - y)Dv^{-}) \cdot D\z$ for every $y \in {\mathbf R}$ and every $\z \in C^{1}_{c}(B_{3/2}(0)).$
By replacing $\z$ with $\z^{2}$ in this and using the Cauchy-Schwarz inequality on the right hand side, 
we get that $\int(|Dv^{+}|^{2} + |Dv^{-}|^{2}) \z^{2} \leq C \int((v^{+}-y)^{2} + (v^{-}-y)^{2})|D\z|^{2}$ 
for any $y \in {\mathbf R}.$\\

\item[(4)] $\int_{B_{\sigma}(z)} (|Dv^{+}|^{2} + |Dv^{-}|^{2}) = \int_{\partial B_{\sigma}(z)} (v^{+}-y)
\frac{\partial v^{+}}{\partial R} + (v^{-}-y) \frac{\partial v^{-}}{\partial R}$ for each $z \in B_{3/2}(0)$ 
and almost every $\s \in (0, 3/2 - |z|).$\\

\item[(5)] $\sum_{i,j=1}^{n} \int_{B_{\sigma}(z)} 
\left((|Dv^{+}|^{2} + |Dv^{-}|^{2})\delta_{ij} - 2D_{i}v^{+}D_{j}v^{+} -2D_{i}v^{-}D_{j}v^{-}\right)D_{i}\z^{j} = 0$ 
for each ball $B_{\s}(z)$ with $\overline B_{\s}(z) \subset B_{3/2}(0)$ and 
each vector field $\z = (\z^{1}, \z^{2} \, \ldots, \z^{n})$ with 
$\z^{j} \in C^{1}_{c}(B_{\sigma}(z))$ for $j = 1,2,3, \ldots,n.$\\ 

\end{itemize}
\end{proposition}

\begin{proof} Part $(2)$ is a direct consequence of the definition (\ref{coarseexcess}) and the 
estimate (\ref{badsetest}). The proofs of parts $(1)$, $(3)$, $(4)$ and $(5)$ are contained in \cite{WN}, Section 4;
part $(1)$ follows from the identity (4.30) of 
\cite{WN}; parts $(3)$, $(4)$ and $(5)$ follow from exactly the arguments of Lemma 4.7, part $(i)$; 
Lemma 4.7, part $(ii)$ and Lemma 4.8 of \cite{WN} respectively. 
\end{proof}

\medskip

\noindent
{\bf Definition:} Let $v \in {\mathcal F}_{\delta}$, $z \in B_{3/2}(0)$ and $y \in {\mathbf R}.$ 
Define the frequency function $N_{v, z, y}(\cdot)$ by 

\begin{equation}\label{frequencyfn}
N_{v, z, y}(\r) = \frac{\r\int_{B_{\r}(z)}|Dv|^{2}}{\int_{\partial \, B_{\r}(z)} (v^{+} - y)^{2} + (v^{-} - y)^{2}}
\end{equation}

\noindent
provided $\r \in (0, 3/2 - |z|)$ and 
$\int_{\partial \, B_{\r}(z)} (v^{+} - y)^{2} + (v^{-} - y)^{2} \neq 0.$\\ 

Whenever $z \in B_{3/2}(0)$ is a Lebesgue point of 
both $v^{+}$ and $v^{-}$, and $v^{+}(z) = v^{-}(z) = y$ (as will be the case in most of our applications of 
the frequency function), we shall let ${\mathcal N}_{v, z}(\r)  = {\mathcal N}_{v, z, y}(\r).$

\medskip

\begin{proposition}\label{freqfn}
Suppose $v \in {\mathcal F}_{\delta},$ $0 \leq \r_{1} < \r_{2},$ $B_{\r_{2}}(z) \subseteq B_{3/2}(0),$
$y \in {\mathbf R}$ and 
$\int_{\partial \, B_{\r}(z)} (v^{+} - y)^{2} + (v^{-} - y)^{2} \neq 0$ for all $\r \in (\r_{1}, \r_{2}).$ 
Then $N_{v, z, y}(\cdot)$ is monotonically non-decreasing in $(\r_{1}, \r_{2}).$
\end{proposition}

\begin{proof}
The argument is the same as in the proof of Lemma 5.13, \cite{WN}. We reproduce it here for 
the reader's convenience. Note first that the identity of Proposition~\ref{properties}, part (4) implies that

\begin{equation} \label{monotonicityofscaledenergy}
\frac{d}{d \rho} \, \left(\rho^{2-n} \int_{B_{\rho}(z)}
|Dv|^2\right) = 2\rho^{2-n} \int_{\partial B_{\rho}(z)}\left| 
\frac{\partial v}{\partial R} \right|^2
\end{equation}

\noindent
for almost all $\rho \in (0, {\rm dist} \, (z, \partial \, B_{3/2}(0))))$, where $\frac{\partial v}{\partial R}(x) = Dv(x) 
\cdot \frac{x - z}{|x - z|}$ is the radial derivative. This follows by taking 
$(x^{j}-z^{j}) \, \zeta_{l}$ in place of $\zeta^{j}$
in the identity of Proposition~\ref{properties}, part (4) and letting $l \to \infty$, where $\zeta_{l}$ 
is a sequence of $C^{\infty}_{c}(B_{\rho}(z))$ functions
converging to the characteristic function of the ball $B_{\rho}(z)$. (We omit the 
details here. This is exactly the argument used to derive the standard monotonicity formula 
for stationary harmonic maps, and can be found e.g. in \cite{S3}, Chapter 2.) Note also that by 
Proposition~\ref{properties}, part (3), 

\begin{equation}\label{freq-1}
\int_{B_{\r}(z)}|Dv|^{2} = \frac{1}{2}\int_{\partial \, B_{\r}(z)} \frac{\partial}{\partial R} \, ((v^{+} - y)^{2} + 
(v^{-} - y)^{2}) 
\end{equation}

\noindent
for almost every $\r \in (0, 3/2 - |z|).$\\
 
Now by a change of variables in the denominator of (\ref{frequencyfn}), we have that
 
$$N_{v,z, y}(\rho) = \frac{\rho^{2-n}\int_{B_{\rho}(z)}\,|Dv|^2}
{\int_{{\bf S}^{n-1}} ({\hat v}^{+}_{z, \r, y})^{2} + ({\hat v}^{-}_{z, \r, y})^{2}}$$

\noindent
where ${\hat v}^{\pm}_{z, \r, y} (\omega) = v^{\pm}(z + \r \,\omega)  - y.$ Using this and the identities 
(\ref{monotonicityofscaledenergy}), (\ref{freq-1}), we have that for 
a.e. $\r \in (\r_{1}, \r_{2})$,\\

\begin{eqnarray} \label{derivativeoffreqfn}
\hspace{.5in} \lefteqn{\frac{d}{d\rho}\,N_{v,z, y}(\rho)} \nonumber\\ 
& =&  \frac{\frac{d}{d\rho} \, \left(\rho^{2-n} \int_{B_{\rho}(z)}
|Dv|^2\right)}{\int_{{\bf S}^{n-1}} ({\hat v}^{+}_{z, \r, y})^{2} + ({\hat v}^{-}_{z, \r, y})^{2}} - 
\frac{\rho^{2-n}\int_{B_{\rho}(z)}|Dv|^2 
\, \frac{d}{d\r} \, \int_{{\bf S}^{n-1}} ({\hat v}^{+}_{z, \r, y})^{2} + ({\hat v}^{-}_{z, \r, y})^{2}}
{\left(\int_{{\bf S}^{n-1}} ({\hat v}^{+}_{z, \r, y})^{2} + ({\hat v}^{-}_{z, \r, y})^{2}\right)^2} \nonumber\\
&=&  \frac{2 \rho^{2-n} \int_{\partial B_{\rho}(z)} \left|\frac{\partial v}{\partial R} \right|^2}
{\int_{{\bf S}^{n-1}} ({\hat v}^{+}_{z, \r, y})^{2} + ({\hat v}^{-}_{z, \r, y})^{2}} -\nonumber\\
&& - \,\,\frac{\frac{\r^{2-n}}{2}
\int_{\partial B_{\rho}(z)} \frac{\partial}{\partial \, R} \, \left((v^{+} - y)^{2} + 
(v^{-} - y)^{2}\right) 
\frac{d}{d\r}\int_{{\bf S}^{n-1}} ({\hat v}^{+}_{z, \r, y})^{2} + ({\hat v}^{-}_{z, \r, y})^{2}}
{\left(\int_{{\bf S}^{n-1}} ({\hat v}^{+}_{z, \r, y})^{2} + ({\hat v}^{-}_{z, \r, y})^{2}\right)^{2}} \nonumber\\
&=&  \frac{2 \rho^{-1}\left( \int_{{\bf S}^{n-1}} ({\hat v}^{+}_{z, \r, y})^{2} + ({\hat v}^{-}_{z, \r, y})^{2}\, 
\int_{{\mathbf S}^{n-1}} \left| \frac{\partial {\hat v}_{z, \r, y}}{\partial R} \right|^2
- \left(\int_{{\mathbf S}^{n-1}} {\hat v}^{+}_{z, \r, y}\frac{\partial {\hat v}^{+}_{z, \r, y}}{\partial \, R} + 
{\hat v}^{-}_{z, \r, y} \frac{\partial {\hat v}^{-}_{z, \r, y}}{\partial \, R} \right)^2\right)}
{\left(\int_{{\bf S}^{n-1}} ({\hat v}^{+}_{z, \r, y})^{2} + ({\hat v}^{-}_{z, \r, y})^{2}\right)^2} \nonumber\\
&\geq& 0.\nonumber\\
\end{eqnarray}

\noindent
The inequality above follows from the Cauchy-Schwarz inequality. This completes the proof.
\end{proof}

\medskip

\noindent
{\bf Remark:} By the definition (\ref{frequencyfn}) of frequency function and the identity (\ref{freq-1}), 
it follows that for $z \in B_{3/2}(0)$ and $\r \in (0, 3/2 - |z|)$,

\begin{equation}\label{alternate}
N_{v, z, y}(\r) = \frac{\r\int_{\partial \, B_{\r}(z)}\frac{\partial}{\partial R} 
((v^{+} - y)^{2} + (v^{-} - y)^{2})}{2\int_{\partial \, B_{\r}(z)} (v^{+} - y)^{2} + (v^{-} - y)^{2}} 
\end{equation}

\noindent
whenever $N_{v, z, y}(\r)$ is defined.\\ 

\begin{lemma} \label{coarseexcessbounds}
Let $v = (v^{+}, v^{-}) \in {\mathcal F}_{\delta}$, $z \in B_{3/2}(0)$ and suppose that 
$\int_{\partial \, B_{\s_{0}}(z)} (v^{+} - y)^{2} + (v^{-} - y)^{2} > 0$ for 
some $\s_{0} \in (0, {\rm dist} \, (z, 3/2 - |z|).$ Then

\begin{itemize}
\item[(a)] $\int_{\partial \, B_{\r}(z)} (v^{+} - y)^{2} + (v^{-} - y)^{2} > 0$ for 
all $\r \in (0, 3/2 - |z|)$ and hence $N_{v, z, y}(\r)$ is defined for all $\r \in (0, 3/2 - |z|).$ 

\item[(b)] For each $\r \in (0, 3/2 - |z|)$ and each $\th \in (0, 1]$, 
$$\frac{\E^{2}_{z, \, y, \, \th\r}}{\E^{2}_{z, \, y, \, \r}} \geq \th^{2(N_{v, z, y}(\r)- 1)}$$ 
\noindent
where $\E_{z, \, y, \, \r} = \left(\r^{-n-2}\int_{B_{\r}(z)} (v^{+} - y)^{2} + (v^{-} - y)^{2}\right)^{1/2}.$ 
\end{itemize}
\end{lemma}

\begin{proof}
Since $\int_{\partial \, B_{\s}(z)} (v^{+} - y)^{2} + (v^{-} - y)^{2}$ is 
(absolutely) continuous as a function of $\s$ and 
$\int_{\partial \, B_{\s_{0}}(z)} (v^{+} - y)^{2} + (v^{-} - y)^{2} > 0$ by hypothesis, 
there exist $\s_{1} \in (0, 3/2 - |z|)$ with 
$\s_{1} < \s_{0}$ such that $\int_{\partial \, B_{\s}(z)} (v^{+} - y)^{2} + (v^{-} - y)^{2} > 0$ 
for all $\s \in (\s_{1}, \s_{0}].$ Hence the frequency function $N_{v, z, y}(\s)$ 
is well defined for all $\s \in (\s_{1}, \s_{0}]$ and by the monotonicity of $N_{v,z, y}(\s)$ and the 
identity (\ref{alternate}), we have that for all $\s \in (\s_{1}, \s_{0}],$

\begin{equation*}
N_{v, z, y}(\s) = \frac{\s \frac{d}{d\s}\int_{{\mathbf S}^{n-1}} (v^{(z,\s) \, +} - y)^{2} + (v^{(z,\s) \, -} -y)^{2}}
{2\int_{{\mathbf S}^{n-1}}(v^{(z,\s) \, +} - y)^{2} + (v^{(z, \s) \, -} - y)^{2}} \leq N_{v,z, y}(\s_{0}) = N_{0}
\end{equation*}

\noindent
where $v^{(z, \s) \, \pm}(\omega) = v^{\pm}(z + \s \, \omega).$ This is equivalent to  

\begin{equation*}\label{cb0}
\frac{d}{d\s} \log \left(\frac{\s^{1-n}\int_{\partial B_{\s}(z)} (v^{+} - y)^{2} + (v^{-} - y)^{2}}{\s^{2N_{0}}}\right) \leq 0
\end{equation*}

\noindent 
and integrating this differential inequality with respect to $\s$ from $\s_{1}$ to $\s_{0}$, we have that  

\begin{equation*} \label{cb1}
\frac{\int_{\partial B_{\s_{0}}(z)} (v^{+} - y)^{2} + (v^{-} - y)^{2}}{\s_{0}^{2N_{0} + n-1}} \leq  
\frac{\int_{\partial B_{\s_{1}}(z)} (v^{+} - y)^{2} + (v^{-} - y)^{2}}{\s_{1}^{2N_{0} + n-1}}.
\end{equation*}

\noindent
This readily implies that $\int_{\partial \, B_{\s_{1}}(z)} (v^{+} - y)^{2} + (v^{-} - y)^{2} > 0.$ Thus 
$\int_{\partial \, B_{\r}(z)} (v^{+} - y)^{2} + (v^{-} - y)^{2} > 0$ for all $\r \in (0, \s_{0}].$ Since 
by Proposition~\ref{properties}, part (2), the function 
$(v^{+} - y)^{2} + (v^{-}  - y)^{2}$ is weakly subharmonic in $B_{3/2}(0),$
it follows from the maximum principle that $\int_{\partial \, B_{\r}(z)} (v^{+} - y)^{2} + (v^{-} - y)^{2} >0$ 
for all $\r \in (\s_{0}, 3/2 - |z|).$ Thus part $(a)$ of the lemma holds.\\

To prove part $(b)$, fix $\r \in (0, 3/2 - |z|).$ Using part $(a)$ and arguing as above, we have that 

\begin{equation*}\label{cb0-1}
\frac{d}{d\s} \log \left(\frac{\s^{1-n}\int_{\partial B_{\s}(z)} (v^{+} - y)^{2} + (v^{-} - y)^{2}}{\s^{2N}}\right) \leq 0
\end{equation*}

\noindent
for all $\s \in (0, \r)$, where $N = N_{v, z, y}(\r)$, and by integrating this from $\s_{1}$ to 
$\s_{2}$, we obtain that for every $\s_{1}, \s_{2} \in (0, 3/2 - |z|)$ with 
$0 < \s_{1} < \s_{2} \leq \r$, 

\begin{equation} \label{cb1-0}
\frac{\int_{\partial B_{\s_{1}}(z)} (v^{+} - y)^{2} + (v^{-} - y)^{2}}{\s_{1}^{2N + n-1}} \leq  
\frac{\int_{\partial B_{\s_{2}}(z)} (v^{+} - y)^{2} + (v^{-} - y)^{2}}{\s_{2}^{2N + n-1}}.
\end{equation}

\noindent
Hold $\s_{1}$ fixed with $0 < \s_{1} < \th\r$, multiply inequality (\ref{cb1-0}) by $\s_{2}^{2N + n -1}$ and integrate 
with respect to $\s_{2}$ from $\th\r$ to $\r$ to obtain, for each $\s_{1} \in (0, \th\r)$, that 

\begin{equation}\label{cb1-1}
\int_{B_{\r}(z) \setminus B_{\th\r}(z)} (v^{+} - y)^{2} + (v^{-} - y)^{2}\leq  
\frac{1}{2N + n} \left(\r^{2N + n} - (\th\r)^{2N + n}\right)
\frac{\int_{\partial B_{\s_{1}}(z)} (v^{+} - y)^{2} + (v^{-} - y)^{2}}{\s_{1}^{2N + n-1}}.
\end{equation}

\noindent
Now multiply both sides of (\ref{cb1-1}) 
by $\s_{1}^{2N + n -1}$ and integrate with respect to $\s_{1}$ from $0$ to $\th\r$. This
gives

\begin{equation*}
(\th\r)^{2N+n}\int_{B_{\r}(z) \setminus B_{\th\r}(z)} (v^{+} - y)^{2} + (v^{-} - y)^{2}\leq  
\left(\r^{2N + n} - (\th\r)^{2N + n}\right)
\int_{B_{\th\r}(z)} (v^{+} - y)^{2} + (v^{-} - y)^{2}
\end{equation*}

\noindent
which, upon rearrangement of terms, gives the desired estimate.
\end{proof}

\medskip

\noindent
{\bf Definition:} For $v \in {\mathcal F}_{\delta}$, $z \in B_{3/2}(0)$ and $y\in {\mathbf R}$
with $\int_{\partial \, B_{\s_{0}}(z)} (v^{+} - y)^{2}  + (v^{-} - y)^{2}>0$ for some $\s_{0} \in (0, 3/2 - |z|),$ 
define ${\mathcal N}_{v, y}(z) = \lim_{\r \downarrow 0} \, N_{v, z, y}(\r).$ Note that $N_{v, z, y}(\r)$ 
is well defined for all $\r \in (0, 3/2 - |z|)$ and this limit exists
by Lemmas~\ref{coarseexcessbounds} and ~\ref{freqfn} above.\\

Whenever $z \in B_{3/2}(0)$ is a Lebesgue point of both $v^{+}$ and $v^{-},$ and $v^{+}(z) = v^{-}(z) = y$, we shall
let ${\mathcal N}_{v}(z) = {\mathcal N}_{v, y}(z).$\\

\medskip

\begin{lemma}\label{constfreq}
Let $v \in {\mathcal F}_{\delta}$ and $z \in B_{3/2}(0).$ Suppose that 
$\int_{\partial \, B_{\s_{0}}(z)} (v^{+} - y)^{2} + (v^{-} - y)^{2} > 0$ for some 
$\s_{0} \in (0, 3/2 - |z|).$ Then $N_{v,z,y}(\rho)$ is constant  for $\r \in (0, 3/2 - |z|)$ {\rm (}with value 
${\mathcal N}_{v,y}(z)${\rm )} 
if and only if $\sqrt{(v^{+} - y)^{2} + (v^{-} - y)^{2}}$ is 
homogeneous of degree ${\mathcal N}_{v, y}(z)$ from the point $z$ in $B_{3/2 - |z|}(z);$
i.e. if and only if 
$$(v^{+}(z + \rho \, \omega)  - y)^{2} + (v^{-}(z + \r \, \omega) - y)^{2}
= \left(\frac{\r}{\r^{\prime}}\right)^{2{\mathcal N}_{v,y}(z)}((v^{+}(z + 
\r^{\prime} \, \omega) - y)^{2} + (v^{-}(z + \r^{\prime} \, \omega) -y)^{2})$$

\noindent 
for each $\r, \, \r^{\prime} \in (0, 3/2 - |z|)$ and $\omega \in {\mathbf S}^{n-1}.$   
\end{lemma}

\begin{proof}
Note first that $N_{v, z, y}(\r)$ is well defined for 
$\r \in (0, 3/2 - |z|)$ by Lemma~\ref{coarseexcessbounds}, part $(a).$ 
If $\sqrt{(v^{+} - y)^{2} + (v^{-} - y)^{2}}$ is homogeneous of some degree 
$\alpha$ from $z$ in $B_{3/2 - |z|}(z)$, it is easy to see using 
the identity (\ref{alternate}) that $N_{v,z, y}(\r) = \alpha (={\mathcal N}_{v,y}(z))$ for $\r \in (0, 3/2 - |z|).$ Conversely, 
suppose $N_{v,z, y}(\r)$ is constant in the interval $(0, 3/2 - |z|).$ Then by (\ref{derivativeoffreqfn}), 

$$\frac{\partial}{\partial R} \, {\hat v}^{\pm}_{z, \r, y}(\omega) 
= \alpha {\hat v}^{\pm}_{z, \r, y}(\omega)$$

\noindent
for some constant $\alpha$, almost all $\r \in (0, 3/2 - |z|)$ and almost all $\omega \in {\mathbf S}^{n-1}$, 
where ${\hat v}^{\pm}_{z, \r, y}(\omega) = v^{\pm}(z +\r \omega) - y$. 
(This just follows from the condition under which equality holds in Cauchy-Schwarz inequality.) 
This is equivalent to the differential identities $\r\frac{d}{d\r} \, (v^{\pm}(z + \r \, \omega)  - y) 
= \a (v^{\pm}(z + \r \, \omega) - y)$ 
which imply that $(v^{+}(z + \r \, \omega)  - y)^{2} 
+ (v^{-}(z + \r \, \omega) - y)^{2} =  \left(\frac{\r}{\r^{\prime}}\right)^{2\a}
((v^{+}(z + \r^{\prime} \, \omega) - y)^{2} + (v^{-}(z + \r \, \omega) - y)^{2})$ for 
$\r, \, \r^{\prime} \in (0, 3/2 - |z|)$ and $\omega \in {\mathbf S}^{n-1}.$ It then follows from (\ref{alternate}) that 
$\a = {\mathcal N}_{v, y}(z).$ 
\end{proof}

\medskip

The estimate in Lemma~\ref{radialexcessbd} below, essentially due to Hardt and Simon \cite{HS}, will play a 
very important role first in our 
proof of continuity of functions in ${\mathcal F}_{\delta}$ (Lemma~\ref{lipschitz} below) and 
later in establishing crucial asymptotic decay properties (Theorem~\ref{regblowupthm}) of these functions.
In the proof of this estimate, we shall need the following:\\
 
\begin{lemma}\label{heightbound}
Let $\s \in (0, 3/2).$ There exist $\e = \e(n, \s) \in (0, 1)$ and $C = C(n,\s) \in (0, \infty)$ such that 
if $M \in {\mathcal I}_{b}$, $\frac{{\mathcal H}^{n} \, (M \cap B_{2}^{n+1}(0))}{\omega_{n}2^{n}} \leq 
3-\d,$ $L$ is an affine hyperplane of ${\mathbf R}^{n+1}$ with 
${\rm dist} \, (L \cap (B_{1}(0) \times {\mathbf R}), B_{1}(0)) < \e$ and  
${\hat E}  =  {\hat E}_{M}(3/2, L) < \e$, then for each $Z \in M \cap (B_{\s}(0) \times {\mathbf R})$ with 
$\Theta_{M} \, (Z) \geq 2$ we have that
$${\rm dist} \, (Z, L) \leq C \, {\hat E}.$$  
\end{lemma}

\begin{proof}
By translating, scaling and rotating, we may assume without loss of generality 
that $L = {\mathbf R}^{n} \times \{0\}.$ Let 
$Z$ be as in the statement of the lemma and write $Z = (z^{\prime}, z^{n+1}).$ 
Set $\s_{0} = 3/2 - \s.$ The monotonicity formula for $M$ (\cite{S1}, Section 17) says that 

\begin{equation}\label{ht-1}
\int_{M \cap B_{\s_{0}/2}^{n+1}(Z)} \frac{(\nu \cdot (X - Z))^{2}}{|X - Z|^{n+2}} 
= \frac{{\mathcal H}^{n} \, (M \cap B_{\s_{0}/2}^{n+1}(Z))}{\omega_{n}(\s_{0}/2)^{n}} - \Theta_{M}(Z)
\end{equation}

\noindent
where $\nu$ denotes the unit normal to $M$. 
Writing $\nu = (\nu^{\prime}, \nu^{n+1})$ where $\nu^{n+1} = \nu \cdot e^{n+1}$, we have that 

\begin{eqnarray}\label{ht-2}
\int_{M \cap B_{\s_{0}/2}^{n+1}(Z)} \frac{(\nu \cdot (X - Z))^{2}}{|X - Z|^{n+2}} &\geq&
(\s_{0}/2)^{-n-2}\int_{M \cap B_{\s_{0}/2}^{n+1}(Z)} (\nu^{\prime} \cdot (x^{\prime} - z^{\prime}) + 
\nu^{n+1} (x^{n+1} - z^{n+1}))^{2}\nonumber\\
&\geq& \frac{1}{2}(\s_{0}/2)^{-n-2}\int_{M \cap B_{\s_{0}/2}^{n+1}(Z)} |x^{n+1} - z^{n+1}|^{2} -\nonumber\\
&& - (\s_{0}/2)^{-n}\int_{M \cap B_{\s_{0}/2}^{n+1}(Z)} 1 - (\nu^{n+1})^{2}\nonumber\\
&\geq&   \frac{1}{2}(\s_{0}/2)^{-n-2}\int_{M \cap B_{\s_{0}/2}^{n+1}(Z)} |x^{n+1} - z^{n+1}|^{2} -\nonumber\\ 
&&- c\, \s_{0}^{-n-2} \int_{M \cap (B_{\s_{0}}(z^{\prime}) \times {\mathbf R})} |x^{n+1}|^{2}\nonumber\\
&\geq&   \frac{1}{2}(\s_{0}/2)^{-n-2}\int_{M \cap B_{\s_{0}/2}^{n+1}(Z)} |x^{n+1} - z^{n+1}|^{2} 
- c\, \s_{0}^{-n-2} {\hat E}^{2}
\end{eqnarray}

\noindent
where $c = c(n)$, and for the second of the inequalities in the above, we have used  
$(a + b)^{2} \geq a^{2}/2 - b^{2}$ with $a = \nu^{n+1}(x^{n+1} - z^{n+1})$, 
$b = \nu^{\prime} \cdot (x^{\prime} - z^{\prime})$ and the fact that $|\nu^{\prime}|^{2} 
= 1 -  (\nu^{n+1})^{2},$ and the third inequality is standard and 
is analogous to (\ref{tiltconvergence}).\\

On the other hand, provided $\e = \e(n, \s)$ is sufficiently small, we have that

\begin{eqnarray}\label{ht-3}
\frac{{\mathcal H}^{n} \, (M \cap B_{\s_{0}/2}^{n+1}(Z))}{\omega_{n}(\s_{0}/2)^{n}} - \Theta_{M}(Z)
&\leq& \frac{{\mathcal H}^{n} \, (M \cap B_{\s_{0}/2}^{n+1}(Z))}{\omega_{n}(\s_{0}/2)^{n}} - 2 \nonumber\\
&\leq& C \s_{0}^{-n}\int_{\Omega \cap B_{\s_{0}/2}(z^{\prime})} 
\sqrt{1 + |Du^{+}|^{2}} - 1 + \nonumber\\
&&+C \s_{0}^{-n}\int_{\Omega \cap B_{\s_{0}/2}(z^{\prime})} \sqrt{1 + |Du^{-}|^{2}} - 1
+ C\s_{0}^{-n} {\hat E}^{2 + \mu}\nonumber\\
&=&C \s_{0}^{-n}\int_{\Omega \cap B_{\s_{0}/2}(z^{\prime})} \frac{|Du^{+}|^{2}}
{1 + \sqrt{1 + |Du^{+}|^{2}}} + \nonumber\\
&& + C \s_{0}^{-n}\int_{\Omega \cap B_{\s_{0}/2}(z^{\prime})} \frac{|Du^{-}|^{2}}
{1 + \sqrt{1 + |Du^{-}|^{2}}} + C\s_{0}^{-n} {\hat E}^{2 + \mu}\nonumber\\
&\leq& C \s_{0}^{-n}\int_{M \cap (B_{\s_{0}/2}(z^{\prime}) \times {\mathbf R})}  1 - (\nu^{n+1})^{2}
+ C\s_{0}^{-n} {\hat E}^{2}\nonumber\\
&\leq&C\s_{0}^{-n-2}\int_{M \cap (B_{\s_{0}}(z^{\prime}) \times {\mathbf R})} |x^{n+1}|^{2} + C\s_{0}^{-n}{\hat E}^{2}
\nonumber\\
&\leq&C\s_{0}^{-n-2}{\hat E}^{2}
\end{eqnarray}

\noindent
where $C = C(n)$ and $\Omega$, $u^{+}$ and $u^{-}$ correspond, respectively, to $\Omega_{k}$, $u^{+}_{k}$ and 
$u^{-}_{k}$ of item (2) of the discussion (with
$M$ in place of $M_{k}$ and $L_{k} \equiv {\mathbf R}^{n} \times \{0\}$) at the beginning of 
Section~\ref{blowupprop}. Note that we have used the estimate (\ref{badsetest}) here.\\

Combining the estimates (\ref{ht-2}) and  (\ref{ht-3}), we have

\begin{equation}\label{ht-4}
\int_{M \cap B_{\s_{0}/2}^{n+1}(Z)} |x^{n+1} - z^{n+1}|^{2} \leq C{\hat E}^{2}.
\end{equation} 

\noindent
Since by monotonicty $\frac{{\mathcal H}^{n} \, (M \cap 
B_{\s_{0}/2}^{n+1}(Z))}{\omega_{n} (\s_{0}/2)^{n}} \geq \Th_{M} \, (Z) \geq 2$, it follows 
readily from the estimate (\ref{ht-4}) and the triangle inequality that

\begin{equation*}
|z^{n+1}|^{2} \leq C\s_{0}^{-n}{\hat E}^{2}
\end{equation*}

\noindent
where $C = C(n).$ This is the required estimate.
\end{proof}

\medskip
 
\begin{lemma} \label{radialexcessbd}
Suppose $v = (v^{+}, v^{-}) \in {\mathcal F}_{\delta}$. If $z \in B_{3/2}(0)$ is a Lebesgue point 
of both $v^{+}$ and $v^{-}$, $v^{+}(z) = v^{-}(z)$
and $v^{+} \not\equiv v^{-}$ (as $L^{2}$ functions) in any ball centered at $z$, then for each 
$\rho \in (0, 3/2 - |z|)$, we have that 
\begin{equation*}
\int_{B_{\rho/2}(z)} R^{2-n}\left(\frac{\partial}{\partial R} \left(\frac{v^{+} - y}{R}\right)\right)^{2} 
+ R^{2-n}\left(\frac{\partial}{\partial R} \left(\frac{v^{-} - y}{R}\right)\right)^{2} 
\leq C\rho^{-n-2} \int_{B_{\rho}(z)} (v^{+} - y)^{2} + (v^{-} - y)^{2}
\end{equation*}

\noindent
where $y = v^{+}(z) = v^{-}(z).$ Here $C = C(n).$\\
\end{lemma}

\begin{proof}
Suppose the hypotheses of the lemma are satisfied for some $z \in B_{3/2}(0).$
Let $\{ M_{k}\} \subset {\mathcal I}_{b}$ be a sequence of hypersurfaces
whose blow-up is $v$.  
First we claim that for each $\t \in (0, 3/2 - |z|)$, there exist infinitely many $k$ such that $M_{k} \cap 
(B_{\t}(z) \times{\mathbf R})$ contains a point $Z_{k}$ with $\Theta_{M_{k}}(Z_{k}) \geq 2.$
For if not,  $M_{k} \cap (B_{\t}(z) \times {\mathbf R})$ would be embedded for some $\t \in (0, 3/2 - |z|)$ and 
all sufficiently
large $k$, and hence, by Lemma~\ref{embedded}, $v^{+}$ and 
$v^{-}$ would both be individually harmonic in $B_{\t}(z).$ Since $v^{+} \geq v^{-}$ and $v^{+}(z) = v^{-}(z)$, 
we would then have by the maximum principle that $v^{+} \equiv v^{-}$ in $B_{\t}(z),$ contradicting 
one of the hypotheses of the lemma. Hence the claim must be true.\\
 
Now take an arbitrary sequence of numbers $\t_{j} \searrow 0$
and apply this claim with $\t_{j}$ in place of $\t.$ This gives a subsequence of $\{k\},$ which we 
continue to denote $\{k\},$ such that
$M_{k} \cap (B_{3/2}(0) \times {\mathbf R})$ contains a point $Z_{k} = (Z_{k}^{\prime}, Z_{k}^{n+1})$ 
with $\Theta_{M_{k}}(Z_{k}) \geq 2,$ satisfying $Z_{k}^{\prime} \to z.$
By the usual monotonicity identity for minimal submanifolds (\cite{S1}, Section 17), 
we have that, for $\r \in (0, 3/2 - |z|)$,

\begin{eqnarray*}
\int_{M_{k} \cap B_{\r/2}^{n+1}(Z_{k})} \frac{((X - Z_{k}) \cdot \nu_{k})^{2}}{|X - Z_{k}|^{n+2}} &=& 
\frac{{\mathcal H}^{n}(M_{k} \cap B_{\r/2}^{n+1}(Z_{k}))}{\omega_{n}(\r/2)^{n}} - \Theta_{M_{k}}(Z_{k}) 
\nonumber\\
&\leq& \frac{{\mathcal H}^{n}(M_{k} \cap B_{\r/2}^{n+1}(Z_{k}))}{\omega_{n}(\r/2)^{n}} - 2 
\end{eqnarray*}

Estimating as in (\ref{ht-3}), we have 

\begin{eqnarray*}
\frac{{\mathcal H}^{n}(M_{k} \cap B_{\r/2}^{n+1}(Z_{k}))}{\omega_{n}(\r/2)^{n}} - 2 &=& 
\frac{{\mathcal H}^{n}((G_{k}^{+} \cup G_{k}^{-}) \cap B_{\r/2}^{n+1}(Z_{k}))}{\omega_{n}(\r/2)^{n}} - 2 + \nonumber\\
&&+ \frac{{\mathcal H}^{n}((M_{k} \setminus G_{k}) \cap B_{\r/2}^{n+1}(Z_{k}))}{\omega_{n}(\r/2)^{n}} \nonumber\\
&\leq& \frac{1}{\omega_{n}(\r/2)^{n}} \int_{B_{\r/2}(Z_{k}^{\prime})} 
\left(\sqrt{1 + |D(\overline{\psi}_{k}u_{k}^{+})|^{2}} -1 \right) + \nonumber\\
&&+ \frac{1}{\omega_{n}(\r/2)^{n}} \int_{B_{\r/2}(Z_{k}^{\prime})} 
\left(\sqrt{1 + |D(\overline{\psi}_{k}u_{k}^{-})|^{2}} -1 \right)
+ \frac{C{\hat E}_{k}^{2 + \mu}}{\omega_{n}(\r/2)^{n}} \nonumber\\
&=& \frac{1}{\omega_{n}(\r/2)^{n}} \int_{B_{\r/2}(Z_{k}^{\prime})} \frac{|D(\overline{\psi}_{k}u_{k}^{+})|^{2}}
{\sqrt{1 + |D(\overline{\psi}_{k}u_{k}^{+})|^{2}} +1} + \nonumber\\
&&+ \frac{1}{\omega_{n}(\r/2)^{n}} \int_{B_{\r/2}(Z_{k}^{\prime})} \frac{|D(\overline{\psi}_{k}u_{k}^{-})|^{2}}
{\sqrt{1 + |D(\overline{\psi}_{k}u_{k}^{-})|^{2}} +1} + \frac{C{\hat E}_{k}^{2 + \mu}}{\omega_{n}(\r/2)^{n}} 
\end{eqnarray*}

\noindent
which implies that 
\begin{equation}\label{rad-1}
\limsup_{k \to \infty} \frac{1}{{\hat E}_{k}^{2}} \left(\frac{{\mathcal H}^{n}(M_{k} \cap B_{\r/2}(Z_{k}))}
{\omega_{n}(\r/2)^{n}} - 2\right) \leq \frac{1}{2\omega_{n} (\r/2)^{n}} \int_{B_{\r/2}(z)} |Dv^{+}|^{2} 
+ |Dv^{-}|^{2}
\end{equation}

\noindent
On the other hand,

\begin{eqnarray*}
\int_{M_{k} \cap B_{\r/2}^{n+1}(Z_{k})} \frac{((X - Z_{k}) \cdot \nu_{k})^{2}}{|X - Z_{k}|^{n+2}} &\geq&
\int_{G_{k}^{+} \cap B_{\r/2}^{n+1}(Z_{k})} \frac{((X - Z_{k}) \cdot \nu_{k})^{2}}{|X - Z_{k}|^{n+2}} +
\int_{G_{k}^{-} \cap B_{\r/2}^{n+1}(Z_{k})} \frac{((X - Z_{k}) \cdot \nu_{k})^{2}}{|X - Z_{k}|^{n+2}} \nonumber\\
&\geq& \int_{B_{\r/2}^{n}(Z_{k}^{\prime})} 
\frac{\left(-(X^{\prime} - Z_{k}^{\prime}) \cdot D(\overline{\psi}_{k}u_{k}^{+})
+(\overline{\psi}_{k}u_{k}^{+} - Z_{k}^{n+1})\right)^{2}}{\left((\overline{\psi}_{k}u_{k}^{+} - Z_{k}^{n+1})^{2}
+ |X^{\prime} - Z_{k}^{\prime}|^{2}\right)^{\frac{n+2}{2}}} + \nonumber\\
&&+ \int_{B_{\r/2}^{n}(Z_{k}^{\prime})} 
\frac{\left(-(X^{\prime} - Z_{k}^{\prime}) \cdot D(\overline{\psi}_{k}u_{k}^{-})
+(\overline{\psi}_{k}u_{k}^{-} - Z_{k}^{n+1})\right)^{2}}{\left((\overline{\psi}_{k}u_{k}^{-} - Z_{k}^{n+1})^{2}
+ |X^{\prime} - Z_{k}^{\prime}|^{2}\right)^{\frac{n+2}{2}}}. 
\end{eqnarray*}

This implies by Fatou's lemma and (\ref{rad-1}) that

\begin{eqnarray*}\label{rad-2}
C\r^{-n} \int_{B_{\r/2}(z)} |Dv^{+}|^{2} 
+ |Dv^{-}|^{2} &\geq&
\liminf_{k \to \infty} \frac{1}{{\hat E}_{k}^{2}} \int_{M_{k} \cap B_{\r/2}(Z_{k})} 
\frac{((X - Z_{k}) \cdot \nu_{k})^{2}}{|X - Z_{k}|^{n+2}}\nonumber\\
&\geq& \int_{B_{\r/2}^{n}(z)} \frac{((v^{+} - y) - (X^{\prime} - z) \cdot Dv^{+})^{2}}{|X^{\prime}-z|^{n+2}} + \nonumber\\
 &&+ \int_{B_{\r/2}^{n}(z)} \frac{((v^{-} - y) - (X^{\prime}-z) \cdot Dv^{-})^{2}}{|X^{\prime}-z|^{n+2}} \nonumber\\
&=& \int_{B_{\r/2}(z)} R^{2-n}\left(\frac{\partial}{\partial R}\left(\frac{v^{+} - y}{R}\right)\right)^{2} + \nonumber\\
&&+ \int_{B_{\r/2}(z)} R^{2-n}\left(\frac{\partial}{\partial R}\left(\frac{v^{+} - y}{R}\right)\right)^{2}
\end{eqnarray*}

\noindent
where $C = C(n)$ and 
$y  = \lim_{k \to \infty} \, \frac{Z_{k}^{n+1}}{{\hat E}_{k}}$, possibly after passing to a subsequence of $\{k\}$.
Note that the existence of $y < \infty$ follows from Lemma~\ref{heightbound}. 
The required estimate follows by combining the inequalities (\ref{rad-1}) and 
(\ref{rad-2}), and using Proposition~\ref{properties}, part (2). Observe that the estimate automatically implies 
that $y = v^{+}(z) = v^{-}(z)$ for if not, the integral on the left hand side would not be finite.
\end{proof}

\medskip

\noindent
{\bf Remark}: Note that the proof of the preceding lemma shows the following:  If $v = (v^{+}, v^{-}) \in 
{\mathcal F}_{\delta}$, $z \in B_{3/2}(0)$ is a Lebesgue point of both $v^{+}$ and $v^{-}$, $v^{+}(z) = v^{-}(z) = y$, 
$v^{+} \not\equiv v^{-}$ (as $L^{2}$ functions) in any ball centered at $z$,  and 
if $\{M_{k}\}$ is a sequence of hypersurfaces in ${\mathcal I}_{b}$ whose blow-up is $v$, then there 
exist a subsequence $\{k_{j}\}$ of $\{k\}$ and points $Z_{k_{j}} = (Z_{k_{j}}^{\prime}, Z_{k_{j}}^{n+1}) 
\in M_{k_{j}} \cap (B_{3/2}(0) \times {\mathbf R})$ such that
$\Th_{M_{k_{j}}}(Z_{k_{j}}) \geq 2$ and 
$\lim_{j \to \infty} \, \left(Z_{k_{j}}^{\prime}, \frac{Z_{k_{j}}^{n+1}}{{\hat E}_{k_{j}}}\right) = (z, y).$\\

\begin{lemma} \label{crossings}
Let $(v^{+}, v^{-}) \in {\mathcal F}_{\delta}$ and $\{M_{k}\}$ be a sequence of 
hypersurfaces in ${\mathcal I}_{b}$ whose blow-up is $(v^{+}, v^{-}).$ 
If $z \in B_{3/2}(0)$ is a Lebesgue point of both $v^{+}$ and $v^{-}$, and if 
$v^{+}(z) > v^{-}(z),$ then there exists $\b > 0$ such that  
$M_{k} \cap (B_{\b}(z) \times {\mathbf R})$ are embedded for all sufficiently large $k,$ and 
hence $v^{+}$ and $v^{-}$ are individually harmonic in $B_{\b}(z).$\\ 
\end{lemma}

\begin{proof}
Suppose that $z \in B_{3/2}(0)$ is a Lebesgue point of both $v^{+}$ and $v^{-}$, $v^{+}(z) > v^{-}(z)$ but 
for no $\b >0$, 
$M_{k} \cap (B_{\b}(z) \times {\mathbf R})$ are embedded for all sufficiently large $k$. Then, taking 
$\b = 1/j$, we can find a subsequence $\{k_{j}\}$ of $\{k\}$ such that  
there exists $Z_{k_{j}}  = (Z_{k_{j}}^{\prime}, Z_{k_{j}}^{n+1}) \in M_{k_{j}}\cap (B_{1/j}(z) \times {\mathbf R})$ 
with $\Theta_{M_{k_{j}}}(Z_{k_{j}}) \geq 2.$ In particular, 
$Z_{k_{j}}^{\prime} \to z.$ By the argument of  
Lemma~\ref{radialexcessbd} above, we then have that

\begin{equation}
\int_{B_{\r}(z)} R^{2-n}\left(\frac{\partial}{\partial R} \left(\frac{v^{+} - y}{R}\right)\right)^{2} 
+ R^{2-n}\left(\frac{\partial}{\partial R} \left(\frac{v^{-} - y}{R}\right)\right)^{2} < \infty
\end{equation}

\noindent
for any $\r \in (0, 3/2 - |z|)$ and some $y \in {\mathbf R}$, implying that $v^{+}(z) = v^{-}(z) \,(= y)$. This 
contradiction shows that there exists $\b > 0$ such that 
$M_{k} \cap (B_{\b}(z) \times {\mathbf R})$ are embedded for all sufficiently large $k$. 
It then follows from Lemma~\ref{embedded} that $v^{+}$ and $v^{-}$ 
are individually harmonic in $B_{\b}(z)$. The lemma is thus proved.
\end{proof}

\medskip

\noindent
{\bf Remark:} Note that the proof of the above lemma shows the following:
If for some $\b \in (0, 1)$ there is no $z \in {\overline B}_{\b}(0)$ such that $v^{+}(z) = v^{-}(z)$, then 
$M_{k} \cap (B_{\b}(0) \times {\mathbf R})$ are embedded for all sufficiently large $k$.\\

\medskip

In the next lemma and subsequently, we shall use the following notation: 
for any $v = (v^{+}, v^{-}) \in {\mathcal F}_{\delta}$ and any $\r \in (0, 3/2)$, 

\begin{equation}\label{blowupnotation-0}
\widetilde{v}_{\r} = (\widetilde{v}^{+}_{\r}, \widetilde{v}^{-}_{\r}) \equiv 
\left(\frac{v^{+}_{\r}}{\E_{\r}}, 
\frac{v^{-}_{\r}}{\E_{\r}}\right)
\end{equation}

\noindent
where $v^{\pm}_{\r}(x) = \frac{v^{\pm}(\frac{2\r x}{3})}{\frac{2\r}{3}}$  and 
$\E^{2}_{\r} = \r^{-n-2} \int_{B_{\r}(0)} (v^{+})^{2} + (v^{-})^{2}.$ More generally, if 
$v = (v^{+}, v^{-}) \in {\mathcal F}_{\delta}$, $z \in B_{3/2}(0)$ and $y \in {\mathbf R}$, we let, for
$\r \in (0, 3/2 - |z|)$, 

\begin{equation}\label{blowupnotation-1}
\widetilde{v}_{z, \, \r, \, y} = (\widetilde{v}^{+}_{z, \, \r, \, y}, \widetilde{v}^{-}_{z, \, \r, \, y}) \equiv 
\left(\frac{v^{+}_{z, \, \r, \, y}}{\E_{z, \, \r, \, y}}, 
\frac{v^{-}_{z, \, \r, \, y}}{\E_{z, \, \r, \, y}}\right)
\end{equation}

\noindent
where $v^{\pm}_{z, \, \r, \, y}(x) = \frac{v^{\pm}(z + \frac{2\r x}{3}) - y}{\frac{2\r}{3}}$  and 
$\E^{2}_{z, \, \r, \, y} = \r^{-n-2} \int_{B_{\r}(z)} (v^{+} - y)^{2} + (v^{-} - y)^{2}.$ \\ 

Note that if $v \in {\mathcal F}_{\delta}$, $z \in B_{3/2}(0)$, $\r \in (0, 3/2 - |z|)$ and $y \in {\mathbf R}$,
then ${\widetilde v}_{z, \, \r, \, y} \in {\mathcal F}_{\delta}.$ In fact if $v$ is the blow-up 
(in the sense of Section~\ref{blowupprop}) of the sequence
of hypersurfaces $\{M_{k}\} \subset {\mathcal  I}_{b}$ off the sequence $\{L_{k}\}$ of affine hyperplanes 
converging to ${\mathbf R}^{n} \times \{0\}$, then ${\widetilde v}_{z, \, \r, \, y}$ 
is the blow-up of the sequence 
${\widetilde M}_{k} \equiv \eta_{(z,{\hat E}_{k}y), \, \frac{2}{3}\r} \, M_{k}$  off the sequence 
${\widetilde L}_{k} \equiv \left(\frac{2}{3}\r\right)^{-1}\left( L_{k} + {\hat E}_{k}\nu^{L_{k}} 
- (z, {\hat E}_{k}y)\right)$ of affine hyperplanes, where ${\hat E}_{k}$ and $\nu^{L_{k}}$ are as defined in 
Section~\ref{blowupprop}. The fact that $\frac{{\mathcal H}^{n} \, ({\widetilde M}_{k} \cap B_{2}^{n+1}(0))}
{\omega_{n}2^{n}} \leq 3 - \d$ for sufficiently large $k$ is easily checked using the approximate graphical 
decomposition (as given by the method of \cite{SS} and explained in the discussion of 
item (2) at the beginning of the present section) 
of $M_{k} \cap (B_{2 - \e}(0) \times {\mathbf R})$ for a suitably small fixed positive $\e$ 
independent of $k$.\\

Finally, if $v \in {\mathcal F}_{\delta}$ and $z \in B_{3/2}(0)$ is a Lebesgue point of both 
$v^{+}$ and $v^{-}$ with $v^{+}(z) = v^{-}(z) = y$, we let, for $\r \in (0, 3/2 - |z|)$,  

\begin{equation}\label{blowupnotation-2}
{\mathcal E}_{z, \r} = {\mathcal E}_{z, \r, y} \hspace{.2in} \mbox{and}  \hspace{.2in} 
{\widetilde v}^{\pm}_{z, \, \r} = {\widetilde v}^{\pm}_{z, \, \r, \, y}.
\end{equation}

\medskip

\begin{proposition}\label{lipschitz}
\begin{itemize}
\item[(a)] If $v \in {\mathcal F}_{\delta}$, then $v$ is (a. e. equal to) a continuous function on $B_{3/2}(0).$
\item[(b)]  For each $\s, \s^{\prime} \in (0, 3/2)$ with $\s^{\prime} < \s$, there exists a finite number 
$C = C(n, \s, \s^{\prime})$ such that if $v \in {\mathcal F}_{\delta}$ and $v^{+}(z) = v^{-}(z)$ 
for some point $z \in B_{\s^{\prime}}(0)$, then 
$$|v(x) - v(z)| \leq C |x - z|\left(\int_{B_{\s}(0)}|v|^{2} \right)^{1/2}$$ 
\noindent
for all $x \in B_{\s^{\prime}}(0).$\\ 
\end{itemize}
\end{proposition}

\begin{proof} 
Let $v \in {\mathcal F}_{\delta}.$ Denote by $\G$ the set of points $z \in B_{3/2}(0)$ 
with the property that there exists $y  = y_{z} \in {\mathbf R}$ 
satisfying 

\begin{equation}\label{lip1}
\int_{B_{\rho/2}(z)} R^{2-n}\left(\frac{\partial}{\partial R} \left(\frac{v^{+} - y}{R}\right)\right)^{2} 
+ R^{2-n}\left(\frac{\partial}{\partial R} \left(\frac{v^{-} - y}{R}\right)\right)^{2} 
\leq C\rho^{-n-2} \int_{B_{\rho}(z)} (v^{+} - y)^{2} + (v^{-} - y)^{2}
\end{equation}

\noindent
for all $\r \in (0, 3/2 - |z|),$ where the constant $C = C(n)$ is as in Lemma~\ref{radialexcessbd}. We claim 
that any $z \in \G$ must be a Lebesgue point of both $v^{+}$ and $v^{-}$ with $v^{+}(z) = v^{-}(z) = y$ 
and that for $z \in \G$, a local Lipschitz estimate 

\begin{equation}\label{lip2}
|v(x) - v(z)|^{2} \leq {\widetilde C}|x - z|^{2} \left(\r_{z}^{-n-2}\int_{B_{\r_{z}}(z)} |v|^{2}\right)
\end{equation}

\noindent
must hold for some $\r_{z} \in (0, 3/2 - |z|)$ and a.e.  $x \in B_{\r_{z}/2}(z),$ where ${\widetilde C} 
= {\widetilde C}(n).$  
In order to prove these claims, fix $z \in \G$ and first note that we may suppose that at least one of 
$v^{+}$ or $v^{-}$ is non-constant in every 
ball $B_{\r}(z)$, $0 < \r < 3/2 - |z|,$ for if both $v^{+}$ and $v^{-}$ were constant 
in some ball $B_{\r^{\prime}}(z)$, $\r^{\prime} \in (0, 3/2 - |z|),$ then by (\ref{lip1}) the value of the 
constant must be $y$, and hence we have the claims trivially with 
(\ref{lip2}) holding for $\r_{z} = \r^{\prime}$ and ${\widetilde C} = 1.$ Then we must have that

\begin{equation}\label{lip2-1}
\int_{\partial \, B_{\r}(z)} (v^{+} - y)^{2}  + (v^{-} - y)^{2}> 0 \hspace{.2in} \mbox{for} \hspace{.1in} 
\r \in (0, 3/2 - |z|)
\end{equation}

\noindent
because otherwise, since $(v^{+} - y)^{2} + (v^{-} - y)^{2}$ is subharmonic 
in $B_{3/2}(0)$ (by Proposition~\ref{properties}, part 2), we would 
have by the maximum principle a $\r > 0$ such that $(v^{+}(x) - y)^{2} + (v^{-}(x) - y)^{2} = 0$ for 
a.e. $x \in B_{\r}(z),$ contrary to the preceding assumption. Hence (\ref{lip2-1}) must hold, so that the frequency 
function $N_{v, z, y}(\r)$ is defined for $\r \in (0, 3/2 - |z|)$ and is monotonically non-decreasing. 
We claim that 

\begin{equation}\label{lip0}
{\mathcal N}_{v}(z) \geq 1.
\end{equation} 

\noindent
To see this, note that by (\ref{lip1}) for each $\r \in (0, 3/2 - |z|)$, 

\begin{equation}\label{lip3}
\int_{B_{1/2}(0)} R^{2-n}\left(\frac{\partial \, \left({\widetilde v}^{+}_{z, \r, y}/R\right)}{\partial R}\right)^{2}
+ R^{2-n}\left(\frac{\partial \, \left({\widetilde v}^{-}_{z, \r, y}/R\right)}{\partial R}\right)^{2} \leq C \int_{B_{1}(0)} 
\left({\widetilde v}^{+}_{z, \r, y}\right)^{2} + \left({\widetilde v}^{-}_{z, \r, y}\right)^{2}
\end{equation}

\noindent
where the notation is as in (\ref{blowupnotation-1}). Since 
${\widetilde v}_{z, \r, y} \in {\mathcal F}_{\delta},$ we have by Lemma~\ref{compactness} that for an 
arbitrary sequence 
$\r_{j} \downarrow 0^{+}$, after passing to a subsequence
which we continue to denote $\{j\}$, that 
${\widetilde v}_{z, \r_{j}, y} \to {\widetilde v} \in {\mathcal F}_{\delta}$, where the 
convergence is in $W^{1, 2}(B_{\s}(0))$ for every $\s \in (0, 3/2).$ By (\ref{lip3}), 

\begin{equation}\label{lip3-1}
\int_{B_{1/2}(0)} R^{2-n} \left(\frac{\partial \, \left({\widetilde v}^{+}/R\right)}{\partial R}\right)^{2}
+ R^{2-n}\left(\frac{\partial \, \left({\widetilde v}^{-}/R\right)}{\partial R}\right)^{2} \leq C\int_{B_{1}(0)} ({\widetilde v}^{+})^{2} +
({\widetilde v}^{-})^{2} < \infty
\end{equation} 

\noindent
and we also have by Lemma~\ref{coarseexcessbounds}  that for each 
$\r \in (0, 1]$, $\int_{B_{\r}(0)} |{\widetilde v}_{z, \r_{j}, y}|^{2} 
\geq \left(\frac{2\r}{3}\right)^{2(N_{v, z, y}(\frac{3}{2} - |z|) - 1)}$, and hence that 
${\widetilde v} \not\equiv 0$ in any ball $B_{\r}(0),$ $\r \in (0, 1].$ 
Consequently, since ${\widetilde v}^{2}$ is subharmonic (by Proposition~\ref{properties}, part (2)), 
we have that $\int_{\partial \, B_{\r}(0)} {\widetilde v}^{2} > 0$ for $\r \in (0, 1]$, and 
therefore the frequency function $N_{{\widetilde v}, 0, 0}(\r)$ is defined for 
$\r \in (0, 1].$ But then 

\begin{equation}\label{lip3-2}
N_{{\widetilde v}, 0, 0} (\r) = \frac{\r\int_{B_{\r}(0)} |D{\widetilde v}|^{2}}
{\int_{\partial \, B_{\r}(0)}|{\widetilde v}|^{2}} = \lim_{j \to \infty}\frac{\frac{2}{3}\r\r_{j}
\int_{B_{\frac{2}{3}\r\r_{j}}(z)} |Dv|^{2}}{\int_{\partial \, B_{\frac{2}{3}\r\r_{j}}(z)}(v^{+} - y)^{2} + 
(v^{-} - y)^{2}}
= {\mathcal N}_{v, y}(z)
\end{equation}

\noindent
for $\r \in (0, 1]$ and hence by Lemma~\ref{constfreq},  
${\widetilde v}$ is homogeneous of degree ${\mathcal N}_{v, y}(z)$ from the origin. 
It then follows directly from the finiteness condition (\ref{lip3-1}) that 
${\mathcal N}_{v, y}(z) \geq 1.$\\

With $z \in \G$ and $y = y_{z},$ we next claim that $\r^{-n-2}\int_{B_{\r}(z)} (v^{+} - y)^{2} + (v^{-} - y)^{2}$ is monotonically non-decreasing 
for $\r \in (0, 3/2 - |z|).$ To see this, 
we use the abbreviation 
$d_{v, z}(x) = \sqrt{(v^{+}(x) - y)^{2} + (v^{-}(x) - y)^{2}},$ and compute as follows:\\

\begin{eqnarray}\label{lip4}
\frac{d}{d\r} \, \rho^{-n-2} \int_{B_{\rho}(z)} d_{v, z}^{2} &=& \frac{d}{d\r} \int_{B_{1}(0)} 
\frac{d^{2}_{v, z} \, (z + \r x)}{\r^{2}} \, dx \nonumber\\
&=& \int_{B_{1}(0)} \frac{2d_{v, z} \, (z + \r x) D \, d_{v, z} \, (z + \r x) \cdot x}{\r^{2}} 
- \frac{2d_{v, z}^{2} \, (z + \r x)}{\r^{3}} \nonumber\\
&=& \frac{2}{\r^{3}} \int_{B_{1}(0)} d_{v, z} \, (z + \r x)\left(D \, d_{v, z} \, (z + \r x) \cdot \r x 
- d_{v, z} \, (z + \r x)\right) \nonumber\\ 
&=&2\rho^{-n-3} \int_{B_{\r}(z)} d_{v, z}({\widetilde x})\left(D \, d_{v, z}({\widetilde x}) \cdot ({\widetilde x}-z) 
- d_{v, z}({\widetilde x})\right) \, d{\widetilde x} \nonumber\\
&=& 2\r^{-n-3} \int_{0}^{\r} \int_{\partial B_{\t}(z)} \left(d_{v, z}({\widetilde x})D \, d_{v, z}({\widetilde x}) 
\cdot({\widetilde x}-z) - d_{v, z}^{2}({\widetilde x})\right) \, 
d{\widetilde x} \, d\t \nonumber\\
&=& 2\r^{-n-3} \int_{0}^{\r}\left(\frac{1}{2}\t \int_{\partial \, B_{\t}(z)}
\frac{\partial}{\partial \, R}d_{v, z}^{2} - \int_{\partial \, B_{\t}(z)} d_{v, z}^{2} \right) \,d\t \nonumber\\
&\geq & 0
\end{eqnarray}

\noindent
for almost every $\r \in (0, 3/2 - |z|).$ The last inequality holds since 
$1 \leq {\mathcal N}_{v, y}(z) \leq N_{z, v, y}(\t) = \frac{\t\int_{\partial \, B_{\t}(z)}\frac{\partial}{\partial R} 
d_{v, z}^{2}}{2\int_{\partial \, B_{\t}(z)}d_{v,z}^{2}}, \, \, $ by (\ref{lip0}) and (\ref{alternate}).
Thus in particular, $\r^{-n-2}\int_{B_{\r}(z)} d^{2}_{v,z}$ remains bounded from above as $\r \to 0$ 
and consequently $z$ must be a Lebesgue point of both $v^{+}$ and $v^{-}$ with 
$v^{+}(z) = v^{-}(z) = y.$\\
 
Now, $|v|^{2}$ is subharmonic in $B_{3/2}(0)$ by Proposition~\ref{properties}, part (2), and hence by 
the mean value property  

\begin{equation}\label{lip4-1-1}
|v(z)|^{2} \leq \omega_{n}^{-1}\r_{z}^{-n}\int_{B_{\r_{z}}(z)} |v|^{2}
\end{equation}

\noindent
where $\r_{z} = \frac{1}{2}(\frac{3}{2} - |z|).$ Also, since $d^{2}_{v, z}$ is subharmonic, again 
by the mean value property we have that for a.e. $x \in B_{\r_{z}/2}(z)$,

\begin{eqnarray}\label{lip5-1}
d^{2}_{v, z}(x) &\leq& \omega_{n}^{-1}(|x - z|)^{-n}\int_{B_{|x-z|}(x)}d^{2}_{v, z}\nonumber\\
&\leq& \omega_{n}^{-1}(|x-z|)^{-n}\int_{B_{2|x-z|}(z)}d^{2}_{v, z}\nonumber\\
&=&\omega_{n}^{-1} \, 2^{n+2}|x -z|^{2} (2|x-z|)^{-n-2}\int_{B_{2|x-z|}(z)} d^{2}_{v, z} \nonumber\\
&\leq& \omega_{n}^{-1} \, 2^{n+2} |x-z|^{2} 
\r_{z}^{-n-2}\int_{B_{\r_{z}}(z)} d^{2}_{v, z}\nonumber\\
&\leq&C|x-z|^{2} \r_{z}^{-n-2}\int_{B_{\r_{z}}(z)}|v|^{2}
\end{eqnarray} 

\noindent
where $C = C(n).$ Here we have used the monotonicity of $\r^{-n-2}\int_{B_{\r}(z)} d_{v, z}^{2}$ and 
the estimate (\ref{lip4-1-1}). This is the required estimate (\ref{lip2}).\\ 

We have thus shown that every $z \in \G$ is a Lebesgue point of both $v^{+}$ and 
$v^{-}$ with
$v^{+}(z) = v^{-}(z) = y_{z},$ and that the local Lipschitz estimate (\ref{lip2}) holds at such $z$.\\ 

Now consider a point $z \in B_{3/2} \setminus \G.$ We claim that there exists $\s_{z} \in (0, 3/2 - |z|)$ 
such that $\left. v^{+}\right|_{B_{\s_{z}}(z)}$ and $\left. v^{-}\right|_{B_{\s_{z}}(z)}$ are 
respectively a. e. equal to harmonic functions ${v}^{z \, +}$ and $v^{z \, -}$ on $B_{\s_{z}}(z).$ To see this, 
consider a sequence of hypersurfaces $\{M_{k}\} \subset {\mathcal I}_{b}$ whose blow-up is $v$. 
There must exist $\s_{z} \in (0, 3/2-|z|)$ such that for all sufficiently large $k$, 
$M_{k} \cap (B_{\s_{z}}(z) \times {\mathbf R})$ must be embedded. For if not, there 
exists a subsequence $\{k_{j}\}, \, j = 1, 2, \ldots$ of $\{k\}$ and points $Z_{k_{j}}  = 
(Z^{\prime}_{k_{j}}, Z^{n+1}_{k_{j}}) \in M_{k_{j}} \cap (B_{1/j}(z) \times {\mathbf R})$ with 
$\Theta_{M_{k_{j}}}(Z_{k_{j}}) \geq 2$ and by exactly the argument of Lemma~\ref{radialexcessbd}, 
this implies that (\ref{lip1}) holds for some $y \in {\mathbf R}$ and all $\r \in (0, 3/2 - |z|)$,
contradicting the fact that $z \in B_{3/2}(0) \setminus \G.$  The claim now follows 
from Lemma~\ref{embedded}. Now define ${\overline v}^{\pm} \, : \, B_{3/2}(0) \to {\mathbf R}$ 
by setting ${\overline v}^{+}(z) = v^{z \, +}(z)$, ${\overline v}^{-}(z) = v^{z \, -}(z)$ 
if $z \in B_{3/2}(0) \setminus \G$ and ${\overline v}^{+}(z) = {\overline v}^{-}(z) = y_{z}$ 
if $z \in \G.$ Since $\G$ is relatively closed in $B_{3/2}(0)$ (which follows 
directly from the definition of $\G$), it follows 
by unique continuation for harmonic functions and the continuity 
estimate (\ref{lip2}) for points $z \in \G$ that ${\overline v}^{\pm}$ 
are well defined and are continuous in $B_{3/2}(0).$ Furthermore, 
$v^{\pm}$ are a. e. equal to ${\overline v}^{\pm}.$ This concludes the proof of 
part $(a)$ of the lemma.\\ 

To prove part $(b)$, let $v \in {\mathcal F}_{\delta}$ ($v$ now assumed to be continuous),
$z \in B_{3/2}(0)$ and suppose that $v^{+}(z) = v^{-}(z) = y.$ Note first 
that we must have that either $v^{+} \equiv v^{-} \equiv y$ in $B_{3/2}(0)$
or that $\int_{\partial \, B_{\s}(z)} (v^{+} - y)^{2} + (v^{-} - y)^{2} >0$ 
for all $\s \in (0, 3/2 - |z|).$ To see this, first note that if 
$\int_{\partial \, B_{\s_{0}}(z)} (v^{+} - y)^{2} + (v^{-} - y)^{2} > 0$ for 
some $\s_{0} \in (0, 3/2 - |z|)$, then by continuity, there exists $\s_{1} \in (0, \s_{0})$
such that  $\int_{\partial \, B_{\s}(z)} (v^{+} - y)^{2} + (v^{-} - y)^{2} > 0$ for 
all $\s \in (\s_{1}, \s_{0}].$ Hence the frequency function $N_{v, z}(\s)$ is 
defined for $\s \in (\s_{1}, \s_{0}]$ and by exactly the argument  leading to 
(\ref{cb1}), we have the estimate

\begin{equation} \label{lip5-1-1}
\frac{\int_{\partial B_{\s_{0}}(z)} (v^{+} - y)^{2} + (v^{-} - y)^{2}}{\s_{0}^{2N + n-1}} \leq  
\frac{\int_{\partial B_{\s}(z)} (v^{+} - y)^{2} + (v^{-} - y)^{2}}{\s^{2N + n-1}}
\end{equation}

\noindent
for each $\s \in (\s_{1}, \s_{0}],$ where $N = N_{v,z}(\s_{0}).$ Letting $\s \to \s_{1}$ in this, 
we see that $\int_{\partial \, B_{\s_{1}}(z)} (v^{+} - y)^{2} + (v^{-} - y)^{2} > 0.$ This argument 
shows that  if $\int_{\partial \, B_{\s_{0}}(z)} (v^{+} - y)^{2} + (v^{-} - y)^{2} > 0$ for 
some $\s_{0} \in (0, 3/2 - |z|)$ then $\int_{\partial \, B_{\s}(z)} (v^{+} - y)^{2} + (v^{-} - y)^{2} > 0$ for 
all $\s \in (0, \s_{0}].$  On the other hand, since $(v^{+} - y)^{2} + (v^{-} - y)^{2}$ is 
subharmonic, if 
$\int_{\partial \, B_{\s}(z)} (v^{+} - y)^{2} + (v^{-} - y)^{2} = 0$ for some 
$\s \in (0, 3/2 - |z|)$, then by the maximum principle we must have that 
$v^{+}(x) = v^{-}(x) = y$ for all $x \in {\overline B_{\s}(z)}.$ Hence, either 
$\int_{\partial \, B_{\s}(z)} (v^{+} - y)^{2} + (v^{-} - y)^{2} > 0$ for all 
$\s \in (0, 3/2 - |z|)$ or $v^{+}(x) = v^{-}(x) = y$ for all $x \in B_{3/2 - |z|}(z).$ If the latter
were the case, it is easy to see using the estimate (\ref{lip5-1-1}) repeatedly 
with suitably chosen center points in place of $z$ that we must have 
$v^{+}(x) = v^{-}(x) = y$ for all $x \in B_{3/2}(0).$\\

If $v^{+} \equiv v^{-} \equiv y$ in $B_{3/2}(0)$, the estimate in part 
$(b)$ holds trivially. Otherwise, we have by the above argument that
the frequency function $N_{v, z}(\s)$ is well defined for $\s \in (0, 3/2 - |z|)$ and we 
claim that ${\mathcal N}_{v}(z) \geq 1.$  
This is easy to see if $\left. v^{+}\right|_{B_{\s}(z)} \equiv \left. v^{-}\right|_{B_{\s}(z)}$ for 
some $\s \in (0, 3/2 - |z|)$, because then $v^{+} = v^{-} = h$ in $B_{\s}(z)$ (where 
$h = \frac{1}{2}(v^{+} + v^{-})$) and hence, since $h$ is 
harmonic (everywhere in $B_{3/2}(0)$),  it follows in this case that
${\mathcal N}_{v}(z) = {\mathcal N}_{h - h(z)}(z) \geq 1.$ 
Else, by Lemma \ref{radialexcessbd}, we have the estimate 
(\ref{lip1}) for each $\r \in (0, 3/2 - |z|),$ and we may then argue exactly 
as in the proof of (\ref{lip0}) above to conclude that 
${\mathcal N}_{v}(z) \geq 1.$ Consequently, we also have the monotonicity estimate (\ref{lip4}),
by the same computation, for each $\r \in (0, 3/2 - |z|).$\\ 

To complete the proof of part $(b)$, let $\s, \s^{\prime} \in (0, 3/2)$ with 
$\s^{\prime} < \s,$ and suppose that $z \in B_{\s^{\prime}}(0)$ and that 
$v^{+}(z) = v^{-}(z).$ Since $|v|^{2}$ is subharmonic, we have by the mean value 
property that  

\begin{equation}\label{lip4-1}
\sup_{B_{\s^{\prime}}(0)} \, |v|^{2} \leq C\int_{B_{\s}(0)} |v|^{2}
\end{equation}

\noindent
where $C = C(n,\s, \s^{\prime}).$ Also, since $d^{2}_{v, z} \equiv (v^{+} - y)^{2} 
+ (v^{-} - y)^{2}$ is subharmonic, again 
by the mean value property we have that for every $x \in B_{\s^{\prime}}(0)$,

\begin{eqnarray}\label{lip5}
d^{2}_{v, z}(x) &\leq& \omega_{n}^{-1}(\lambda|x - z|)^{-n}\int_{B_{\lambda|x-z|}(x)}d^{2}_{v, z}\nonumber\\
&\leq& \omega_{n}^{-1}(\lambda|x-z|)^{-n}\int_{B_{2\lambda|x-z|}(z)}d^{2}_{v, z}\nonumber\\
&=&C|x -z|^{2} (2\lambda|x-z|)^{-n-2}\int_{B_{2\lambda|x-z|}(z)} d^{2}_{v, z} \nonumber\\
&\leq&C|x-z|^{2} \int_{B_{\s - \s^{\prime}}(z)} (v^{+} - y)^{2} + (v^{-} - y)^{2}\nonumber\\
&\leq&C|x-z|^{2}\int_{B_{\s}(0)}|v|^{2}
\end{eqnarray} 

\noindent
where $\lambda = \frac{\s - \s^{\prime}}{2(\s + \s^{\prime})}$ and $C = C(n , \s, \s^{\prime}).$ Here we have used 
the monotonicity of $\r^{-n-2}\int_{B_{\r}(z)} d_{v, z}^{2}$ and the estimate
(\ref{lip4-1}). This completes the proof of part $(b)$ and the lemma. 
\end{proof}

\medskip

We next establish several important properties of $w$:\\

\begin{proposition}\label{wproperties}
Suppose $v = (v^{+}, v^{-}) \in {\mathcal F}_{\delta}$ and recall the notation $w = \frac{1}{2}(v^{+} - v^{-}).$ 
We have the following:
\begin{itemize}
\item[(1)] $w \geq 0.$\\

\item[(2)] $\int |Dw|^{2} \z = -\int wDw \cdot D\z$ for every $\z \in C^{1}_{c}(B_{3/2}(0)).$\\

\item[(3)] $\int_{B_{\sigma}(z)} |Dw|^{2} = \int_{\partial B_{\sigma}(z)} w 
\frac{\partial w}{\partial R}$ for each ball $B_{\s}(z)$ with $\overline B_{\s}(z) \subset B_{3/2}(0).$\\

\item[(4)] $\sum_{i,j=1}^{n} \int_{B_{\sigma}(z)} 
\left(|Dw|^{2} \delta_{ij} - 2D_{i}wD_{j}w \right)D_{i}\z^{j} = 0$
for every ball $B_{\s}(z)$ with $\overline B_{\s}(z) \subset B_{3/2}(0)$ and 
every $\z^{j} \in C^{1}_{c}(B_{\sigma}(z)),$ $j = 1,2,3,\ldots,n.$\\

\item[(5)] $\Delta w = 0$ in $B_{3/2}(0) \setminus Z_{w}$ where $Z_{w}$ is the zero set of $w$.\\

\item[(6)] either $Z_{w}  = \emptyset$ or ${\mathcal H}^{n-2}(Z_{w}) = \infty.$\\

\item[(7)] if $\int_{\partial \, B_{\r_{1}}(z_{1})} w^{2} > 0$ for some $z_{1} \in B_{3/2}(0)$ and 
$\r_{1} \in (0, 3/2 - |z_{1}|)$, then $\int_{\partial \, B_{\r}(z_{1})} w^{2} >0$ for all $\r \in (0, \r_{1}].$\\
 
\item[(8)] Either $w \equiv 0$ in $B_{3/2}(0)$ or $\int_{\partial \, B_{\r}(z)} w^{2} > 0$ for each 
$z \in B_{3/2}(0)$ and each $\r \in (0, 3/2 - |z|).$\\
 
\item[(9)] Either $w \equiv 0$ in $B_{3/2}(0)$ or the frequency function 
$N_{w, z}(\rho) \equiv \frac{\rho\int_{B_{\rho}(z)} |Dw|^{2}}
{\int_{\partial B_{\rho}(z)} w^{2}}$ is defined for each $z \in B_{3/2}(0)$ and each
$\r \in (0, 3/2 - |z|)$ and is 
monotonically non-decreasing as a function of $\rho.$
Hence ${\mathcal N}_{w}(z) \equiv \lim_{\rho \downarrow 0} \, N_{w, z}(\rho)$ exists for 
each $z \in B_{3/2}(0)$ unless $w \equiv 0.$\\
\end{itemize}
\end{proposition}

\noindent
\begin{proof} Part (1) follows from the definition of $w$. 
Part  (2) follows directly by substituting $v^{+} = h + w$, $v^{-} = h-w$ 
in the identity of 
part (2) of Proposition~\ref{properties}, and observing that $h$, being harmonic, 
satisfies the identity $\int |Dh|^{2} \z = -\int hDh \cdot D\z.$ Similarly, part
(4) follows by substituting $v^{+} = h + w$, $v^{-} = h-w$ in the identity of part (4) of 
Proposition~\ref{properties} and observing that 
$\sum_{i,j=1}^{n} \int_{B_{\sigma}(z)} 
\left(|Dh|^{2} \delta_{ij} - 2D_{i}hD_{j}h \right)D_{i}\z^{j} = 0$. Part
(3) follows from part (2) by taking a smooth approximation to the characteristic function of the ball 
$B_{\s}(z)$. Part (5) follows from Lemma~\ref{crossings}.\\

To see part (6), note first that it suffices to show that for each given 
$\s \in (0, 3/2)$, either $Z_{w} \cap B_{\s}(0) = \emptyset$ or 
${\mathcal H}^{n-2} \, (Z_{w} \cap B_{\s}(0)) = \infty$. So fix $\s \in (0, 3/2)$ and 
suppose that ${\mathcal H}^{n-2} \, (Z_{w} \cap B_{\s}(0)) < \infty.$  
By continuity of $w$ (Lemma~\ref{lipschitz}),  $Z_{w}$ is closed, so that by exactly the same 
construction as in (\ref{schoenest0}), we have for each 
$\t \in (0, 3/2 - \s)$ a sequence of Lipschitz functions $\b_{\ell} \, : \, B_{3/2}(0) \to {\mathbf R},$ 
$\ell = 1, \, 2, \, 3, \, \ldots$, with $\b_{\ell}(x) \equiv 1$ for each $\ell$ and each $x$ with ${\rm dist} \, (x, Z_{w} \cap B_{\s}(0)) > 
\t$, $\b_{\ell} \equiv 0$ in some neighborhood of $Z_{w} \cap B_{\s}(0),$ 
$0 \leq \b_{\ell} \leq 1$ everywhere and $\int_{B_{3/2}(0)} |D \b_{\ell}|^{2} \to 0$ as $\ell \to \infty.$ Now, given  
$\varphi \in C^{\infty}_{c}(B_{\s}(0))$, we have that $\b_{\ell}\varphi$ is 
Lipschitz with compact support in $B_{\s}(0) \setminus Z_{w}$, and hence, 
since $w$ is harmonic in $B_{3/2}(0) \setminus Z_{w}$, 

$$\int_{B_{\s}(0)} Dw \cdot D(\b_{\ell}\varphi) = 0$$

\noindent
which implies that 

$$\int_{B_{\s}(0) \setminus (Z_{w})_{\t}} Dw \cdot D\varphi
 = -\int_{B_{\s}(0) \cap (Z_{w})_{\t}} \b_{\ell} Dw \cdot D\varphi + \int_{B_{\s}(0)} \varphi Dw \cdot D\b_{\ell}$$

\noindent
where $(Z_{w})_{\t}$ denotes the $\t$ neighborhood of $Z_{w}.$ Hence

\begin{eqnarray}\label{w-prop0}
\left| \int_{B_{\s}(0) \setminus (Z_{w})_{\t}} Dw \cdot D\varphi \right| 
&\leq& {\rm sup} \, |D\varphi| \left(\int_{B_{\s}(0)} |Dw|^{2}\right)^{1/2} \, 
\left({\mathcal H}^{n}(B_{\s}(0) \cap (Z_{w})_{\t})\right)^{1/2} + \nonumber\\ 
&& + \, {\rm sup} \, |\varphi| \left(\int_{B_{\s}(0)} |Dw|^{2}\right)^{1/2}\, 
\left(\int_{B_{\s}(0)}|D\b_{\ell}|^{2}\right)^{1/2}.
\end{eqnarray}

\noindent
Letting first $\ell \to \infty$ and then $\t \to 0$ in this, we conclude that $w$ is harmonic 
in $B_{\s}(0).$ Since $w \geq 0$ and ${\mathcal H}^{n-2} (Z_{w} \cap B_{\s}(0)) < \infty,$
it follows from the maximum principle that $Z_{w} \cap B_{\s}(0) = \emptyset.$ 
This proves the assertion in part (6).\\ 

To see part (7), first note that it follows from the identity of part (4) that 

\begin{equation}\label{w-prop1}
\frac{d}{d\r} \, \left(\r^{2-n}\int_{B_{\r}(z)} |Dw|^{2} \right) 
= 2\r^{2-n}\int_{\partial \, B_{\r}(z)} \left|\frac{\partial \, w}{\partial \, R}\right|^{2}.
\end{equation}

\noindent
(See \cite{S3}, p. 24 for the details of this claim.) Also, the identity of part (3) 
and the definition of $N_{w, z}(\r)$ directly imply that

\begin{equation}\label{w-prop2}
N_{w,z}(\r) = \frac{\r\frac{d}{d\r}\left(\r^{1-n}\int_{\partial \, B_{\r}(z)} w^{2}\right)}
{2\r^{1-n}\int_{\partial \, B_{\r}(z)} w^{2}}
\end{equation}

\noindent
whenever $N_{w,z}(\r)$ is defined. To prove (7), suppose $\int_{\partial \, B_{\r_{1}}(z_{1})} w^{2} > 0$ for 
some $z_{1} \in B_{3/2}(0)$ and $\r_{1} > 0$.
Then by continuity, there exist $\r_{0}$ with $0 < \r_{0} < \r_{1}$ such that 
$\int_{\partial \, B_{\r}(z_{1})} w^{2} > 0$ for all $\r \in (\r_{0}, \r_{1}],$ and hence $N_{w, z}(\r)$ is defined 
for all $\r \in (\r_{0}, \r_{1}].$ A computation similar to that of (\ref{derivativeoffreqfn})
using the identity (\ref{w-prop1}), the identity of part
(3) of the present lemma and the Cauchy-Schwarz inequality  then implies that 

\begin{equation}\label{w-prop3}
\frac{d}{d\r} \, N_{w, z}(\r) \geq 0
\end{equation} 
 
\noindent
for $\r \in (\r_{0}, \r_{1}].$ Thus in particular, $N_{w,z}(\r) \leq N_{2} \equiv N_{w, z}(\r_{1})$ for 
$\r \in (\r_{0}, \r_{1}].$ Using the expression (\ref{w-prop2}) in this last inequality and 
integrating the resulting differential inequality then gives

\begin{equation}\label{w-prop4}
\frac{\s^{1-n}\int_{\partial \, B_{\s}(z)} w^{2}}{\s^{2N_{2}}} \geq 
\frac{\t^{1-n}\int_{\partial \, B_{\t}(z)} w^{2}}{\t^{2N_{2}}}
\end{equation}
 
\noindent
for all $\s, \t$ with $\r_{0} < \s \leq \t \leq \r_{1}.$ Using this with $\t = \r_{1}$ and 
$\s = \s_{j}$ where $\s_{j} \downarrow \r_{0}$, we conclude that $\int_{\partial \, B_{\r_{0}}(z)}  w^{2} > 0.$
It follows that $\int_{\partial \, B_{\r}(z)} w^{2} > 0$ for all $\r \in (0, \r_{1}]$ as required.\\ 

To see parts (8) and (9), let ${\mathcal O} = \{z \in B_{3/2}(0) \, : \, 
\int_{\partial \,B_{\r}(z)} w^{2} > 0 \,\, \mbox{for each} \, \, \r \in (0, 3/2 - |z|)\}.$ Since $w^{2}$ is subharmonic 
(by (2)), it follows from
the maximum principle that if $w(z) \neq 0$ for some $z \in B_{3/2}(0)$, then   
$z \in {\mathcal O}.$ Thus if $w \not\equiv 0$, then ${\mathcal O} \neq \emptyset.$ 
We argue that ${\mathcal O}$ is open as follows. Suppose $z \in {\mathcal O}$ and 
consider $z^{\prime} \in B_{3/2}(0)$ with $|z^{\prime}  - z| < \frac{1}{4}(\frac{3}{2} - |z|).$ 
By the maximum principle and the fact that $z \in {\mathcal O}$, it follows 
that $\int_{\partial \, B_{\r}(z^{\prime})} w^{2} > 0$ for each $\r$ with $|z^{\prime} - z| < \r 
< 3/2 - |z|.$ On the other hand, it follows 
from part (7) that $\int_{\partial \, B_{\r}(z^{\prime})} w^{2} > 0$ for 
each $\r \in (0, |z^{\prime} - z|]$, giving that $z^{\prime} \in {\mathcal O}.$
Thus ${\mathcal O}$ is open. It is easy to see by the maximum principle again that 
${\mathcal O}$ is relatively closed in $B_{3/2}(0)$. Thus, we conclude that either $w \equiv 0$ in 
$B_{3/2}(0)$ or that $N_{w, z}(\r)$ is defined for all $z \in B_{3/2}(0)$ and all 
$\r \in (0, 3/2 - |z|)$ with (\ref{w-prop3}) satisfied.
\end{proof}

\medskip

\noindent
{\bf Remark}: Although we shall not need it anywhere in the present paper, 
we point out here that $w$ is weakly subharmonic in $B_{3/2}(0)$. To see this, choose a small positive constant 
$\epsilon$, and let $\gamma_{\epsilon} \, : \, {\mathbf R} \to {\mathbf R}$ 
be a smooth cut-off function with $\gamma_{\epsilon}(t) = 0$ if $t \leq \epsilon$, 
$\gamma_{\epsilon}(t) = 1$ if $t > 2\epsilon$, $\gamma_{\epsilon}(t) \geq 0$, and 
$0 \leq \gamma_{\epsilon}^{\prime}(t) \leq 2/\epsilon$ for all $t$.
Then, since $w$ is harmonic in $B_{3/2}(0) \setminus Z_{w}$, 
we have that for any smooth, 
non-negative function $\varphi$ with compact support in $B_{3/2}(0)$, 

\begin{equation} \label{wsubharm1}
\int_{B_{3/2}(0)}\varphi \gamma_{\epsilon}(w) \Delta w = 0.
\end{equation}

Integrating by parts in this we get

\begin{equation} \label{wsubharm2}
\int_{B_{3/2}(0)} \gamma_{\epsilon}(w)D\varphi \cdot Dw  
= -\int_{B_{3/2}(0)} \varphi \gamma_{\epsilon}^{\prime}(w)|Dw|^{2}. 
\end{equation}

Since the right hand side of the above is non-positive, we have that 
$\int_{B_{3/2}(0)}  \gamma_{\epsilon}(w)D\varphi \cdot Dw \leq 0.$ The assertion follows 
by letting $\epsilon \to 0$ in this.\\

\begin{lemma}\label{homogeneous}
Let $v =(v^{+}, v^{-}) \in {\mathcal F}_{\delta}$ with $v^{+}(0) = v^{-}(0) = 0.$ If $v$ is homogeneous 
of degree 1 from the origin, then $\mbox{\rm graph} \, v^{+} \cup \mbox{\rm graph} \, v^{-} = 
P_{1} \cup P_{2},$ where $P_{1}, P_{2}$ are hyperplanes of ${\mathbf R}^{n+1}$, possibly with 
$P_{1} \equiv P_{2}.$
\end{lemma}

\medskip

\begin{proof}
Since $h = \frac{1}{2}(v^{+} + v^{-})$ is harmonic,
and homogeneous of degree 1 by hypothesis, $h$ must be a linear 
function. Hence, if $v^{+} \equiv v^{-}$, the lemma holds with $P_{1} \equiv P_{2}.$
So suppose $v^{+} \not\equiv v^{-}.$ By rotating coordinates, we may and we shall assume that $h \equiv 0.$ 
Let $w = \frac{1}{2}(v^{+} - v^{-}).$ By Proposition~\ref{wproperties}, part (6)
${\mathcal H}^{n-2}(Z_{w} \cap B_{1}(0))  = \infty.$ Choose an arbitrary  
point $z \in (Z_{w} \setminus \{0\}) \cap B_{1}(0)$ and blow up 
$(v^{+}, v^{-})$ at $z$. This gives

\begin{equation}\label{homog-1}
{\widetilde v}\equiv
({\widetilde v}^{+}, {\widetilde v}^{-}) = \lim_{j \to \infty} \, {\widetilde v}_{z, \, \s_{j}} 
\end{equation}

\noindent
for some sequence of numbers $\s_{j} \searrow 0$, where ${\widetilde v}_{z, \, \s_{j}}$
is as in (\ref{blowupnotation-2}) with $y = v^{+}(z) = v^{-}(z) = 0.$ Note that since 
${\widetilde v}_{z, \, \s_{j}} \in {\mathcal F}_{\delta}$, the convergence in (\ref{homog-1}) is,
by Lemma~\ref{compactness}, in 
$W^{1,2}(B_{\s}(0))$ for each $\s \in (0, 3/2).$ Setting $\r = \s_{j}$ and 
$\th = 2\r/3$ in Lemma~\ref{coarseexcessbounds} and letting $j \to \infty$ 
it follows that $\int_{B_{\r}(0)} |{\widetilde v}|^{2} \geq 
\left(\frac{2\r}{3}\right)^{2(N_{v, z}(1) - 1)}$ for each $\r \in (0, 1]$ so that ${\widetilde v}$ is not 
identically zero in any ball $B_{\r}(0).$ Hence we have the assertions (\ref{lip3-1}) and 
(\ref{lip3-2}), by exactly the same reasoning. Thus, ${\widetilde v}$ is 
homogeneous of degree ${\mathcal N}_{v}(z)$ from the origin, and consequently by the finiteness of the left 
hand side of (\ref{lip3-1}), we immediately have that ${\mathcal N}_{v}(z) \geq 1.$  
On the other hand, by homogeneity of $v$ it follows that ${\mathcal N}_{v}(z) \leq {\mathcal N}_{v}(0) = 1,$ 
and hence we conclude that

\begin{equation}\label{homog-2}
{\mathcal N}_{v}(z) = {\mathcal N}_{v}(0)
\end{equation}
 
\noindent
for any $z \in Z_{w}$, and therefore, that 
$v$ is invariant under translations in the direction of any element of $Z_{w}.$
(See \cite{WN}, Lemma 5.17.) Since $w \not\equiv 0$ by assumption and 
${\mathcal H}^{n-2}(Z_{w} \cap B_{1}(0)) = \infty$, this means that 
$v$ is invariant under translations precisely by 
the elements of an $(n-1)$-dimensional linear subspace, and hence each of $v^{+}$ and $v^{-}$ must be 
a function of a single variable. Since by proposition~\ref{wproperties}, part (5) 
$v^{\pm}$ are harmonic in $B_{3/2} \setminus Z_{w}$, 
it follows that the union of the graphs of $v^{+}$ and $v^{-}$ must be equal to the union of 
four distinct closed, $n$-dimensional 
half spaces of ${\mathbf R}^{n+1}$ meeting along a common $(n-1)$-dimensional subspace. \\

To complete the proof, note that since $v^{+} + v^{-} \equiv 0$, it suffices to show that 
the two half spaces that make up $\mbox{graph} \, w$ make equal angles with ${\mathbf R}^{n} 
\times \{0\}.$ This follows from the identity of Proposition~\ref{wproperties}, part (4). 
Specifically, suppose without loss of generality that $Z_{w} = {\mathbf R}^{n-1} \times \{(0, 0)\}$ 
and $w(x) = \overline w (x^{1}).$ Setting $\z^{2} = \z^{3} =  \ldots \z^{n} \equiv 0$ in the identity 
of conclusion (4) of Proposition~\ref{wproperties}, 
we get the statement that $\int \left(\frac{d\overline w}{d x^{1}}\right)^{2} \frac{\partial \z^{1}}{\partial x^{1}} = 0$ 
for every $\z^{1} \in C_{c}^{1}(B_{1}(0)).$ If $\a^{+}$ and $\a^{-}$ are the angles that 
$\mbox{graph} \, \overline w$ makes with the positive and negative $x^{1}$-axes respectively, 
this identity says that $\tan^{2} \a^{-}\int_{B_{1}(0) \cap \{x^{1} < 0\}} \frac{\partial \z^{1}}{\partial x^{1}}
+ \tan^{2}\a^{+}\int_{B_{1}(0) \cap \{x^{1} > 0\}} \frac{\partial \z^{1}}{\partial x^{1}} = 0$ for 
every $\z^{1} \in C_{c}^{1}(B_{1}(0)).$ Taking a standard cut-off function for $\z^{1}$ in this 
yields $\a^{-} = \a^{+}.$ The lemma is thus proved. 
\end{proof}

\medskip

The argument of the preceding lemma shows the following.\\

\begin{lemma}\label{freqlowerbd}
Suppose $v = (v^{+}, v^{-}) \in {\mathcal F}_{\delta}$, $v^{+}(z) = v^{-}(z)$ and that 
$v$ is non constant in $B_{3/2}(0).$ Then ${\mathcal N}_{v}(z) \geq 1.$ 
\end{lemma}

\medskip
We conclude this section by mentioning the following upper semi-continuity result, which follows directly 
from the monotonicity of $N_{v, z}(\cdot).$\

\begin{lemma}\label{usc}
Suppose ${v}_{k} \in {\mathcal F}_{\d}$ for $k = 1, 2, 3, \ldots$, $z \in B_{3/2}(0),$ $v_{k} \to v$ 
in $W^{1,2}_{\rm loc} \, (B_{3/2}(0))$ and that $v_{k}$, $v$ are not identically equal to $0$ in 
$B_{3/2}(0)$ for all $k = 1, 2, 3, \ldots.$ Then ${\mathcal N}_{v}(z) \geq \limsup_{k \to \infty} \, 
{\mathcal N}_{v_{k}}(z).$ 
\end{lemma}

\section{A transverse picture}\label{transpicture}
\setcounter{equation}{0}

In this section, we analyze the situation where a hypersurface $M \in {\mathcal I}_{b}$ 
is weakly close to a multiplicity 2 hyperplane but when it is 
scaled ``vertically'' (i.e. blown up) by its height excess relative to this hyperplane, it becomes close to a 
transversely intersecting pair of hyperplanes. The geometric 
meaning of this is of course that $M$ is in fact significantly closer, in a weak sense, to a transverse 
pair of hyperplanes (with a small angle) than it is to the multiplicity 2 hyperplane; i.e. the ``fine excess''
of $M$ measured relative to a suitably chosen transverse pair of hyperplanes is significantly  smaller 
than the ``coarse excess'' of $M$ relative to the multiplicity 2 hyperplane. We obtain in this case 
(in Lemma~\ref{transverse} and its variant Lemma~\ref{transverse-1} below)
improvement of the fine excess at a fixed smaller scale.  The arguments
used to prove excess improvement here are in part variants of those developed by 
L. Simon in \cite{S}, and are in fact carried out 
in detail in \cite{WN}, although the results are not presented there in the form below. Here we state 
the lemmas in the form needed for the purposes of the present paper, and outline their proof, 
referring the reader to \cite{WN} and \cite{S} for details.\\

The lemmas have two applications; we shall need Lemma~\ref{transverse}  to handle one case of the main 
excess decay lemma (Lemma~\ref{mainlemma1}) of the paper, and we shall apply Lemma~\ref{transverse-1} 
in Sections~\ref{regblowup} to 
prove regularity of functions in ${\mathcal F}_{\delta}$ whenever their graphs are close to 
transversely intersecting pairs of hyperplanes. (See Lemmas~\ref{localdecomp} and \ref{blowupregularity}.)\\ 

\begin{lemma} \label{transverse}
Let $\theta \in (0, 1/8)$, $\delta \in (0, 1)$ and $\t \in (0, 1).$ There exists 
a number $\epsilon_{0} = \epsilon_{0}(n, \theta, \delta, \t)>0$ 
such that the following holds. Suppose $M \in {\mathcal I}_{b}$ and\\ 

\begin{itemize}
\item[(1)] $\frac{{\mathcal H}^{n}(M \cap B_{2}^{n+1}(0))}{\omega_{n}2^{n}} \leq 3 - \delta$\\
\item[(2)] ${\hat E}_{M}^{2}(3/2, L)  \equiv \left(\frac{3}{2}\right)^{-n-2}\int_{M \cap (B_{3/2}(0) \times {\mathbf R})} 
{\rm dist}^{2} \, (x, L) \leq \epsilon_{0}$ for some affine hyperplane $L$ with 
$d_{\mathcal H} \, (L \cap (B_{1}(0) \times {\mathbf R}), B_{1}(0)) \leq \e_{0}$ and\\
\item[(3)] $\int_{M \cap (B_{1}(0) \times {\mathbf R})} {\rm dist}^{2} \, (x, P) 
\leq \epsilon_{0}{\hat E}^{2}_{M}(3/2, L)$ for some pair of affine hyperplanes $P = P^{+} \cup P^{-}$ with 
$P^{+} \cap P^{-} \cap (B_{\th/4}(0) \times {\mathbf R}) \neq \emptyset$.\\
\end{itemize}

Then, either 

\begin{itemize}
\item[$(a)$] there exists an affine hyperplane ${\widetilde L}$ with $d_{\mathcal H} \, ({\widetilde L} 
\cap (B_{1}(0) \times {\mathbf R}), L \cap (B_{1}(0) \times {\mathbf R})) \leq C{\hat E}_{M}(3/2, L)$, 
$C = C(n)$ such that 
$$\left(\frac{1}{2}\right)^{-n-2}\int_{M \cap (B_{1/2}(0) \times {\mathbf R})} 
{\rm dist}^{2} \, (x, {\widetilde L}) 
\leq \t {\hat E}_{M}^{2}(3/2, L) \hspace{.2in} {\rm or}$$

\item[$(b)$] there exists a pair of affine hyperplanes ${\widetilde P} = {\widetilde P}^{+} \cup {\widetilde P}^{-}$ 
with ${\widetilde P}^{+} \cap {\widetilde P}^{-} \cap (B_{\th}(0) \times {\mathbf R}) \neq \emptyset$ such that 
\begin{itemize}
\item[(i)] $$\th^{-2}d_{\mathcal H}^{2}({\widetilde P} \cap (B_{\th}(0) \times {\mathbf R}), P \cap (B_{\th}(0) \times 
{\mathbf R})) \leq C \int_{M \cap (B_{1}(0) \times {\mathbf R})} {\rm dist}^{2} \, (x, P),$$
\item[(ii)] $$\theta^{-n-2} \int_{M \cap (B_{\theta}(0) \times {\mathbf R})} {\rm dist}^{2} \, (x, \widetilde P) \leq 
C \theta^{2} \int_{M \cap (B_{1}(0) \times {\mathbf R})} {\rm dist}^{2} \, (x, P) \hspace{.2in} 
{\rm and}$$
\item[(iii)] $M \cap ((B_{\theta}(0) \setminus S_{\widetilde P}({\th}^{2}/16)) \times {\mathbf R}) = {\rm graph} \, u^{+}
\cup {\rm graph} \, u^{-} $ where, for $\s \in (0, 1)$, 
$$S_{\widetilde P}({\s}) = \{x \in {\mathbf R}^{n} \times \{0\}\, : \, {\rm dist} \, (x, 
\pi \, ({\widetilde P}^{+} \cap {\widetilde P}^{-})) \leq \s \},$$ 
$u^{\pm}  \in C^{2}(B_{\th}(0) \setminus S_{\widetilde P}({\th}^{2}/16)) $ with $u^{+} > u^{-}$ and, for 
$x \in B_{\th}(0) \setminus S_{\widetilde P}({\th}^{2}/16)$, ${\rm dist} \, ((x, u^{+}(x)), {\widetilde P}) 
= {\rm dist} \, ((x, u^{+}(x)), {\widetilde P}^{+})$ and ${\rm dist} \, ((x, u^{-}(x)), {\widetilde P}) 
= {\rm dist} \, ((x, u^{-}(x)), {\widetilde P}^{-}).$
\end{itemize}
\noindent
Here $C = C(n)>0$ and $\pi \, : \, {\mathbf R}^{n+1} \to {\mathbf R}^{n} \times 
\{0\}$ is the orthogonal projection.\\
\end{itemize}
\end{lemma}  

\begin{proof} 
We argue by contradiction, so consider a sequence $\{M_{k}\} \subset {\mathcal I}_{b}$ satisfying   

\begin{itemize}
\item[(1)] $\frac{{\mathcal H}^{n}(M_{k} \cap B_{2}(0))}{\omega_{n}2^{n}} \leq 3 - \delta$\\
\item[(2)] ${\hat E}_{k}^{2} \equiv {\hat E}^{2}_{M_{k}}(3/2, L_{k}) = \left(\frac{3}{2}\right)^{-n-2}
\int_{M_{k} \cap (B_{3/2}(0) \times {\mathbf R})} {\rm dist}^{2} \, (x, L_{k})
\leq \frac{1}{k}$ for some affine hyperplane $L_{k}$ with $d_{\mathcal H} \, (L_{k} \cap (B_{1}(0) \times {\mathbf R}), 
B_{1}(0)) \leq \frac{1}{k}$ and\\
\item[(3)] $\int_{M_{k} \cap (B_{1}(0) \times {\mathbf R})} {\rm dist}^{2} \, (x, P_{k}) 
\leq \frac{1}{k}{\hat E}_{k}^{2}$ for some pair of affine hyperplanes $P_{k}= P_{k}^{+} \cup P_{k}^{-}$ with 
$P_{k}^{+} \cap P_{k}^{-} \cap (B_{\th/4}(0) \times {\mathbf R}) \neq \emptyset.$\\
\end{itemize}

\noindent
Write $P_{k} = P_{k}^{(1)} \cup P_{k}^{(2)}$ where $P_{k}^{(1)}$, $P_{k}^{(2)}$ are affine hyperplanes. It follows 
from $(2)$ and $(3)$ above that 

\begin{eqnarray}\label{trans-1}
&{\rm either} \hspace{.2in} {\rm dist}_{\mathcal H} \, (L_{k} \cap (B_{1}(0) \times {\mathbf R}),
P_{k}^{(1)} \cap (B_{1}(0) \times {\mathbf R})) \leq C{\hat E}_{k}\\ \nonumber 
&{\rm or} \hspace{.2in}
{\rm dist}_{\mathcal H} \, (L_{k} \cap (B_{1}(0) \times {\mathbf R}),
P_{k}^{(2)} \cap (B_{1}(0) \times {\mathbf R})) \leq C{\hat E}_{k}
\end{eqnarray}

\noindent
where $C = C(n).$ For $i=1, 2$, define 
$p_{k}^{(i)} \, : \, L_{k} \to L_{k}^{\perp}$ by $P_{k}^{(i)} = {\rm graph} \, p_{k}^{(i)} 
\equiv \{ x + p_{k}^{(i)}(x) \,  : \, x \in L_{k}\}$ (if $P_{k}^{(i)}$ is perpendicular to $L_{k}$, 
tilt $P_{k}^{(i)}$ slightly) and set 

\begin{equation}\label{trans-2}
p^{(i)} = \lim_{k \to \infty} \, ({\hat E}_{k})^{-1} p_{k}^{(i)} \circ \varphi_{k}
\end{equation}
 
\noindent
and $P^{(i)}= {\rm graph} \, p^{(i)}$, where $\varphi_{k} \, : \, {\mathbf R}^{n} \times \{0\} \to {\mathbf R}$ is such that 
${\rm graph} \, \varphi_{k} = L_{k}.$  The limit exists, possibly after passing to a subsequence.
Let $P = P^{(1)} \cup P^{(2)}.$  Note that by (\ref{trans-1}), at most 
one of $P^{(1)}$ and  $P^{(2)}$ can be perpendicular to ${\mathbf R}^{n} \times \{0\}.$\\

Now blow up the $M_{k}$'s by ${\hat E}_{k}$, 
to produce $v^{+}, v^{-} \, : \, B_{3/2}(0) \to {\mathbf R}$ as described in Section~\ref{blowupprop}. 
 Condition (3) says that 
$\mbox{graph} \, \left.v^{+}\right|_{B_{1}(0)} \, \cup \, \mbox{graph} \, \left. v^{-}\right|_{B_{1}(0)}  
\subseteq P.$\\ 

Suppose $\left. v^{+}\right|_{B_{1}(0)} \equiv \left. v^{-}\right|_{B_{1}(0)}$. 
Then $w = \frac{1}{2}(v^{+} - v^{-}) \equiv 0$ on $\overline B_{1}(0)$ and hence by part (8) of 
Lemma~\ref{wproperties}, $w \equiv 0$ on $B_{3/2}(0).$ 
It follows from this and the fact that $\frac{1}{2}(v^{+} + v^{-})$ is harmonic everywhere that 
$\mbox{graph} \, \left.v^{+}\right|_{B_{3/2}(0)} =  \mbox{graph} 
\, \left.v^{-}\right|_{B_{3/2}(0)} 
= L \, \cap \, (B_{3/2}(0) \times {\mathbf R})$ for some affine hyperplane $L$ (in fact 
$L = P^{(1)}$ or $L = P^{(2)}$), so that in this case, 
for sufficiently large $k$, 
option $(a)$ of the conclusion of the lemma holds with $M_{k}$ in place of 
$M$ and ${\widetilde L}_{k} = {\rm graph} \, (\varphi_{k} + {\hat E}_{k}\varphi)$ in place of ${\widetilde L}$
where $\varphi \, : \, {\mathbf R}^{n} \times \{0\} \to {\mathbf R}$ is such that 
$L = {\rm graph} \, \varphi.$\\

If on the other hand $\left. v^{+}\right|_{B_{1}(0)} \not\equiv \left. v^{-}\right|_{B_{1}(0)}$, then 
$P$ must be the union of distinct affine hyperplanes and
$\mbox{graph} \, \left.v^{+}\right|_{B_{1}(0)} \, \cup \, \mbox{graph} \, \left. v^{-}\right|_{B_{1}(0)} 
= P \, \cap \, (B_{1}(0) \times {\mathbf R}).$ Note that by Lemma~\ref{properties}, part (2), 
$P \cap (B_{1}(0) \times {\mathbf R}) \subset \{(x^{\prime}, x^{n+1}) \in {\mathbf R}^{n+1} \, : \, 
|x^{n+1}| \leq C\}$ where $C = C(n).$ If $\sup_{B_{1}(0)} \, |v^{+} - v^{-}| < \t/2$, we again 
have option $(a)$ of the conclusion of the lemma with $M_{k}$ in place of $M$ and 
${\widetilde L}_{k} = {\rm graph} \, (\varphi_{k} + {\hat E}_{k}\varphi)$ in place of 
${\widetilde L}$, where $\varphi \, : \, {\mathbf R}^{n} \times \{0\} \to {\mathbf R}$ 
is the affine function such that $\left. \varphi\right|_{B_{3/2}(0)} = \frac{1}{2}(v^{+} + v^{-}).$ So suppose 

\begin{equation}\label{trans-3}
\sup_{B_{1}(0)} \, |v^{+} - v^{-}| \geq \t/2.
\end{equation}

\noindent
Denote by  $\G$ the axis of $P$ (thus $\G = P^{+} \cap P^{-}$) 
and for $\s \in (0, 1)$, let $N(\s)$ be the 
tubular neighborhood of radius $\s$ around $\G.$ (Thus 
$N(\s) = \{ X \in {\mathbf R}^{n+1} \, : \, {\rm dist} \, (X, \G) \leq \s \}.$) 
We claim that for each given $\t \in (0, 1/2)$, $M_{k} \cap (B_{1}(0) \times {\mathbf R})$ must be
embedded outside $N(\t)$ for all sufficiently large $k.$ For if not, we would have a number $\t \in (0, 1/2)$ and 
a subsequence of $\{M_{k}\}$ which we continue to denote $\{M_{k}\}$ such that 
$(M_{k} \setminus N(\t)) 
\cap (B_{1}(0) \times {\mathbf R})$ contains a point 
$Z_{k} = (Z_{k}^{\prime}, Z_{k}^{n+1})$ with $\Theta_{M_{k}}(Z_{k}) \geq 2.$
The argument of the proof of Lemma~\ref{radialexcessbd} (with $\r = 1/2$) then gives that

$$\int_{B_{1/4}(z)} R^{2-n}\left(\frac{\partial}{\partial R}\left(\frac{v^{+}(x) - y}{R}\right)\right)^{2} 
+R^{2-n}\left(\frac{\partial}{\partial R}\left(\frac{v^{-}(x) - y}{R}\right)\right)^{2}  \, dx< \infty$$

\noindent
where $z = \lim_{k \to \infty} \, Z_{k}^{\prime},$ $y = \lim_{k \to \infty} \, \frac{Z_{k}^{n+1}}{\hat E_{k}}$
(both limits exist after possibly passing to a subsequence, the latter by Lemma~\ref{heightbound}), $R = |x- z|$
and $\frac{\partial}{\partial \, R}$ denotes radial differentiation. 
This implies that $v^{+}(z) = v^{-}(z) $ ($= y$), which is impossible since 
$z \in \overline {B_{1}(0)} \setminus {\pi} \, (N(\t))$ while any point ${\widetilde z}$ 
with $v^{+}({\widetilde z}) = v^{-}({\widetilde z})$ must be contained in 
${\pi} \, (\G) \cap B_{1}(0).$ Thus, if $\{\t_{k}\}$ is any
sequence of numbers with $\t_{k} \searrow 0,$ we can find a subsequence of $\{M_{k}\}$ (which we 
again denote $\{M_{k}\}$) such that $M_{k} \cap (B_{1}(0) \times {\mathbf R})$ is embedded outside
$N(\t_{k}).$\\
   
Now blow up $M_{k} \cap (B_{1}(0) \times {\mathbf R})$ by the fine excess 
$E_{k} = \sqrt{\int_{M_{k} \cap (B_{1}(0) \times {\mathbf R})} \mbox{dist}^{2} \, (x, P_{k})}$  
exactly as described in Section 6 of \cite{WN}, and outlined in the paragraph below. 
Note that although in \cite{WN} $M_{k}$ are assumed
to be free of singularities, this assumption is not necessary for the blow up argument of 
Section 6 of \cite{WN}.\\ 

Thus, let ${\widetilde q}_{k}$ be a rigid motion of ${\mathbf R}^{n+1}$ 
such that ${\widetilde q}_{k} \, (\mbox{axis of} \, \, P_{k}) 
= {\mathbf R}^{n-1} \times \{0\} \times \{0\}$, ${\widetilde q}_{k} (a_{k}) = 0$, where 
$a_{k}$ is the nearest point of the axis of $P_{k}$ to the origin of ${\mathbf R}^{n+1}$, and 
${\widetilde q}_{k} \, {\widetilde L}_{k} = {\mathbf R}^{n} \times \{0\}$, where 
${\widetilde L}_{k} = {\rm graph} \, \frac{1}{2}(p_{k}^{+} + p_{k}^{-}).$ Following the notation of  
\cite{WN}, Section 6, let ${\mathbf H}_{k}^{(1)} = {\widetilde q}_{k} \, P_{k}^{+} \cap \{x^{1} > 0\}$, 
${\mathbf H}_{k}^{(2)} = {\widetilde q}_{k} \, P_{k}^{+} \cap \{x^{1} < 0\}$, 
${\mathbf H}_{k}^{(3)} = {\widetilde q}_{k} \, P_{k}^{-} \cap \{x^{1} < 0\}$ and
${\mathbf H}_{k}^{(4)} = {\widetilde q}_{k} \, P_{k}^{-} \cap \{x^{1} > 0\}.$ (Note 
that strictly speaking, in Section 6 of \cite{WN}, the definitions of ${\mathbf H}_{k}^{(i)}$ are 
in terms of the blow-up $(v^{+}, v^{-}) \equiv (p^{+}, p^{-}),$ and the fine excess $E_{k}$ (which 
is denoted $\b_{k}$ in \cite{WN}) is defined relative to the pair of affine hyperplanes 
$P_{k}^{(0)} \equiv {\rm graph} \, {\hat E}_{k} p^{+} \cup {\rm graph} \, {\hat E}_{k} p^{-}.$  Since here we need to 
prove improvement of the excess $E_{k}$ defined relative to $P_{k}$---and not the improvement 
of excess relative to $P_{k}^{(0)}$---the above are the correct definitions of the half-spaces 
${\mathbf H}_{k}^{(i)}$ to adopt.) Now, exactly as 
in \cite{WN}, Section 6, we may express, by Allard's regularity theorem, 
${\widetilde q}_{k} \, M_{k} \cap (B_{1}^{n+1}(0) \setminus T_{k}) = \cup_{i=1}^{4} {\rm graph} \, g_{k}^{(i)}$, 
where $g_{k}^{(i)} \in C^{2}(U_{k}^{(i)}, {\mathbf H}_{k}^{(i) \, \perp})$, 
$i = 1, \ldots, 4$ satisfy the estimates as in \cite{WN}, Section 6 and $T_{k}$,  $U_{k}^{(i)}$ are as defined 
there. Defining ${\widetilde g}_{k}^{(i)}$, $i = 1, \ldots, 4$ as in \cite{WN}, Section 6, 
we obtain,  as in \cite{WN}, Section 6, functions (the blow-up) 
$w^{(1)}, w^{(4)} \in C^{2}({\mathbf R}^{n \, +} \cap B_{1}(0))$, 
and $w^{(2)}, w^{(3)} \in C^{2}({\mathbf R}^{n \, -} \cap B_{1}(0))$, where 
${\mathbf R}^{n \, +} \equiv \{ x \in {\mathbf R}^{n} \times \{0\} \, : \, x^{1} > 0\}$ and 
${\mathbf R}^{n \, -} \equiv \{ x \in {\mathbf R}^{n} \times \{0\} \, : \, x^{1} < 0\}$, such that 
$E_{k}^{-1} {\widetilde g}_{k}^{(i)} \to w^{(i)}$ for $i=1, \ldots, 4$, where for each $i$, the convergence is in 
the $C^{2}$-norm on each compact subset of the domain of $w^{(i)}$ and also in 
the $L^{2}$-norm on the domain of $w^{(i)}$.        
By Lemma 6.23 of \cite{WN}, the blow-up $\{w^{(i)}\}_{i=1}^{4}$ (restricted to a suitably 
smaller ball, say $B_{1/2}(0)$) consists of two 
harmonic functions $w^{(13)}$ and $w^{(24)}$  in the sense that the union of the closures 
of the graphs of $w^{(1)}$, $w^{(3)}$ in $B_{1/2}(0) \times {\mathbf R}$ 
is the graph of a harmonic function $w^{(13)}$ over $B_{1/2}(0)$ and 
similarly the union of the closures of the graphs of $w^{(2)}$, $w^{(4)}$ in 
$B_{1/2}(0) \times {\mathbf R}$ is the graph of 
a harmonic function $w^{(24)}$ over $B_{1/2}(0).$ For $x \in {\mathbf R}^{n} \times \{0\}$, 
let $l^{(13)}(x) = w^{(13)}(0) + Dw^{(13)}(0) \cdot x$, 
$l^{(24)}(x) = w^{(24)}(0) + Dw^{(24)} \cdot x$ and let the affine functions 
$h_{k}^{(13)}, h_{k}^{(24)} \, : \, {\mathbf R}^{n} \times \{0\} \to {\mathbf R}$ 
be defined by ${\rm closure} \, {\mathbf H}_{k}^{(1)} \cup {\rm closure} \, {\mathbf H}_{k}^{(3)} = 
{\rm graph} \, h_{k}^{(13)}$ and  ${\rm closure} \, {\mathbf H}_{k}^{(2)} \cup {\rm closure} \,  {\mathbf H}_{k}^{(2)} = 
{\rm graph} \, h_{k}^{(24)}.$ Set ${\widetilde P}_{k} = {\widetilde q}_{k}^{-1} \, ({\rm graph} \, (h_{k}^{(13)} + E_{k}l^{(13)}) 
\cup {\rm graph} \, (h_{k}^{(24)} + E_{k}l^{(24)})).$ Then, using standard estimates for harmonic functions, 
and the ``non-concentration of excess'' estimate of part (ii) of Lemma 6.22, \cite{WN}, we conclude that

\begin{equation}\label{transrefined-1}
\th^{-n-2}\int_{M_{k} \cap (B_{\th}(0) \times {\mathbf R})} \mbox{dist}^{2} \, (X, 
{\widetilde P}_{k}) \leq C\th^{2}E_{k}^{2}
\end{equation}

\noindent
for sufficiently large $k$, where $C = C(n).$ If we write, using our usual notation, 
${\widetilde P}_{k} = {\widetilde P}_{k}^{+} \cup {\widetilde P}_{k}^{-}$, 
then, since $\sup_{B_{1}(0)} \, |h_{k}^{(13)} - h_{k}^{(24)}| \geq \frac{1}{4}\t{\hat E}_{k}$
(by (\ref{trans-2}) and (\ref{trans-3})) and $E_{k}/{\hat E}_{k} \to 0$, we must have that 
${\widetilde P}_{k}^{+} \cap {\widetilde P}_{k}^{-} \cap (B_{\th}(0) \times {\mathbf R}) 
\neq \emptyset$ for all sufficiently large $k$.\\ 
 
 Finally, note that conclusion $(b)(i)$ of the lemma with $M_{k}$, ${\widetilde P}_{k}$, $P_{k}$ in 
place of $M$, ${\widetilde P}$, $P$ follows directly from the definition of ${\widetilde P}_{k}$ 
and conclusion $(b)(iii)$ with $M_{k}$, ${\widetilde P}_{k}$ in place of $M$, ${\widetilde P}$ and 
appropriate functions $u_{k}^{\pm} \in C^{2}(B_{\th}(0) \setminus S_{{\widetilde P}_{k}}(\th^{2}/16))$ 
in place of $u^{\pm}$ follows from Allard's regularity theorem and the fact that 
$E_{k}/{\hat E}_{k} \to 0.$ 
\end{proof}

\medskip

In addition to the hypotheses of Lemma~\ref{transverse}, if we also assume that 
$0 \in {\overline M},$ $\Theta_{M}(0) \geq 2$ and that $P = P^{+} \cup P^{-}$ is a pair of hyperplanes
(so that $0 \in P^{+} \cap P^{-}$),
then the conclusions of the lemma hold with ${\widetilde P} = {\widetilde P}^{+} \cup 
{\widetilde P}^{-}$ equal to a pair of hyperplanes 
(so that $0 \in {\widetilde P}^{+} \cap {\widetilde P}^{-}$). 
This follows from the fact that under these additional hypotheses, 
we have for the fine blow-up the estimate 

\begin{equation}\label{ddr}
\int_{B_{1/2}(0)}R^{2-n}\left(\frac{\partial (w^{(13)}/R)}{\partial R}\right)^{2} + 
R^{2-n}\left(\frac{\partial (w^{(24)}/R)}{\partial R}\right)^{2} < C < \infty
\end{equation}

\noindent
where $C = C(n)$, $R = |x|$ and $\frac{\partial}{\partial \, R}$ denotes 
the radial derivative, and $w^{(13)}$, $w^{(24)}$ are as in the proof of Lemma~\ref{transverse} above. 
This estimate says in particular that $w^{(13)}(0) = w^{(24)}(0) = 0$. Since we have, by hypothesis, 
that $0 \in P_{k}^{+} \cap P_{k}^{-}$ for each $k$, we immediately conclude 
that $0 \in {\widetilde P}_{k}^{+} \cap {\widetilde P}_{k}^{-}.$ 
(Notation as in the proof of Lemma~\ref{transverse} above.) 
The estimate (\ref{ddr}) was first proved in \cite{S} (see \cite{S}, Lemma 3.4 and 
\cite{S}, Section 5.1, inequality (12)) and in view of Lemmas 6.21 and 6.22 of \cite{WN}, the same 
proof as in \cite{S} yields it here as well.\\ 

Thus we have the following variant of Lemma~\ref{transverse}.\\

\begin{lemma} \label{transverse-1}
Let $\theta \in (0, 1/8)$, $\delta \in (0, 1)$ and $\t \in (0, 1).$ There exists 
$\epsilon_{0} = \epsilon_{0}(n, \theta, \delta, \t)>0$ 
such that the following holds. Suppose $M \in {\mathcal I}_{b},$ $0 \in {\overline M}$ and\\ 

\begin{itemize}
\item[(1)] $\Theta_{M}(0) \geq 2$\\ 
\item[(2)] $\frac{{\mathcal H}^{n}(M \cap B_{2}^{n+1}(0))}{\omega_{n}2^{n}} \leq 3 - \delta$\\
\item[(3)] ${\hat E}_{M}^{2}(3/2, L) \equiv \left(\frac{3}{2}\right)^{-n-2}\int_{M \cap (B_{3/2}(0) \times {\mathbf R})} 
{\rm dist}^{2} \, (x, L) \leq \epsilon_{0}$ for some affine hyperplane $L$ with 
$d_{\mathcal H} \, (L \cap (B_{1}(0) \times {\mathbf R}), B_{1}(0)) \leq \e_{0}$ and\\
\item[(4)] $\int_{M \cap (B_{1}(0) \times {\mathbf R})} {\rm dist}^{2} \, (x, P) 
\leq \epsilon_{0}{\hat E}^{2}_{M}(3/2, L)$ for some pair of hyperplanes $P = P^{+} \cup P^{-}$.\\
\end{itemize}

Then, either 

\begin{itemize}
\item[$(a)$] there exists an affine hyperplane ${\widetilde L}$ with 
$d_{\mathcal H} \, ({\widetilde L} 
\cap (B_{1}(0) \times {\mathbf R}), L \cap (B_{1}(0) \times {\mathbf R})) \leq C{\hat E}_{M}(3/2, L)$, $C = C(n),$ such that
$$\left(\frac{1}{2}\right)^{-n-2}\int_{M \cap (B_{1/2}(0) \times {\mathbf R})} 
{\rm dist}^{2} \, (x, {\widetilde L}) 
\leq \t {\hat E}_{M}^{2}(3/2, L) \hspace{.2in} {\rm or}$$

\item[$(b)$] there exists a pair of hyperplanes ${\widetilde P} = {\widetilde P}^{+} \cup {\widetilde P}^{-}$ 
such that 
\begin{itemize}
\item[(i)] $$d_{\mathcal H}^{2}({\widetilde P} \cap (B_{1}(0) \times {\mathbf R}), P \cap (B_{1}(0) \times 
{\mathbf R})) \leq C \int_{M \cap B_{1}(0) \times {\mathbf R}} {\rm dist}^{2} \, (x, P),$$
\item[(ii)] $$\theta^{-n-2} \int_{M \cap (B_{\theta}(0) \times {\mathbf R})} {\rm dist}^{2} \, (x, \widetilde P) \leq 
C \theta^{2} \int_{M \cap (B_{1}(0) \times {\mathbf R})} {\rm dist}^{2} \, (x, P) \hspace{.2in} 
{\rm and}$$
\item[(iii)] $M \cap ((B_{\theta}(0) \setminus S_{\widetilde P} ({\th}^{2}/16)) \times {\mathbf R}) = {\rm graph} \, u^{+}
\cup {\rm graph} \, u^{-} $ where, for $\s \in (0, 1)$, 
$$S_{\widetilde P}({\s}) = \{x \in {\mathbf R}^{n} \times \{0\}\, : \, {\rm dist} \, (x, 
{\pi} \, ({\widetilde P}^{+} \cap {\widetilde P}^{-})) \leq \s \},$$ 
$u^{\pm}  \in C^{2}(B_{\th}(0) \setminus S_{\widetilde P}({\th}^{2}/16)) $ with $u^{+} > u^{-}$ and, for 
$x \in B_{\th}(0) \setminus S_{\widetilde  P}({\th}^{2}/16)$, ${\rm dist} \, ((x, u^{+}(x)), {\widetilde P}) 
= {\rm dist} \, ((x, u^{+}(x)), {\widetilde P}^{+})$ and ${\rm dist} \, ((x, u^{-}(x)), {\widetilde P}) 
= {\rm dist} \, ((x, u^{-}(x)), {\widetilde P}^{-}).$
\end{itemize}
\noindent
Here $C = C(n)>0$ and ${\pi} \, : \, {\mathbf R}^{n+1} \to {\mathbf R}^{n} \times 
\{0\}$ is the orthogonal projection.\\
\end{itemize}
\end{lemma}

\section{Regularity of blow-ups off affine hyperplanes}\label{regblowup}
\setcounter{equation}{0}
In order to handle one case of the proof of the main excess decay lemma (Lemma~\ref{mainlemma1})---namely, 
the case in which the ``fine excess''  of a hypersurface $M \in {\mathcal I}_{b}$ 
(i.e. the height excess of $M$ measured relative to a pair of affine hyperplanes) 
is of the same order as the ``coarse excess''  of $M$ (i.e. the excess of $M$ 
relative to a single affine hyperplane), which 
geometrically corresponds to the situation where $M$ has ``lots'' of self intersections distributed 
more or less evenly---it is necessary to understand, in sufficient detail, the asymptotic behavior of 
the 2-valued functions 
belonging to the class ${\mathcal F}_{\delta}.$  Our goal in this section is to do that. At the end of this section, 
we prove the following regularity theorem for any $v \in {\mathcal F}_{\delta}$:\\

\begin{theorem}\label{regblowupthm}
Let $v = (v^{+}, v^{-}) \in {\mathcal F}_{\delta}.$ There exists a relatively closed (possibly empty) 
subset $S_{v}$ of $B_{3/2}(0)$ (the branch set of $v$) such that  
\begin{itemize}
\item[$(a)$] if $\Omega \subset B_{3/2}(0) \setminus S_{v}$ is open and simply connected, then there exist two 
harmonic functions $v^{1}, v^{2} \, : \, \Omega \to {\mathbf R}$ such that
$$(\mbox{\rm graph} \, v^{+} \cup \mbox{\rm graph} \, v^{-}) \cap (\Omega \times {\mathbf R})  = 
\mbox{\rm graph} \, v^{1} \cup \mbox{\rm graph} \, v^{2}$$ and

\item[$(b)$] for each $z \in S_{v} \cap B_{1}(0)$, there 
exists an affine function $l_{z}\, : \, {\mathbf R}^{n} \times \{0\} \to {\mathbf R}$
such that 

$$\r^{-n-2} \int_{B_{\r}(z)} (v^{+} - l_{z})^{2} + (v^{-} - l_{z})^{2} \leq C\r^{\lambda}\int_{B_{5/4}(0)} 
(v^{+})^{2} + (v^{-})^{2}$$

\noindent
for all $\r \in (0, 1/64)$, where $C$, $\lambda$ are positive constants depending only on $n$ and $\d$. In fact, 
$l_{z}(x) = h(z) + Dh(z) \cdot (x - z)$ where $h = \frac{1}{2}(v^{+} + v^{-}).$ (Recall that 
$h$ is harmonic in $B_{3/2}(0).$)\\
\end{itemize}
\end{theorem}

We begin with a series of lemmas.\\

\begin{lemma}\label{monotonicity}
If $(v^{+}, v^{-}) \in {\mathcal F}_{\delta}$ and $v^{+}(z) = v^{-}(z) = y$, then 
$\E_{z, \r} \equiv \rho^{-n-2}\int_{B_{\rho}(z)} (v^{+} - y)^{2} 
+ (v^{-} - y)^{2}$ is monotonically increasing as a function of $\rho.$ Therefore,
$\rho^{-n-2}\int_{B_{\rho}(z)} (v^{+} - y)^{2} + (v^{-} - y)^{2} \leq C\int_{B_{1}(0)} (v^{+})^{2} + 
(v^{-})^{2} \leq C$ for all $z \in B_{1/4}(0) \cap Z_{w}$ and all $\rho \in (0,1/4)$ where 
$C = C(n).$  
\end{lemma}

\medskip

\begin{proof}
The first assertion follows directly from Lemma~\ref{freqlowerbd} and the estimate (\ref{lip4}). 
The second assertion follows from the first and the estimate $|y|^{2} \leq C\int_{B_{1}(0)} h^{2}$ 
which holds since $y = h(z)$ and
$h = \frac{1}{2}(v^{+} + v^{-})$ is harmonic in $B_{3/2}(0).$ 
\end{proof}
  
\medskip

\begin{lemma} \label{Nbound}
Let $\alpha_{0} \in (0, \pi/2)$, $\delta_{0} \in (0, 1)$. There exists 
$\epsilon_{1} = \epsilon_{1}(n, \alpha_{0}, \delta_{0}) 
\in (0, 1)$ such that 
if $P_{0} = P_{0}^{+} \cup P_{0}^{-}$ is a pair of hyperplanes with $\a_{0} \leq \angle \, P_{0} < \pi,$ 
$(v^{+}, v^{-}) \in {\mathcal F}_{\delta},$ $v^{+}(0) = v^{-}(0)$ and 
$$\int_{B_{1}(0)} (v^{+} - p_{0}^{+})^{2} + (v^{-} - p_{0}^{-})^{2} \leq \epsilon_{1} $$

\noindent
then $1 \leq N_{v, 0}(1) \leq 1 + \delta_{0}.$
\end{lemma}

\medskip

\begin{proof}
Since $v^{+}(0) = v^{-}(0)$, the lower bound $ N_{v, 0}(1) \geq 1$ follows from 
the monotonicity of $N_{v,0}$ and Lemma~\ref{freqlowerbd}.
If the upper bound fails to hold for some $\d \in (0, 1)$, there exists a sequence
$v_{k} = (v_{k}^{+}, v_{k}^{-}) \in {\mathcal F}_{\delta},$ $k = 1, 2, \ldots$, with $v^{+}_{k}(0) = v^{-}_{k}(0)$, 
and  a sequence of pairs of hyperplanes 
$P_{k} = P^{+}_{k} \cup P^{-}_{k}$ with $\a_{0} \leq \angle \, P_{k} < \pi$ satisfying

\begin{equation} \label{L2conv}
\int_{B_{1}(0)} (v_{k}^{+} - p_{k}^{+})^{2} + (v_{k}^{-} - p_{k}^{-})^{2} \leq  \frac{1}{k},
\end{equation}

\noindent
and yet $N_{v_{k}, 0}(1) > 1 + \d$ for all $k$. In view of Proposition~\ref{properties}, 
part (2), the inequality (\ref{L2conv}) implies that 
$\int_{B_{1}(0)} (p_{k}^{+})^{2} + (p_{k}^{-})^{2} \leq C$ for each $k$ where $C = C(n)$. Passing to a subsequence, 
$P_{k} \to P$ for some pair of hyperplanes $P = P^{+} \cup P^{-}$ with $\a_{0} \leq \angle \, P <\pi.$ 
By (\ref{L2conv}) again, the sequence $\{v_{k}\}$ 
converges to $p \equiv (p^{+}, p^{-})$ in $L^{2}(B_{1}(0))$,  and by Lemma~\ref{compactness}, the convergence is 
in $W^{1,2}(B_{1}(0)).$ Thus $N_{v_{k}, 0}(1) \to N_{p, 0}(1).$ 
But since $p$ is homogeneous of degree 1, $N_{p, 0}(1) = 1.$ This proves the lemma.
\end{proof} 

\medskip

\begin{lemma} \label{localdecomp}
Let $\overline\th \in (0, 1/8),$ $\d \in (0, 1)$ and $\overline\alpha \in (0, \pi).$ There exists a number  
$\overline\epsilon = \overline\epsilon(n, \overline\th, \d, \overline\alpha) \in (0, 1)$ 
such that the following holds. 
If ${\overline P} = {\overline P}^{+} \cup {\overline P}^{-}$ is a  
pair of hyperplanes of ${\mathbf R}^{n+1}$ with $\pi > \angle \, \overline P \geq \overline\alpha,$ 
$v = (v^{+}, v^{-}) \in {\mathcal F}_{\delta},$  $v^{+}(0) = v^{-}(0) = 0$, 

$$\int_{B_{1}(0)} {\rm dist}^{2}\, ((x, v^{+}(x)), \overline P) + {\rm dist}^{2}\, ((x, v^{-}(x)), \overline P) 
\leq \overline\epsilon \hspace{.2in} \mbox{and if}$$ 

$$\int_{B_{1}(0) \setminus S_{\overline P}(\overline\th/16)} (v^{+} - {\overline p}^{+})^{2} + (v^{-} - {\overline p}^{-})^{2} 
\leq \overline\e,$$  
\noindent
where $S_{\overline P}(\s) = \left\{ x \in {\mathbf R}^{n} \times \{0\} \, : \, {\rm dist} \, (x, {\pi} \, ({\overline P}^{+} 
\cap {\overline P}^{-})) \leq \s\right\}$, 
then there exists a pair of hyperplanes ${\widetilde P} = {\widetilde P}^{+} \cup  {\widetilde P}^{-}$ 
of ${\mathbf R}^{n+1}$ such that 
\begin{eqnarray*}
(\overline{\th})^{-n -2}\int_{B_{\overline{\th}}(0)} & {\rm dist}^{2} \, ((x,v^{+}(x)), {\widetilde P})  + 
{\rm dist}^{2} \, ((x, v^{-}(x)), {\widetilde P}) \nonumber\\
&\leq {\overline C}{\overline\th}^{2} 
\int_{B_{1}(0)} {\rm dist}^{2} \, ((x, v^{+}(x)), \overline P) + {\rm dist}^{2} \, ((x, v^{-}(x)), \overline P)
\end{eqnarray*}

\noindent
and

\begin{eqnarray*}
(\overline \th)^{-n-2}\int_{B_{\overline\th}(0) \setminus S_{\widetilde P}({\overline\th^{2}/16})} 
&(v^{+} - {\widetilde p}^{+})^{2} + (v^{-} - {\widetilde p}^{-})^{2} \nonumber\\
&\leq {\overline C}{\overline\th}^{2} \int_{B_{1}(0)} {\rm dist}^{2} \, ((x, v^{+}(x)), \overline P) 
+ {\rm dist}^{2} \, ((x, v^{-}(x)), \overline P)
\end{eqnarray*}

\noindent
Here ${\overline C} = C(n) \in (0, \infty).$\\
\end{lemma}

\begin{proof} 
By the definition of ${\mathcal F}_{\delta}$, there exists a sequence $M_{k}$ of hypersurfaces in ${\mathcal I}_{b}$ 
with $\frac{{\mathcal H}^{n} \, (M_{k} \cap B_{2}^{n+1}(0))}{\omega_{n}2^{n}} \leq 3 -\d$ 
and a sequence $L_{k}$ of affine hyperplanes converging to ${\mathbf R}^{n} \times \{0\}$ 
such that ${\hat E}_{M_{k}}(3/2, L_{k}) \to 0$ and the blow-up of $\{M_{k}\}$
off $\{L_{k}\}$ (as described in Section~\ref{blowupprop}) is $(v^{+}, v^{-})$.
Since $\int_{B_{1}(0) \setminus S_{\overline P}({\overline \th}/16)} (v^{+} - {\overline p}^{+})^{2} + 
(v^{-} - {\overline p}^{-})^{2} \leq {\overline \e}$, 
if ${\overline \e} = {\overline \e}(n, {\overline \a})$ is sufficiently 
small, it follows from Lemma~\ref{wproperties} part (8) that $v^{+} \not\equiv v^{-}$ in any 
ball $B_{\s}(0)$, $0< \s \leq 1.$  Thus, since $v^{+}(0) = v^{-}(0)$,  we have by the remark following
Lemma~\ref{radialexcessbd} that possibly after taking a subsequence of $\{k\}$ which we continue to 
denote $\{k\}$,  there exists $Z_{k} = (Z_{k}^{\prime}, Z_{k}^{n+1}) 
\in M_{k} \cap (B_{1}(0) \times {\mathbf R})$  such that $\Theta_{M_{k}}(Z_{k}) \geq 2$ and 
$\left(Z_{k}^{\prime},  \frac{Z_{k}^{n+1}}{{\hat E}_{k}}\right) \to (0,0).$ Let ${\widetilde M}_{k} \equiv 
\eta_{Z_{k}, \, 1 - |Z_{k}|} \, M_{k}.$ By the monotonicity of mass ratio, for sufficiently large $k$,  
$\frac{{\mathcal H}^{n} \, ({\widetilde M}_{k} \cap B_{2}^{n+1}(0))}{\omega_{n}2^{n}} \leq 3 - \d/2$ 
and the blow-up, 
as in Section~\ref{blowupprop}, of the sequence of hypersurfaces ${\widetilde M}_{k}$ off the sequence 
$(1 - |Z_{k}|)^{-1}(L_{k} - Z_{k})$ of affine hyperlanes is also $(v^{+}, v^{-}).$
Thus, by replacing the original sequence $M_{k}$ with ${\widetilde M}_{k}$,  we may assume that 
$0 \in M_{k}$ and $\Th_{M_{k}}(0) \geq 2$ for all $k$ so that the hypotheses $(1)$ and $(2)$ of 
Lemma~\ref{transverse-1} are satisfied with $M_{k}$ in place of $M$ and $\d/2$ in place of $\d.$\\

By hypothesis we have

\begin{equation*}\label{d-0} 
\int_{B_{1}(0)} {\rm dist}^{2} \, ((x, v^{+}(x)), \overline P) + {\rm dist}^{2} \, ((x, v^{-}(x)), \overline P) \leq 
\overline\epsilon
\end{equation*}

\noindent
which together with the squared triangle inequality 
${\rm dist}^{2} \, (X, {\overline P}) \leq 2{\rm dist}^{2} \, (Y, {\overline P}) + 
2|X - Y|^{2}$ implies, for sufficiently large $k$, that

\begin{equation}\label{d}
\int_{B_{1}(0)} {\rm dist}^{2} \, \left(\left(x, \frac{\overline{\psi}_{k}(x)u_{k}^{+}(x)}{{\hat E}_{k}}\right), \overline P\right)
+ {\rm dist}^{2} \, \left(\left(x, \frac{\overline{\psi}_{k}(x)u_{k}^{-}(x)}{{\hat E}_{k}}\right), \overline P\right) 
\leq 4\overline\epsilon
\end{equation}

\noindent
where the notation is as in (\ref{supconvergence}) and (\ref{blowups}). Let
${\overline P}_{k} = {\rm graph} \, {\hat E}_{k} {\overline p}^{+} \cup {\rm graph} \, {\hat E}_{k}{\overline p}^{-}.$ Then

\begin{eqnarray}\label{e}
\int_{M_{k} \cap (B_{1}(0) \times {\mathbf R})} \mbox{dist}^{2} \, (X, {\overline P}_{k}) &=& 
\int_{G_{k}^{+} \cap (B_{1}(0) \times {\mathbf R})} \mbox{dist}^{2} \, (X, {\overline P}_{k}) \nonumber\\
& &+ \int_{G_{k}^{-} \cap (B_{1}(0) \times {\mathbf R}) } \mbox{dist}^{2} \, (X, {\overline P}_{k}) \nonumber\\
&& + \int_{(M_{k} \setminus G_{k}) \cap (B_{1}(0) \times {\mathbf R})} 
\mbox{dist}^{2} \, (X, {\overline P}_{k}) \nonumber\\
&\leq& C\int_{B_{1}(0)} {\rm dist}^{2} \,  ((x, \overline{\psi}_{k}(x)u_{k}^{+}(x)) , {\overline P}_{k}) 
+ {\rm dist}^{2} \, ((x, \overline{\psi}_{k}(x)u_{k}^{-}(x)), {\overline P}_{k}) \nonumber\\
&&+ C{\hat E}_{k}^{2 + \mu}
\end{eqnarray}

\noindent 
where $C = C(n).$ The inequality in the above follows from the estimates (\ref{badsetest}), (\ref{cutoffmeasure})
together with the general fact that if $L = {\rm graph} \, \ell$ is a hyperplane of ${\mathbf R}^{n+1}$, 
where $\ell \, : \, {\mathbf R}^{n} \times \{0\} \to {\mathbf R}$ is 
given by $\ell(x^{\prime}) = a\cdot x^{\prime}$ for some $a \in {\mathbf R}^{n} \times \{0\}$, then 
for any point $(x^{\prime}, x^{n+1}) \in {\mathbf R}^{n+1}$ and any number $\lambda > 0$, 

\begin{equation}\label{reg5-0-0}
{\rm dist}^{2} \, ((x^{\prime}, \lambda x^{n+1}), L^{\lambda}) = \frac{\lambda^{2}(1 + |a|^{2})}{1 + \lambda^{2}|a|^{2}}
{\rm dist}^{2} \, ((x^{\prime}, x^{n+1}), L),
\end{equation}

\noindent
where $L^{\lambda} = {\rm graph} \, \lambda \ell.$ By the inequalities (\ref{d}) and (\ref{e}), we have that for all 
sufficiently large $k$, 

\begin{equation}\label{g}
\int_{M_{k} \cap (B_{1}(0) \times {\mathbf R})} \mbox{dist}^{2} \, (X, {\overline P}_{k}) \leq 
C{\overline \e}{\hat E}_{k}^{2}.
\end{equation}

Now note that there exists a constant $C_{1} \in (0, 1)$ depending only on $\overline\a$ such that 

\begin{equation}\label{h-1}
d_{\mathcal H} \, ({\overline P} \cap (B_{1/2}(0) \times {\mathbf R}), L \cap (B_{1/2}(0) \times {\mathbf R})) \geq C_{1}
\end{equation}
 
\noindent
for any affine hyperplane $L$. In view of (\ref{g}), 
given any $\t \in (0, 1)$, if ${\overline \e} = \overline \e (n, \overline\th, \d, \t)$ is sufficiently 
small, we may apply Lemma~\ref{transverse-1} with $\th = {\overline \th}$, 
${\overline P}$ in place of $P^{0}$, $\d/2$ in place of $\d$ and
$M_{k}$ in place of $M.$ Lemma~\ref{transverse-1} then gives for each $k$ either a pair of hyperplanes 
$\widetilde P_{k} = {\widetilde P}_{k}^{+} \cup {\widetilde P}_{k}^{-}$ with 

\begin{equation}\label{h-2}
d^{2}_{\mathcal H} \, ({\widetilde P}_{k} \cap (B_{1}(0) \times {\mathbf R}), {\overline P}_{k} \cap 
(B_{1}(0) \times {\mathbf R}) \leq C\int_{M_{k} \cap (B_{1}(0) \times {\mathbf R})} 
{\rm dist}^{2} \, (x, {\overline P}_{k}) 
\end{equation}

\noindent
such that

\begin{equation}\label{i}
\overline\theta^{-n-2} \int_{M_{k} \cap (B_{\overline\theta}(0) \times {\mathbf R})} 
{\rm dist}^{2} \, (x, \widetilde P_{k}) \leq C \overline\theta^{2} 
\int_{M_{k} \cap (B_{1}(0) \times {\mathbf R})} {\rm dist}^{2} \, (x, {\overline P}_{k})
\end{equation}
 
\noindent
where $C = C(n),$ or an affine hyperplane ${\widetilde L}_{k}$ with
$d_{\mathcal H} \, ({\widetilde L}_{k} \cap (B_{1}(0) \times {\mathbf R}), L_{k} \cap (B_{1}(0) \times 
{\mathbf R})) \leq C{\hat E}_{k}$, $C = C(n)$, satisfying

\begin{equation}\label{i-1}
\int_{M_{k} \cap (B_{1/2}(0) \times {\mathbf R})} {\rm dist}^{2} \, (x, {\widetilde L}_{k}) 
\leq \t {\hat E}_{k}^{2}. 
\end{equation}

\noindent
However, if (\ref{i-1}) holds for infinitely many $k$, 
we see by dividing (\ref{i-1}) by ${\hat E}_{k}^{2}$ and passing to the limit 
as $k \to \infty$ that $\int_{B_{1}(0)} (v^{+} - \ell)^{2} + (v^{-} - \ell)^{2} 
\leq \t$ for some affine function $\ell$, which, in view of (\ref{h-1}), contradicts the 
hypothesis  $\int_{B_{1}(0) \setminus S_{\overline P}({\overline \th}/16)} (v^{+} - {\overline p}^{+})^{2} 
+ (v^{-} - {\overline p}^{-})^{2} \leq {\overline \e}$ provided 
$\t = \t(n, C_{1}) \in (0, 1)$ (hence $\t = \t(n, {\overline\a})$) is chosen sufficiently small. (Here 
$C_{1}$ is as in (\ref{h-1}.)) Thus if 
${\overline \e} = {\overline \e} (n, \overline \th, \d, \overline \a)$ is 
chosen sufficiently small, option (\ref{i-1}) cannot occur for infinitely many $k$, and hence we must 
have (\ref{i}) for all sufficiently large $k$.  It follows, upon dividing the inequality (\ref{i}) 
by ${\hat E}_{k}^{2}$ and letting $k \to 
\infty$ after possibly passing to a subsequence, (and using the estimates (\ref{e}), (\ref{h-2}) and 
${\mathcal H}^{n}((M_{k} \setminus G_{k}) \cap (B_{1}(0) \times {\mathbf R})) 
\leq C\, {\hat E}_{k}^{2 + \mu}$) that for 
some pair of hyperplanes ${\widetilde P}$, 

\begin{eqnarray}
\overline\th^{-n-2}\int_{B_{\overline\th}(0)} & \mbox{dist}^{2} \, ((x, v^{+}(x)) , {\widetilde P}) + 
\mbox{dist}^{2} \,((x, v^{-}(x)) , {\widetilde P}) \nonumber\\
& \leq C\overline\th^{2} \int_{B_{1}(0)} {\rm dist}^{2} \, ((x, v^{+}(x)) ,{\overline P}) 
+ {\rm dist}^{2} \, ((x, v^{-}(x)), {\overline P})
\end{eqnarray}

\noindent
where $C = C(n).$ 
The remaining claim follows directly from conclusion $(b)(iii)$ of Lemma~\ref{transverse-1}. The lemma is thus proved.
\end{proof}

\medskip

The next lemma says that if the graph of $w = \frac{1}{2}(v^{+} - v^{-})$ stays close, in $B_{1}(0) \times {\mathbf R},$  
to a pair of $n$-dimensional half-spaces of ${\mathbf R}^{n+1}$ 
meeting at an angle $< \pi$ along an $(n-1)$-dimensional axis, and if $Z_{w}$ is the zero set of $w$, then 
$Z_{w} \cap B_{1/2}(0)$ cannot have too large a gap.\\

\begin{lemma}\label{gapcontrol}
Let $(v^{+}, v^{-}) \in {\mathcal F}_{\delta}$, $w = \frac{1}{2}(v^{+} - v^{-})$ and 
$\gamma \in (0, 1/2).$ Suppose that 
$\int_{B_{1}(0)} (w - L)^{2} \leq \g$ where $L \, : \, {\mathbf R}^{n} \times \{0\}\to {\mathbf R}^{+} \cup \{0\}$ 
is such that $\mbox{\rm graph} \, L$ is equal to the union of two $n$-dimensional half-spaces of 
${\mathbf R}^{n+1}$ meeting along ${\mathbf R}^{n-1} \times \{(0,0)\},$ each making the same 
angle $\b \in (0, \pi/2)$ with ${\mathbf R}^{n} \times \{0\}.$ If $B_{r}(q) \cap Z_{w} = \emptyset$ for 
some $q \in ({\mathbf R}^{n-1} \times \{(0,0)\}) \cap B_{1/2}^{n+1}(0)$ and $r > 0$, 
then $r \leq C\g^{1/2n}$ where $C$ depends only on $n$ and $\b$.\\ 
\end{lemma}

\begin{proof}
Let $Q = \{x \in B_{1}(0) \, : \, |w(x) - L(x)| \geq \g^{1/4} \}.$ Since $\int_{B_{1}(0)} (w - L)^{2} \leq \g$, 
it follows that 

\begin{equation}\label{gapcontrol-0}
{\mathcal L}^{n}(Q) \leq \g^{1/2}.
\end{equation}

Suppose $B_{r}(q) \cap Z_{w} = \emptyset$
for some $q \in ({\mathbf R}^{n-1} \times \{(0,0)\}) \cap B_{1/2}(0)$  and $r > 0.$ Then by  
Proposition~\ref{wproperties}, part (6), $w$ is harmonic 
(and positive) in $B_{r}(q)$, so that by the Harnack inequality we have that 

\begin{equation}\label{gapcontrol-1}
\mbox{sup}_{B_{r/2}(q)} \, w \leq 3^{n} \, \mbox{inf}_{B_{r/2}(q)} \, w. 
\end{equation}

With $\mu = \mu(n)  \in (0, 1/2)$ to 
be chosen, let $\Lambda = {\mathcal L}^{n}(B_{1}(0) \cap ({\mathbf R}^{n-1} \times [-\mu, \mu])).$ 
If $r$ is such that $\Lambda \left(\frac{r}{2}\right)^{n} > \g^{1/2}$, 
then in view of (\ref{gapcontrol-0}), there must exist a point 
$x_{0} \in B_{r/2}(q)  \cap ({\mathbf R}^{n-1} \times [-\mu r/2, \mu r/2])$ 
with $|w(x_{0}) - L(x_{0})| < \g^{1/4}.$  Then, $w(x_{0}) \leq \g^{1/4} +  C_{1}\mu r$ where $C_{1} = C_{1}(\b),$ 
so that

\begin{equation}\label{gapcontrol-2}
\mbox{inf}_{B_{r/2}(q)} \, w \leq \g^{1/4} + C_{1}\mu r.
\end{equation}

On the other hand, choosing $\mu^{\prime} = \mu^{\prime}(n) \in (0, 1/2)$ such that 
$\Lambda^{\prime} \equiv {\mathcal L}^{n}(B_{1}(0) \cap ({\mathbf R}^{n-1} \times 
[-\mu^{\prime}, \mu^{\prime}])) < \frac{\omega_{n}}{4}$, 
if $r$ also satisfies $(\omega_{n} - \Lambda^{\prime})\left(\frac{r}{2}\right)^{n} > \g^{1/2}$, 
then, again in view of (\ref{gapcontrol-0}),  
there must exist a point 
$x_{1} \in B_{r/2}(q) \setminus ({\mathbf R}^{n-1} \times [-\mu^{\prime} r/2, \mu^{\prime} r/2])$ such that 
$|w(x_{1}) - L(x_{1})| < \g^{1/4}.$ Then $w(x_{1}) > L(x_{1}) - \g^{1/4} \geq C_{1}\mu^{\prime} r  - \g^{1/4}$
and hence 

\begin{equation}\label{gapcontrol-3}
\mbox{sup}_{B_{r/2}(q)} \, w \geq C_{1}\mu^{\prime} r - \g^{1/4}.  
\end{equation}

Taking $\mu= \frac{\mu^{\prime}}{2 \cdot 3^{n}}$ and combining the 
inequalities (\ref{gapcontrol-1}),  (\ref{gapcontrol-2}) and (\ref{gapcontrol-3}), we then have that 
$r \leq  C\g^{1/4}$ where $C = C(\b, n).$ Thus in all cases, $r \leq C\g^{1/2n}.$
\end{proof}

\bigskip

\begin{lemma}\label{blowupregularity}
Let $\a \in (0, \pi)$ and $\d \in (0, 1).$ There exist numbers $\epsilon = \epsilon(n, \d, \a) \in (0, 1)$ and 
$\k = \k(n, \a) \in (0, 1)$ such that 
the following is true. If ${\widetilde P}_{0} = {\widetilde P}_{0}^{+} \cup {\widetilde P}_{0}^{-}$ 
is a pair of hyperplanes with $\a \leq \angle {\widetilde P}_{0} < \pi,$ ${\widetilde p}_{0}^{+} + 
{\widetilde p}_{0}^{-} \equiv 0,$ and
if $(v^{+}, v^{-}) \in {\mathcal F}_{\delta}$ satisfies $v^{+}(0) = v^{-}(0) =0,$  
and $\int_{B_{1}(0)} ({v}^{+} - {\widetilde p}_{0}^{+})^{2} 
+ ({v}^{-} - {\widetilde p}_{0}^{-})^{2} \leq \epsilon$, 
then there exist two harmonic functions 
$v_{1}, v_{2} \, : \,B_{\k}(0) \to {\mathbf R}$ such that 
$\left. v^{+}\right|_{B_{\k}(0)} = \mbox{max} \, \{v_{1}, v_{2} \}$ and  
$\left. v^{-}\right|_{B_{\k}(0)} = \mbox{min} \, \{v_{1}, v_{2} \}.$ Furthermore, the vanishing order of 
$v_{1} - v_{2}$ at any point $z \in B_{\k}(0)$ where $v_{1}(z) = v_{2}(z)$ is equal to 1.\\
\end{lemma}

\begin{proof} 

The hypotheses

\begin{equation}\label{reg2-1}
\int_{B_{1}(0)} ({v}^{+} - {\widetilde p}_{0}^{+})^{2} 
+ ({v}^{-} - {\widetilde p}_{0}^{-})^{2} \leq \epsilon
\end{equation}

\noindent
and $\a \leq \angle{\widetilde P}_{0}$  together with the fact that $\E_{1}^{2} = \int_{B_{1}(0)} (v^{+})^{2} +
(v^{-})^{2} \leq \left(\frac{3}{2}\right)^{n+2}$ (Proposition~\ref{properties}, part $(2)$) imply that

\begin{equation}\label{reg2-1-1}
\Lambda \leq \int_{B_{1}(0)} ({\widetilde p}_{0}^{+})^{2} + ({\widetilde p}_{0}^{-})^{2} \leq 
2\left(\frac{3}{2}\right)^{n+2} + 2\e
\end{equation}

\noindent
for some $\Lambda = \Lambda(n, \a)>0$, and consequently that

\begin{equation}\label{reg2-2}
\E_{1}^{2} = \int_{B_{1}(0)} (v^{+})^{2} + (v^{-})^{2} \geq \frac{\Lambda}{2} - \e \geq \frac{\Lambda}{4}, 
\end{equation}

\noindent
provided $\e = \e(n, \a) < \Lambda/4$.\\ 

Set ${\overline P}_{0} = {\rm graph} \, \frac{1}{\E_{1}}{\widetilde p}_{0}^{+} \cup 
{\rm graph} \, \frac{1}{\E_{1}}{\widetilde p}_{0}^{-}$ and 
$S^{(0)} = \{x \in {\mathbf R}^{n} \times \{0\} \, : \, {\rm dist} \, (x, {\pi} \, ({\overline P}^{+}_{0} \cap 
{\overline P}^{-}_{0})) \leq \th/16\}.$ Note that inequality (\ref{reg2-1}) implies that

\begin{eqnarray}\label{reg2-2-1}
\int_{B_{1}(0)} {\rm dist}^{2} \, ((x, {\widetilde v}_{1}^{+}(x)), {\overline P}_{0}) &+& 
{\rm dist}^{2} \, ((x, {\widetilde v}_{1}^{-}(x)), {\overline P}_{0}) \\ \nonumber
&\leq& \int_{B_{1}(0)} ({\widetilde v}_{1}^{+} - {\overline p}_{0}^{+})^{2} + 
({\widetilde v}_{1}^{-} - {\overline p}_{0}^{-})^{2} \\ \nonumber
&\leq& \left(\frac{2}{3}\right)^{-n-2}
\frac{4\e}{\Lambda} 
\end{eqnarray}

\noindent
(notation as in (\ref{blowupnotation-0})) which of course in particular says that

\begin{equation}\label{reg2-2-2}
\int_{B_{1}(0) \setminus S^{(0)}} ({\widetilde v}^{+}_{1} - {\overline p}_{0}^{+})^{2} + 
({\widetilde v}^{-}_{1} - {\overline p}_{0}^{-})^{2} \leq \left(\frac{2}{3}\right)^{-n-2}\frac{4\e}{\Lambda}.
\end{equation}

Since $\angle {\widetilde P}_{0} \in [\a, \pi)$ and $\E_{1}^{2} \leq \left(\frac{3}{2}\right)^{n+2}$, we have
that 

\begin{equation}\label{reg3-1}
{\overline \a}_{0} \leq \angle \, {\overline P}_{0} < \pi
\end{equation}

\noindent
for some ${\overline \a}_{0} = {\overline \a}_{0}(n, \a)>0.$ Now choose $\th = \th(n) \in (0, 1)$ such that 

$${\overline C}\th <\frac{1}{4}$$

\noindent
where ${\overline C} = {\overline C}(n)$ is as in Lemma~\ref{localdecomp}. 
If we then choose $\epsilon = \epsilon(n, \d, \a)$ so that 

\begin{equation}\label{reg4}
\left(\frac{2}{3}\right)^{-n-2}\frac{4\e}{\Lambda} < \overline\epsilon(n, \th, \d, {\overline \a}_{0})
\end{equation}

\noindent
where $\overline\e$ is as in Lemma~\ref{localdecomp},  we may apply Lemma~\ref{localdecomp}
with ${\overline P}_{0}$ in place of $\overline{P},$  
${\overline \a}_{0}$ in place of ${\overline \a}$, $\th$ in place of $\overline \th$ and ${\widetilde v}_{1}$ in place 
of $v$ to conclude that there exists a pair of hyperplanes ${\widetilde P}_{1}$ such that

\begin{eqnarray}\label{reg4-1}
\th^{-n-2}\int_{B_{\th}(0)} & 
{\rm dist}^{2}((x, {\widetilde v}_{1}^{+}(x)) , {\widetilde P}_{1}) 
+ {\rm dist}^{2} \, ((x, {\widetilde v}_{1}^{-}(x)), {\widetilde P}_{1}) \nonumber\\
&\leq {\overline C}\th^{2}\int_{B_{1}(0)} 
{\rm dist}^{2} \, ((x, {\widetilde v}_{1}^{+}(x)), {\overline P}_{0}) + {\rm dist}^{2} \, ((x, {\widetilde v}_{1}^{-}(x)), 
{\overline P}_{0}) \hspace{.2in} {\rm and}
\end{eqnarray}

\begin{eqnarray}\label{reg4-1-1}
\th^{-n-2}\int_{B_{\th}(0) \setminus S_{{\widetilde P}_{1}}(\th^{2}/16)} & 
({\widetilde v}_{1}^{+} - {\widetilde p}_{1}^{+})^{2} 
+ ({\widetilde v}_{1}^{-} - {\widetilde p}_{1}^{-})^{2} \nonumber\\
&\leq {\overline C}\th^{2}\int_{B_{1}(0)} 
{\rm dist}^{2} \, ((x, {\widetilde v}_{1}^{+}(x)), {\overline P}_{0}) + {\rm dist}^{2} \, ((x, {\widetilde v}_{1}^{-}(x)), 
{\overline P}_{0})
\end{eqnarray}

\noindent
where ${\overline C} = {\overline C}(n)$ is as in Lemma~\ref{localdecomp} and 
$S_{{\widetilde P}_{1}}(\s) = \{x  \in {\mathbf R}^{n} \times \{0\} \, : \, {\rm dist} \, (x, 
{\pi}\, ({\widetilde P}_{1}^{+} \cap {\widetilde P}_{1}^{-})) \leq \s\}.$\\   

Now, if $\e \leq \e_{1}(n,\a,1/2)$ where $\e_{1}$ is as in Lemma~\ref{Nbound}, we 
have by Lemmas~\ref{monotonicity}, ~\ref{coarseexcessbounds} and~\ref{Nbound} that

\begin{equation}\label{reg4-1-2}
1 \geq \frac{\E_{\th}^{2}}{\E_{1}^{2}} \geq \th^{2(N_{v}(1) - 1)} \geq \th.
\end{equation}

\noindent
Setting 
$${\overline P}_{1} = {\rm graph} \, \frac{\E_{1}}{\E_{\th}} 
{\widetilde p}^{+}_{1} \cup {\rm graph} \, \frac{\E_{1}}{\e_{\th}}{\widetilde p}_{1}^{-},$$ 
we conclude from (\ref{reg4-1}), (\ref{reg4-1-1}), (\ref{reg4-1-2}) and (\ref{reg5-0-0}) (with 
$\lambda = \frac{\E_{1}}{\E_{\th}} \in [1, \th^{-1/2}]$, so that 
${\rm dist}^{2} \, ((x^{\prime}, \lambda x^{n+1}), L^{\lambda}) \leq \th^{-1}
{\rm dist}^{2} \, ((x^{\prime}, x^{n+1}), L)$, where $L$, $L^{\lambda}$ are as in 
(\ref{reg5-0-0})) that

\begin{eqnarray}\label{reg5}
\int_{B_{1}(0)} & {\rm dist}^{2} \, ((x, {\widetilde v}_{\th}^{+}(x)), {\overline P}_{1}) + 
{\rm dist}^{2} \, ((x, {\widetilde v}_{\th}^{-}(x)),  {\overline P}_{1}) \nonumber\\
& \leq {\overline C}\th\int_{B_{1}(0)} 
{\rm dist}^{2} \, ((x, {\widetilde v}_{1}^{+}(x)),  {\overline P}_{0}) 
+ {\rm dist}^{2} \, ((x, {\widetilde v}_{1}^{-}(x)),  {\overline P}_{0})\nonumber\\
& \leq 4^{-1} {\e}_{2}\hspace{.2in} {\rm and}
\end{eqnarray}

\begin{eqnarray}\label{reg5-0}
\int_{B_{1}(0) \setminus S^{(1)}} & 
({\widetilde v}_{\th}^{+} - {\overline p}_{1}^{+})^{2} 
+ ({\widetilde v}_{\th}^{-} - {\overline p}_{1}^{-})^{2} \nonumber\\
&\leq {\overline C}\th\int_{B_{1}(0)} 
{\rm dist}^{2} \, ((x, {\widetilde v}_{1}^{+}(x)), {\overline P}_{0}) + {\rm dist}^{2} \, ((x, {\widetilde v}_{1}^{-}(x)), 
{\overline P}_{0})\nonumber\\
&\leq 4^{-1}{\e}_{2}
\end{eqnarray}

\noindent
where $S^{(1)} = \{x \in {\mathbf R}^{n} \times \{0\} \, : \, {\rm dist} \, (x, {\pi} \, ({\overline P}^{+}_{1} 
\cap {\overline P}^{-}_{1})) \leq \th/16\}$ and $\e_{2} = \left(\frac{2}{3}\right)^{-2n-2} \frac{4\e}{\Lambda}.$\\  

We claim that for each $j = 1, 2, \ldots,$ we can find a pair of hyperplanes ${\overline P}_{j}$ 
such that 

\begin{eqnarray}\label{reg5-1}
\int_{B_{1}(0)} && {\rm dist}^{2} \, ((x, {\widetilde v}_{\th^{j}}^{+}(x)),  {\overline P}_{j}) + 
{\rm dist}^{2} \, ((x, {\widetilde v}_{\th^{j}}^{-}(x)),  {\overline P}_{j}) \nonumber\\
&\leq& {\overline C}\th\int_{B_{1}(0)} {\rm dist}^{2} \, 
((x, {\widetilde v}_{\th^{j-1}}^{+}(x)), {\overline P}_{j-1})
+ {\rm dist}^{2} \, ((x, {\widetilde v}_{\th^{j-1}}^{-}(x)) ,{\overline P}_{j-1}) \hspace{.2in} {\rm and}
\end{eqnarray}

\begin{eqnarray}\label{reg5-1-1}
\int_{B_{1}(0) \setminus S^{(j)}} && ({\widetilde v}_{\th^{j}}^{+} - {\overline p}_{j}^{+})^{2} + 
({\widetilde v}_{\th^{j}}^{-} -  {\overline p}_{j}^{-})^{2} \nonumber\\
&\leq& {\overline C}\th\int_{B_{1}(0)} {\rm dist}^{2} \, 
((x, {\widetilde v}_{\th^{j-1}}^{+}(x)), {\overline P}_{j-1})
+ {\rm dist}^{2} \, ((x, {\widetilde v}_{\th^{j-1}}^{-}(x)) ,{\overline P}_{j-1})
\end{eqnarray}

\noindent
where $S^{(j)} = \{x \in {\mathbf R}^{n} \times \{0\} \, : \, {\rm dist} \, (x, {\pi} \, ({\overline P}_{j}^{+} \cap 
{\overline P}_{j}^{-})) \leq \th/16\}.$  We prove this by induction. Note that by (\ref{reg5}) and (\ref{reg5-0}), 
the assertion is true for $j=1.$ Suppose that it holds for all $j=1, 3, \ldots, i$ for some $i$. Thus

\begin{eqnarray}\label{reg6}
\int_{B_{1}(0)} && {\rm dist}^{2} \, ((x, {\widetilde v}_{\th^{i}}^{+}(x)),  {\overline P}_{i}) + 
{\rm dist}^{2} \, ((x, {\widetilde v}_{\th^{i}}^{-}(x)),  {\overline P}_{i}) \nonumber\\
&\leq& {\overline C}\th\int_{B_{1}(0)} {\rm dist}^{2} \, 
((x, {\widetilde v}_{\th^{i-1}}^{+}(x)), {\overline P}_{i-1})
+ {\rm dist}^{2} \, ((x, {\widetilde v}_{\th^{i-1}}^{-}(x)) ,{\overline P}_{i-1})\nonumber\\
 &\leq&\left({\overline C}\th\right)^{i} 
\int_{B_{1}(0)} {\rm dist}^{2} \, ((x, {\widetilde v}_{1}^{+}(x)),  {\overline P}_{0}) 
+ {\rm dist}^{2} \, ((x, {\widetilde v}_{1}^{-}(x)),  {\overline P}_{0})\nonumber\\
&\leq& 4^{-i} {\e}_{2}, 
\end{eqnarray}

\begin{eqnarray}\label{reg6-0}
\int_{B_{1}(0) \setminus S^{(j)}} && ({\widetilde v}_{\th^{j}}^{+} - {\overline p}_{j}^{+})^{2} + 
({\widetilde v}_{\th^{j}}^{-} -  {\overline p}_{j}^{-})^{2} \nonumber\\
&\leq& {\overline C}\th \int_{B_{1}(0)}  {\rm dist}^{2} \, ((x, {\widetilde v}^{+}_{\th^{j-1}}(x)), {\overline P}_{j-1}) 
+ {\rm dist}^{2} \, ((x, {\widetilde v}^{-}_{\th^{j-1}}(x)), {\overline P}_{j-1})\nonumber\\
&& \leq 4^{-j}{\e}_{2}\hspace{.2in} \mbox{and}
\end{eqnarray}

\begin{equation}\label{reg6-0-0}
\int_{B_{1}(0) \setminus S^{(j-1)}} ({\widetilde v}_{\th^{j-1}}^{+} - {\overline p}_{j-1}^{+})^{2} + 
({\widetilde v}_{\th^{j-1}}^{-} -  {\overline p}_{j-1}^{-})^{2} \leq 4^{-(j-1)}{\e}_{2} 
\end{equation}

\noindent
for $j = 1, 3, \ldots, i.$ Writing ${\widetilde P}_{j} = {\rm graph} \, 
\frac{\E_{\th^{j}}}{\E_{\th^{j-1}}} {\overline p}_{j}^{+} \cup {\rm graph} \, 
\frac{\E_{\th^{j}}}{\E_{\th^{j-1}}} {\overline p}_{j}^{-}$ and 
using the fact that $\E_{\th^{j}} \leq \E_{\th^{j-1}}$ (by Lemma~\ref{monotonicity}), 
we see from the inequality (\ref{reg6-0}) that 

\begin{equation}\label{reg6-0-1}
\th^{-n-2}\int_{B_{\th}(0) \setminus S_{{\widetilde P}_{j}}(\th^{2}/16)} 
({\widetilde v}_{\th^{j-1}}^{+} - {\widetilde p}_{j}^{+})^{2} + 
({\widetilde v}_{\th^{j-1}}^{-} -  {\widetilde p}_{j}^{-})^{2} \leq 4^{-j}{\e}_{2}
\end{equation}

\noindent
for $j = 1, 2, \ldots, i$, 
which together with the inequality (\ref{reg6-0-0}) implies, by the triangle inequality and homogeneity 
of ${\widetilde P}_{j}$, ${\overline P}_{j-1}$, that

\begin{equation}\label{reg6-2}
\int_{B_{1}(0)} ({\widetilde p}_{j}^{+} - {\overline p}_{j-1}^{+})^{2} + 
({\widetilde p}_{j}^{-} - {\overline p}_{j-1}^{-})^{2} \leq {\widetilde C}_{1}4^{-(j-1)}{\e}_{2}
\end{equation}
  
\noindent
for $j=1, 3, \ldots, i$, where ${\widetilde C}_{1} = {\widetilde C}_{1}(n, \a).$ Therefore,

\begin{equation}\label{reg6-3}
\|{\widetilde p}_{j}^{+} - {\widetilde p}_{j}^{-}\|_{L^{2}(B_{1}(0))} \geq 
\|{\overline p}_{j-1}^{+} - {\overline p}_{j-1}^{-}\|_{L^{2}(B_{1}(0))} 
- 2\sqrt{{\widetilde C}_{1}{\e}_{2}} \, 2^{-(j-1)}
\end{equation}

\noindent
and hence, by the definition of ${\widetilde p}_{j}^{\pm}$ and the fact that 
$\E_{\th^{j}} \leq \E_{\th^{j-1}}$, 

\begin{equation}\label{reg6-4}
\|{\overline p}_{j}^{+} - {\overline p}_{j}^{-}\|_{L^{2}(B_{1}(0))} \geq 
\|{\overline p}_{j-1}^{+} - {\overline p}_{j-1}^{-}\|_{L^{2}(B_{1}(0))} 
- 2\sqrt{{\widetilde C}_{1} {\e}_{2}} \, 2^{-(j-1)}.
\end{equation}

\noindent
Summing over $j$, we conclude from this that 

\begin{equation}\label{reg6-5}
\|{\overline p}_{i}^{+} - {\overline p}_{i}^{-}\|_{L^{2}(B_{1}(0))} \geq 
\|{\overline p}_{0}^{+} - {\overline p}_{0}^{-}\|_{L^{2}(B_{1}(0))} - 4\sqrt{{\widetilde C}_{1}{\e}_{2}}.
\end{equation}

\noindent
By inequality (\ref{reg6-0}), Proposition~\ref{properties}, part $(2)$ and homogeneity of 
${\overline p}_{j}^{\pm}$, it follows that 
$\int_{B_{1}(0)} ({\overline p}_{j}^{+} + {\overline p}_{j}^{-})^{2} \leq C$ for 
some fixed constant $C = C(n) \in (0, \infty),$ and hence, provided $\e = \e(n, \a)$ is sufficiently small, 
we have from the estimate (\ref{reg6-5}) that 
$$\pi > \angle \, {\overline P}_{i} \geq \b$$ 

\noindent
where $\b = \b(n, \a) \in (0, \pi/2)$ is a fixed angle. Thus, since
$({\widetilde v}^{+}_{\th^{j}}, {\widetilde v}_{\th^{j}}^{-}) \in {\mathcal F}_{\delta}$, we may apply 
Lemma~\ref{localdecomp} with $\th$ in place of $\overline\th$, $b$ in place of 
$\overline\a$, 
$({\widetilde v}^{+}_{\th^{i}}, {\widetilde v}^{-}_{\th^{i}})$ in place of $(v^{+}, v^{-})$ and 
${\overline P}_{i}$ in place of ${\overline P}$ to conclude that 
there exists a pair of hyperplanes ${\widetilde P}_{i+1}$ such that

\begin{eqnarray}\label{reg6-6}
\th^{-n-2}\int_{B_{\th}(0)} & 
{\rm dist}^{2}((x, {\widetilde v}_{\th^{i}}^{+}(x)) , {\widetilde P}_{i+1}) 
+ {\rm dist}^{2} \, ((x, {\widetilde v}_{\th^{i}}^{-}(x)), {\widetilde P}_{i+1}) \nonumber\\
&\leq {\overline C}\th^{2} \int_{B_{1}(0)} 
{\rm dist}^{2} \, ((x, {\widetilde v}_{\th^{i}}^{+}(x)), {\overline P}_{i}) 
+ {\rm dist}^{2} \, ((x, {\widetilde v}_{\th^{i}}^{-}(x)), {\overline P}_{i}) \hspace{.2in} {\rm and}
\end{eqnarray}

\begin{eqnarray}\label{reg6-6-1}
\th^{-n-2}\int_{B_{\th}(0) \setminus S_{{\widetilde P}_{i+1}}({\th}^{2}/16)} & 
({\widetilde v}_{\th^{i}}^{+} - {\widetilde p}_{i+1}^{+})^{2} 
+ ({\widetilde v}_{\th^{i}}^{-} - {\widetilde p}_{i+1}^{-})^{2} \nonumber\\
&\leq {\overline C}\th^{2} \int_{B_{1}(0)} 
{\rm dist}^{2} \, ((x, {\widetilde v}_{\th^{i}}^{+}(x)), {\overline P}_{i}) 
+ {\rm dist}^{2} \, ((x, {\widetilde v}_{\th^{i}}^{-}(x)), {\overline P}_{i}).
\end{eqnarray}

\noindent
It follows from the triangle inequality, the inequalities (\ref{reg6-0}), (\ref{reg6-6-1}) and 
homogeneity of ${\widetilde P}_{i+1}$, ${\overline P}_{i}$ that

\begin{equation}\label{reg6-7}
\int_{B_{1}(0)} ({\widetilde p}_{i+1}^{+} - {\overline p}_{i}^{+})^{2} + 
({\widetilde p}_{i+1}^{-} - {\overline p}_{i}^{-})^{2} \leq {\widetilde C}_{1}4^{-i}{\e}_{2}
\end{equation}
  
\noindent
where ${\widetilde C}_{1} = {\widetilde C}_{1}(n, \a)$ is as in (\ref{reg6-2}).\\ 

Note again that by Lemmas~\ref{monotonicity}, ~\ref{coarseexcessbounds}, the monotonicity of the frequency function 
$N_{v}(\cdot)$ and Lemma~\ref{Nbound}, we have

\begin{equation}\label{reg6-11}
1 \geq \frac{\E^{2}_{\th^{i+1}}}{\E^{2}_{\th^{i}}} \geq \th^{2(N_{v}(1) - 1)} \geq \th
\end{equation}

\noindent
so setting ${\overline P}_{i+1} = {\rm graph} \, \frac{\E_{\th^{i}}}{\E_{\th^{i+1}}} 
{\widetilde p}^{+}_{i+1} \cup   {\rm graph} \, \frac{\E_{\th^{i}}}{\E_{\th^{i+1}}} {\widetilde p}_{i+1}^{-}$ 
and using the bound (\ref{reg6-11}), we obtain from (\ref{reg6-6}), (\ref{reg6-6-1}) 
and (\ref{reg5-0-0}) that

\begin{eqnarray}\label{reg6-12}
\int_{B_{1}(0)} & {\rm dist}^{2} \, ((x, {\widetilde v}_{\th^{i+1}}^{+}(x)), {\overline P}_{i+1}) + 
{\rm dist}^{2} \, ((x, {\widetilde v}_{\th^{i+1}}^{-}(x)),  {\overline P}_{i+1}) \nonumber\\
& \leq {\overline C}\th\int_{B_{1}(0)} {\rm dist}^{2} \, 
((x, {\widetilde v}_{\th^{j}}^{+}(x)),  {\overline P}_{j}) 
+ {\rm dist}^{2} \, ((x, {\widetilde v}_{\th^{j}}^{-}(x)),  {\overline P}_{j}) \hspace{.2in} {\rm and}
\end{eqnarray}

\begin{eqnarray}\label{reg6-13}
\int_{B_{1}(0) \setminus S^{(i+1)}} & 
({\widetilde v}_{\th^{i+1}}^{+} - {\overline p}_{i+1}^{+})^{2} 
+ ({\widetilde v}_{\th^{i+1}}^{-} - {\overline p}_{i+1}^{-})^{2} \nonumber\\
&\leq {\overline C}\th \int_{B_{1}(0)} 
{\rm dist}^{2} \, ((x, {\widetilde v}_{\th^{i}}^{+}(x)), {\overline P}_{i}) 
+ {\rm dist}^{2} \, ((x, {\widetilde v}_{\th^{i}}^{-}(x)), {\overline P}_{i})
\end{eqnarray}

\noindent
where $S^{(i+1)} = \{x \in {\mathbf R}^{n} \times \{0\} \, : \, {\rm dist} \, (x, \pi \, ({\overline P}_{i}^{+} \cap 
{\overline P}_{i}^{-})) \leq \th/16\}.$ This completes the induction.\\
 
We thus obtain a sequence of pairs of hyperplanes 
${\overline P}_{j}$, $j=1, 2, 3, \ldots$ satisfying  (\ref{reg6}) and 
(\ref{reg6-0}). Now let 
$P^{\pm}_{j} = {\rm graph} \, \E_{\th^{j}}{\overline p}_{j}^{\pm}.$ Then (\ref{reg6}), (\ref{reg6-0}) 
and (\ref{reg5-0-0}) say that 

\begin{equation}\label{reg6-13-1}
\left(\frac{2}{3}\th^{j}\right)^{-n-2}\int_{B_{\frac{2}{3}\th^{j}}(0)} {\rm dist}^{2} \, ((x, {v}^{+}(x)), {P}_{j}) 
+ {\rm dist}^{2} \, ((x, {v}^{-}(x)), {P}_{j}) \leq 4^{-j}\left(\frac{3}{2}\right)^{n+2}{\e}_{2}\hspace{.2in} \mbox{and}
\end{equation}

\begin{equation}\label{reg6-13-2}
\left(\frac{2}{3}\th^{j}\right)^{-n-2}\int_{B_{\frac{2}{3}\th^{j}}(0) \setminus S_{P_{j}}(\th^{j+1}/24)} ({v}^{+} - {p}_{j}^{+})^{2} 
+ ({v}^{-} -  {p}_{j}^{-})^{2} \leq 4^{-j}\left(\frac{2}{3}\right)^{n+2}{\e}_{2}
\end{equation}
 
\noindent
for all $j = 0, 1, 2, \ldots$, where we have used the fact that 
$\E_{\th^{j}} \leq \E_{1} \leq \left(\frac{3}{2}\right)^{n+2}.$ 
By the triangle inequality and the homogeneity of $P_{j}$, $P_{j-1}$, (\ref{reg6-13-2}) implies that 

\begin{equation}\label{reg6-14}
\|(p_{j}^{+}, p_{j}^{-}) - (p_{j-1}^{+}, p_{j-1}^{-})\|_{L^{2}(B_{1}(0))} \leq C4^{-(j-1)}{\e}_{2}
\end{equation}

\noindent
where $C = C(n,\a).$ i.e. that $(p_{j}^{+}, p_{j}^{-})$ is a Cauchy sequence. Hence 
there exists a pair of hyperplanes $P$ such that $P_{j} \to P.$ 
We then have by the triangle inequality and the inequalities (\ref{reg6-13-1}), (\ref{reg6-13-2}) and 
(\ref{reg6-14}) that 

\begin{equation}\label{reg8}
(\frac{2}{3}\th^{j})^{-n-2}\int_{B_{\frac{2}{3}\th^{j}}(0)} {\rm dist}^{2} \, ((x, {v}^{+}(x)),  {P}) + 
{\rm dist}^{2} \, ((x, v^{-}(x)),  {P}) \leq C4^{-j}{\e},
\end{equation}

\begin{equation}\label{reg8-0}
\left(\frac{2}{3}\th^{j}\right)^{-n-2}\int_{B_{\frac{2}{3}\th^{j}}(0) \setminus S_{P_{j}}(\th^{j+1}/24)} ({v}^{+} - {p}^{+})^{2} 
+ ({v}^{-} -  {p}^{-})^{2} \leq C4^{-j}{\e} \hspace{.2in} \mbox{and}
\end{equation}

\begin{equation}\label{reg8-1}
\|(p_{j}^{+}, p_{j}^{-}) - (p^{+}, p^{-})\|_{L^{2}(B_{1}(0))} \leq C4^{-j}{\e}
\end{equation}

\noindent
for all $j=0, 1, 2, \ldots$, where $C = C(n, \a)$. Now, given any $\r \in (0, 1/4)$, 
there exists a unique non-negative integer $j^{\star}$ such that $\frac{2}{3}\th^{j^{\star}+1} \leq \r < 
\frac{2}{3}\th^{j^{\star}}.$ Using the estimates (\ref{reg8}), (\ref{reg8-0}) and (\ref{reg8-1}) 
with $j = j^{\star}$, we obtain that

\begin{equation}\label{reg9}
\r^{-n-2}\int_{B_{\r}(0)} {\rm dist}^{2} \, ((x, {v}^{+}(x)), {P}) + 
{\rm dist}^{2} \, ((x, v^{-}(x)), {P}) \leq C\r^{\mu}{\e},
\end{equation} 

\begin{equation}\label{reg9-0}
\r^{-n-2}\int_{B_{\r}(0) \setminus S_{P_{j^{\star}}}(\r/16)} ({v}^{+} - {p}^{+})^{2} 
+ ({v}^{-} -  {p}^{-})^{2} \leq C\r^{\mu}{\e}\hspace{.2in} \mbox{and}
\end{equation}

\begin{equation}\label{reg9-0-1}
\|(p_{j^{\star}}^{+}, p_{j^{\star}}^{-}) - (p^{+}, p^{-})\|_{L^{2}(B_{1}(0))} \leq C\r^{\mu}{\e}
\end{equation}

\noindent
where $C = C(n, \a)>0$ and $\mu = \mu(n, \a)>0.$ Since (\ref{reg9-0-1}) implies 

\begin{equation}\label{reg9-0-1-1}
d_{\mathcal H} \, (T_{P_{j^{\star}}} \cap B_{\r}(0), T_{P} \cap B_{\r}(0)) \leq C\r^{1+\mu}{\e}
\end{equation}

\noindent
where $C = C(n, \a)$ and $T_{P}$ denotes the orthogonal projection of the 
axis $P^{+} \cap P^{-}$ of $P$ onto ${\mathbf R}^{n} \times \{0\},$ we deduce from (\ref{reg9-0}) that 

\begin{equation}\label{reg9-0-2}
\r^{-n-2}\int_{B_{\r}(0) \setminus S_{P}(\r/8)} ({v}^{+} - {p}^{+})^{2} 
+ ({v}^{-} -  {p}^{-})^{2} \leq C\r^{\mu}{\e}\hspace{.2in} \mbox{and}
\end{equation}

\noindent
provided $C\e \leq 1/16,$ where $C$ is as in (\ref{reg9-0-1-1}). Thus, we have the estimates 
(\ref{reg9}) and (\ref{reg9-0-2}) for all $\r \in (0, 1/4]$ provided 
$\e = \e(n, \a) \in (0, 1)$ is sufficiently small. Note also that (\ref{reg8-1}) in particular says 

\begin{equation}\label{reg9-1}
\|(p^{+}, p^{-}) - ({\widetilde p}_{0}^{+}, {\widetilde p}_{0}^{-})\|_{L^{2}(B_{1}(0))} \leq C{\e}
\end{equation}

\noindent
where $C = C(n, \a),$ which implies that if $\e = \e(n, \a)$ is sufficiently small, $P$ must be 
a transverse pair of hyperplanes with $\a/2 \leq \angle \, P < \pi$. Hence, provided 
$\e = \e(n, \a) \in (0, 1)$ is chosen sufficiently small, the estimate of Lemma~\ref{lipschitz}, 
part $(b)$ together with the estimate (\ref{reg9-0-2}) implies that 

\begin{equation}\label{reg9-1-1}
Z_{w} \cap (B_{\r}(0) \setminus S_{P}(\r/8)) = \emptyset \hspace{.2in} \mbox{for each $\r \in (0, 1/4]$}
\end{equation}

\noindent
where $Z_{w} = \{z \, : \, v^{+}(z) = v^{-}(z)\}.$ i.e. that $Z_{w}  \cap B_{1/4}(0)$ is contained in a cone with vertex at the origin, axis the orthogonal projection 
of the axis of $P$ onto ${\mathbf R}^{n} \times \{0\}$ and with a fixed cone angle depending only on $n$.\\

Next we argue that provided $\e = \e(n, \a)$ is sufficiently small, 
the decay estimates (\ref{reg9}), (\ref{reg9-0-2}) and the cone condition 
(\ref{reg9-1-1}) hold uniformly for each ``base point'' 
$z \in Z_{w}$ sufficiently close to the origin, with a unique choice of a pair of affine hyperplanes 
$P_{z}$ depending on $z$. So let $z \in B_{1/4}(0)$ be such that 
$v^{+}(z) = v^{-}(z).$ Set 
$V^{(z) \, \pm}(x) = {\widetilde v}^{\pm}_{z, 1/2}(x)$ for $x \in B_{1}(0)$ where the notation 
is as in (\ref{blowupnotation-1}). Then $(V^{(z) \, +}, V^{(z) \, -}) \in {\mathcal F}_{\delta}$ and  
$V^{(z) \, \pm}(0) = 0.$  Note that by the standard estimates for harmonic functions we have that, 
since $y = h(z)$, 

\begin{equation}\label{reg10-3}
|y|, \,\, |Dh(z)| \leq C|z|\left(\int_{B_{1}(0)}(v^{+})^{2} + (v^{-})^{2}\right)^{1/2} \leq C|z|
\end{equation}
 
\noindent
for all $z \in B_{1/4}(0)$, where $C = C(n).$ Also note that it follows from the inequality (\ref{reg2-1}) that
provided $\e = \e(n \a) \in (0, 1)$ is sufficiently small, 
\begin{equation}
{\widetilde C} \geq \E_{z, 1/2}^{2} \geq C >0
\end{equation}

\noindent
where ${\widetilde C} = {\widetilde C}(n)$ and $C = C(n, \a).$\\ 

Now set ${\widetilde p}^{(z) \, \pm}(x) = \frac{1}{2\E_{z, 1/2}}{\widetilde p}_{0}^{\pm}(x).$
Then $\pi > \angle {\widetilde P}^{(z)} \geq {\widetilde \a},$ where 
${\widetilde a} = {\widetilde a}(n, \a)>0.$  
It is then easy to see directly from the definition of $V^{(z) \, \pm}$ and the estimates (\ref{reg10-3}) that 
there exists $\g = \g(n, \a)>0$ and $\k = \k(n, \a)>0$ such that for all $z \in B_{\k}(0)$ with 
$v^{+}(z) = v^{-}(z)$, 

\begin{eqnarray}\label{reg10-6}
\int_{B_{1}(0)} \left(V^{(z) \, +} - {\widetilde p}^{(z) \, +}\right)^{2} 
&+& \left(V^{(z) \, -} - {\widetilde p}^{(z) \, -}\right)^{2} \nonumber\\ 
&\leq&  
\frac{3^{n+2}}{\E_{z, 1/2}^{2}}\int_{B_{1/3}(z)} (v^{+}(x) - \frac{1}{3}{\widetilde p}_{0}^{+}(3x - z))^{2} \nonumber\\ 
&& \hspace{1in} + \,\,(v^{-}(x) - \frac{1}{3}{\widetilde p}_{0}^{-}(3x - z) )^{2} \, dx\nonumber\\
&\leq& \frac{2 \cdot 3^{n+2}}{\E_{z, 1/2}^{2}}\int_{B_{1}(0)}(v^{+}(x) - {\widetilde p}_{0}^{+}(x))^{2} + 
(v^{-}(x) - {\widetilde p}_{0}^{-}(x))^{2} \nonumber\\
&&\hspace{.5in} + \,\,\left({\widetilde p}_{0}^{+}(x) - {\widetilde p}_{0}^{+}\left(x - \frac{z}{3}\right)\right)^{2} 
+ \left({\widetilde p}_{0}^{-}(x) - {\widetilde p}_{0}^{-}\left(x - \frac{z}{3}\right)\right)^{2} \, dx\nonumber\\
&\leq& \frac{2\cdot3^{n+2}}{\E_{z, 1/2}^{2}}\int_{B_{1}(0)}(v^{+} - {\widetilde p}_{0}^{+})^{2} + 
(v^{-} - {\widetilde p}_{0}^{-})^{2}  + C|z|^{2}\nonumber\\
&\leq&\e
\end{eqnarray}

\noindent
where $C = C(n, \a)$ and $\e = \e(n, {\widetilde \a})$ is as in the argument (with ${\widetilde \a}$ in place of $\a$) 
leading to the estimates (\ref{reg9}) and (\ref{reg9-0-2}), 
provided $\int_{B_{1}(0)}({v}^{+} - {\widetilde p}_{0}^{+})^{2} + ({v}^{-} - {\widetilde p}_{0}^{-})^{2} 
\leq \g\e.$\\  
  
Therefore, if the hypotheses of the lemma are satisfied with $\g\e$ in place of $\e$, 
we may repeat the argument leading to the estimates (\ref{reg9}), (\ref{reg9-0-2}), (\ref{reg9-1}) and the cone 
condition (\ref{reg9-1-1}) with 
$V^{(z) \, \pm}$  in place of $v^{\pm}$ and  ${\widetilde p}^{(z) \, \pm}$ in 
place of ${\widetilde p}_{0}^{\pm}.$ This will yield for each $z \in B_{\k}(0)$ with $v^{+}(z) = v^{-}(z)$ a pair
of transverse hyperplanes $P_{z} = P_{z}^{+} \cup P_{z}^{-}$ satisfying

\begin{equation}\label{reg11}
\r^{-n-2}\int_{B_{\r}(z)} {\rm dist}^{2} \, ((x, {v}^{+}(x)),  (z, y) + {P}_{z}) + 
{\rm dist}^{2} \, ((x, v^{-}(x)), (z, y) + {P}_{z})\leq C\r^{\mu} {\e},
\end{equation} 

\begin{equation}\label{reg11-0}
\r^{-n-2}\int_{B_{\r}(z) \setminus S_{(z, y) + P_{z}}(\r/8)} ({v}^{+}(x) - (y + {p}_{z}^{+}(x - z)))^{2} 
+ ({v}^{-}(x) -  (y + {p}_{z}^{-}(x - z)))^{2} \leq C\r^{\mu}{\e},
\end{equation}

\begin{equation}\label{reg11-0-0}
\|(p_{z}^{+}, p_{z}^{-}) - ({\widetilde p}_{0}^{+}, {\widetilde p}_{0}^{-})\|_{L^{2}(B_{1}(0))} \leq C \e \hspace{.2in}
\mbox{and}
\end{equation}

\begin{equation}\label{reg11-0-1}
Z_{w} \cap (B_{\r}(z) \setminus S_{(z, y) + P_{z}}(\r/8)) = \emptyset
\end{equation}

\noindent
for all $\r \in (0, 1/12).$ Here $y = v^{+}(z) = v^{-}(z)$ and $C = C(n, \a).$\\ 

Note that by the estimates (\ref{reg11-0}), (\ref{reg11-0-0}), Lemma~\ref{monotonicity} and the triangle inequality, 
it follows that for each $z \in B_{\k}(0)$ and $\r \in (0, 1/12)$, 

\begin{equation}\label{reg11-0-2}
C \leq \E_{z, \r} \leq {\widetilde C}
\end{equation}

\noindent
for fixed $C = C(n \a)>0$ and ${\widetilde C} = {\widetilde C}(n) < \infty.$\\ 

Next we assert that $Z_{w} \cap B_{\k/2}(0)$ projects fully 
onto the axis ${\widetilde P}_{0}^{+} \cap {\widetilde P}_{0}^{-} \cap B_{\k/2}(0).$  
To see this, first note that since ${\widetilde p}_{0}^{+} + {\widetilde p}_{0}^{-} \equiv 0$ by 
hypothesis, we have that ${\widetilde P}_{0}^{+} \cap {\widetilde P}_{0}^{-} \subset {\mathbf R}^{n} \times 
\{0\}.$ For notational convenience, 
(and without loss of generality, by making an orthogonal rotation of ${\mathbf R}^{n} \times \{0\}$) 
let us assume that ${\widetilde P}_{0}^{+} \cap {\widetilde P}_{0}^{-} = {\mathbf R}^{n-1} \times \{(0,0)\}.$ 
If there is a point 
$(\xi, 0,0) \in ({\mathbf R}^{n-1} \times \{(0,0)\}) \cap B_{\k/2}(0)$ with ${\mathbf p}^{-1}(\xi, 0,0) 
\cap Z_{w} = \emptyset$, 
where ${\mathbf p} \, : \, {\mathbf R}^{n} \times \{0\} \to {\mathbf R}^{n-1} \times \{(0,0)\}$ is the orthogonal 
projection, then, since $Z_{w}$ is a closed set, there must exist $r > 0$ such that 

\begin{equation}\label{reg12}
(B_{r}^{n-1}(\xi,0,0) \times {\mathbf R} \times \{0\}) \cap Z_{w} = \emptyset \hspace{.2in} {\rm and}
\end{equation}

$$(\overline B_{r}^{n-1}(\xi, 0,0) \times {\mathbf R} \times \{0\}) \cap Z_{w} \neq \emptyset.$$  

\noindent
Choose $z \in (\overline B_{r}^{n-1}(\xi, 0,0) \times {\mathbf R} \times \{0\}) \cap Z_{w}.$\\ 

Note next the following fact: Let $\a_{1} \in (0, \pi).$ Then for any given $\eta$, there exists $\z =\z(\a_{1}, \eta)$ with 
$\z \downarrow 0$ as $\eta \downarrow 0$ such that 
if $v = (v^{+}, v^{-}) \in {\mathcal F}_{\delta}$ satisfies $\int_{B_{1}(0)} {\rm dist}^{2} \,((x, v^{+}(x)), {P}_{1}) + 
{\rm dist}^{2} \, ((x, v^{-}(x)), {P}_{1}) \leq \z$ and 
$\int_{B_{1}(0) \setminus S_{{P}_{1}}(1/8)} (v^{+} - {p}_{1}^{+})^{2} + 
(v^{-} - {p}_{1}^{-})^{2} \leq \z$ for some pair of hyperplanes ${P}_{1}$ 
with $\a_{1} \leq \angle \, {P}_{1} < \pi$, then $\int_{B_{1}(0)} (w - {L}_{1})^{2} \leq \eta$ 
where $w = \frac{1}{2}(v^{+} - v^{-})$ and $L_{1} = \frac{1}{2}({p}_{1}^{+} - {p}_{1}^{-}).$
(This can easily be seen by arguing by contradiction.) Since the estimates (\ref{reg11}) and 
(\ref{reg11-0}) say that for each $\r \in (0, 1/8)$,     
$\int_{B_{1}(0)} {\rm dist}^{2}\,((x, {\widetilde v}_{z,\r}^{+}(x)), {P}_{z}^{(\r)}) + 
{\rm dist}^{2} \, ((x, {\widetilde v}_{z, \r}^{-}(x)), {P}_{z}^{(\r)}) \leq C\r^{\mu}\e$ 
and $\int_{B_{1}(0) \setminus S_{P_{z}^{(\r)}}(1/8)} ({\widetilde v}_{z, \r}^{+} - {p}_{z}^{(\r) \, +})^{2} 
+ ({\widetilde v}_{z, \r}^{-} - {p}_{z}^{(\r) \, -})^{2} \leq C\r^{\mu}\e$
where $p_{z}^{(\r) \, \pm} = \frac{1}{\E_{z,\r}} p_{z}^{\pm}$ and the estimates 
(\ref{reg11-0-2}) say that $P_{z}^{(\r)}$ satisfies $\a_{1} \leq \angle \, P_{z}^{(\r)} < \pi$ for 
some $\a_{1} = \a_{1}(n, \a) > 0$, 
it follows that for any given $\eta \in (0, 1/2)$, there exists 
$\r = \r(n, \a, \eta) \in ), 1/2)$ such that $\int_{B_{1}(0)} ({\widetilde w}_{z,\r} - L_{z}^{(\r)})^{2} 
\leq \eta$ where ${\widetilde w}_{z,\r} = \frac{1}{2}({\widetilde v}_{z,\r}^{+} - {\widetilde v}_{z, \r}^{-})$ and 
$L_{z}^{(\r)} = \frac{1}{2}(P_{z}^{(\r) \, +} - P_{z}^{(\r) \, -}).$ Thus, we may  apply
Lemma~\ref{gapcontrol} with ${\widetilde v}^{\pm}_{z, \r}$ in place of $v^{\pm}$ for a 
suitable choice of sufficiently small $\r \in (0, r/4)$ to arrive at a contradiction of (\ref{reg12}). 
(Note that here we have also used the fact that $\pi \, (P_{z}^{(\r) \, +} \cap P_{z}^{(\r) \, -})$ remains 
close to ${\widetilde P}_{0}^{+} \cap {\widetilde P}_{0}^{-}$ as $\r \downarrow 0,$ which follows 
from the estimate (\ref{reg11-0-0}).) 
Hence $Z_{w} \cap B_{\k/2}(0)$ must have full 
projection onto ${\widetilde P}_{0}^{+} \cap {\widetilde P}_{0}^{-} \cap B_{\k/2}(0).$\\ 

It then follows first from the estimates (\ref{reg11}), 
(\ref{reg11-0-0}) and (\ref{reg11-0-1}) that $Z_{w} \cap B_{\k/2}(0)$ 
is equal to a Lipschitz graph (over ${\widetilde P}_{0}^{+} \cap {\widetilde P}_{0}^{-} \cap B_{\k/2}(0)$) and then by 
the estimate (\ref{reg11-0}) that this graph is $C^{1,\mu}.$ This implies directly
that the union of the graphs of $v^{+}$,  $v^{-}$ over $B_{\k/2}(0)$ is equal to the 
union of the graphs of two harmonic functions $v^{1}, v^{2} \, : B_{\k/2}(0) \to {\mathbf R}$. 
Specifically, if we let $\Omega^{\pm}$ denote the 
two components of $B_{\k/2}(0) \setminus Z_{w}$ and define a function $v^{1}$ 
on $B_{\k/2}(0)$ by setting $v^{1}(x) = v^{+}(x)$ if $x \in \overline\Omega^{+}$,  
and $v^{1}(x) = v^{-}(x)$ if $x \in \Omega^{-}$, we see first by 
(\ref{reg11}) that $v^{1} \in C^{1}(B_{\k/2}(0))$ and then by integration by 
parts that $\int_{B_{\k/2}(0)} Dv^{1} \cdot D\z = \int_{B_{\k/2}(0) \cap {\Omega^{+}}} 
Dv^{1}\cdot D\z + \int_{B_{\k/2}(0) \cap \Omega^{-}} Dv^{1} \cdot D\z = 0$ for every 
$\z \in C^{1}_{c}(B_{\k/2}(0)).$ Thus $v^{1}$ is harmonic. Similarly, we may 
define $v^{2} \, : \, B_{\k/2}(0) \to {\mathbf R}$ by setting $v^{2}(x) = v^{-}(x)$ if $x \in \overline\Omega^{+}$,  
and $v^{2}(x) = v^{+}(x)$ if $x \in \Omega^{-}$, and check that $v^{2}$ is also harmonic.\\

Finally, since by (\ref{reg11}) and (\ref{reg11-0}) the tangent planes to the graphs of $v^{1}$ and $v^{2}$ at any point 
$(z, y)$ where $v^{1}(z) = v^{2}(z)=y$ are transversely intersecting, 
it follows that the vanishing order of $v^{1} - v^{2}$ at such a point must be equal to 1. Thus
the lemma holds with $\k/2$ in place of $\k$ and $\g\e$ in place of $\e$.
\end{proof}

\medskip

\noindent
{\bf Definition}: Given $v = (v^{+}, v^{-}) \in {\mathcal F}_{\delta}$, we shall call a point $z \in B_{3/2}(0)$ 
a {\em branch point} of $v$ if there exists no $\s > 0$ such that 
$(\mbox{graph} \, v^{+} \cup \mbox{graph} \, v^{-}) \cap (B_{\s}(z) \times {\mathbf R})$ is 
equal to the union of the graphs of two harmonic functions over $B_{\s}(z).$\\

\noindent
{\bf Remark}: It follows directly from Proposition~\ref{properties}(2) and 
Proposition~\ref{wproperties} (5) 
that if $z$ is a branch point of $v = (v^{+}, v^{-})$, then $z \in Z_{w}$, i.e. that $v^{+}(z) = v^{-}(z).$
Furthermore, if $v^{+} \equiv v^{-}$, then $v^{\pm}$ are each harmonic, so no point 
$z \in B_{3/2}(0)$ is a branch point in this case.\\

Using Lemma~\ref{blowupregularity}  and adapting techniques due to L. Simon \cite{S}, we  
establish in the next two lemmas crucial uniform asymptotic decay estimates for any function 
$v \in {\mathcal F}_{\delta}$ at a branch point.\\ 

\begin{lemma} \label{spherical}
Let $\d \in (0, 1).$ There exists a constant $c = c(n, \d)>0$ such that the following is true. If 
$v=(v^{+}, v^{-}) \in {\mathcal F}_{\delta}$, 
$v^{+}(0) = v^{-}(0)=0$, $Dh(0) = 0$, where $h = \frac{1}{2}(v^{+} + v^{-})$, 
and if either 

\begin{itemize}
\item[$(a)$] the origin is a branch point of $v$ or 
\item[$(b)$] $w \not\equiv 0$ and ${\mathcal N}_{w}(0) > 1$, where $w = \frac{1}{2}(v^{+} - v^{-})$,
\end{itemize}
\noindent
then  
$$\int_{B_{1}(0) \setminus B_{1/2}(0)} \left(\frac{\partial (v^{+}/R)}{\partial R}\right)^{2}
+ \left(\frac{\partial (v^{-}/R)}{\partial R}\right)^{2} \geq c \int_{B_{1}(0)} (v^{+})^{2} + (v^{-})^{2}.$$
\end{lemma}

\begin{proof} If the lemma is not true, there exists a sequence of functions 
$v_{k} = (v^{+}_{k}, v^{-}_{k}) \in 
{\mathcal F}_{\delta}$ satisfying $v_{k}^{+}(0) = v_{k}^{-}(0) = 0$ and
$Dh_{k}(0) = 0$ where $h_{k} = \frac{1}{2}(v_{k}^{+} + v_{k}^{-})$, such that 
for each $k$, either the origin is a branch point of $v_{k}$ (in which case $w_{k} \not\equiv 0,$ where 
$w_{k} = \frac{1}{2}(v_{k}^{+} - v_{k}^{-}$)  or  
$w_{k} \not\equiv 0$ and ${\mathcal N}_{w_{k}}(0) > 1$, and

\begin{equation} \label{sp1}
\int_{B_{1}(0) \setminus B_{1/2}(0)} \left(\frac{\partial (v_{k}^{+}/R)}{\partial R}\right)^{2}
+ \left(\frac{\partial (v_{k}^{-}/R)}{\partial R}\right)^{2} \leq \frac{1}{k}\int_{B_{1}(0)} (v^{+}_{k})^{2} + (v^{-}_{k})^{2}.
\end{equation}

Let ${\widetilde v}_{k}^{\pm}(x) 
= \frac{3}{2}\frac{v_{k}^{\pm}(2x/3)}{\left(\int_{B_{1}(0)}(v_{k}^{+})^{2} + (v_{k}^{-})^{2}\right)^{1/2}}.$  
Then ${\widetilde v}_{k} \equiv ({\widetilde v}_{k}^{+}, {\widetilde v}_{k}^{-}) \in {\mathcal F}_{\delta},$ 
and by Lemma~\ref{compactness}, 
after passing to a subsequence which we continue to denote 
$\{k\}$, $({\widetilde v}_{k}^{+}, {\widetilde v}_{k}^{-}) \to {v} = ({v}^{+}, {v}^{-}) 
\in {\mathcal F}_{\delta}$ where the convergence is in $W^{1,2}(B_{\s}(0))$  for 
every $\s \in (0, 3/2).$ By (\ref{sp1}), 

\begin{equation}\label{sp2}
\int_{B_{3/2}(0) \setminus B_{3/4}(0)} \left(\frac{\partial ({\widetilde v}_{k}^{+}/R)}{\partial R}\right)^{2}
+ \left(\frac{\partial ({\widetilde v}_{k}^{-}/R)}{\partial R}\right)^{2} 
\leq \frac{1}{k} \left(\frac{3}{2}\right)^{n-2}.
\end{equation}

We claim that ${v}$ cannot be identically equal to zero on $B_{1}(0)$. To see this, first note that 
for any $r, s \in (3/4, 3/2)$ and $\omega \in {\mathbf S}^{n-1}$, we have that  

\begin{equation*}
\left|\frac{{\widetilde v}_{k}(r \, \omega)}{r} - \frac{{\widetilde v}_{k}(s \, \omega)}{s}\right| 
=\left|\int_{s}^{r} \frac{\partial \, ({\widetilde v}_{k}(R \, \omega)/R)}{\partial \, R} dR \right|
\leq \int_{3/4}^{3/2}\left| \frac{\partial \, ({\widetilde v}_{k}(R \, \omega)/R)}{\partial \, R}\right| dR
\end{equation*}

\noindent
which implies, by the triangle inequality, Cauchy-Schwarz inequality and the 
fact that $r, s \in (3/4, 3/2)$ that

\begin{equation*}
|{\widetilde v}_{k}(r \, \omega)|^{2} \leq c\left(|{\widetilde v}_{k}(s \, \omega)|^{2}
+ \int_{3/4}^{3/2} R^{n-1}\left| \frac{\partial \, ({\widetilde v}_{k}(R \, \omega)/R)}{\partial \, R}\right|^{2} dR\right)
\end{equation*}

\noindent
where $c = c(n) \in [1, \infty).$ Integrating this with respect to $\omega$ yields

\begin{equation*}
\int_{{\mathbf S}^{n-1}} |{\widetilde v}_{k}(r \, \omega)|^{2} d\omega 
\leq c \left(\int_{{\mathbf S}^{n-1}} |{\widetilde v}_{k}(s \, \omega)|^{2} d\omega 
+ \int_{B_{3/2}(0) \setminus B_{3/4}(0)} \left|\frac{\partial \, ({\widetilde v}_{k}/R)}{\partial \, R}\right|^{2}\right)
\end{equation*}

\noindent
where $c = c(n) \in [1, \infty).$ First multiplying both sides of the above by $r^{n-1}$ and integrating with 
respect to $r$ over
the interval $(3/4, 3/2),$ and then multiplying both sides of the resulting inequality by $s^{n-1}$ and 
integrating it with respect to $s$ over the interval $(3/4, 1)$ gives

\begin{equation*}
\int_{B_{3/2}(0) \setminus B_{3/4}(0)}|{\widetilde v}_{k}|^{2}  
\leq c \left(\int_{B_{1}(0) \setminus B_{3/4}(0)} |{\widetilde v}_{k}|^{2} 
+ \int_{B_{3/2}(0) \setminus B_{3/4}(0)} \left|\frac{\partial \, ({\widetilde v}_{k}/R)}{\partial \, R}\right|^{2}\right)
\end{equation*}

\noindent
where $c = c(n) \in [1, \infty)$. Since 
$\|{\widetilde v}_{k}\|^{2}_{L^{2}(B_{3/2}(0))} = \left(\frac{3}{2}\right)^{n+2}$, this implies that

\begin{equation*}
\left(\frac{3}{2}\right)^{n+2} \leq c \left(\int_{B_{1}(0)} |{\widetilde v}_{k}|^{2}  
+ \int_{B_{3/2}(0) \setminus B_{3/4}(0)} \left|\frac{\partial \, ({\widetilde v}_{k}/R)}{\partial \, R}\right|^{2}\right)
\end{equation*}

\noindent
which in view of (\ref{sp2}) immediately implies that $v \not\equiv 0$ in $B_{1}(0).$ Hence, 
by Lemma~\ref{wproperties}, part (8), $\int_{\partial \, B_{\r}(0)} |v|^{2} > 0$ for all
$\r \in (0, 3/2).$ By (\ref{sp2}) again, 
${v}$ is homogeneous of degree one in the region $B_{1} \setminus B_{3/4}$ which implies
that $N_{{v},0}(\rho) = 1$ for $3/4 \leq \rho \leq 1.$ (This can be seen easily by the expression 
$N_{{v}, 0}(\r) 
= \frac{\r\frac{d}{d\r}\int_{{\bf S}^{n-1}}|{\hat v}_{\r}|^{2}}{2\int_{{\bf S}^{n-1}}|{\hat v}_{\r}|^{2}}$ 
where ${\hat v}_{\r}(x) = {v}(\r x)$).  By Lemma~\ref{freqlowerbd}, ${\mathcal N}_{{\widetilde v}_{k}}(0) 
\geq 1$, and hence by Lemma~\ref{usc}, we have that ${\mathcal N}_{v}(0) \geq 1.$ Hence by monotonicity of 
$N_{{v}, 0}(\cdot),$ it follows that 
$N_{{v}, 0}(\r) = 1$ for every $\r \in (0, 1).$ By Lemma~\ref{constfreq}, 
this means that $v$ is homogeneous of degree 1 from the origin, and hence 
by lemma~\ref{homogeneous}, ${\rm graph} \, {v}^{+} \cup 
{\rm graph} \, {v}^{-} = P_{1} \cup P_{2}$ for hyperplanes $P_{1}, P_{2}.$  
Thus, if ${v}^{+}$ is not identically equal to ${v}^{-}$, by Lemma~\ref{blowupregularity}, 
for sufficiently large $k$, $v_{k}^{+} = \mbox{max} \, \{v_{k}^{1}, v_{k}^{2}\}$ and 
$v_{k}^{-} = \mbox{min} \, \{v_{k}^{1}, v_{k}^{2}\}$ in $B_{\k}(0)$ for some $\k >0$, where 
$v_{k}^{1}, \, v_{k}^{2}$ are 
harmonic functions in $B_{\k}(0)$, each equal to zero at the origin, and with the 
difference $v_{k}^{1} - v_{k}^{2}$ having vanishing order at the origin equal to 1. But this contradicts 
either of the hypotheses that $v_{k}$ has a branch point at $0$ or that ${\mathcal N}_{w_{k}}(0) >1$. 
Thus we must have that ${v}^{+} \equiv {v}^{-} \equiv L$ for some linear function $L$. But then 
since $D{\widetilde h}_{k}(0) = 0$ for 
every $k$, where ${\widetilde h}_{k} = \frac{1}{2}({\widetilde v}^{+}_{k} + 
{\widetilde v}^{-}_{k}),$ and ${\widetilde h}_{k} \to \frac{1}{2}(v^{+} + v^{-})$ smoothly
in $B_{1}(0)$ (since ${\widetilde h}_{k}$ are harmonic with uniformly bounded 
$L^{2}(B_{3/2}(0))$ norm), $L$ would have to be identically zero, which is impossible. 
The lemma is thus proved.
\end{proof}

\medskip

\begin{lemma} \label{branchptdecay}
Let $v = (v^{+}, v^{-}) \in {\mathcal F}_{\delta}$, $v^{+}(0) = v^{-}(0) = 0$, and suppose either that the origin  
is a branch point of $v$ or that ${\mathcal N}_{w}(0) > 1$. Then

$$\r^{-n-2}\int_{B_{\r}(0)} (v^{+} - l)^{2} + (v^{-} - l)^{2} \leq C\r^{\nu}\int_{B_{1}(0)} (v^{+})^{2} +
(v^{-})^{2}$$

\noindent
for some linear function $l \, : \, {\mathbf R}^{n} \times \{0\}\to {\mathbf R}$ and all $\r \in (0, 1/16)$. In fact, 
$l(x) = Dh(0) \cdot x$ where $h = \frac{1}{2}(v^{+}+ v^{-}).$ Here 
$C = C(n, \d) \in (0, \infty)$ and $\nu = \nu(n, \d) \in (0, 1).$\\
\end{lemma}

\begin{proof}
Let $l(x) = Dh(0)\cdot x.$ Note that there exists $C = C(n)$ such that 

\begin{equation}\label{branchpt0}
|Dh(0)| \leq C\left(\int_{B_{1}(0)} (v^{+})^{2} + (v^{-})^{2}\right)^{1/2} \leq C.
\end{equation}
 
By the definition of ${\mathcal F}_{\delta}$, there exists a sequence of hypersurfaces $M_{k} \in {\mathcal I}_{b}$
and a sequence of affine hyperplanes $L_{k} \to {\mathbf R}^{n} \times \{0\}$ such that 
the blow-up of $\{M_{k}\}$ by the height excesses 
${\hat E}_{k}$ of $M_{k}$ relative to $L_{k}$ (as in Section~\ref{blowupprop}) is $(v^{+}, v^{-}).$ For each $k$, let 
$l_{k} \, : \, {\mathbf R}^{n} \times \{0\} \to {\mathbf R}$ be the affine function such that 
$L_{k} = {\rm graph} \, l_{k}$. Let 
$\left(v^{+}_{(l)}, v^{-}_{(l)}\right)$ be the blow-up produced by blowing up the $M_{k}$'s  by their 
height excesses ${\hat E}_{k}^{(l)}$ relative to the affine hyperplanes given by 
${\rm graph} \, (l_{k} + {\hat E}_{k}l).$ Since by (\ref{branchpt0}),  $\frac{{\hat E}_{k}^{(l)}}{{\hat E}_{k}} \leq C,$
where $C = C(n) < \infty$, we have that

\begin{equation*}\label{branchpt00}
C_{l}\left(v^{+}_{(l)}, v^{-}_{(l)}\right) = (v^{+} - l, v^{-} - l)
\end{equation*}

\noindent
where $0 < C_{l} \leq C = C(n) < \infty$. (Note that here we are assuming 
that not both $v^{+}$, $v^{-}$ are identical to $l$; if this were the case, the lemma is trivially 
true.) It then follows that since $\left(v^{+}_{(l)}, v^{-}_{(l)}\right) \in {\mathcal F}_{\delta}$ 
(by the definition of ${\mathcal F}_{\delta}$), 
all the properties and estimates we have established for $(v^{+}, v^{-})$ will hold with 
$v^{\pm} - l$ in place of $v^{\pm}.$ In particular, 
Lemma~\ref{radialexcessbd} (with $z=0$, $y=0$) holds with $v^{+} - l$, $v^{-} - l$ 
in place of $v^{+}$, $v^{-}$. Thus 

\begin{equation}\label{branchpt1}
\int_{B_{\rho/2}(0)} R^{2-n}\left(\frac{\partial\left(v^{+}/R\right)}{\partial R}\right)^{2} 
+ R^{2-n}\left(\frac{\partial\left(v^{-}/R\right)}{\partial R}\right)^{2} 
\leq C\rho^{-n-2} \int_{B_{\rho}(0)} (v^{+}  - l)^{2} + (v^{-} - l)^{2}
\end{equation}

\noindent
for all $\r \in (0, 1/8)$, where $C = C(n).$ On the other hand, applying Lemma~\ref{spherical} with 
$({\widetilde V}^{+}, {\widetilde V}^{-}) \equiv 
\frac{(v^{+}_{\r} - l, \, v^{-}_{\r} - l)}{\left(\r^{-n-2}\int_{B_{\r}(0)}(v^{+} - l)^{2} 
+ (v^{-} - l)^{2}\right)^{1/2}} \in {\mathcal F}_{\delta}$ 
in place of $(v^{+}, v^{-})$, where $v^{\pm}_{\r}(x) = \frac{1}{(2\r/3)}v^{\pm}(\frac{2}{3}\r x)$  (noting 
that, by definition of $l$, $D{\widetilde H}(0) = 0$ 
where ${\widetilde H} = \frac{1}{2}({\widetilde V}^{+} + {\widetilde V}^{-})$), 
we have that 

\begin{equation*} \label{branchpt2}
\int_{B_{2\r/3}(0) \setminus B_{\r/3}(0)} \left(\frac{\partial (v^{+}/R)}{\partial R}\right)^{2}
+ \left(\frac{\partial (v^{-}/R)}{\partial R}\right)^{2} \geq c\r^{-4}\int_{B_{2\r/3}(0)} (v^{+} - l)^{2} + 
(v^{-} - l)^{2}
\end{equation*}

\noindent
where $c = c(n, \d) > 0$, which gives, since $R = |X| \leq 2\r/3$ for $X \in B_{2\r/3}(0)$, that 

\begin{equation} \label{branchpt3}
\int_{B_{2\r/3}(0) \setminus B_{\r/3}(0)} R^{2-n}\left(\frac{\partial (v^{+}/R)}{\partial R}\right)^{2}
+ R^{2-n}\left(\frac{\partial (v^{-}/R)}{\partial R}\right)^{2} \geq c\r^{-n-2}\int_{B_{2\r/3}(0)} (v^{+} - l)^{2} + 
(v^{-} - l)^{2}.\\
\end{equation}

Replacing $\r$ with $3\r/2$ in the inequalities (\ref{branchpt1}) and (\ref{branchpt3}), and combining them 
gives 

\begin{eqnarray*} \label{branchpt4} 
\int_{B_{\r}(0) \setminus B_{\r/2}(0)} R^{2-n}\left(\frac{\partial (v^{+}/R)}{\partial R}\right)^{2}
&+& R^{2-n}\left(\frac{\partial (v^{-}/R)}{\partial R}\right)^{2} \nonumber\\ 
&\geq& 
\frac{c}{C}\int_{B_{\r/2}(0)} R^{2-n}\left(\frac{\partial (v^{+}/R)}{\partial R}\right)^{2}
+ R^{2-n}\left(\frac{\partial (v^{-}/R)}{\partial R}\right)^{2} 
\end{eqnarray*}
 
\noindent
for $\r \in (0, 1/12).$ This implies that 

\begin{equation*} \label{branchpt5}
\int_{B_{\r/2}(0)} R^{2-n}\left(\frac{\partial (v^{+}/R)}{\partial R}\right)^{2}
+ R^{2-n}\left(\frac{\partial (v^{-}/R)}{\partial R}\right)^{2} \leq 
\kappa \int_{B_{\r}(0)} R^{2-n}\left(\frac{\partial (v^{+}/R)}{\partial R}\right)^{2}
+ R^{2-n}\left(\frac{\partial (v^{-}/R)}{\partial R}\right)^{2}
\end{equation*}

\noindent
where $\kappa =  \kappa(n, \d) = \frac{1}{1 + \frac{c}{C}} \in (0, 1)$
(here $C$ and $c$ are as in inequalities (\ref{branchpt1}) and (\ref{branchpt3})). By iterating this 
starting with $\r = 1/16$, we obtain that 

\begin{eqnarray*} \label{branchpt6}
\int_{B_{\frac{2^{-j}}{16}}(0)} R^{2-n}\left(\frac{\partial (v^{+}/R)}{\partial R}\right)^{2}
&+& R^{2-n}\left(\frac{\partial (v^{-}/R)}{\partial R}\right)^{2} \nonumber\\
&\leq&  \kappa^{j} \int_{B_{1/16}(0)} R^{2-n}\left(\frac{\partial (v^{+}/R)}{\partial R}\right)^{2}
+ R^{2-n}\left(\frac{\partial (v^{-}/R)}{\partial R}\right)^{2}
\end{eqnarray*}

\noindent
for every $j = 1, 2 ,\ldots.$  Combining this with inequalities (\ref{branchpt1}) and 
(\ref{branchpt3}), we have

\begin{equation} \label{branchpt7}
\left(\frac{2^{-j}}{16}\right)^{-n-2}\int_{B_{\frac{2^{-j}}{16}}(0)} (v^{+} - l)^{2} + 
(v^{-} - l)^{2} \leq C\k^{j}\int_{B_{1/8}(0)} (v^{+}  - l)^{2} + (v^{-} - l)^{2}
\end{equation}

\noindent 
for all $j$. Now given any $\r \in (0, 1/16)$, there exists a unique non-negative integer $j$ such that 
$\frac{2^{-j-1}}{16} \leq \r < \frac{2^{-j}}{16}$, and using (\ref{branchpt7}) with this $j$ gives

\begin{eqnarray*} \label{branchpt9}
\r^{-n-2}\int_{B_{\r}(0)} (v^{+} - l)^{2} + (v^{-} - l)^{2} &\leq& 
C\r^{\nu} \int_{B_{1/8}(0)} (v^{+} - l)^{2} + (v^{-} - l)^{2} \nonumber\\
&\leq& C\r^{\nu} \int_{B_{1}(0)} (v^{+})^{2} + (v^{-})^{2} 
\end{eqnarray*}

\noindent
for all $\r \in (0, 1/16),$ where $C = C(n, \d) \in (0, \infty)$ and $\n = \n(n, \d) \in (0, 1).$  (The last inequality 
in (\ref{branchpt9}) follows from (\ref{branchpt0}).) This is the desired estimate.
\end{proof}

\medskip

The preceding two lemmas imply the existence of a fixed positive ``frequency gap'' 
for the functions in ${\mathcal F}_{\delta}.$ Specifically,
we have the following:\\
  
\begin{lemma}\label{freqgap}
Let $\d \in (0, 1).$ There exists a fixed constant $\nu_{0}>0$ depending only on $n$ and $\d$ such that 
if $v = (v^{+}, v^{-}) \in {\mathcal F}_{\delta},$ 
$w = \frac{1}{2}(v^{+} - v^{-}),$ $z \in Z_{w} \cap B_{1/2}(0)$ and either 
${\mathcal N}_{w}(z) > 1$ or $z$ is a branch point of $v$, then ${\mathcal N}_{w}(z) \geq 1 + \nu_{0}.$\\ 
\end{lemma}

\medskip

\begin{proof}
Recall that $N_{w,z}(\r) = \frac{\r\frac{d}{d\r}\int_{{\bf S}^{n-1}} w^{2}_{z, \r}}{2\int_{{\bf S}^{n-1}}w^{2}_{z, \r}}.$ 
Fix any $\r \in (0, 1/2).$ Then we have by the monotonicity of $N_{w, z}(\cdot)$ that 
$\frac{\s\frac{d}{d\s}\int_{{\bf S}^{n-1}} w^{2}_{z, \s}}{2\int_{{\bf S}^{n-1}}w^{2}_{z, \s}} \leq N_{w, z}(\r)$
for all $\s \in (0, \r].$ Integrating this differential inequality (cf. Lemma~\ref{coarseexcessbounds}) gives

\begin{equation}\label{freqgap-1}
\s^{-n}\int_{B_{\s}(z)}w^{2} \geq \left(\r^{1-n-2N_{w, z}(\r)}\int_{\partial B_{\r}(z)} w^{2}\right)\s^{2N_{w,z}(\r)}
\end{equation}

\noindent
for all $\s \in (0, \r].$ On the other hand, Lemma~\ref{branchptdecay}, applied with 
$v^{\pm}(z + (\cdot))$ in place of $v^{\pm}(\cdot)$, implies that 

\begin{equation} \label{freqgap-2}
\s^{-n-2}\int_{B_{\s}(z)}w^{2} \leq C\s^{\n}
\end{equation}

\noindent
for all $\s \in (0, 1/8).$ The estimates (\ref{freqgap-1}) and (\ref{freqgap-2}) readily imply that 

\begin{equation*}\label{freqgap-3}
N_{w, z}(\r) \geq 1 + \frac{\n}{2}
\end{equation*}

\noindent
for all $\r \in (0, 1/2).$ This gives ${\mathcal N}_{w}(z) \geq 1 + \frac{\nu}{2}.$
\end{proof}

\medskip

\begin{lemma}\label{freq1}
Let $(v^{+}, v^{-}) \in {\mathcal F}_{\delta}$ and $w = \frac{1}{2} (v^{+} - v^{-}).$ Suppose $z \in B_{1}(0)$ and 
$v^{+}(z) = v^{-}(z).$  If  
${\mathcal N}_{w}(z) = 1$, then there exists $\s = \s(z) > 0$ and two harmonic functions 
$v^{1}, \, v^{2} \, : \, B_{\s}(z) \to {\mathbf R}$ such that 
$\left. v^{+}\right|_{B_{\s}(z)} = \mbox{max} \, \{v^{1}, v^{2}\}$ and 
$\left. v^{-}\right|_{B_{\s}(z)} = \mbox{min} \, \{v^{1}, v^{2}\}.$\\
\end{lemma}  

\begin{proof}
This follows immediately from Lemma~\ref{freqgap} and the definition of branch point.
\end{proof} 

\bigskip

\begin{proof}[\bf Proof of Theorem~\ref{regblowupthm}]
Let $v = (v^{+}, v^{-}) \in {\mathcal F}_{\delta}$ and $w = \frac{1}{2}(v^{+} - v^{-}).$ 
Let $S_{v} = \{z \in B_{1}(0) \, : \, z \hspace{.1in} \mbox{is a branch point of $v$}\}$. Then $S_{v}$ is a relatively closed 
subset of $B_{1}(0)$ by definition. Also by the definition of $S_{v}$, if $z \in B_{1}(0) \setminus S_{v}$, then 
the graphs of $v^{\pm}$ decompose, locally near $z$, as the union of the graphs of 
two harmonic functions, and hence, the same is true
over any open, simply connected subset $\Omega \subseteq B_{1}(0) \setminus S_{v}.$ This proves part $(a)$ of the 
lemma.\\ 

Part $(b)$ follows by applying Lemma~\ref{branchptdecay} to the function 
${\widetilde v}_{z, \, \frac{3}{8}}$ (notation  as in \ref{blowupnotation-1})
and changing variables. Note that ${\widetilde v}_{z, \, \frac{3}{8}} \in {\mathcal F}_{\delta}.$
\end{proof}

\bigskip

\section{Improvement of excess relative to pairs of hyperplanes}\label{excessdecay}
\setcounter{equation}{0}

In this section, we prove the main excess decay lemma (Lemma~\ref{mainlemma1} below) needed for the proof of 
Theorem~\ref{maintheorem}. Roughly speaking, this lemma says that 
whenever a hypersurface $M \in {\mathcal I}_{b}$ satisfying $0 \in {\overline M}$ and 
$\frac{{\mathcal H}^{n} \, (M \cap (B_{1}(0) \times {\mathbf R}))}{\omega_{n}} \leq 3- \d$ 
for some fixed $\d \in (0, 1)$ is sufficiently $L^{2}$-close, in the cylinder $B_{1}(0) \times {\mathbf R},$  
to a pair of affine hyperlanes of ${\mathbf R}^{n+1}$---i.e. has small height excess relative to 
a pair of affine hyperplanes---then, at one of three possible smaller scales, 
it is closer by a fixed factor to a new pair of affine hyperplanes; i.e. the height excess improves. 
By iterating this result, we shall prove in the next section our main regularity theorem, Theorem~\ref{maintheorem}.
The principal quantity we are interested in keeping track of that measures the height excess of $M$ at 
scale $\r \in (0, 1)$ and that is improving 
is $E_{M} \, (\r, P) \equiv \sqrt{\r^{-n-2}\int_{M \cap (B_{\r}(0) \times {\mathbf R})} {\rm dist}^{2} \, (x, P)},$ where $P$ 
denotes a 
pair of affine hyperplanes. However, in the proof of Lemma~\ref{mainlemma1} (see case $(a)$ of the proof), 
we need to make sure that the ``sheets'' of $M$  separate whenever this excess is significantly
smaller than a certain ``coarse excess,'' which measures the $L^{2}$ deviation of $M$ from 
a {\it single} affine hyperplane. In order to achieve this, it is necessary to modify the definition of 
the improving quantity and consider the sum of $E_{M}^{2} \, (\r, P)$ and a quantity that 
measures the squared $L^{2}$-distance {\it of $P$ from $M$} (see the statement of Lemma~\ref{mainlemma1} for 
the precise definition of this quantity). The main point that necessitates this is simply that smallness of 
$E_{M} \, (\r, P)$ alone need not imply separate closeness of the ``individual sheets'' of $M$ to 
each of the two affine hyperplanes that make up $P$; $M$ may consist of two sheets both of 
which are close to the same single affine hyperplane of $P$.\\ 

In the proof of Lemma~\ref{mainlemma1}, we shall need the elementary facts 
asserted in Lemmas~\ref{lowerupperbounds} and \ref{separation} below. But first we need to recall/introduce 
some notation we shall use in this section and the next. The purpose of the items (1) through (6) below 
is to fix notation that will enable us to define in a convenient way 
the ``second term'' of the improving quantity of 
Lemma~\ref{mainlemma1} referred to in the preceding paragraph, and facilitate statement and 
proof of Lemma~\ref{mainlemma1}.\\

Fix $\d \in (0, 1).$ Let $\r \in (0, 1]$, $M \in {\mathcal I}_{b}$ and suppose that $0 \in {\overline M}$ 
and $\frac{{\mathcal H}^{n} \, (M \cap (B_{\r}(0) \times {\mathbf R})}{\omega_{n}\r^{n}}
\leq 3 -\d.$\\ 

\begin{itemize}
\item[(1)] ${\mathcal A}(M, \r)$ denotes the set of 
affine hyperplanes $L$ of ${\mathbf R}^{n+1}$ 
satisfying $L \cap (B_{1}(0) \times {\mathbf R}) \subset \{ (x^{\prime},x^{n+1}) \in {\mathbf R}^{n+1} \, : \, 
|x^{n+1}| \leq 1/8\}$ and

$${\hat E}_{M}(\r, L) \equiv \r^{-n-2}\int_{M \cap (B_{\r}(0) \times {\mathbf R})}
{\rm dist}^{2} \, (x, L) \leq \frac{3}{2}\inf_{L^{\prime}} \, \r^{-n-2}\int_{M \cap (B_{\r}(0) \times {\mathbf R})}
{\rm dist}^{2} \, (x, {L^{\prime}}),$$ 

\noindent
where the ${\rm inf}$ is taken over all affine hyperplanes $L^{\prime}$ of ${\mathbf R}^{n+1}$ 
satisfying $L^{\prime} \cap (B_{1}(0) \times {\mathbf R}) \subset 
\{(x^{\prime}, x^{n+1}) \in {\mathbf R}^{n+1} \, : \, |x^{n+1}| \leq 1/8\}.$\\

\item[(2)] Given an affine hyperplane $L$ of ${\mathbf R}^{n+1}$ with 
 $L \cap (B_{1}(0) \times {\mathbf R}) \subset \{|x^{n+1}| \leq 1/8\},$ let ${\mathcal R} \, (M, L, \r)$ 
denote the set of regular values $t \in (1/4, 1/2)$ of the 
function $g(X) \equiv 1 - (\nu(X) \cdot \nu^{L})^{2}$ on $M$, satisfying

$${\mathcal H}^{n-1} \, (M \cap (B_{3\r/4}(0) \times {\mathbf R}) \cap \{X \, : \, g(X) = t\}) 
\leq C{\hat E}_{M}^{2} \, (\r, L)$$

\noindent
where $\nu, \nu^{L}$ are the unit normals to $M$, $L$ respectively, and 
$C = C(n)$ is the constant as in inequality (\ref{schoenest}). Note that ${\mathcal R}(M, L, \r)$
contains infinitely many numbers (see the argument of \cite{SS}, p. 753.)\\

\item[(3)] Given affine hyperplane $L$ of ${\mathbf R}^{n+1}$ with 
$L \cap (B_{1}(0) \times {\mathbf R}) \subset \{|x^{n+1}| \leq 1/8\}$ and $t \in {\mathcal R}(M, L, \r)$, and assuming 
${\hat E}_{M} \, (\r, L) \leq \overline\e$ where $\overline\e  = \e(n) \in (0, 1)$ is a sufficiently small fixed constant depending
only on $n$, let $G_{M}^{(L, \, t)}(\r)$ denote the graphical part, relative to $L$, of 
$M \cap q_{L}^{-1} \, (B_{3\r/4}(0) \times {\mathbf R})$  chosen in the sense of \cite{SS}. 
(See item $(2)$ of the discussion at the beginning of Section~\ref{blowupprop}.) Here 
$q_{L}$ denotes a rigid motion of ${\mathbf R}^{n+1}$ with $q_{L}(a_{L}) = 0$ and 
$q_{L}(L) = {\mathbf R}^{n} \times \{0\}$ where $a_{L}$ is the nearest point of $L$ to $0 \in {\mathbf R}^{n+1}.$ 
Thus, for any given radius $\r \in (0, 1]$ and 
choices of an affine hyperplne $L$ with $L \cap (B_{1}(0) \times {\mathbf R}) \subset 
\{|x^{n+1}| \leq 1/8\}$ and $t \in {\mathcal R}(M, L, \r)$, provided 
${\hat E}_{M} \, (\r, L) \leq \overline\e$, $G_{M}^{(L, \, t)}(\r)$ is uniquely determined, and is the union of two 
Lipschitz graphs over a domain $\subset L$ with Lipschitz constants $\leq 1$, and moreover, 

\begin{equation}\label{prelim}
{\mathcal H}^{n} \, ((M \setminus G_{M}^{(L, \, t)}(\r)) \cap q_{L}^{-1} \, (B_{3\r/4}(0) \times {\mathbf R})) 
\leq C({\hat E}_{M}(\r, L))^{2 + \mu},
\end{equation}

\noindent
where $C = C(n)$ and $\mu = \mu(n)$ are fixed positive constants depending only on $n$. (However, we remark here
that in the proof of Lemma~\ref{mainlemma1}, we do not need such precise control of the size of the complement 
of $G_{M}^{(L, \, t)}(\r)$ as is given by the estimate (\ref{prelim}); all we need is that $G_{M}^{(L, \, t)}(\r)$ 
has $n$-dimensional measure larger than a fixed fraction of the measure of $B_{\r/2}(0).$ See 
Lemma~\ref{lowerupperbounds} below.)\\ 

\item[(4)] Given affine hyperplanes $L$, $U$ of ${\mathbf R}^{n+1}$ with
$L \cap (B_{1}(0) \times {\mathbf R}) \subset 
\{|x^{n+1}| \leq 1/8\}$, $U \cap (B_{1}(0) \times {\mathbf R}) \subset \{|x^{n+1}| \leq 1/8\}$ such 
that ${\hat E}_{M} \, (\r, L) \leq \overline\e$ ($\overline\e$ as in $(3)$ above), and 
$t \in {\mathcal R}(M, L, \r)$, let 
$$U^{\star} \, (M, L, t, \r) = {U} \cap \pi^{-1} \, (\pi \, G_{M}^{(L, t)}(\r)).$$ 

\noindent
Given $L$ as above and a pair of affine hyperplanes $P = P_{1} \cup P_{2}$ of ${\mathbf R}^{n+1}$ 
(with $P_{1}$, $P_{2}$ affine hyperplanes of ${\mathbf R}^{n+1}$) such that 
$P \cap (B_{1}(0) \times {\mathbf R}) \subset \{|x^{n+1}| \leq 1/8\}$, define 

$$P^{\star}(M, L, t, \r)  = P_{1}^{\star} \, (M, L, t, \r) \cup P_{2}^{\star} \, (M, L, t, \r).$$\\

\noindent
\item[(5)] If $P = P_{1} \cup P_{2}$ is a pair of affine hyperplanes of ${\mathbf R}^{n+1}$ (with 
$P_{1}$, $P_{2}$ affine hyperplanes) such that 
$P \cap (B_{1}(0) \times {\mathbf R}) \subset \{|x^{n+1}|<1/8\}$,
we set, for $\t \in (0, 1/2)$, $S_{P}(\t) = \{x \in {\mathbf R}^{n} \times \{0\} 
\, : \, {\rm dist}\,(x, {\pi} \, (P_{1} \cap P_{2})) \leq \t\}$ if $P_{1}$ and $P_{2}$ 
are distinct with ${\pi} \, (P_{1} \cap P_{2}) \cap 
B_{1/8}(0) \neq \emptyset,$ and $S_{P}(\t) = \emptyset$ otherwise.\\

\noindent
\item[(6)] If $U$ is an affine hyperplane of ${\mathbf R}^{n+1}$, we shall denote by $U^{T}$ the hyperplane 
obtained by translating $U$ parallel to itself. If $P = P_{1} \cup P_{2}$ is a pair of affine hyperplanes, with 
$P_{1}$, $P_{2}$ affine hyperplanes, then we shall let $P^{T} = P_{1}^{T} \cup P_{2}^{T}.$\\  
\end{itemize}

\begin{lemma}\label{lowerupperbounds}
Let $\d \in (0, 1)$. There exist constants $c_{1} = c_{1}(n, \d) \in (0, \infty),$ $c_{2} = c_{2}(n, \d) \in (0, \infty)$
and $\z_{0} = \z_{0}(n, \d) \in (0, 1)$ such that the following is true.
If $M \in {\mathcal I}_{b}$, $\frac{{\mathcal H}^{n} \, (M \cap (B_{1}(0) \times {\mathbf R}))}{\omega_{n}} \leq 
3 - \d$, $L \in {\mathcal A}(M, 1),$ $t \in {\mathcal R}(M, L, 1)$, 
${\hat E}_{M} \, (1, L) \leq 1,$ $P = P^{+} \cup P^{-}$ is a pair of affine hyperplanes with 
${\rm dist}_{\mathcal H} \, (P \cap (B_{1}(0) \times {\mathbf R}), B_{1}(0)) \leq \z_{0}$ and

$$\int_{M \cap (B_{1}(0) \times {\mathbf R})} {\rm dist}^{2} \, (x, P) 
+ \int_{P^{\star} \cap ((B_{1/2}(0) \setminus S_{P}(1/16)) \times {\mathbf R})} {\rm dist}^{2} \, (x, G_{M}^{(L, \, t)}(1)) 
\leq \z_{0} {\hat E}_{M}^{2} \, (1, L)$$

\noindent
where $P^{\star} \equiv P^{\star}(M, L, t, 1)$, then 

$$c_{1}{\hat E}_{M} \, (1, L) \leq {\rm sup}_{B_{1}(0)} \, |p^{+} - p^{-}| \leq c_{2} \, {\hat E}_{M} \, (1, L).$$\\ 
\end{lemma}

\begin{proof}
Note first that it follows from the conditions 

$${\rm dist}_{\mathcal H} \, (P \cap (B_{1}(0) \times {\mathbf R}), B_{1}(0)) 
\leq \z_{0},$$

 $$\int_{M \cap (B_{1}(0) \times {\mathbf R})} {\rm dist}^{2} \, (x, P) \leq \z_{0}$$
 
 \noindent
 and the triangle inequality that 
$\int_{M \cap (B_{1}(0) \times {\mathbf R})} |x^{n+1}|^{2} \leq C\z{0}$, 
where $C = C(n)$, so that by the definition of 
${\mathcal A}(M, 1)$, we have that 

\begin{equation}\label{lower-0}
{\hat E}_{M} \, (1, L) \leq \frac{3}{2}C\z_{0}. 
\end{equation}

\noindent
Thus, if $\z_{0} = \z_{0}(n)$ is sufficiently small, $G_{M}^{(L, \, t)}(1) \neq \emptyset,$ and in fact by the estimate 
(\ref{prelim}), 

\begin{equation}\label{lower-0-1}
{\mathcal H}^{n} \, (G_{M}^{(L, \, t)}(1) \cap (B_{1/2}(0) \times {\mathbf R})) \geq 
\frac{1}{2}\omega_{n} \left(\frac{1}{2}\right)^{n}.
\end{equation}
 
To see the lower bound of the asserted inequalities in the conclusion of the lemma, 
let $U = {\rm graph} \, \frac{1}{2}(p^{+} + p^{-}).$ Then, 
by the definition of ${\hat E}_{M}(1, L)$ and the triangle inequality, we have that 

\begin{eqnarray}\label{lower-1}
\frac{2}{3}{\hat E}_{M}(1, L) &\leq& \int_{M \cap (B_{1}(0) \times {\mathbf R})} {\rm dist}^{2} \, (x, {U})\noindent\\
&\leq& 2\int_{M \cap (B_{1}(0) \times {\mathbf R})} {\rm dist}^{2} \, (x, P) + c \,{\rm sup} \, |p^{+} - p^{-}|
\end{eqnarray}  

\noindent
where $c = c(n).$ Provided we take $\z_{0} < 1/4$, the lower bound follows directly from this since 
$\int_{M \cap (B_{1}(0) \times {\mathbf R})} {\rm dist}^{2} \, (x, P) \leq \z_{0}{\hat E}_{M}(1, L)$ by hypothesis.\\

To see the upper bound, we argue by contradiction. If the assertion is not true, then there exist a sequence 
of hypersurface $M_{k} \in {\mathcal I}_{b},$ $k = 1, 2, 3, \ldots,$ 
with $\frac{{\mathcal H}^{n} \, (M_{k} \cap (B_{1}(0) \times {\mathbf R}))}{\omega_{n}}
\leq 3 - \d$, a sequence of affine hyperplanes $L_{k}$ with $L_{k} \cap (B_{1}(0) \times {\mathbf R}) 
\subset \{|x^{n+1}| \leq 1/8\}$ and

\begin{equation}\label{upper-1}
{\hat E}_{k}^{2}  \equiv \int_{M_{k} \cap (B_{1}(0) \times {\mathbf R})} {\rm dist}^{2} \, (x, L_{k}) 
\leq \frac{3}{2} {\rm inf}_{L^{\prime}} \, \int_{M_{k} \cap (B_{1}(0) \times {\mathbf R})}
{\rm dist} \, (x, {L^{\prime}})
\end{equation}

\noindent
where for each $k$, the ${\rm inf}$ is taken over all affine hyperplanes $L^{\prime}$ satisfying
$L^{\prime} \cap (B_{1}(0) \times {\mathbf R}) \subset \{|x^{n+1}| \leq 1/8\}$, 
a sequence of numbers $t_{k} \in {\mathcal R}(M_{k}, L_{k}, 1)$, 
a sequence $P_{k} = P_{k}^{+} \cup P_{k}^{-}$ of pairs of affine hyperplanes with 

\begin{equation}\label{upper-1-1}
{\rm dist}_{\mathcal H} \, (P_{k} \cap (B_{1}(0) \times {\mathbf R}), B_{1}(0)) \to 0 \hspace{.2in} 
\mbox{as $k \to \infty$ and} 
\end{equation}

\begin{equation}\label{upper-2}
\int_{M_{k} \cap (B_{1}(0) \times {\mathbf R})} {\rm dist}^{2} \, (x, P_{k})
+ \int_{P_{k}^{\star} \cap ((B_{1/2}(0) \setminus S_{P_{k}}(1/16)) \times {\mathbf R})} {\rm dist}^{2} \, 
(x, G_{k}) \leq \frac{1}{k}{\hat E}_{k}^{2},
\end{equation}

\noindent
and yet, for each $k$,  

\begin{equation}\label{upper-3}
{\rm sup}_{B_{1}(0)} \, |p_{k}^{+} - p_{k}^{-}| \geq k {\hat E}_{k}.
\end{equation}

\noindent
Here we are using the abbreviations $G_{k} = G_{M_{k}}^{(L_{k}, t_{k})}(1)$ and 
$P_{k}^{\star} = P_{k}^{\star}(M_{k},L_{k},t_{k},1).$ Note then by 
(\ref{lower-0}), ${\hat E}_{k} \to 0$,  and by (\ref{upper-1-1})
and (\ref{upper-2}), $M_{k} \cap (B_{1}(0) \times {\mathbf R}) \to B_{1}(0) \times \{0\}$ in Hausdorff 
distance. Consequently, $L_{k} \to {\mathbf R}^{n} \times \{0\}.$ Note also that by 
(\ref{lower-0-1}), 

\begin{equation}\label{upper-3-1}
{\mathcal H}^{n} \, (G_{k} \cap (B_{1/2}(0) \times {\mathbf R})) \geq \frac{1}{2}\omega_{n} \left(\frac{1}{2}\right)^{n}
\end{equation}

\noindent
for all sufficiently large $k.$ Let $v \in L^{2}(B_{1}(0); {\mathbf R}^{2}) \cap 
W^{1, 2}_{\rm loc} (B_{1}(0); {\mathbf R}^{2})$ be the blow-up, in the 
sense of Section~\ref{blowupprop}, of $M_{k}$ by ${\hat E}_{k}$. In view of Proposition~\ref{properties}, part (2), it follows
from (the bound on the first term on the left hand side of) (\ref{upper-2}) and (\ref{upper-3}) that 
$v^{+} \equiv v^{-} \equiv l$ for some affine function $l$. 
Indeed, if we write $P_{k} = P_{k}^{(1)} \cup P_{k}^{(2)}$ where $P_{k}^{(1)}, P_{k}^{(2)}$ are affine 
hyperplanes, and define functions $p_{k}^{(1)}, p_{k}^{(2)} \, : \, {\mathbf R}^{n} \times \{0\} \to {\mathbf R}$ by 
$P_{k}^{(i)} = {\rm graph} \, p_{k}^{(i)}$, $i = 1, 2$, then, after 
possibly passing to a subsequence, $l = \lim_{k \to \infty} \,(p_{k}^{(1)} - \phi_{k})/{\hat E}_{k}$ 
or $l = \lim_{k \to \infty} \, (p_{k}^{(2)} - \phi_{k})/{\hat E}_{k}$ where 
$\phi_{k} \, : \, {\mathbf R}^{n} \times \{0\} \to {\mathbf R}$ is such that 
$L_{k}  = {\rm graph} \, \phi_{k}.$  
(The existence of one of these two limits, is guaranteed by Lemma~\ref{properties}, 
part (2) and the bound on the first term on the left hand side of (\ref{upper-2}).) By relabeling if necessary, we 
assume that 
$l = \lim_{k \to \infty} \, (p_{k}^{(1)} - \phi_{k})/{\hat E}_{k}.$ Note then that by (\ref{upper-3}), 

\begin{equation}\label{upper-3-1-1}
\lim_{\k \to \infty} \, (p_{k}^{(2)} - \phi_{k})(x)/{\hat E}_{k}  = \infty
\end{equation}

\noindent
for each $x \in \cup_{j=1}^{\infty} 
\left({\mathbf R}^{n} \times \{0\} \setminus \cap _{k=j}^{\infty} \{x \, : \, p_{k}^{(1)}(x) = p_{k}^{(2)}(x)\}\right)$ 
and that (\ref{upper-2}) in particular says that 

\begin{equation}\label{upper-3-2}
\int_{P_{k}^{(2) \, {\star}} \cap ((B_{1/2}(0) \setminus S_{P_{k}}(1/16)) \times {\mathbf R})} {\rm dist}^{2} \, 
(x, G_{k}) \leq \frac{1}{k}{\hat E}_{k}^{2}.
\end{equation}

\noindent
If we let ${\widetilde L}_{k} = {\rm graph} \, (\phi_{k} + {\hat E}_{k} l)$, we have

\begin{equation}\label{uppper-4}
\int_{M_{k} \cap (B_{1/2}(0) \times {\mathbf R})} {\rm dist}^{2} \, (x, {\widetilde L}_{k}) \leq 
\frac{1}{16}{\hat E}_{k}^{2}
\end{equation}

\noindent
for infinitely many $k$, which implies by the triangle inequality that 

\begin{equation}\label{upper-5}
\int_{G_{k} \cap (B_{1/2}(0) \times {\mathbf R})} {\rm dist}^{2} \, (x, {P}^{(1)}_{k}) \leq \frac{1}{8}{\hat E}_{k}^{2}
\end{equation}

\noindent
for infinitely many $k$. Now let ${\widetilde G}_{k} = \{x \in G_{k} \cap (B_{1/2}(0) \times {\mathbf R}) 
\, : \, {\rm dist} \, (x, P^{(1)}_{k}) \leq \sqrt{\frac{2^{n-1}}{\omega_{n}}} {\hat E}_{k}\}.$ Then by 
(\ref{upper-5}) and (\ref{upper-3-1}), 

\begin{equation}\label{upper-6}
{\mathcal H}^{n} \, ({\widetilde G}_{k}) \geq  \frac{1}{4} \omega_{n} \left(\frac{1}{2}\right)^{n}.
\end{equation}

\noindent
Since $G_{k}$ is the union of two Lipschitz graphs with Lipschitz constants $\leq 1$, 
for any $x = (x^{\prime}, x^{n+1}) \in \pi \, G_{k} \times {\mathbf R}$, 
${\rm dist} \, (x, G_{k})$ is bounded below by a fixed positive constant times the 
``vertical distance'' $\min \, \{|x^{n+1}  - y_{1}^{n+1}|, |x^{n+1} - y_{2}^{n+1}| \, : 
\, (x^{\prime}, y_{1}^{n+1}), (x^{\prime}, y_{2}^{n+1}) \in G_{k}\}$. Moreover, by (\ref{upper-6}),  
${\mathcal H}^{n} \, (P_{k}^{(2) \, {\star}} \cap ((B_{1/2}(0) \setminus S_{P_{k}}(1/16)) \times {\mathbf R})) 
\geq C = C(n) >0$ which, in view of (\ref{upper-3}), contradicts (\ref{upper-3-2}). This completes the proof 
of the lemma.   
\end{proof}

\medskip

\begin{lemma}\label{separation}
Let $\d \in (0, 1)$, $\eta \in (0, 1)$ and $c_{1} \in (0, \infty)$ be given. 
There exists a number $\z = \z(n, \d, \eta, c_{1}) \in (0, 1)$ such that the following holds. If 
$P = P^{+}  \cup P^{-}$ is a pair of affine hyperplanes of ${\mathbf R}^{n+1}$ with 
$\sup_{B_{1}(0)} |p^{+} - p^{-}| \geq c_{1}$ and $v = (v^{+}, v^{-}) \in 
{\mathcal F}_{\delta}$ satisfies 

\begin{eqnarray*}
&&\int_{B_{1}(0)} {\rm dist}^{2} \, ((x, v^{+}(x)), P) + {\rm dist}^{2} \, ((x, v^{-}(x)), P)\nonumber\\
&& \hspace{.5in} + \int_{B_{1/2}(0) \setminus S_{P}(1/8)} {\rm dist}^{2} \, ((x, p^{+}(x)), V) + 
{\rm dist}^{2} \, ((x, p^{-}(x)), V) \leq \z
\end{eqnarray*}

\noindent
where $V = {\rm graph} \, v^{+} \cup {\rm graph} \, v^{-}$, then 

$$\int_{B_{1}(0)} (v^{+} - p^{+})^{2} + (v^{-} - p^{-})^{2} \leq \eta.$$  
\end{lemma}

\begin{proof}
If the assertion is false, then there exist numbers $\d \in (0, 1)$, $\eta \in (0, 1)$, $c_{1} \in (0, \infty),$   
a sequence of functions $v_{k} = (v_{k}^{+}, v_{k}^{-}) \in {\mathcal F}_{\delta}$ and 
a sequence of affine hyperplanes $P_{k} = (P_{k}^{+}, P_{k}^{-})$ of ${\mathbf R}^{n+1}$ such that 

\begin{equation}\label{sep-1}
\sup_{B_{1}(0)} |p_{k}^{+} - p_{k}^{-}| \geq c_{1} \hspace{.3in} \mbox{and}
\end{equation}

\begin{eqnarray}\label{sep-2}
&&\int_{B_{1}(0)} {\rm dist}^{2} \, ((x, v_{k}^{+}(x)), P_{k}) + {\rm dist}^{2} \, ((x, v_{k}^{-}(x)), P_{k}) \nonumber\\
&&\hspace{.5in} + \int_{B_{1/2}(0) \setminus S_{P_{k}}(1/8)} {\rm dist}^{2} \, ((x, p_{k}^{+}(x)), V_{k}) + 
{\rm dist}^{2} \, ((x, p_{k}^{-}(x)), V_{k}) \leq \frac{1}{k}
\end{eqnarray}

\noindent
where $V_{k} = {\rm graph} \, v_{k}^{+} \cup {\rm graph} \, v_{k}^{-},$ and yet 

\begin{equation}\label{sep-3}
\int_{B_{1}(0)} (v_{k}^{+} - p_{k}^{+})^{2} + (v_{k}^{-} - p_{k}^{-})^{2} \geq \eta
\end{equation}

\noindent
for all $k = 1, 2, 3, \ldots.$ After passing to a subsequence, we have by Lemma~\ref{compactness} 
that $v_{k} \to v$ for some $v \in {\mathcal F}_{\delta}$, where the convergence is in 
$W^{1, 2}(B_{1}(0);{\mathbf R}^{2}),$
and that $P_{k} \to P$ for some affine pair of hyperplanes of ${\mathbf R}^{n+1}$ satisfying 
 $\sup_{B_{1}(0)} |p^{+} - p^{-}| \geq c_{1}.$ Note that since $v^{\pm}$ are bounded in $B_{1}(0)$
 (by Proposition~\ref{properties}; part $(3)$ says $|v|^{2}$ is subharmonic in $B_{3/2}(0)$, and the 
 mean value property and part $(2)$ say $|v|^{2}$ is bounded in $B_{1}(0)$) and continuous
 (by Proposition~\ref{lipschitz}), 
(\ref{sep-2}) says that $v^{+} \equiv p^{+}$ and $v^{-} \equiv p^{-}$ on $B_{1}(0).$ This immediately contradicts 
(\ref{sep-3}) for sufficiently large $k$.
\end{proof}

\medskip

\begin{lemma}\label{mainlemma1}
Let $\th \in (0, 1/16),$ $\b \in (0, \th/16)$ and $\g \in (0, \b/16).$ Let $\delta \in (0, 1).$ There exist numbers 
$\e_{0} = \e_{0}(n, \d, \th, \b, \g) \in (0, 1/2)$ and $\lambda = \lambda(n, \d) \in (0, 1)$
such that the following is true. Suppose $M \in {\mathcal I}_{b}$, $0 \in \overline M,$ $\r \in (0, 1],$ 
$$\frac{{\mathcal H}^{n}(M \cap (B_{\r}(0) \times {\mathbf R}))}{\omega_{n}\r^{n}} \leq 3 - \delta \hspace{.3in}
\mbox{and}$$ 

$$\r^{-n-2}\int_{M \cap (B_{\r}(0) \times {\mathbf R})} {\rm dist}^{2}(x, P) 
+\r^{-n-2}\int_{P^{\star} \cap ((B_{\r/2}(0) \setminus S_{P}(\r/16)) \times {\mathbf R})} 
{\rm dist}^{2} \, (x, G_{M}^{(L, \, t)}(\r)) \leq \e_{0}$$

\noindent
for some affine hyperplane $L \in {\mathcal A}(M, \r)$, number 
$t \in {\mathcal R}(M, L, \r)$ and some pair of affine hyperplanes $P$ of ${\mathbf R}^{n+1}$ satisfying
${\rm dist}_{\mathcal H} \, (P \cap (B_{1}(0) \times {\mathbf R}), B_{1}(0)) \leq \e_{0}.$
Here we have used the notation $P^{\star} \equiv P^{\star}(M, L, t, \r).$ 
Then there exists a pair of affine hyperplanes ${\widetilde P},$ an affine hyperplane ${\widetilde L}$
and a number ${\widetilde t} \in (1/4, 1/2)$ such that  

\begin{eqnarray*}
(1) &&\r^{-2}d_{\mathcal H}^{2} \, ({\widetilde P} \cap (B_{\r}(0) \times {\mathbf R}), 
P \cap (B_{\r}(0) \times {\mathbf R})) \nonumber\\
&&\hspace{.3in}\leq C\left(\r^{-n-2}\int_{M \cap (B_{\r}(0) \times {\mathbf R})} {\rm dist}^{2} \, (x, P)
+ \r^{-n-2}\int_{P^{\star} \cap ((B_{\r/2}(0) \setminus S_{P}(\r/16)) \times {\mathbf R})} 
{\rm dist}^{2} \, (x, G_{M}^{(L, \, t)}(\r)) \right),
\end{eqnarray*}

\begin{eqnarray*}
(2) && d_{\mathcal H}^{2} \, ({\widetilde P}^{T} \cap (B_{1}(0) \times {\mathbf R}), 
P^{T} \cap (B_{1}(0) \times {\mathbf R})) \nonumber\\
&&\hspace{.3in}\leq C\left(\r^{-n-2}\int_{M \cap (B_{\r}(0) \times {\mathbf R})} {\rm dist}^{2} \, (x, P)
+\r^{-n-2}\int_{P^{\star} \cap ((B_{\r/2}(0) \setminus S_{P}(\r/16)) \times {\mathbf R})} 
{\rm dist}^{2} \, (x, G_{M}^{(L, \, t)}(\r)) \right)
\end{eqnarray*}

\noindent
and\\

\noindent
{\rm (3)} one of the following options $(A)$, $(B)$ or $(C)$ holds:\\

\begin{itemize}
\item[$(A)$] ${\widetilde L} \in {\mathcal A}(M, \th\r),$ ${\widetilde t} \in {\mathcal R}(M, {\widetilde L}, \th\r),$ 

\begin{eqnarray*}
&&(\th\r)^{-n-2}\int_{M \cap (B_{\th\r}(0) \times {\mathbf R})} {\rm dist}^{2} \, (x, {\widetilde P}) 
 + (\th\r)^{-n-2}\int_{{\widetilde P}^{\star} \cap ((B_{\th\r/2}(0) \setminus S_{\widetilde P}(\th\r/16)) \times {\mathbf R})} 
{\rm dist}^{2} \, (x, G_{M}^{({\widetilde L}, \, {\widetilde t})}(\th\r)) \nonumber\\
&&\hspace{.3in} \leq C_{1}\th^{\lambda} \left(\r^{-n-2}\int_{M \cap (B_{\r}(0) \times {\mathbf R})} {\rm dist}^{2} \, (x, P)
+ \r^{-n-2}\int_{P^{\star} \cap ((B_{\r/2}(0) \setminus S_{P}(\r/16)) \times {\mathbf R})} {\rm dist}^{2} \, (x, G_{M}^{(L, \t)}
(\r))\right),\\
\end{eqnarray*}
\noindent
where ${\widetilde P}^{\star} \equiv {\widetilde P}^{\star}(M, {\widetilde L}, 
{\widetilde t}, \th\r)$, and\\

$$\frac{{\mathcal H}^{n}(M \cap (B_{\th\r}(0) \times {\mathbf R}))}{\omega_{n}(\th\r)^{n}} \leq 3 - \delta.$$ 

\medskip

\item[$(B)$]  ${\widetilde L} \in {\mathcal A}(M, \b\r),$ ${\widetilde t} \in {\mathcal R}(M, {\widetilde L}, \b\r),$ 

\begin{eqnarray*}
&& (\b\r)^{-n-2}\int_{M \cap (B_{\b\r}(0) \times {\mathbf R})} {\rm dist}^{2} \, (x, {\widetilde P}) 
+ (\b\r)^{-n-2}\int_{{\widetilde P}^{\star} \cap ((B_{\b\r/2}(0) \setminus S_{\widetilde P}(\b\r/16)) \times {\mathbf R})} 
{\rm dist}^{2} \, (x, G_{M}^{({\widetilde L}, \, {\widetilde t})}(\b\r)) \nonumber\\
&& \hspace{.3in} \leq C_{2}\b^{\lambda}\left(\r^{-n-2}\int_{M \cap (B_{\r}(0) \times {\mathbf R})} {\rm dist}^{2} \, (x, P)
+ \r^{-n-2}\int_{P^{\star} \cap ((B_{\r/2}(0) \setminus S_{P}(\r/16)) \times {\mathbf R})} 
{\rm dist}^{2} \, (x, G_{M}^{(L, \, t)}(\r))\right),\\
\end{eqnarray*}
\noindent
where ${\widetilde P}^{\star} \equiv {\widetilde P}^{\star}(M, {\widetilde L}, 
{\widetilde t}, \b\r)$, and\\

$$\frac{{\mathcal H}^{n}(M \cap (B_{\b\r}(0) \times {\mathbf R}))}{\omega_{n}(\b\r)^{n}} \leq 3 - \delta.$$ 

\medskip

\item[$(C)$] ${\widetilde L} \in {\mathcal A}(M, \g\r),$ ${\widetilde t} \in {\mathcal R}(M, {\widetilde L}, \g\r),$ 
\begin{eqnarray*}
&&(\g\r)^{-n-2}\int_{M \cap (B_{\g\r}(0) \times {\mathbf R})} {\rm dist}^{2} \, (x, {\widetilde P}) 
+ (\g\r)^{-n-2}\int_{{\widetilde P}^{\star} \cap ((B_{\g\r/2}(0) \setminus S_{\widetilde P}(\g\r/16)) \times {\mathbf R})} 
{\rm dist}^{2} \, (x, G_{M}^{({\widetilde L}, \, {\widetilde t})}(\g\r)) \nonumber\\
&&\hspace{.3in}\leq C_{3}\g^{\lambda}\left(\r^{-n-2}\int_{M \cap (B_{\r}(0) \times {\mathbf R})} {\rm dist}^{2} \, (x, P)
+ \r^{-n-2}\int_{P^{\star} \cap ((B_{\r/2}(0) \setminus S_{P}(\r/16)) \times {\mathbf R})} 
{\rm dist}^{2} \, (x, G_{M}^{(L, \, t)}(\r))\right),\\
\end{eqnarray*}
\noindent
where ${\widetilde P}^{\star} \equiv {\widetilde P}^{\star}(M, {\widetilde L}, 
{\widetilde t}, \g\r)$, and\\

$$\frac{{\mathcal H}^{n}(M \cap (B_{\g\r}(0) \times {\mathbf R}))}{\omega_{n}(\g\r)^{n}} \leq 3 - \delta.$$ 
\end{itemize}

\noindent
Here the dependence of the constants $C, C_{i}$, $i = 1, 2, 3$ on the parameters is as follows: 
$C = C(n, \d, \th, \b, \g)$, $C_{1} = C_{1}(n, \d)$, $C_{2} = C_{2}(n, \d, \th)$ and $C_{3} = C_{3}(n, \d, \th, \b).$\\
\end{lemma}

\begin{proof}
Note first that conclusion (2) follows from conclusion (1). Since the hypotheses and the conclusions of the 
lemma are scale invariant, it suffices to prove the lemma 
assuming $\r = 1$, and we shall make this assumption in what follows.  
Let $\{M_{k}\} \subset {\mathcal I}_{b}$ be an arbitrary sequence of hypersurfaces with 
$0 \in \overline M_{k}$,

\begin{equation}\label{main0-1}
\frac{{\mathcal H}^{n}(M_{k} \cap (B_{1}(0) \times {\mathbf R}))}{\omega_{n}} \leq 3 - \delta,
\end{equation}

\begin{equation}\label{main1}
\int_{M_{k} \cap (B_{1}(0) \times {\mathbf R})} {\rm dist}^{2} \, (x, P_{k}) 
+ \int_{P_{k}^{\star} \cap ((B_{1/2}(0) \setminus S_{P_{k}}(1/16)) \times {\mathbf R})} {\rm dist}^{2} \, (x, G_{k}) 
\searrow 0
\end{equation}

\noindent
for a sequence of affine hyperplanes $L_{k} \in {\mathcal A}(M_{k}, 1),$ a sequence 
of numbers $t_{k} \in {\mathcal R}(M_{k}, L_{k}, 1)$ and 
a sequence of pairs of affine hyperplanes $P_{k} = P_{k}^{1} \cup P_{k}^{2}$ (where 
$P_{k}^{1}, P_{k}^{2}$ are affine hyperplanes, possibly with $P_{k}^{1} \equiv P_{k}^{2}$), 
satisfying

\begin{equation}\label{main0-11}
d_{\mathcal H} \, (P_{k} \cap (B_{1}(0) \times {\mathbf R}), B_{1}(0)) \searrow 0.
\end{equation}

\noindent
Here we use the notation $G_{k} \equiv G_{M_{k}}^{(L_{k}, t_{k})}(1)$ and 
$P_{k}^{\star} \equiv P_{k}^{\star}(M_{k}, L_{k}, t_{k}, 1).$ Note that (\ref{main1}) and 
(\ref{main0-11}) imply that $M_{k} \cap (B_{1}(0) \times {\mathbf R}) \to B_{1}(0) \times \{0\}$ 
in Hausdorff distance, and hence by the triangle inequality that 
$\int_{M_{k} \cap (B_{1}(0) \times {\mathbf R})} |x^{n+1}|^{2} 
\to 0.$ By the definition of ${\mathcal A}(M_{k}, 1)$, it then follows that 

\begin{equation}\label{main-0-2}
{\hat E}_{k} \to 0,
\end{equation}

\noindent
where we use the notation ${\hat E}_{k} \equiv {\hat E}_{M_{k}} \, (1, L_{k})$.
This in turn says that ${\rm dist}_{\mathcal H} \, (M_{k} \cap (B_{1}(0) \times {\mathbf R}), 
L_{k} \cap (B_{1}(0) \times {\mathbf R})) \to 0$, so that $L_{k} \to {\mathbf R}^{n} \times \{0\}.$ 
Note also that (\ref{main-0-2}) in particular implies that for all sufficiently large $k$, 

\begin{equation}\label{main-0-3}
{\mathcal H}^{n} \, (G_{k}) \geq \frac{1}{2} \omega_{n} \left(\frac{1}{2}\right)^{n}
\end{equation}

\noindent
and hence that

\begin{equation}\label{main-0-4}
{\mathcal H}^{n} \, (P_{k}^{\star} \cap (B_{1/2}(0) \setminus S_{P_{k}}(1/16)) \times {\mathbf R})) \geq 
\frac{1}{4} \omega_{n} \left(\frac{1}{2}\right)^{n}.
\end{equation}

We show that for infinitely many $k$, we can find pairs of affine hyperplanes ${\widetilde P}_{k},$ 
affine hyperplanes ${\widetilde L}_{k}$ and numbers ${\widetilde t}_{k} \in (1/4, 1/2)$ such that the conclusions 
of the lemma hold with $M_{k}, P_{k}, L_{k}, t_{k}, {\widetilde P}_{k}, {\widetilde L}_{k}, 
{\widetilde t}_{k}$ in place of $M$, $P$, $L$, $t$, ${\widetilde P}$,  ${\widetilde L}$ and 
${\widetilde t}$ respectively, and with the constants $C, C_{i}$, $i = 1, 2, 3$ and 
$\lambda$ fixed depending only on the specified parameters as in the statement of the lemma. In view 
of the arbitrariness of $\{M_{k}\}$, this will prove the lemma.\\ 

First notice that for any given $\t>0$, we have, since ${\hat E}_{k} \to 0$,  that for all sufficiently large $k$ depending on $\t$, 
${\mathcal H}^{n} \, (G_{k} \cap (B_{\t}(0) \times {\mathbf R})) \to 2\omega_{n}\t^{n}$ and 
${\mathcal H}^{n} \, ((M_{k} \setminus G_{k}) \cap (B_{\t}(0) \times {\mathbf R})) \to 0$, so that the last 
of the conclusions in each of the options $(3)(A)$, $(3)(B)$  and $(3)(C)$ hold with $M_{k}$ in place of $M$ 
for all sufficiently large $k$. It only remains to show that the other conclusions hold with $M_{k},$ $P_{k}$,
$L_{k}$, $t_{k}$
in place of $M$, $P$, $L$, $t$ respectively and with suitable choices 
of ${\widetilde P}_{k}$, ${\widetilde L}_{k}$ and ${\widetilde t}_{k}$, 
in place of ${\widetilde P}$, ${\widetilde L}$ and ${\widetilde t}$ respectively.\\ 

Let $\z = \z(n, \th, \d) \in (0, 1/8)$ be a small number to be determined depending only 
on $n$, $\th$ and $\d$. We divide the rest of the proof of the lemma into two cases according 
to the following two possibilities, one of which must hold for infinitely many $k$:\\

\begin{itemize}
\item[$(a)$] $\int_{M_{k} \cap (B_{1}(0) \times {\mathbf R})} {\rm dist}^{2} \, (x, P_{k}) 
+ \int_{P_{k}^{\star} \cap ((B_{1/2}(0) \setminus S_{P_{k}}(1/16)) \times {\mathbf R})} {\rm dist}^{2} \, (x, G_{k}) 
< \z \, {\hat E}_{k}^{2}.$\\

\item[$(b)$] $\int_{M_{k} \cap (B_{1}(0) \times {\mathbf R})} {\rm dist}^{2} \, (x, P_{k}) 
+ \int_{P_{k}^{\star} \cap ((B_{1/2}(0) \setminus S_{P_{k}}(1/16)) \times {\mathbf R})} {\rm dist}^{2} \, (x, G_{k}) 
\geq \z \, \hat{E}_{k}^{2}.$\\
\end{itemize}

Suppose first that possibility $(a)$ occurs. By Lemma~\ref{lowerupperbounds}, provided we choose
$\z \leq \z_{0} = \z_{0}$, where $\z_{0} = \z_{0}(n, \d)$ is as in Lemma~\ref{lowerupperbounds}, we have
in this case that 

\begin{equation} \label{main5}
P_{k} = {\rm graph} \, {\hat E}_{k} p_{k}^{0 \, +} \cup {\rm graph} \, {\hat E}_{k} p_{k}^{0 \, -} 
\end{equation}

\noindent 
for infinitely many $k$, with $P_{k}^{0} = P_{k}^{0 \, +} \cup P_{k}^{0 \, -}$ ($P_{k}^{0 \, \pm} 
= {\rm graph} \, p_{k}^{0 \, \pm}$) equal to a pair of affine hyperplanes 
satisfying 

\begin{equation}\label{main5-1}
c_{1} \leq \mbox{sup}_{B_{1}(0)} \, |p_{k}^{0 \, +} - p_{k}^{0 \, -}| \leq c_{2},
\end{equation}

\noindent
where $c_{1} = c_{1}(n, \d)$, $c_{2} = c_{2}(n, \d)$ are the positive constants given by Lemma~\ref{lowerupperbounds}. 
Note that (\ref{main5}) and (\ref{main5-1}) say that the blow-up by ${\hat E}_{k}$ 
of a subsequence of the sequence $\{P_{k}\}$ is a transverse pair of planes. 
So let $P^{0} = P^{0 \, +} \cup P^{0 \, -}$ be a subsequential limit of $\{P_{k}^{0}\}$
and consider the blow-up $v = (v^{+}, v^{-})$ of 
$M_{k}$ by ${\hat E}_{k}.$  We have directly from the defining condition of case $(a)$ and the identity (\ref{reg5-0-0}) 
that 

\begin{eqnarray}\label{main7-0}
&& \int_{B_{2/3}(0)} {\rm dist}^{2} \, ((x, v^{+}(x)), P^{0}) + {\rm dist}^{2} \, ((x, v^{-}(x)), P^{0}) \nonumber\\
&& \hspace{.3in} + \int_{B_{1/2}(0) \setminus S_{P^{0}}(1/8)} 
{\rm dist}^{2} \, ((x, p^{0 \, +}(x)), V) + {\rm dist}^{2} \, ((x, p^{0 \, -}(x)), V) \leq C\z
\end{eqnarray}

\noindent
where $V = {\rm graph} \, v^{+} \cup {\rm graph} \, v^{-}$ and $C = C(n, \d).$ In view of the lower bound of 
(\ref{main5-1}), we then have by Lemma~\ref{separation} that for 
any given $\eta \in (0, 1)$, 

\begin{equation}\label{main7-1}
\int_{B_{2/3}(0)} (v^{+} - p^{0 \, +})^{2} + (v^{-} - p^{0 \, -})^{2} \leq \eta
\end{equation}

\noindent
provided $\z = \z(n, \d, \eta) \in (0, 1)$ is sufficiently small.\\ 

We now separate the analysis of case $(a)$ into two further 
possibilities depending on the nature of $P^{0}$. Precisely one of the following must hold:\\ 

\begin{itemize}
\item[$(a)(i)$] $P^{0 \, +} \cap P^{0 \, -} \cap (B_{\th}(0) \times {\mathbf R}) = \emptyset$ or\\ 

\item[$(a)(ii)$] $P^{0 \, +} \cap P^{0 \, -} \cap (B_{\th}(0) \times {\mathbf R}) \neq \emptyset.$\\ 
\end{itemize}

\noindent
Suppose first that $(a)(i)$ holds. Taking $\eta =  \eta(n,\d, \th) >0$ in (\ref{main7-1}) sufficiently
small, we see by the estimate of Proposition~\ref{lipschitz}, part $(b)$ and the 
fact that $P^{0 \, +} \cap P^{0 \, -} \cap (B_{\th}(0) \times {\mathbf R}) = \emptyset$ that 
(\ref{main7-1}) implies, provided only that $\z = \z(n, \d, \th)$ is chosen sufficiently small, that 
we have $Z_{w} \cap B_{3\th/4}(0) = \emptyset,$ where $w = \frac{1}{2}(v^{+} - v^{-})$ 
and $Z_{w}$ is the zero set of $w$. 
By the remark following Lemma~\ref{crossings}, this means that $M_{k} \cap (B_{\th/2}(0) \times {\mathbf R})$ 
are embedded for all sufficiently large $k$, and hence by Schoen-Simon 
regularity theorem (\cite{SS}, Theorem 1), $M_{k} \cap (B_{\th/4}(0) \times {\mathbf R})$ decomposes as the 
disjoint union of minimal graphs ${\mathcal U}_{k}^{(1)}$, ${\mathcal U}_{k}^{(2)}$ (over the affine hyperplanes 
$P_{k}^{1}$ and $P_{k}^{2}$). By standard elliptic estimates, we then have that

\begin{equation}\label{main7-7}
(\s\th)^{-n-2}\int_{M_{k} \cap (B_{\s\th}(0) \times {\mathbf R})} 
\mbox{dist}^{2} \, (x, {\widetilde P}_{k}) \leq C\s^{2}
\th^{-n-2}\int_{M_{k} \cap (B_{\th}(0) \times {\mathbf R})} \mbox{dist}^{2} \, (x, P_{k}) 
\end{equation}

\noindent
for all $\s \in (0, 1/4)$, where $C = C(n)$ and ${\widetilde P}_{k}$ is the union of the tangent planes 
${\widetilde P}_{k}^{1},$ ${\widetilde P}_{k}^{2}$ to ${\mathcal U}_{k}^{(1)}$, ${\mathcal U}_{k}^{(2)}$ 
respectively at points $Z_{k}^{(1)} \in 
{\mathcal U}_{k}^{(1)}$, $Z_{k}^{(2)} \in {\mathcal U}_{k}^{(2)}$ with 
$\pi \, (Z_{k}^{(i)}) = 0$ for $i=1, 2.$  Taking $\s = \b/\th$ in this, we conclude that 

\begin{equation}\label{main7-8}
\b^{-n-2}\int_{M_{k} \cap (B_{\b}(0) \times {\mathbf R})} 
\mbox{dist}^{2} \, (x, {\widetilde P}_{k}) \leq C_{2}\b^{2}
\int_{M_{k} \cap (B_{1}(0) \times {\mathbf R})} \mbox{dist}^{2} \, (x, P_{k}) 
\end{equation}

\noindent
where $C_{2} = C_{2}(n, \th).$ Note that by the definition of ${\widetilde P}_{k}^{i}$ and elliptic estimates
again, it follows that $(\th/8)^{-2}{\rm dist}^{2}_{\mathcal H} \, ({\widetilde P}_{k} \cap (B_{\th/8}(0) \times {\mathbf R}), 
P_{k} \cap (B_{\th/8}(0) \times {\mathbf R})) \leq C\th^{-n-2}\int_{M_{k} \cap (B_{\th/4}(0) \times {\mathbf R})} 
{\rm dist}^{2} \, (x, P_{k})$ where $C = C(n),$ which implies that 
${\rm dist}^{2}_{\mathcal H} \, ({\widetilde P}_{k} \cap (B_{1}(0) \times {\mathbf R}), 
P_{k} \cap (B_{1}(0) \times {\mathbf R})) \leq C\int_{M_{k} \cap (B_{1}(0) \times {\mathbf R})} 
{\rm dist}^{2} \, (x, P_{k})$ where $C = C(n, \th).$\\  

Now for each $k$, take any ${\widetilde L}_{k} \in {\mathcal A}(M_{k}, \b)$ and any 
${\widetilde t}_{k} \in {\mathcal R}(M_{k}, {\widetilde L}_{k}, \b).$ Since $\b^{-n-2}\int_{M_{k} \cap (B_{\b}(0) \times 
{\mathbf R})}|x^{n+1}|^{2} \to 0$ as $k \to \infty$ (by Hausdorff convergence), it follows from the definition of 
${\mathcal A}(M_{k}, \b)$ that ${\hat E}_{M_{k}} \, (\b, {\widetilde L}_{k}) \to 0$ as $ k \to \infty,$ which in turn 
implies that ${\rm dist}_{\mathcal H} \, ({\widetilde L}_{k} \cap (B_{1}(0) \times {\mathbf R}), B_{1}(0) \times \{0\}) 
\to 0$ as $k \to \infty.$ Thus, we have in the present case (i.e. case $(a)(i)$) 
that $G_{M_{k}}^{({\widetilde L}_{k}, {\widetilde t}_{k})}(\b) \cap 
(B_{\b}(0) \times {\mathbf R}) = ({\mathcal U}_{k}^{(1)} \cup {\mathcal U}_{k}^{(2)}) 
\cap (B_{\b}(0) \times {\mathbf R}).$  If we write ${\mathcal U}_{k}^{(i)} \cap (B_{\b}(0) \times {\mathbf R})= 
{\rm graph} \, {\widetilde u}_{k}^{i}$, where ${\widetilde u}_{k}^{i} \, : \, 
B_{\b}(0) \to {\mathbf R}$, we have by (\ref{main7-1}) that provided 
$\z = \z(n, \d, \th)$ is suffciently small, for 
each $x = (x^{\prime}, x^{n+1}) \in ({\widetilde P}_{k}^{i})^{\star}$ ($={\widetilde P}_{k}^{i} \cap (B_{\b}(0) 
\times {\mathbf R})$), 

\begin{eqnarray}\label{main7-9}
{\rm dist} \, (x, G_{M_{k}}^{({\widetilde L}_{k}, {\widetilde t}_{k})}(\b)) \leq 
|x^{n+1} - {\widetilde u}_{k}^{i}(x^{\prime})| &\leq& 2 \, {\rm dist} \, ((x^{\prime}, u_{k}^{i}(x^{\prime})), 
{\widetilde P}_{k}^{i}) \nonumber\\
&=&2{\rm dist} \, ((x^{\prime}, u_{k}^{i}(x^{\prime})), {\widetilde P}_{k})
\end{eqnarray}

\noindent
for $i=1, 2$, which implies that 

\begin{equation}\label{main7-10}
\int_{{\widetilde P}_{k}^{\star} \cap (B_{\b}(0) \times {\mathbf R})}  
{\rm dist}^{2} \, (x, G_{M_{k}}^{({\widetilde L}_{k}, {\widetilde t}_{k})}(\b)) \leq 4\int_{M_{k} \cap 
(B_{\b}(0) \times {\mathbf R})} {\rm dist}^{2} \, (x, {\widetilde P}_{k}).
\end{equation}
 
\noindent
Thus, we conclude in case $(a)(i)$ that for infinitely many $k,$ the conclusions of the lemma 
hold with option $(3) (B)$, with $M_{k}$, $P_{k}$, $L_{k}$, $t_{k}$, ${\widetilde P}_{k}$, 
${\widetilde L}_{k}$ and ${\widetilde t}_{k}$ in place of 
$M$, $P$, $L$, $t$, ${\widetilde P}$, ${\widetilde L}$ and ${\widetilde t}$ respectively and with 
$\lambda = 2.$\\

If $(a)(ii)$ holds for infinitely many $k$, then we have 

$$\int_{M_{k} \cap (B_{1}(0) \times {\mathbf R})} {\rm dist}^{2} \, (x, P_{k}) < \z \, {\hat E}_{k}^{2}$$

for infinitely many $k$, where $P_{k}$ is as in (\ref{main5}) with $P^{0}_{k} = {\rm graph} \, p_{k}^{0 \, +} 
\cup {\rm graph} \, p_{k}^{0 \, -}$ equal to a transverse pair of 
affine hyperplanes satisfying (\ref{main5-1}) and $P_{k}^{0 \, +} \cap P_{k}^{0 \, -} \cap (B_{3\th/2}(0) \times {\mathbf R})
\neq \emptyset$.
Thus $\pi - \overline{\a} > \angle P_{k}^{0} > \overline{\a}$ for some fixed angle $\overline{\a} = \overline{\a}(n, \d) \in 
(0, \pi)$. Let $\t \in (0, 1)$ be arbitrary for the moment. Choosing the constant $\z  = \z(n, \th, \d, \t)> 0$ so that, 
in addition to the restrictions already imposed upon $\z$, we also have

\begin{equation}\label{main8}
\z \leq \left(\frac{2}{3}\right)^{n+2}\e_{0}
\end{equation}

\noindent
where $\e_{0} = \e_{0}(n, \d, \overline\a, 6\th, \t)$ is as in Lemma~\ref{transverse}, we have by 
Lemma~\ref{transverse} (with ${\overline \a}$ in place of $\a_{0}$, $6\th$ in place of $\th$ and
$\eta_{0, \, 2/3} \, M_{k}$ in place of $M$) that for infinitely many $k$, either there exists a pair 
of hyperplanes ${\widetilde P}_{k}$ with 

$$d_{\mathcal H} \, ({\widetilde P}_{k} \cap (B_{1}(0) \times {\mathbf R}), P_{k} \cap (B_{1}(0) \times {\mathbf R})
\leq C \int_{M_{k} \cap (B_{1}(0) \times {\mathbf R})} {\rm dist}^{2} \, (x, P_{k})$$

\noindent
satisfying 

\begin{equation}\label{main9}
\th^{-n-2} \int _{M_{k} \cap (B_{\th}(0) \times {\mathbf R})} \mbox{dist}^{2}\,(x, {\widetilde P}_{k}) 
\leq C\th^{2}\int_{M_{k} \cap (B_{1}(0) \times {\mathbf R})} 
\mbox{dist}^{2} \, (x, P_{k}),\\
\end{equation}

\noindent
where $C = C(n),$ or that 

\begin{equation}\label{main9-1}
\int_{M_{k} \cap (B_{1/2}(0) \times {\mathbf R})} {\rm dist}^{2} \, (x,L^{\prime}_{k}) 
\leq \t {\hat E}_{k}^{2}
\end{equation}

\noindent
for some affine hyperplane $L_{k}^{\prime}$  with 
$d_{\mathcal H} \, (L_{k}^{\prime} \cap (B_{1}(0) \times {\mathbf R}), L_{k} \cap (B_{1}(0) \times {\mathbf R})) 
\leq C{\hat E}_{k}$, $C = Cn).$ 
However, if (\ref{main9-1}) holds for infinitely many $k$, we must have that 

\begin{equation}\label{main9-2}
\int_{B_{1/2}(0)} (v^{+} - \ell^{\prime})^{2}  + (v^{-} - \ell^{\prime})^{2}
\leq \t
\end{equation}

\noindent
for some affine function $\ell^{\prime} \, : \, {\mathbf R}^{n} \times \{0\} \to {\mathbf R}$, which contradicts (\ref{main7-1}) 
provided we choose  $\eta = \eta(c_{1}) \in (0, 1)$ and $\t = \t(c_{1}) \in (0, 1)$ sufficiently small depending only 
on $c_{1}$ (hence only on $n$and $\d$), where 
$c_{1}$ is as in (\ref{main5-1}). Thus, provided
$\z = \z(n,\th, \d) \in (0, 1)$ is sufficiently small depending only on $n$, $\th$ and $\d$, we must have the option
(\ref{main9}) for infinitely many $k.$ 

Next in this case, we check that 

\begin{equation}\label{main9-2-1}
\int_{{\widetilde P}_{k}^{\star} \cap ((B_{\th/2}(0) \setminus S_{{\widetilde P}_{k}}(\th/16)) \times {\mathbf R})}
{\rm dist}^{2} \, (x, G_{M_{k}}^{({\widetilde L}_{k}, {\widetilde t}_{k})}(\th)) \leq 
4\int_{M_{k} \cap (B_{\th}(0) \times {\mathbf R})} {\rm dist}^{2} \, (x, {\widetilde P}_{k})
\end{equation}

\noindent
for arbitrary choices of ${\widetilde L}_{k} \in {\mathcal A}(M_{k}, \th)$ and 
${\widetilde t}_{k} \in {\mathcal R}(M_{k}, {\widetilde L}_{k}, \th).$ Reasoning as in 
case $(a)(i)$ (see paragraph preceding inequalities (\ref{main7-9})), we see that 
${\hat E}_{M_{k}} \, (\th, {\widetilde L}_{k}) \to 0$ as $k \to \infty,$ and by 
Lemma~\ref{transverse}, part $(b)(iii)$, that 

$$G_{M_{k}}^{({\widetilde L}_{k}, {\widetilde t}_{k})}(\th) \cap ((B_{\th/2}(0) \setminus 
S_{{\widetilde P}_{k}}(\th/16)) \times {\mathbf R})  = {\rm graph} \, u_{k}^{+} \cup {\rm graph} \, 
u_{k}^{-}$$

\noindent
where $u_{k}^{\pm} \in C^{2}(B_{\th/2}(0) \setminus S_{{\widetilde P}_{k}}(\th/16))$ (in fact 
$u_{k}^{\pm}$ solve the minimal surface equation), $u_{k}^{+} > u_{k}^{-}$, and 
${\rm dist} \, ((x^{\prime}, u_{k}^{\pm}(x^{\prime})), {\widetilde P}_{k}) 
= {\rm dist} \, ((x^{\prime}, u_{k}^{\pm}(x^{\prime})), {\widetilde P}_{k}^{\pm}) 
\geq \frac{1}{2}|u_{k}^{\pm}(x^{\prime}) - {\widetilde p}_{k}^{\pm}(x^{\prime})|$ for every
$x^{\prime} \in B_{\th/2}(0) \setminus S_{{\widetilde P}_{k}}(\th/16),$ where 
${\widetilde P}_{k}^{\pm} = {\rm graph}\, {\widetilde p}_{k}^{\pm}$. Hence we have in this case  
for any $x = (x^{\prime}, x^{n+1}) \in 
{\widetilde P}_{k}^{\star} \cap ((B_{\th/2}(0) \setminus S_{{\widetilde P}_{k}}(\th/16))
\times {\mathbf R}) ( = {\widetilde P}_{k} \cap ((B_{\th/2}(0) \setminus S_{{\widetilde P}_{k}}(\th/16) \times 
{\mathbf R}),$ provided $\z = \z(n, \th, \d) \in (0, 1)$ is chosen sufficiently small (so as to ensure that 
${\rm dist} \, (x, {\rm graph} \, u_{k}^{\pm}) 
\leq {\rm dist} \, (x, {\rm graph} \, u_{k}^{\mp})$ whenever 
$x \in {\widetilde P}_{k}^{\star \, \pm} \cap ((B_{\th/2}(0) \setminus S_{{\widetilde P}_{k}}(\th/16))
\times {\mathbf R}))$, that

\begin{equation}\label{main9-2-2}
{\rm dist} \, (x, G_{M_{k}}^{({\widetilde L}_{k}, {\widetilde t}_{k})}(\th))
\leq 2 {\rm dist} \, ((x^{\prime}, u_{k}^{\pm}(x^{\prime})), {\widetilde P}_{k})
\end{equation}

\noindent
for $x \in {\widetilde P}_{k}^{\star} \cap ((B_{\th/2}(0) \setminus S_{{\widetilde P}_{k}}(\th/16)) \times {\mathbf R}),$ 
where the sign $\pm$ is chosen according to whether $x \in {\widetilde P}^{\pm}_{k}.$ This of course 
implies (\ref{main9-2-1}). We thus have in case $(a)(ii),$ for infinitely many $k$, the conclusions of the lemma 
with option $(3) (A)$, with $M_{k}$, $P_{k}$, $L_{k}$, $t_{k}$, ${\widetilde P}_{k}$, 
${\widetilde L}_{k}$, ${\widetilde t}_{k}$ in place of 
$M$, $P$, $L$, $t$, ${\widetilde P}$, ${\widetilde L}$, ${\widetilde t}$ respectively and with 
$\lambda = 2.$\\ 

It now remains to analyze possibility $(b).$ We shall take $\z = \z(n, \th, \d) \in (0, 1)$ 
to be fixed for the remainder of the proof. If possibility $(b)$ holds for infinitely many $k$, consider the blow-up 
$v = (v^{+}, v^{-})$ of $\{M_{k}\}$ 
by the excess $\hat{E}_{k}$ off $L_{k}$, as described in Section~\ref{blowupprop}. (To be precise, 
since the excess 
${\hat E}_{k}$ is at scale 1 here, we are in fact applying the analysis of section~\ref{blowupprop} with 
$\eta_{0, \, 2/3} \, M_{k}$  in place of $M_{k}.$) 
Thus $v^{+}, v^{-} \in L^{2}(B_{1}(0)) \cap W^{1, \, 2} _{\rm loc} \, (B_{1}(0))$  satisfy the asymptotic decay properties as given by 
Theorem~\ref{regblowupthm}. Let $w = \frac{1}{2}(v^{+} - v^{-}),$ and $Z_{w}$ be the zero set of $w.$ 
One of the following 2 possibilities must occur:\\

\begin{itemize}
\item[$(b)(i)$] Either $Z_{w} \cap B_{2\b}(0) = \emptyset$ or $w \equiv 0$ or  
${\mathcal N}_{w}(z) = 1$ for every $z \in Z_{w} \cap B_{2\b}(0)$; i.e. $v$ has no branch 
point in $B_{2\b}(0).$\\

\item[$(b)(ii)$] $Z_{w} \cap B_{2\b}(0) \neq \emptyset$,  $w \not\equiv 0$ and there exists 
a point $z \in Z_{w} \cap B_{2\b}(0)$ with ${\mathcal N}_{w}(z) > 1$; i.e. $v$ has a branch 
point $z \in B_{2\b}(0)$.\\
\end{itemize}

If $(b)(i)$ occurs, then by Lemma~\ref{localdecomp}, the union of the graphs of $v^{+}$, $v^{-}$ over 
$B_{2\b}(0)$ is, locally 
near every point of $B_{2\b}(0)$, the union of the graphs of two harmonic functions, and hence, 
since $B_{2\b}(0)$ is 
simply connected, the union of the graphs of $v^{+}$, $v^{-}$ over $B_{2\b}(0)$ is globally the union of 
the graphs of two harmonic functions $v^{1}, \, v^{2} \, : B_{2\b}(0) \to {\mathbf R}$. 
Let $l^{i}$, $i=1,2$ be the affine part of the Taylor series of $v^{i}$ around $0$
(i.e. $l^{i}(x) = v^{i}(0) + x \cdot Dv^{i}(0)$ for $x \in B_{2\b}(0)$),
let $P_{k}^{(i)} = {\rm graph} \, (\varphi_{k} + {\hat E}_{k} l^{i})$ where
$L_{k} = {\rm graph} \, \varphi_{k}$ and 
set ${\widetilde P}_{k} = P_{k}^{(1)} \cup P_{k}^{(2)}$. Then

\begin{eqnarray}\label{decay1} 
\g^{-n-2}\int_{M_{k} \cap (B_{\g}(0) \times {\mathbf R})} \mbox{dist}^{2} \, (X,{\widetilde P}_{k})    
&=&\g^{-n-2}\int_{G_{k}^{+} \cap (B_{\g}(0) \times {\mathbf R})} \mbox{dist}^{2} \, (X,{\widetilde P}_{k}) +\nonumber\\
&&+ \, \g^{-n-2}\int_{G_{k}^{-} \cap (B_{\g}(0) \times {\mathbf R})} \mbox{dist}^{2} \, (X,{\widetilde P}_{k}) 
+ \nonumber\\
&& + \, \g^{-n-2}\int_{(\eta_{0, \, 2/3} \, M_{k} \setminus (G_{k}^{+} \cup G_{k}^{-})) 
\cap (B_{\g}(0) \times {\mathbf R})} \mbox{dist}^{2} \, (X,{\widetilde P}_{k}) \nonumber\\
&\leq& c\,\g^{-n-2}\int_{B_{\g}(0)}{\rm dist}^{2} \,((x, \varphi_{k}(x) + \overline\psi_{k} u_{k}^{+}(x)),
{\widetilde P}_{k}) 
+ \nonumber\\
&&+ \, c\g^{-n-2}\int_{B_{\g}(0)}{\rm dist}^{2} \, ((x, \varphi_{k}(x) + \overline\psi_{k} u_{k}^{-}(x)), 
{\widetilde P}_{k}) +\nonumber\\
&& + c \, \g^{-n-2} {\hat E}_{k}^{2 + \mu} \nonumber\\
&\leq& c \, \g^{-n-2}\int_{B_{\g}(0)}(\overline\psi_{k}u_{k}^{+} - {\hat E}_{k}v^{+})^{2} + \nonumber\\ 
&&+ c \, \g^{-n-2}\int_{B_{\g}(0)}(\overline\psi_{k}u_{k}^{-} - {\hat E}_{k}v^{-})^{2} + \nonumber\\
&&+ \, c \, \g^{-n-2}\int_{B_{\g}(0)}{\rm dist}^{2}\,((x,\varphi_{k}(x) + {\hat E}_{k}v^{+}(x)),{\widetilde P}_{k}) 
+ \nonumber\\ 
&&+ \, c \, \g^{-n-2}\int_{B_{\g}(0)}{\rm dist}^{2} \,((x,\varphi_{k}(x) + {\hat E}_{k}v^{-}(x)),{\widetilde P}_{k})
\nonumber\\
&&+ \,c \, \g^{-n-2}{\hat E}_{k}^{2+\mu} \nonumber\\
&=& c \, \g^{-n-2}\int_{B_{\g}(0)}(\overline\psi_{k} u_{k}^{+} - {\hat E}_{k}v^{+})^{2} + \nonumber\\ 
&&+ c \, \g^{-n-2}\int_{B_{\g}(0)}(\overline\psi_{k}u_{k}^{-} - {\hat E}_{k}v^{-})^{2} + \nonumber\\
&&+ \, c\, \g^{-n-2}\int_{B_{\g}(0)}{\rm dist}^{2} \, ((x, \varphi_{k}(x) + {\hat E}_{k}v^{1}(x)),{\widetilde P}_{k}) 
+ \nonumber\\ 
&&+ \, c \, \g^{-n-2}\int_{B_{\g}(0)}{\rm dist}^{2} \, ((x, \varphi_{k}(x) + {\hat E}_{k}v^{2}(x)), {\widetilde P}_{k})
\nonumber\\ 
&&+ \, c \, \g^{-n-2}{\hat E}_{k}^{2+\mu} \nonumber\\
&\leq& c \, \g^{-n-2} q({\hat E}_{k}) {\hat E}_{k}^{2} 
+ c\, \g^{-n-2}{\hat E}_{k}^{2}\int_{B_{\g}(0)} (v^{1} - l^{1})^{2} + \nonumber\\
&&+\, c\,\g^{-n-2}{\hat E}_{k}^{2}\int_{B_{\g}(0)}(v^{2} - l^{2})^{2}
+ \, c \, \g^{-n-2}{\hat E}_{k}^{2 + \mu} \nonumber\\
&\leq& c \, \g^{-n-2} q({\hat E}_{k}){\hat E}_{k}^{2} + c\g^{2}\b^{-n-4}{\hat E}_{k}^{2} 
\left(\int_{B_{1/2}(0)}(v^{+})^{2} + (v^{-})^{2}\right) + \nonumber\\
&& + \, c \, \g^{-n-2}{\hat E}_{k}^{2 + \mu} \nonumber\\
&\leq& c \, {\hat E}_{k}^{2}\left(\g^{-n-2} q({\hat E}_{k}) + \g^{2}\b^{-n-4} + 
\g^{-n-2}{\hat E}_{k}^{\mu}\right) \nonumber\\
\end{eqnarray}

\noindent
where $q(t) \to 0$ as $t \to 0$ and $c$ depends only on $n$ and $\d$. It follows from this that for all sufficiently large $k$, 

\begin{eqnarray}\label{main9-3}
&&\g^{-n-2}\int_{M_{k} \cap (B_{\g}(0) \times {\mathbf R})} \mbox{dist}^{2} \, (X,{\widetilde P}_{k}) 
\leq C\g^{2}\b^{-n-4}{\hat E}_{k}^{2} \nonumber\\
&&\hspace{.5in}\leq C_{3}\g^{2} \left(\int_{M_{k} \cap (B_{1}(0) \times {\mathbf R})} {\rm dist}^{2} \, (x, P_{k}) 
+ \int_{P_{k}^{\star} \cap ((B_{1/2}(0) \setminus S_{P_{k}}(1/16)) \times {\mathbf R})} {\rm dist}^{2} \, (x, G_{k})\right)
\end{eqnarray}

\noindent
where $C= C(n) > 0$ and we have set $C_{3} = \frac{C\b^{-n-4}}{\z}$, with $\z = \z(n, \th, \d)$  as in 
the definitions of cases $(a)$ and $(b)$, so that $C_{3} = C_{3}(n, \d, \th, \b).$\\

 Notice next that by the definition of ${\widetilde P}_{k}$, we see that 

$$d^{2}_{\mathcal H} \, ({\widetilde P}_{k} \cap (B_{1}(0) \times {\mathbf R}), 
L_{k} \cap (B_{1}(0) \times {\mathbf R})) \leq  C{\hat E}_{k}^{2}$$

\noindent
where $C = C(n).$ On the other hand, it follows from the inequality 
$$\z{\hat E}^{2}_{k} \leq 
\int_{M_{k} \cap (B_{1}(0) \times {\mathbf R})} {\rm dist}^{2} \, (x, P_{k}) + 
\int_{P_{k}^{\star} \cap ((B_{1/2}(0) \setminus S_{P_{k}}(1/16)) \times {\mathbf R})} {\rm dist}^{2} \, (x, G_{k})$$

\noindent
and the triangle inequality that  

\begin{eqnarray*}
&&d^{2}_{\mathcal H} \, (L_{k} \cap (B_{1}(0) \times {\mathbf R}), 
P_{k} \cap (B_{1}(0) \times {\mathbf R})) \nonumber\\
&&\hspace{.5in}\leq 
C\left( \int_{M_{k} \cap (B_{1}(0) \times {\mathbf R})} {\rm dist}^{2} \, (x, P_{k}) + 
\int_{P_{k}^{\star} \cap ((B_{1/2}(0) \setminus S_{P_{k}}(1/16)) \times {\mathbf R})} {\rm dist}^{2} \, (x, G_{k}) \right)
\end{eqnarray*}

\noindent
with $C = C(n, \th, \d)$, and therefore, by the triangle inequality again, we have that 

\begin{eqnarray}\label{main9-4}
&& d^{2}_{\mathcal H} \, ({\widetilde P}_{k} \cap (B_{1}(0) \times {\mathbf R}), 
P_{k} \cap (B_{1}(0) \times {\mathbf R})) \nonumber\\
&&\hspace{.5in} \leq C\left(\int_{M_{k} \cap (B_{1}(0) \times {\mathbf R})} {\rm dist}^{2} \, (x, P_{k})
+ \int_{P_{k}^{\star} \cap ((B_{1/2}(0) \setminus S_{P_{k}}(1/16)) \times {\mathbf R})} 
{\rm dist}^{2} \, (x, G_{k})\right)
\end{eqnarray}

\noindent
where $C = C(n, \th, \d).$\\

We next show in case $(b)(i)$ that for any choice of ${\widetilde L}_{k} \in {\mathcal A}(M_{k}, \g)$ and 
${\widetilde t}_{k} \in {\mathcal R}(M_{k}, {\widetilde L}_{k}, \g)$,

\begin{eqnarray}\label{main9-4-0-1}
&&\g^{-n-2}\int_{{\widetilde P}_{k}^{\star} \cap (B_{\g/2}(0) \setminus S_{{\widetilde P}_{k}}(\g/16)) \times {\mathbf R})} 
{\rm dist}^{2} \, (x, G_{M_{k}}^{({\widetilde L}_{k}, {\widetilde t}_{k})}(\g)) \nonumber\\
&& \hspace{.5in}\leq C_{3}\g^{2}\left(\int_{M_{k} \cap (B_{1}(0) \times {\mathbf R})} {\rm dist}^{2} \, (x, P_{k}) 
+ \int_{P_{k}^{\star} \cap ((B_{1/2}(0) \setminus S_{P_{k}}(1/16)) \times {\mathbf R})} {\rm dist}^{2} \, (x, G_{k})\right)
\end{eqnarray}

\noindent
where $C_{3} = C_{3}(n, \d, \th ,\b).$ For this, recall first that since $v^{i}$ are harmonic in $B_{2\b}(0)$, we have the 
estimates 
$\sup_{B_{4\g}(0)} \, |v^{i} - l^{i}|^{2} \leq C\g^{4}\b^{-n-4}\int_{B_{2\b}(0)}|v^{i}|^{2},$ $C = C(n)$, so that 
by Proposition~\ref{properties}, part (2), we have that

\begin{equation}\label{main9-4-0-2}
\sup_{B_{4\g}(0)} \, |v^{i} - l^{i}|\leq \G\g^{2}\b^{\frac{-n-4}{2}}
\end{equation}

\noindent
for $i=1, 2$, where $\G = \G(n)$. Consider first the case when 

\begin{equation}\label{main9-4-0-2-1}
\sup_{B_{4\g}(0)} \, |l^{1} - l^{2}| \geq \a\G\g^{2}\b^{\frac{-n-4}{2}}. 
\end{equation}

\noindent
where $\a >1$ is to be chosen depending only on $n$. In this case, if $\a > 68$, the estimates (\ref{main9-4-0-2}) say that for 
each $k$, there is no point 
$x \in B_{4\g}(0) \setminus S_{{\widetilde P}_{k}}(\g/32)$ such that 
$v^{1}(x) = v^{2}(x)$, and hence by the argument of Lemma~\ref{crossings}, it follows that 
for infinitely many $k$, $M_{k} \cap ((B_{3\g}(0) \setminus S_{{\widetilde P}_{k}}(\g/28)) \times {\mathbf R})$ 
must be embedded. But then by Schoen-Simon regularity theorem (\cite{SS}, Theorem 1), 
$M_{k} \cap ((B_{2\g}(0) \setminus S_{{\widetilde P}_{k}}(\g/24)) \times {\mathbf R}) 
= {\rm graph} \, {\widetilde u}_{k}^{+} \cup {\rm graph} \, {\widetilde u}_{k}^{-}$ where 
${\widetilde u}_{k}^{\pm} \, : \, B_{2\g}(0) \setminus S_{{\widetilde P}_{k}}(\g/24) \to {\mathbf R}$
are smooth solutions of the minimal surface equation in their domain, with ${\widetilde u}_{k}^{+} > {\widetilde u}_{k}^{-}.$  
Hence we have by elliptic theory the pointwise estimates 

\begin{equation}\label{main9-4-0-2-2}
\sup_{B_{\g}(0) \setminus S_{{\widetilde P}_{k}}(\g/16)} \, |{\widetilde u}_{k}^{\pm} - {\widetilde p}_{k}^{\pm}|^{2}
\leq C  \g^{-n}\int_{B_{3\g/2}(0) \setminus S_{{\widetilde P}_{k}}(\g/20)} |{\widetilde u}_{k}^{+} - {\widetilde p}_{k}^{+}|^{2} 
+ |{\widetilde u}_{k}^{-} - {\widetilde p}_{k}^{-}|^{2}
\end{equation}

\noindent
where $C = C(n).$ Recall our notation that ${\widetilde p}_{k}^{\pm} \, : \, B_{1}(0) \to {\mathbf R}$ are such that 
${\rm graph} \, {\widetilde p}_{k}^{\pm} = {\widetilde P}_{k}^{\pm}.$ Note also that by elliptic estimates 
again, $\sup_{B_{7\g/4}(0) \setminus S_{{\widetilde P}_{k}}(\g/22)} \, |D{\widetilde u}_{k}^{\pm}| 
\to 0$ as $k \to \infty$ (since $M_{k} \cap (B_{1}(0) \times {\mathbf R})  \to B_{1}(0) \times \{0\}$ 
in Hausdorff distance),  and hence $M_{k} \cap ((B_{3\g/2}(0) \setminus S_{{\widetilde P}_{k}}(\g/20)) \times {\mathbf R})
\subset G_{k}$ for infinitely many $k$. (This follows from the way $G_{k}$ is defined.) Hence, 
${\hat E}_{k}^{-1}({\widetilde u}_{k}^{\pm} - \varphi_{k}) \to v^{\pm}$ in $L^{2}(B_{3\g/2}(0) \setminus 
S^{v}(\g/18))$, where $\varphi_{k} \, : \, {\mathbf R}^{n} \times \{0\} \to {\mathbf R}$, 
${\rm graph} \, \varphi_{k} = L_{k}$ and $S^{v}(\g/18)$ denotes the set 
$\{ x \in B_{1}(0) \, : \, {\rm dist} \, (x, A) \leq \g/18\}$ with 
$A = \{l^{1}(x) = l^{2}(x)\}.$ Hence, by the estimates (\ref{main9-4-0-2}) and (\ref{main9-4-0-2-2}), we have that 

\begin{equation}\label{main9-4-0-2-3}
\sup_{B_{\g}(0) \setminus S_{{\widetilde P}_{k}}(\g/16)} \, |{\widetilde u}_{k}^{\pm} - {\widetilde p}_{k}^{\pm}|^{2}
\leq  2C\omega_{n}\G^{2}\g^{4}\b^{-n-4}{\hat E}_{k}^{2}
\end{equation}

\noindent
where $C = C(n)$ is as in (\ref{main9-4-0-2-2}). Thus, if $\a = \a(n)$ in (\ref{main9-4-0-2-1}) 
is chosen sufficiently large, the estimates (\ref{main9-4-0-2-3}) imply, by exactly the same reasoning used to
justify inequality (\ref{main9-2-2}), that for each $x = (x^{\prime}, x^{n+1}) 
\in {\widetilde P}_{k}^{\star} \cap (B_{\g/2}(0) \setminus S_{{\widetilde P}_{k}}(\g/16)) \times {\mathbf R})$, 

\begin{equation}\label{main9-4-0-3}
{\rm dist} \, (x, G_{M_{k}}^{({\widetilde L}_{k}, {\widetilde t}_{k})}(\g))
\leq 2 {\rm dist} \, ((x^{\prime}, {\widetilde u}_{k}^{\pm}(x^{\prime})), {\widetilde P}_{k})
\end{equation}

\noindent
where the sign $\pm$ is chosen according to whether $x \in {\widetilde P}_{k}^{\pm}.$ In view of the 
estimate (\ref{main9-3}), this gives (\ref{main9-4-0-1}).\\

Suppose the condition (\ref{main9-4-0-2-1}) fails to hold. Note that we have 

\begin{equation}\label{main9-4-0-4}
G_{M_{k}}^{({\widetilde L}_{k}, {\widetilde t}_{k})}(\g) = 
{\rm graph} \, {\widetilde u}_{k}^{+} \cup {\rm graph} \, {\widetilde u}_{k}^{+}
\end{equation}

\noindent
where ${\widetilde u}_{k}^{\pm} \, : \, \pi \, (G_{M_{k}}^{({\widetilde L}_{k}, 
{\widetilde t}_{k})}(\g)) \to {\mathbf R}$ are Lipschitz with Lipschitz constant $\leq 3/2$ and 
${\widetilde u}_{k}^{+} \geq {\widetilde u}_{k}^{-}.$ From this we see that for $x = (x^{\prime}, x^{n+1}) 
\in {\widetilde P}_{k}^{\star} \cap ((B_{\g/2}(0) \setminus S_{{\widetilde P}_{k}}(\g/16)) \times {\mathbf R})$,

\begin{eqnarray}\label{main9-4-0-5}
{\rm dist} \, (x, G_{M_{k}}^{({\widetilde L}_{k}, {\widetilde t}_{k})}(\g))
&\leq& \min \, \{|x^{n+1} - {\widetilde u}_{k}^{+}(x^{\prime})|, |x^{n+1} - {\widetilde u}_{k}^{-}(x^{\prime})|\}\nonumber\\
&\leq& 2 {\rm dist} \, ((x^{\prime}, {\widetilde u}_{k}^{\pm}(x^{\prime})), {\widetilde P}_{k}) 
+ 2{\hat E}_{k}|l^{1}(x^{\prime}) - l^{2}(x^{\prime})|\nonumber\\
&\leq& 2 {\rm dist} \, ((x^{\prime}, {\widetilde u}_{k}^{\pm}(x^{\prime})), {\widetilde P}_{k}) 
+ 2\a\G\g^{2}\b^{\frac{-n-4}{2}}{\hat E}_{k}
\end{eqnarray}

\noindent
where in the second of the inequalities here we have used the fact that 
${\widetilde u}_{k}^{\pm}$ are Lipschitz functions with Lipschitz constants $\leq 3/2$, and 
the sign $\pm$ there is chosen according to whether $x \in {\widetilde P}_{k}^{\pm}$. 
By the estimate (\ref{main9-3}) and the defining property of case $(b)(i)$, we again have 
from this the required estimate (\ref{main9-4-0-1}). We have thus
shown that in case $(b)(i),$ for infinitely many $k$, the conclusions of the lemma 
with option $(3) (C)$ hold, with $M_{k}$, $P_{k}$, $L_{k}$, $t_{k}$, ${\widetilde P}_{k}$, 
${\widetilde L}_{k}$ and ${\widetilde t}_{k}$ in place of 
$M$, $P$, $L$, $t$, ${\widetilde P}$, ${\widetilde L}$ and ${\widetilde t}$ respectively and with 
$\lambda = 2.$\\ 
 
Finally, suppose $(b)(ii)$ occurs. i.e. that there exists a point 
$z \in Z_{w} \cap B_{2\b}$ with ${\mathcal N}_{w}(z) > 1$. 
Then $v^{+}(z) = v^{-}(z)$ and  by Lemma~\ref{branchptdecay}, we have that 

\begin{equation}\label{main9-4-1}
\rho^{-n-2} \int_{B_{\rho}(z)} (v^{+} - l_{z})^{2}  + (v^{-} - l_{z})^{2} \leq 
C\rho^{\nu} \int_{B_{1}(0)} (v^{+})^{2} + (v^{-})^{2}
\end{equation}
 
\noindent
for some affine function $l_{z}$ and 
all $\rho \in (0, 1/64).$ Here $C = C( n, \d) >0$ and $\nu = \nu(n, \d)>0.$ Now fix this 
$z.$ We obtain from (\ref{main9-4-1})  that  

\begin{eqnarray}
\rho^{-n-2}\int_{B_{\rho}(0)} (v^{+} - l_{z})^{2}  + (v^{-} - l_{z})^{2}  &\leq& 
 \left(1 + \frac{|z|}{\rho}\right)^{n+2}
(\rho + |z|)^{-n-2}\int_{B_{\rho + |z|}(z)} (v^{+} - l_{z})^{2}  + (v^{-} - l_{z})^{2} \nonumber\\
&\leq& C\left(1 + \frac{|z|}{\rho}\right)^{n+2} \left(1 + \frac{|z|}{\rho}\right)^{\nu}\rho^{\nu} 
\int_{B_{1}(0)} (v^{+})^{2} + (v^{-})^{2}
\end{eqnarray}
 
\noindent
provided $\rho + |z| \leq 1/64.$ In particular, taking $\rho = \b$ in this and 
using the fact that $z \in B_{2\b}(0)$, (so that $1 + \frac{|z|}{\b} \leq 3$) and 
since $3\b < 1/64$, we have that 

\begin{equation}\label{main9-4-2}
\b^{-n-2}\int_{B_{\b}(0)} (v^{+} - l_{z})^{2}  + (v^{-} - l_{z})^{2}  \leq C\b^{\nu}   
\int_{B_{1}(0)} (v^{+})^{2} + (v^{-})^{2}
\end{equation}

\noindent
where $C = C(n, \d).$ With this, we can estimate as in (\ref{decay1}) 
to conclude that if possibility $(b)(ii)$ occurs, then we must have that

\begin{eqnarray}\label{main9-4-3}
&&\b^{-n-2}\int_{M_{k} \cap (B_{\b}(0) \times {\mathbf R})} \mbox{dist}^{2} \, (X,{\widetilde P}_{k}) \leq 
C\b^{\nu}{\hat E}_{k}^{2}\nonumber\\ 
&&\hspace{.5in}\leq C_{2}\b^{\nu} \left(\int_{M_{k} \cap (B_{1}(0) \times {\mathbf R})} {\rm dist}^{2} \, (x, P_{k}) 
+ \int_{P_{k}^{\star} \cap ((B_{1/2}(0) \setminus S_{P_{k}}(1/16)) \times {\mathbf R})} {\rm dist}^{2} \, (x, G_{k})\right)
\end{eqnarray}

\noindent
for all sufficiently large $k$, where ${\widetilde P}_{k} 
= \mbox{graph} \, (\varphi_{k} + {\hat E}_{k}l_{z})$. Here $C = C(n, \d)$ 
is as in the estimate (\ref{main9-4-2}) and we have set $C_{2} = \frac{C}{\z}$ where 
$\z = \z(n, \d, \th)$ is as in the definition of cases $(a)$ and $(b)$, so that $C_{2} = C_{2}(n, \d, \th).$\\

Arguing exactly as in the proof of the estimate (\ref{main9-4}), we also have in this case that

\begin{eqnarray}\label{main9-4-4}
&& d^{2}_{\mathcal H} \, ({\widetilde P}_{k} \cap (B_{1}(0) \times {\mathbf R}), 
P_{k} \cap (B_{1}(0) \times {\mathbf R})) \nonumber\\
&&\hspace{.5in} \leq C\left(\int_{M_{k} \cap (B_{1}(0) \times {\mathbf R})} {\rm dist}^{2} \, (x, P_{k})
+ \int_{P_{k}^{\star} \cap ((B_{1/2}(0) \setminus S_{P_{k}}(1/16)) \times {\mathbf R})} 
{\rm dist}^{2} \, (x, G_{k}) \right)
\end{eqnarray}

\noindent
where $C = C(n, \th, \d).$\\

To complete the proof of the lemma, we now check in case $(b)(ii)$ that for any choice of 
${\widetilde L}_{k} \in {\mathcal A}(M_{k}, \b)$ and 
${\widetilde t}_{k} \in {\mathcal R}(M_{k}, {\widetilde L}_{k}, \b)$,

\begin{eqnarray}\label{main9-4-5}
&&\b^{-n-2}\int_{{\widetilde P}_{k}^{\star} \cap (B_{\b/2}(0) \setminus S_{{\widetilde P}_{k}}(\b/16)) \times {\mathbf R})} 
{\rm dist}^{2} \, (x, G_{M_{k}}^{({\widetilde L}_{k}, {\widetilde t}_{k})}(\b)) \nonumber\\
&& \hspace{.5in}\leq C_{2}\b^{\nu}\left(\int_{M_{k} \cap (B_{1}(0) \times {\mathbf R})} {\rm dist}^{2} \, (x, P_{k}) 
+ \int_{P_{k}^{\star} \cap ((B_{1/2}(0) \setminus S_{P_{k}}(1/16)) \times {\mathbf R})} {\rm dist}^{2} \, (x, G_{k})\right)
\end{eqnarray}

\noindent
where $C_{2} = C_{2}(n, \d, \th)$ is as in the estimate (\ref{main9-4-3}). But this follows directly from the pointwise 
estimate that for each $x = (x^{\prime}, x^{n+1}) 
\in {\widetilde P}_{k}^{\star} \cap (B_{\b/2}(0) \setminus S_{{\widetilde P}_{k}}(\b/16)) \times {\mathbf R})$, 

\begin{eqnarray}\label{main9-4-6}
{\rm dist} \, (x, G_{M_{k}}^{({\widetilde L}_{k}, {\widetilde t}_{k})}(\b))
&\leq&\min \, \{|x^{n+1} - {\widetilde u}_{k}^{+}(x^{\prime})|, |x^{n+1} - {\widetilde u}_{k}^{-}(x^{\prime})|\}\nonumber\\
&\leq& 2 \min \, \{{\rm dist} \, ((x^{\prime}, {\widetilde u}_{k}^{+}(x^{\prime})), {\widetilde P}_{k}), 
{\rm dist} \, ((x^{\prime}, {\widetilde u}_{k}^{-}(x^{\prime})), {\widetilde P}_{k})\}
\end{eqnarray}

\noindent
where ${\widetilde u}_{k}^{\pm}$ are defined exactly as in (\ref{main9-4-0-4}) with $\b$ in place of $\g.$ In the second 
of the inequalities above, we have used the fact that ${\widetilde u}_{k}^{\pm}$ are Lipschitz functions with 
Lipschitz constants $\leq 3/2$, and that ${\widetilde P}_{k}$ is a single affine hyperplane. The required estimate 
(\ref{main9-4-5}) follows from this and the estimate (\ref{main9-4-3}). We have thus
shown that in case $(b)(ii),$ for infinitely many $k$, the conclusions of the lemma 
with option $(3) (B)$ hold, with $M_{k}$, $P_{k}$, $L_{k}$, $t_{k}$, ${\widetilde P}_{k}$, 
${\widetilde L}_{k}$ and ${\widetilde t}_{k}$ in place of 
$M$, $P$, $L$, $t$, ${\widetilde P}$, ${\widetilde L}$ and ${\widetilde t}$ respectively and with 
$\lambda = \nu.$ This completes the proof of the lemma.
\end{proof}

\bigskip

\section{Main regularity theorems}\label{mainthms}
We are now ready to prove Theorems~\ref{maintheorem}, \ref{maincorollary} and 
\ref{lowdimreg}.\\

\begin{proof}[\bf Proof of Theorem~\ref{maintheorem}] First choose $\th=\th(n, \d) \in (0, 1/16)$ such 
that $C_{1}\th^{\lambda} < 1/4$, then choose 
$\b=\b(n, \d) \in (0, \th/16)$ such that $C_{2}\b^{\lambda} < 1/4$, and finally choose 
$\g = \g (n,\d) \in (0, \b/16)$ such that $C_{3}\g^{\lambda} < 1/4$, where $C_{1}$, $C_{2},$ $C_{3}$ and 
$\lambda$ are as in Lemma~\ref{mainlemma1}.\\

Suppose $M$ satisfies the hypotheses of Theorem~\ref{maintheorem}. Note first that since  
$L^{2}$ closeness of $M$ to a hyperplane implies closeness in Hausdorff distance, the hypothesis 
$\int_{M \cap (B_{1}(0) \times {\mathbf R})} |x^{n+1}|^{2} \leq \e$ implies that 
$d_{\mathcal H} \, (L_{0} \cap (B_{1}(0) \times {\mathbf R}), B_{1}(0)) \leq \t(\e)$ and 
$\int_{B_{1/2}(0)} {\rm dist}^{2} \, (x, G_{M}^{L_{0}, t_{0})}(1)) \leq \t(\e)$ 
for any $L_{0} \in {\mathcal A} \, (M, 1)$ and  any $t_{0} \in {\mathcal R}(M, L_{0}, 1)$, where 
$\t(\e) \downarrow 0$ as $\e \downarrow 0.$ Fix such $L_{0}$ and $t_{0}$.\\   

In what follows, let us use the notation

$$Q(\r, P, L, t) = \r^{-n-2}\int_{M \cap (B_{\r}(0) \times {\mathbf R})} {\rm dist}^{2}(x, P) 
+\r^{-n-2}\int_{P^{\star} \cap ((B_{\r/2}(0) \setminus S_{P}(\r/16)) \times {\mathbf R})} 
{\rm dist}^{2} \, (x, G_{M}^{(L, \, t)}(\r)).$$

If $\e = \e(n, \d) \in (0, 1)$ is sufficiently small, by iterating Lemma~\ref{mainlemma1} starting with 
$P = P_{0} equiv {\mathbf R}^{n} \times \{0\}$, $L=L_{0}$ and $t = t_{0},$ we get a sequence of 
pairs of affine hyperplanes 
${P_{j}}$, a sequence of affine hyperplanes $L_{j} \in {\mathcal A} \, (M, \th^{k_{j}}\b^{l_{j}}\g^{m_{j}})$ and a sequence of 
numbers $t_{j} \in {\mathcal R} \, (M, L_{j}, \th^{k_{j}}\b^{l_{j}}\g^{m_{j}})$ satisfying at the $j$th iteration either

\begin{eqnarray}\label{mt-1}
Q(\th^{k_{j}}\b^{l_{j}}\g^{m_{j}}, P_{j}, L_{j}, t_{j}) &\leq& 4^{-1}Q(\th^{k_{j}-1}\b^{l_{j}}\g^{m_{j}}, P_{j-1}, L_{j-1}, t_{j-1}) \nonumber\\
&\leq& 4^{-j}Q_{1}
\end{eqnarray}  

\noindent
or 

\begin{eqnarray}\label{mt-2}
Q(\th^{k_{j}}\b^{l_{j}}\g^{m_{j}}, P_{j}, L_{j}, t_{j}) &\leq& 4^{-1}Q(\th^{k_{j}}\b^{l_{j}-1}\g^{m_{j}}, P_{j-1}, L_{j-1}, t_{j-1}) \nonumber\\
&\leq& 4^{-j}Q_{1}
\end{eqnarray}

\noindent
or

\begin{eqnarray}\label{mt-2-1}
Q(\th^{k_{j}}\b^{l_{j}}\g^{m_{j}}, P_{j}, L_{j}, t_{j}) &\leq& 4^{-1}Q(\th^{k_{j}}\b^{l_{j}}\g^{m_{j}-1}, P_{j-1}, L_{j-1}, t_{j-1}) \nonumber\\
&\leq& 4^{-j}Q_{1}
\end{eqnarray}

\noindent
where $k_{j}$, $l_{j}$, $m_{j}$ are non-negative integers with  $k_{j}+l_{j}+m_{j} = j,$ 
$P_{0} = {\mathbf R}^{n} \times \{0\}$, 
and 

$$Q_{1} = \int_{M \cap (B_{1}(0) \times {\mathbf R})}|x^{n+1}|^{2} + \int_{B_{1}(0)} {\rm dist}^{2} \,
(x, G_{M}^{({\mathbf R}^{n} \times \{0\}, t)}(1)).$$

\noindent 
Let us denote the sequence of scales so generated $\{s_{j}\}.$ Thus, for each $j = 0, 1, 2, \ldots,$ 
$s_{j} = \th^{k_{j}}\b^{l_{j}}\g^{m_{j}}$ for some non-negative integers $k_{j}, \, l_{j}, \, m_{j}$ with 
$k_{j}  + l_{j} + m_{j} = j,$ and, $s_{j+1} = \th s_{j}$ or $\b s_{j}$ or $\g s_{j}.$ Then (\ref{mt-1})--(\ref{mt-2-1-1}) may 
be rewritten as

\begin{eqnarray}\label{mt-2-1-1}
Q(s_{j}, P_{j}, L_{j}, t_{j}) &\leq& 4^{-1} \, Q(s_{j-1}, P_{j-1}, L_{j-1}, t_{j-1})\nonumber\\
&\leq& 4^{-j}Q_{1}.
\end{eqnarray}

\noindent
The lemma also gives us that\\

\begin{eqnarray}\label{mt-3}
\mbox{dist}^{2}_{\mathcal H} \, (P_{j} \cap (B_{1}(0) \times {\mathbf R}), P_{j-1} \cap (B_{1}(0) 
\times {\mathbf R})) 
&\leq& C \, Q(s_{j-1}, P_{j-1}, L_{j-1}, t_{j-1})\nonumber\\ 
&\leq&  C4^{-j}Q_{1}
\end{eqnarray}

\noindent
and that 
\begin{eqnarray}\label{mt-3-0}
\mbox{dist}^{2}_{\mathcal H} \, (P_{j}^{T} \cap (B_{1}(0) \times {\mathbf R}), P_{j-1}^{T} \cap (B_{1}(0) 
\times {\mathbf R})) 
&\leq& C \, Q(s_{j-1}, P_{j-1}, L_{j-1}, t_{j-1})\nonumber\\ 
&\leq&  C4^{-j}Q_{1}
\end{eqnarray}

\noindent
where $C$ depends only on $n$ and $\d$. Thus, $\{P_{j}\}$ is a Cauchy sequence of pairs of affine hyperplanes, and hence 
there exists a pair of affine hyperplanes $P$ such that 
$P_{j} \to P.$ By (\ref{mt-3}), (\ref{mt-3-0}) and (\ref{mt-1}) respectively, we have that 

\begin{equation}\label{mt-3-1}
\mbox{dist}^{2}_{\mathcal H} \, (P\cap (B_{1}(0) \times {\mathbf R}), P_{j-1} \cap (B_{1}(0) 
\times {\mathbf R})) \leq  C4^{-j}Q_{1},
\end{equation}

\begin{equation}\label{mt-3-1-1}
\mbox{dist}^{2}_{\mathcal H} \, (P^{T}\cap (B_{1}(0) \times {\mathbf R}), P^{T}_{j-1} \cap (B_{1}(0) 
\times {\mathbf R})) \leq  C4^{-j}Q_{1},\hspace{.2in} \mbox{and}
\end{equation}

\begin{equation}\label{mt-4}
s_{j}^{-n-2}\int_{M \cap (B_{s_{j}}(0) \times {\mathbf R})} 
\mbox{dist}^{2} \, (x, P) \leq C4^{-j}Q_{1}
\end{equation}

\noindent
where $C$ depends only on $n$ and $\d$. Note that (\ref{mt-3-1}) and (\ref{mt-3-1-1}) in particular say that 

\begin{equation}\label{mt-5}
\mbox{dist}^{2}_{\mathcal H} \, (P \cap (B_{1}(0) \times {\mathbf R}), B_{1}(0)) \leq CQ_{1} \hspace{.2in} 
\mbox{and}
\end{equation}

\begin{equation}\label{mt-5-1}
\mbox{dist}^{2}_{\mathcal H} \, (P^{T} \cap (B_{1}(0) \times {\mathbf R}), B_{1}(0)) \leq CQ_{1}.
\end{equation}

Now, given any $\r \in (0, 1/8)$, 
there exists a unique $j$ with $s_{j+1} \leq \r < s_{j}.$ 
Since $\g < \b < \th$, this implies that 
$\g^{j+1} \leq \r < \th^{j},$ or, equivalently, that $\frac{\log \r}{\log \th} > j \geq \frac{\log \r}{\log \g} - 1.$ 
Hence, by (\ref{mt-4}), we conclude that 

\begin{eqnarray}\label{mt-6}
\r^{-n-2}\int_{M \cap (B_{\r}(0) \times {\mathbf R})} \mbox{dist}^{2} \, (x, P) & \leq & 
s_{j+1}^{-n-2}\int_{M \cap (B_{s_{j}}(0) \times {\mathbf R})} \mbox{dist}^{2} \, (x, P) \nonumber\\
&=& \left(\frac{s_{j}}{s_{j+1}}\right)^{n+2} \, s_{j}^{-n-2}\int_{M \cap (B_{s_{j}}(0) \times {\mathbf R})} 
\mbox{dist}^{2} \, (x, P) \nonumber\\
&\leq & C\r^{\k}Q_{1}
\end{eqnarray}  

\noindent
for all $\r \in (0, 1/8)$, where $\k = -\log 4/ \log \g$ and $C$ depends only on $n$.\\ 

Next observe that
we can move the base point and repeat the entire argument leading to the estimates (\ref{mt-5}), (\ref{mt-5-1}) and
(\ref{mt-6}). Specifically, for any given $X \in M \cap B_{3/4}^{n+1}(0)$, we have 

\begin{eqnarray}\label{mt-7}
\frac{{\mathcal H}^{n} \, (M \cap B_{7/8}^{n+1}(X))}{\omega_{n}(7/8)^{n}} &=& 
\frac{{\mathcal H}^{n} \, (G \cap B_{7/8}^{n+1}(X)) + {\mathcal H}^{n}((M \setminus G) \cap B_{7/8}^{n+1}(X))}
{\omega_{n}(7/8)^{n}}\nonumber\\
&\leq& \frac{1}{\omega_{n} (7/8)^{n}}\int_{\Omega} \sqrt{1 + |Du^{+}|^{2}} + \sqrt{1 + |Du^{-}|^{2}} 
+ C{\hat E}^{2 + \mu}\nonumber\\
 &\leq& 2\sqrt{2} + C{\hat E}^{2 + \mu}\nonumber\\
&\leq& 3 - \d/16
\end{eqnarray}

\noindent
provided $\e = \e(n, \d) \in (0, 1)$ is sufficiently small. Here $G$ denotes the graphical part of 
$M\cap (B_{7/8}(0) \times {\mathbf R})$ as described in 
Section~\ref{blowupprop}, $\Omega \subset B_{7/8}(0)$, $u^{\pm} \, : \, \Omega \to 
{\mathbf R}$ are such that $G = {\rm graph} \, u^{+} \cup {\rm graph} \, u^{-}$ and 
${\hat E}^{2} = \int_{M \cap (B_{1}(0) \times {\mathbf R})}
|x^{n+1}|^{2}$, $C = C(n, \d)$ and $\mu = \mu(n).$ 

\noindent
Thus, provided $\e$ is sufficiently small, we can repeat the argument leading to the estimates 
(\ref{mt-5}), (\ref{mt-5-1}) and (\ref{mt-6}), iteratively applying Lemma~\ref{mainlemma1} with $\d/16$ 
in place of $\d$ and $\eta_{X, \, 7/16} \, M$ in place of $M$ and starting with 
$P = \eta_{X, 7/16} \, ({\mathbf R}^{n} \times \{0\})$ and arbitrary $L \in {\mathcal A} \, (\eta_{X, 7/16} \, M, 1)$ and 
$t \in {\mathcal R} \, (\eta_{X, 7/16} \, M, L, 1)$ to conclude that for every 
$X \in {\overline M} \cap B_{3/4}^{n+1}(0)$, there exists a pair of affine hyperplanes $P_{X}$ such that

\begin{equation}\label{mt-8}
\mbox{d}^{2}_{\mathcal H} \, (P_{X} - X\cap (B_{1}(0) \times {\mathbf R}), B_{1}(0)) \leq CQ_{1},\\
\end{equation}

\begin{equation}\label{mt-8-1}
\mbox{d}^{2}_{\mathcal H} \, (P_{X}^{T} \cap (B_{1}(0) \times {\mathbf R}), B_{1}(0)) \leq CQ_{1} \hspace{.2in}
\mbox{and}\\
\end{equation}

\begin{equation}\label{mt-9}
\r^{-n-2}\int_{M \cap (B_{\r}(X^{\prime}) \times {\mathbf R})} \mbox{dist}^{2} \, (x, P_{X}) \leq C\r^{\k}Q_{1}
\end{equation}  

\noindent
for all $\r \in (0, 1/8),$ where $X^{\prime} = {\pi} \, (X).$ It follows from this that provided $\e = \e(n, \d)$ 
is sufficiently small, 
${\overline M} \cap (B_{1/2}(0) \times {\mathbf R})$ is the graph of a 2-valued $C^{1,\k}$ function. The
proof of this claim is as follows:\\

First note that by choosing $\e = \e(n, \d)$ sufficiently small, we may assume that 
$\overline M \cap (B_{1/2}(0) \times {\mathbf R}) \subseteq \overline M \cap B_{3/4}^{n+1}(0).$ 
 For $X \in M \cap B_{3/4}^{n+1}(0)$, let $P_{X}$ be as in (\ref{mt-8})---(\ref{mt-9}). Note that by
(\ref{mt-8}),  $P_{X} \cap {\pi}^{-1} \, (X^{\prime})$ consists precisely of two 
(possibly coinciding) points $X$ and $\widetilde X.$
Multiplying the inequality (\ref{mt-9}) by $\r^{2}$ and letting $\r \to 0,$ we see 
that for each $X \in M \cap (B_{1/2}(0) \times {\mathbf R})$, 
$M \cap {\pi}^{-1} \, (X^{\prime}) = P_{X} \cap {\pi}^{-1} \, (X^{\prime})$ so that 
$M \cap {\pi}^{-1} \, (X^{\prime})$ consists of (possibly coinciding) two points.  Furthermore, (\ref{mt-9})
says that the two tangent planes to $M$ at $X$ and $\widetilde X$ are the two hyperplanes whose 
union is $P_{X}.$ Thus, in view of (\ref{mt-8-1}), we have that for each 
$X \in M \cap (B_{1/2}(0) \times {\mathbf R})$, 
$|\nu_{1}(X)  - e_{n+1}|, \, |\nu_{2}(\widetilde X) - e_{n+1}| \leq CQ_{1},$ where $\nu_{1}, \nu_{2}$ denote the 
(locally defined) upward pointing unit normals to $M$. (Thus, in case $X = {\widetilde X}$, $\nu_{1}(X), \nu_{2}(X)$ 
are the two upward pointing unit normals at $X$ to the respective smooth sheets whose 
union is $M \cap B_{\s}^{n+1}(X)$ for some $\s > 0.$) This means that 

\begin{equation}\label{mt-10}
M \cap (B_{1/2}(0) \times {\mathbf R}) = \mbox{graph} \, u^{+} \cup \mbox{graph} \, u^{-}
\end{equation}

\noindent
where $u^{\pm} \, : \, B_{1/2} \setminus {\pi} \, (\mbox{sing} \, M) \to {\mathbf R}$
are Lipschitz functions with $u^{+} \geq u^{-}$ and Lipschitz constants $\leq CQ_{1}.$ The functions 
$u^{+}$, $u^{-}$ then extend uniquely as Lipschitz functions 
${\overline u}^{+}, \, {\overline u}^{-} \, : \, B_{1/2}(0) \to {\mathbf R}$ respectively, 
with the same Lipschitz constants, and we have that

\begin{equation}\label{mt-11}
{\overline M} \cap (B_{1/2}(0) \times {\mathbf R}) = \mbox{graph} \, {\overline u}^{+} \cup \mbox{graph} \, 
{\overline u}^{-}.
\end{equation}

Now note that since $M \cap (B_{1/2}(0) \times {\mathbf R})$ is a Lipschitz graph with Lipschitz 
constant $\leq 1$, it follows that $Q_{1} \leq C{\hat E}^{2}$ for some fixed constant 
$C = C(n),$ where ${\hat E} = \int_{M \cap (B_{1}(0) \times {\mathbf R})} |x^{n+1}|^{2}$, 
and hence we may replace $Q_{1}$ with ${\hat E}^{2}$ in all of the above estimates. 
Note also that since $0 \in {\overline M}$, the estimate for the Lipschitz constant implies the height bound

\begin{equation}\label{mt-11-1}
|u^{+}(x)|, |u^{-}(x)| \leq C{\hat E}.
\end{equation}

It now remains to show that the union of the two Lipschitz graphs in (\ref{mt-11}) is  
the graph of a single 2-valued $C^{1,\k}$ function, with its $C^{1, \k}$ norm bounded by 
a constant times ${\hat E}$. We proceed as follows:\\

Take any two points $X_{1}, X_{2} \in M \cap (B_{1/2}(0) \times {\mathbf R})$ with $X_{1}^{\prime} 
\neq X_{2}^{\prime}$ and let 
$r = |X_{1}^{\prime} - X_{2}^{\prime}|.$ By  (\ref{mt-9}), we have that 

\begin{equation}\label{mt12}
(2r)^{-n-2}\int_{M \cap (B_{2r}(X_{2}^{\prime}) \times {\mathbf R})} {\rm dist}^{2} \, (x, P_{X_{2}}) \leq Cr^{\k}{\hat E}^{2} 
\end{equation}

\noindent
and hence, since $B_{r}(X_{1}^{\prime}) \subset B_{2r}(X_{2}^{\prime})$ it follows that

\begin{equation}\label{mt14}
r^{-n-2}\int_{M \cap (B_{r}(X_{1}^{\prime}) \times {\mathbf R})} {\rm dist}^{2} \, (x, P_{X_{2}}) \leq Cr^{\k}{\hat E}^{2}.
\end{equation}

Also by (\ref{mt-8}) and (\ref{mt-8-1}) we have 

\begin{equation}\label{mt15}
d^{2}_{\mathcal H} \, (P_{X_{2}} \cap (B_{1}(X_{1}^{\prime}) \times {\mathbf R}), B_{1}(X_{1}^{\prime})) \leq C{\hat E}^{2}.
\end{equation}

This means that provided $\e = \e(n, \d)$ is sufficiently small, we may use Lemma~\ref{mainlemma1} exactly as 
it was used in the argument leading to 
(\ref{mt-5}), (\ref{mt-5-1}) and (\ref{mt-6}), with 
$\eta_{X_{1}, \, r} \, M$ in place of $M$, $\eta_{X_{1}, \, r} \, P_{X_{2}}$ in place of 
$P$ of the lemma (which was taken to be the multiplicity 2 hyperplane corresponding to 
${\mathbf R}^{n} \times \{0\}$ in the argument leading to (\ref{mt-5}), (\ref{mt-5-1}) and 
(\ref{mt-6}) above) to conclude that there exists
a pair of affine hyperplanes $P_{X_{1}}^{\prime}$ such that

\begin{equation}\label{mt-16}
(\r r)^{-n-2}\int_{M \cap B_{\r r}(X_{1}^{\prime}) \times {\mathbf R})} {\rm dist}^{2} \, (x, P_{X_{1}}^{\prime})
\leq C\r^{\k}r^{-n-2}\int_{M \cap (B_{r}(X_{1}^{\prime}) \times {\mathbf R})} {\rm dist}^{2} \, (x, P_{X_{2}}) 
\end{equation}

\noindent
for all $\r \in (0, 1/8)$ and\\

\begin{equation}\label{mt-17}
d_{\mathcal H}^{2} \,(P_{X_{1}}^{\prime \, T} \cap (B_{1}(0) \times {\mathbf R}), 
P_{X_{2}}^{T} \cap  (B_{1}(0) \times{\mathbf R})) \leq Cr^{-n-2}\int_{M \cap (B_{r}(X_{1}^{\prime}) 
\times {\mathbf R})} {\rm dist}^{2} \, (x, P_{X_{2}}). 
\end{equation}

In view of (\ref{mt-9}) (with $X = X_{1}$), (\ref{mt-16}) implies that $P^{\prime}_{X_{1}} \equiv P_{X_{1}}$, 
and hence, (\ref{mt-17}) combined with (\ref{mt14}) gives that
 
\begin{equation}\label{mt-18}
d_{\mathcal H}^{2} \,(P_{X_{1}}^{T} \cap (B_{1}(0) \times {\mathbf R}), 
P_{X_{2}}^{T} \cap  (B_{1}(0) \times{\mathbf R})) \leq C{\hat E}^{2}|X_{1}^{\prime} - X_{2}^{\prime}|^{\k}
\end{equation}

\noindent
for all $X_{1}, \, X_{2} \in {\overline M} \cap (B_{\s/4}(0) \times {\mathbf R}).$ This says that, in the notation
introduced in Section~\ref{notation}, the 2-valued function 
$u \, : \, B_{1/2}(0) \to {\mathbf T}({\mathbf R})$ defined by 
$u(x) = \{{\overline u}^{+}(x), {\overline u}^{-}(x)\}$ satisfies 

\begin{equation}
{\mathcal G} \, (Du(x_{1}), Du(x_{2})) \leq C{\hat E}|x_{1} - x_{2}|^{\k/2}
\end{equation}

\noindent
for all $x_{1}, x_{2} \in B_{1/2}(0).$ i.e. that $u$ is a $C^{1,\k/2}(B_{1/2}(0))$ function
with $\|u\|_{C^{1, \k/2}(B_{1/2}(0))} \leq C{\hat E}.$ The theorem is thus proved.
\end{proof}

\medskip
\begin{proof}[\bf Proof of Theorem~\ref{maincorollary}:] By Theorem~\ref{maintheorem}, 
${\overline M} \cap (B_{1/2}(0) \times {\mathbf R})$ is either the graph of a single $C^{1, \a}$ function 
$u^{0}$ or the graph of a 2-valued $C^{1, \a}$ function $u$, with  the appropriate estimate
for the $C^{1, \a}$ norm in either case. In case ${\overline M} \cap (B_{1/2}(0) \times {\mathbf R})$ is 
the graph of a 2-valued function $u$, we have that locally in a neighborhood $\Omega_{x}$ of any point $x$ 
of the open set $B_{1/2}(0) \setminus {\pi}\, ({\rm sing} \, M),$ $u$  is given by two functions, 
each satisfying the minimal surface equation in $\Omega_{x}.$ Since ${\mathcal H}^{n-2} \, ({\rm sing} \, M) = 0$ 
by assumption, $B_{1/2}(0) \setminus {\pi} \, ({\rm sing} \, M)$ is simply connected, and hence 
$M \cap ((B_{1/2}(0) \setminus {\pi} \, ({\rm sing} \, M) ) \times {\mathbf R})$ is equal to the union 
of the graphs of two functions ${\widetilde u}_{1}, \, {\widetilde u}_{2} \, : \, B_{1/2}(0) \setminus 
{\pi} \, ({\rm sing} \, M) \to {\mathbf R}$ each satisfying the minimal surface equation. 
But then by the removable singularity theorem of L. Simon \cite{S4}, 
${\widetilde u}_{1}$, ${\widetilde u}_{2}$ extend as functions 
$u_{1}, \, u_{2} \, : \, B_{1/2 }(0) \to {\mathbf R}$  satisfying the minimal surface equation.  
\end{proof}

\bigskip

\begin{proof}[\bf Proof of Theorem~\ref{lowdimreg}:]
We argue by contradiction. Were the assertion false, there would exist a sequence of hypersurfaces 
$M_{k} \in {\mathcal I}_{b},$ $k = 1, 2, 3, \ldots,$ with $0 \in {\overline M}_{k}$, 

\begin{equation}\label{low-1}
\frac{{\mathcal H}^{n} \, (M_{k})}{\omega_{n}2^{n}} \leq 2 + 1/k
\end{equation}

\noindent
such that for each $k$, 
the conclusion of the theorem fails with $M_{k}$ in place of $M,$ $1/k$ in place of 
$\s$ and with any choice of orthogonal rotation $q$ of ${\mathbf R}^{n+1}$ and any choice of 
pair of hyperplanes $P$. By Allard's integer varifold compactness theorem, we obtain,  
possibly after passing to a subsequence of $\{k\}$ which we continue to 
denote $\{k\}$ that $M_{k} \to V$ as varifolds for some integer multiplicity stationary varifold 
of $B_{2}^{n+1}(0).$\\  

First consider the case when, for a further subsequence of $\{k\}$ which 
we shall continue to denote $\{k\}$, there exist points $Z_{k} \in {\overline M}_{k} \cap B_{1/k}^{n+1}(0)$ 
with $\Theta_{M_{k}} \, (Z_{k}) \geq 2-1/k.$ Then by (\ref{low-1}), upper semicontinuity of 
density and the continuity of mass under varifold convergence, it follows that 
$2 \leq \Theta_{\|V\|} \, (0) \leq \frac{{\mathcal H}^{n} \, ({\rm spt} \, \|V\| \cap B_{2}^{n+1}(0))}{\omega_{n} 2^{n}}
\leq 2$, so that by the monotonicity formula, $V$ must be a cone with $\Theta_{\|V\|} \, (0) = 2.$ 
By Lemma~\ref{conelemma}, part $(b)$ below, $V$ must either be a pair of transverse hyperplanes or 
a hyperplane with multiplicity $2.$ Thus for infinitely many $k$ of the original 
sequence, the conclusions of the theorem hold, by Theorem~\ref{maintheorem} or 
Theorem 1 of \cite{WN2}, with $M_{k}$ in place of $M,$ $1/k$  in place of $\s$
and a suitable choice of an orthogonal transformation $q$ (which carries ${\rm spt}\, \|V\|$ to 
${\mathbf R}^{n} \times \{0\}$ in case $V$ is a single hyperplane with multiplicity $2$, or 
$L \equiv {\rm graph} \, \frac{1}{2}(p^{1} + p^{2})$ to ${\mathbf R}^{n} \times \{0\}$ in case 
${\rm spt} \, \|V\| = P^{1} \cup P^{2}$ where $P^{i} = {\rm graph} \, p^{i}$, $i=1, 2$ are 
transverse hyperplanes.) \\

The remaining alternative is that there is a number $\r >0$ such that for infinitely many $k$, 
$M_{k} \cap B_{\r}^{n+1}(0)$ is an embedded stable minimal hypersurface 
of $B_{\r}^{n+1}(0)$ (with no singularities since $2 \leq n \leq 6$).
But then by Theorem~\ref{curvest} below, there exists a fixed number $\G>0$ such that 

$$\sup_{M_{k} \cap B_{\r/2}^{n+1}(0)} \, |A_{k}| \leq \frac{\G}{\r}$$

\noindent
where $A_{k}$ denotes the second fundamental form of $M_{k}$ and $|A_{k}|$ its length. If $\nu_{k}$ 
denotes an oriented unit normal to $M_{k},$ then for $X \in {\widetilde M}_{k} \cap B_{\r/2}^{n+1}(0)$ 
(where ${\widetilde M}_{k}$ is the connected component of $M_{k}$ that contains the origin) and 
for any unit speed geodesic $\g$ of $M$ from the origin to $X$, we have that 

$$|\nu_{k}(X) - \nu_{k}(0)| \leq \int_{0}^{\ell} \left|\frac{d}{dt} \nu_{k}(\g(t))\right| dt
\leq \frac{\G \ell}{\r}$$ 

\noindent
where $\ell$ is the geodesic distance from $0$ to $X$. Thus, $|\nu_{k}(X) - \nu_{k}(0)| \leq 1/2$ 
for all points $X \in {\widetilde M}_{k} \cap B_{\r/2}^{n+1}(0)$ contained in a geodesic ball of radius 
$\frac{\r}{2\G}$ centered at $0,$ which means that there exists $\r_{1} > 0$ such that for all suffciently large $k$, 
${\widetilde M}_{k} \cap B_{\r_{1}}^{n+1}(0)$ is graphical over the tangent plane to $M_{k}$ at $0.$ Hence for
infinitely many $k$, the conclusions of the theorem hold again, with $M_{k}$ in place of 
$M$, $1/k$ in place of $\s$ and $q_{k}$ in place of $q,$ where $q_{k}$ is the rotation that carries the tangent     
plane of $M_{k}$ at $0$ to ${\mathbf R}^{n} \times \{0\}.$ Note that the $C^{1, \a}$ estimate 
of the conclusion of the theorem holds in this case by standard elliptic estimates.  This completes the proof
of the theorem.   
\end{proof}

\section{Compactness and decomposition theorems}\label{compact}
\setcounter{equation}{0}

In this section we prove Theorems~\ref{compactnessthm} and ~\ref{decomposition}.  
First we need the following lemma in which we shall use the following notation. Given 
a $p$ dimensional rectifiable varifold $V = (\Sigma, \th)$ in ${\mathbf R}^{p+1}$, where 
$\th$ is the multiplicity of $V$ (see \cite{S1}, chapter 4
for an exposition of the the theory of rectifiable varifolds), we let $V \times {\mathbf R}^{n-p}$ 
denote the rectifiable varifold $(\Sigma \times {\mathbf R}^{n-p}, \th_{1})$ of ${\mathbf R}^{n+1}$ where 
$\th_{1}(x,y) = \th(x)$ for $(x,y) \in \Sigma \times {\mathbf R}^{n-p}.$ \\

\begin{lemma} \label{conelemma}
(a) Suppose ${\mathbf C}$ is a cone with $\Theta \, ({\|{\mathbf C}\|}, 0) \leq 3$ belonging to
the varifold closure of immersed, stable minimal hypersurfaces $M$ of $B_{2}^{n+1}(0)$ 
with ${\mathcal H}^{n-2} \, ({\rm sing} \, M) = 0.$ If either $(i)$ $2 \leq n \leq 6$ or $(ii)$ $n \geq 7$ and 
${\mathbf C}$ has the form ${\mathbf C} = {\mathbf C}_{0} \times {\mathbf R}^{n-p}$ for some $p \leq 6$, 
then ${\mathbf C}$ must be the sum of at most $3$ (multiplicity 1 varifolds associated with) 
hyperplanes of ${\mathbf R}^{n+1}$.\\

\noindent
(b) There exists a fixed number $\e \in (0, 1)$ such that if ${\mathbf C}$ is a cone 
with $\Theta \, ({\|{\mathbf C}\|}, 0) \leq 2 + \e$ belonging to the varifold 
closure of immersed stable minimal hypersurfaces $M$ of $B_{2}^{n+1}(0)$ satisfying 
${\mathcal H}^{n-2} \, ({\rm sing} \, M) < \infty,$ and if either $(i)$ $2 \leq n \leq 6$ or 
(ii) $n \geq 7$ and ${\mathbf C}$ has the form ${\mathbf C} = {\mathbf C}_{0} \times {\mathbf R}^{n-p}$ 
for some $p \leq 6$, then ${\mathbf C}$ is equal to the sum of at most $2$ 
(multiplicity 1 varifolds associated with) hyperplanes of ${\mathbf R}^{n+1}.$\\  
\end{lemma}

\medskip

\begin{proof}
To see part $(a)$, first recall the following two standard facts about any stationary cone $W$; namely, 
that $\Theta \,(\|W\|, X) \leq \Theta \,(\|W\|, 0)$ for any $X \in {\rm spt} \, \|W\|$ and,
that if $\Theta \, (\|W\|, X) = \Theta \, (\|W\|, 0)$ for some $X \in {\rm spt} \, \|W\|$, then 
${\rm spt} \, \|W\|$ is invariant under translations by elements of the line $\{tX \, :\, t \in {\mathbf R}\}.$ 
In view of these facts, we may assume without loss of generality that 
$\Theta (\|{\mathbf C}\|, X) < 3$ for every $X \in \mbox{spt} \, \|{\mathbf C}\| \setminus \{0\}.$ 
If $\Theta(\|{\mathbf C}\|, X) \in \{1, 2\}$ for every $X \in \mbox{spt} \, \|{\mathbf C}\| \setminus \{0\}$ then by 
Allard's regularity theorem, Theorem 1 of \cite{WN2} and Theorem~\ref{maintheorem} of the present paper, 
it follows that $\mbox{spt} \, \|{\mathbf C}\| \setminus \{0\}$ is a regular, immersed submanifold, and hence 
J. Simons' theorem (\cite{SJ}, see also \cite{S1}, appendix $B$) concerning the non-existence of non-trivial 
stable minimal hypercones 
of dimension $\leq 6$ (applied to the cross section ${\mathbf C}_{0}$ in case $n \geq 7$ and 
${\mathbf C} = {\mathbf C}_{0} \times {\mathbf R}^{n-p}$ for some $p \leq 6$) implies that
${\mathbf C}$ must be a union of hyperplanes. If there is a point 
$X \in \mbox{spt} \, \|{\mathbf C}\| \setminus \{0\}$ with $\Th \, (\|{\mathbf C}\|, X) \not\in \{1, 2\}$, by taking a 
tangent cone to ${\mathbf C}$ at $X$, 
we produce a cone ${\mathbf C}^{\prime}$ singular at the origin having the form 
(after a rotation) ${\mathbf C}^{\prime} 
= {\mathbf C}^{\prime}_{0} \times {\mathbf R}$ in case $2 \leq n \leq 6$ 
where the dimension of the cross section ${\mathbf C}^{\prime}_{0}$ is one less than the dimension
of ${\mathbf C}$, or, in case $n \geq 7$ and ${\mathbf C} = {\mathbf C}_{0} \times 
{\mathbf R}^{n-p}$ for some $p \leq 6$, having the form ${\mathbf C}^{\prime}_{0} \times {\mathbf R}^{n-p+1}$ where 
the dimension of ${\mathbf C}^{\prime}_{0}$ is one less than the dimension of ${\mathbf C}_{0}.$ 
If ${\mathbf C}^{\prime}$ has  
density $1$ or $2$ everywhere except at $\{0\}$, then the preceding argument tells us that it must be a union of 
hyperplanes, giving a contradiction. So it must have a singular point $X^{\prime} \in {\rm spt} \, 
\|{\mathbf C}^{\prime}\|$ other than the origin. Proceeding inductively, taking a tangent cone to 
${\mathbf C}^{\prime}$ at $X^{\prime}$, we arrive at a contradiction after a finite number of steps.\\

To see part $(b)$, first consider the case when $\Th \, (\|{\mathbf C}\|, 0) \leq 2.$ In this case, we may write, 
after a rotation ${\mathbf C} = {\mathbf C}_{0} \times {\mathbf R}^{n-q}$ for some $q \leq 6$ so that 
${\mathbf C}_{0}$ is a stable cone with $\Th \, (\|{\mathbf C}_{0}\|, 0) \leq 2$ and 
$\Th \, (\|{\mathbf C}_{0}\|, X) < 2$ for every $X \in {\rm spt} \, \|{\mathbf C}_{0}\| \setminus \{0\}.$ 
By taking successive tangent cones at possible singular points away from the origin and 
using J. Simons' theorem as before, we immediately arrive at the conclusion in this case.\\

To show the existence of an $\e$ as asserted in the lemma, we argue by contradiction. If there were no such 
$\e$, then there would exist a sequence of cones ${\mathbf C}_{k},$ $k = 1, 2, \ldots$ in ${\mathbf R}^{n+1}$ 
each occurring as the weak limit 
of a sequence of immersed stable minimal hypersurfaces $M_{k_{j}}$, $j = 1, 2, \ldots$ of $B_{2}^{n+1}(0)$ 
with ${\mathcal H}^{n-2} \, ({\rm sing} \, M_{k_{j}}) < \infty$ for each $k$ and $j$, such that 
$\Th \, (\|{\mathbf C}_{k}\|, 0) \leq 2 + k^{-1}$ and with the additional property in case $n \geq 7$ that 
each ${\mathbf C}_{k}$ has the form ${\mathbf C}_{k} = {\mathbf C}_{0}^{(k)} \times {\mathbf R}^{n-p}$  
for some $p \leq 6$, and yet ${\mathbf C}_{k}$ is not a union of 
hyperplanes for any $k$. In view of the uniform mass bound (implied by the density hypothesis), we may 
extract a subsequence, which we will continue to denote ${\mathbf C}_{k}$ such that 
${\mathbf C}_{k} \to {\mathbf C}$ for some cone ${\mathbf C}$ where the convergence is as varifolds.
By continuity of mass under varifold convergence, we have that $\Th \, (\|C\|, 0) \leq 2.$ Furthermore, 
in case $n \geq 7$, ${\mathbf C}$ has the form ${\mathbf C} = {\mathbf C}_{0} \times {\mathbf R}^{n-p}.$ 
Hence by the discussion of the previous paragraph
we have that ${\mathbf C}$ is either a single multiplicity 1 hyperplane or a pair of hyperplanes. 
Hence by Allard's regularity theorem, Theorem 1.1 of \cite{WN2} or Theorem~\ref{maintheorem}
of the present paper, we must have that for each sufficiently large $k$ and each 
sufficiently large $j$ (depending on $k$)  $M_{k_{j}} \cap B_{1}^{n+1}(0)$ must either be 
a multiplicity 1 $C^{1, \a}$ graph or the union of two multiplicity 1 $C^{1, \a}$ 
graphs or a single 2 valued $C^{1, \a}$ graph, with an interior estimate, in each case, 
for the $C^{1, \a}$ norm of the function(s) defining the graph  
in terms of the $L^{2}$ norm of the function(s) over a larger ball. But this means that for all sufficiently large $k$, 
${\rm spt} \, \|{\mathbf C}_{k}\| \cap B_{1}^{n+1}(0)$ must either be immersed or a equal to a 
2 valued $C^{1, \a}$ graph. In all cases, by taking the tangent cone at the origin (which on the one hand 
must be equal to the tangent plane(s) to the graph at the origin and on the other hand
coincide with ${\mathbf C}_{k}$ since ${\mathbf C}_{k}$ is already a cone), we see that 
${\rm spt} \, \|{\mathbf C}_{k}\|$ must be the union of at most 2 hyperplanes, contrary to the assumption.
The lemma is thus proved.
\end{proof}

\medskip

\begin{proof}[\bf Proof of Theorem~\ref{compactnessthm}]
First note that by Allard's varifold compactness theorem (\cite{AW}, \cite{S1}), we obtain a stationary 
integral varifold $V$
of $B_{2}^{n+1}(0)$ such that for some subsequence of $\{M_{k}\}$ which we continue to 
denote $\{M_{k}\}$, we have $M_{k} \to V$ as varifolds. Next we claim that there exists 
$\s = \s(n, \d) \in (0, 1/2)$ such that  

\begin{equation}\label{cpt-1}
\frac{{\mathcal H}^{n}(M_{k} \cap B_{1}^{n+1}(X))}{\omega_{n}} \leq 3 - \d/2
\end{equation}

\noindent
for all $k$ and all $X \in M_{k} \cap B_{\s}^{n+1}(0).$ To see this, fix any $k$ and suppose that
$X \in M_{k} \cap B_{1/2}^{n+1}(0).$ Then by the monotonicity of mass ratio, we have that 

\begin{eqnarray}\label{cpt-2}
\frac{{\mathcal H}^{n} \, (M_{k} \cap B_{1}^{n+1}(X))}{\omega_{n}} 
&\leq& \frac{{\mathcal H}^{n} \, (M_{k} \cap B_{1 + |X|}^{n+1}(0))}{\omega_{n}}\nonumber\\
 &=& (1 + |X|)^{n}\frac{{\mathcal H}^{n} \, (M_{k} \cap B_{1 + |X|}^{n+1}(0))}{\omega_{n}(1 + |X|)^{n}}\nonumber\\
&\leq& (1 + |X|)^{n} \frac{{\mathcal H}^{n} \, (M_{k} \cap B_{2}^{n+1}(0))}{\omega_{n}2^{n}}\nonumber\\
&\leq& (1 + |X|)^{n}(3 - \d),
\end{eqnarray}

\noindent
which readily implies (\ref{cpt-1}) provided $X \in M_{k} \cap B_{\s}^{n+1}(0)$ for a suitable 
choice of $\s = \s(n, \d) \in (0, 1/2).$ It then follows that   
$\frac{{\mathcal H}^{n}({\rm spt} \, \|V\| \cap B_{1}^{n+1}(X))}{\omega_{n}} \leq 3 - \d/2$
for all $X \in {\rm spt} \, \|V\| \cap B_{\s}^{n+1}(0),$ so that 
$\Theta \, (\|V\|, X) \leq 3 - \d/2$ for all $X \in {\rm spt} \, \|V\| \cap B_{\s}^{n+1}(0).$ 
Hence, if $X \in {\rm spt} \, \|V\| \cap B_{\s}^{n+1}(0)$ is a singular point of 
${\rm spt} \, \|V\| \cap B_{\s}^{n+1}(0)$ and ${\mathbf C}$ is any tangent cone to $V$ at $X$
having, after a rotation, the form ${\mathbf C} = {\mathbf C}_{0} \times {\mathbf R}^{n-p}$ for some 
$p \in \{1, 2, \ldots, n\}$, then by Lemma~\ref{conelemma},
we must have, in case $n \geq 7$, that $p \geq 7;$ 
otherwise, Lemma~\ref{conelemma} says that 
${\mathbf C}$ must be a union of hyperplanes, and since $\Theta \, (\|{\mathbf C}\|, 0) = \Theta\, (\|V\|, X) < 3$, 
it must either be a multiplicity 1 hyperplane, a multiplicity 2 hyperplane or a transverse pair of 
hyperplanes, in all of which cases, by Allard's regularity theorem, Theorem~\ref{maincorollary} of the 
present paper or Theorem 1 of \cite{WN2}, ${\rm spt} \, \|V\|$ would be a regular immersed submanifold near $X$, 
contrary to the hypothesis that $X$ is a singular point. Hence, in case $n \geq 7$, Federer's 
dimension reducing principle implies that ${\rm dim} \, {\rm sing} \, {\rm spt}\, \|V\| 
\cap B_{\s}^{n+1}(0) \leq n-7.$ In case $2 \leq n \leq 6$, Lemma~\ref{conelemma} says that 
any tangent cone at any point $X \in {\rm spt}\,\|V\| \cap B_{\s}^{n+1}(0)$ is either 
a multiplicity 1 hyperplane, a multiplicity 2 hyperplane or a transverse pair of hyperlanes, so that
$X$ must be a regular point of ${\rm spt} \, \|V\|.$\\

It remains to show that when $n = 7,$ ${\rm sing} \, {\rm spt} \, \|V\| \cap B_{\s}^{n+1}(0)$ 
is discrete. This follows by the standard argument. Were it not true, there exist  singular 
points $X$ and $X_{j},$ $j = 1, 2 , \ldots,$ of ${\rm spt} \, \|V\| \cap B_{\s}^{n+1}(0)$ such that      
$X_{j} \neq X$ for all $j$ and $X_{j} \to X.$ Let $\r_{j} = |X - X_{j}|.$ Then after passing to 
a subsequence, $\eta_{X, \, \r_{j}} \, V \to {\mathbf C}$, where ${\mathbf C}$ is a cone with singularities at the origin 
and at a point $Y = \lim_{j \to \infty} \, \frac{X - X_{j}}{\r_{j}} \in {\rm spt} \|{\mathbf C}\| \cap {\mathbf S}^{n-1}.$
(This last claim that $Y$ is a singular point of ${\mathbf C}$ follows from the appropriate regularity 
theorem---i.e. Allard's theorem, Theorem~\ref{maincorollary} of the present paper or Theorem 1 of \cite{WN2}.) 
But then since ${\mathbf C}$ is a cone, this means that the entire ray defined by $Y$ consists of 
singularities of  ${\mathbf C},$ which is impossible since in dimension $n = 7$, the singular set is 
$0$-dimensional. This concludes the proof of the theorem.  
\end{proof}

\bigskip

\begin{proof}[\bf Proof of Theorem~\ref{decomposition}:] Let $\e = \e(n) \in (0, 1)$ be as in 
Lemma~\ref{conelemma}, part $(b)$, and choose $\s = \s(n) \in (0, 1/2)$  as 
in the proof of Theorem~\ref{compactnessthm} (i.e. via the estimate (\ref{cpt-2})), 
so that  $\Th \, (\|V\|, X) \leq 2 + \e/2$ for all $X \in {\rm spt} \, \|V\| \cap B_{\s}^{n+1}(0).$ 
Let $B$ be the set of branch points of ${\rm spt} \, |V\| \cap B_{\s}^{n+1}(0).$ Thus
$$B = \{Z \in {\rm sing} \, V \cap B_{\s}^{n+1}(0) \, : \, 
\mbox{ $V$ has a (unique) multiplicity 2 tangent plane at $Z$}\}.$$ 

\noindent
Set $S = {\rm sing} \, V \cap B_{\s}^{n+1}(0) 
\setminus B.$ Then ${\rm sing} \, V \cap B_{\s}^{n+1}(0) = B \cup S$, $B \cap S = \emptyset$ by definition, 
and by Theorem~\ref{maintheorem}, 
$S$ is relatively closed in ${\rm spt} \, \|V\| \cap B_{\s}^{n+1}(0).$ By Theorem~\ref{maincorollary}, 
it follows readily that if ${\mathcal H}^{n-2} \, (B) = 0$, then $B = \emptyset.$ To estimate the Hausdorff dimension of 
$S$, we proceed as follows. Consider an arbitrary point $Z \in S.$ Let ${\mathbf C}$ be any tangent cone to 
$V$ at $Z$. Then by the definition of $S$ and Theorem 1 of \cite{WN2}, ${\mathbf C}$ cannot be equal to 
a pair of hyperplanes. Hence if $2 \leq n \leq 6$, it follows from Lemma~\ref{conelemma}, part $(b)$
that $S = \emptyset.$ If $n \geq 7,$ Lemma~\ref{conelemma}, part $(b)$ says that, after a rotation, ${\mathbf C} 
= {\mathbf C}_{0} \times {\mathbf R}^{n-p}$ for some $p \geq 7.$ It then follows by the dimension 
reducing principle of Federer that 

\begin{equation}\label{decomposition-1}
{\mathcal H}^{n-7+\g} \, (S) = 0
\end{equation}

\noindent
for every $\g > 0.$\\

It only remains to show that $S$ is finite when $n = 7.$ To see this, suppose $S$ is an infinite set. Then 
there exists a point $Z \in {\rm spt} \, \|V\| \cap {\overline B}_{\s}^{n+1}(0)$ and a sequence of points 
$Z_{j} \in S$ with $Z_{j} \neq z$ for each $j$, such that $Z_{j} \to Z$ as $j \to \infty.$  Let $r_{j} = |Z_{j} - Z|$
and $V_{j} = \eta_{Z, \, r_{j} \, \#} \, V.$ Then after passing to a subsequence, 
$V_{j} \to {\mathbf C}$ as varifolds, where ${\mathbf C}$ is a cone. Let $\z_{j} = r_{j}^{-1}(Z_{j} - z).$ 
Then $\z_{j} \in {\rm sing} \, V_{j} \cap {\mathbf S}^{n},$ and hence, after passing to a further subsequence, 
$\z_{j} \to \z \in {\rm sing} \, {\mathbf C} \cap {\mathbf S}^{n}.$  Now write
${\rm sing} \, {\mathbf C} \cap B_{\s}^{n+1}(0) = B_{{\mathbf C}} \cup S_{{\mathbf C}}$, where 
$B_{\mathbf C}$ is the set of 
branch points of ${\mathbf C}$ in $B_{\s}^{n+1}(0)$ (thus each point of $B_{\mathbf C}$ is a singular point of 
${\mathbf C}$ where ${\mathbf C}$ has a unique multiplicity 2 
tangent plane) and $S_{\mathbf C}$ is the complement of $B_{\mathbf C}$ in ${\rm sing} \, 
{\mathbf C} \cap B_{\s}^{n+1}(0).$ 
Similarly, write ${\rm sing} \, V_{j} \cap B_{\s}^{n+1}(0) = B_{V_{j}} \cup S_{V_{j}}$ with 
$B_{V_{j}}, S_{V_{j}}$ having analogous meaning. Then $\z_{j} \in S_{V_{j}}$ since 
$Z_{j} \in S$. By (\ref{decomposition-1}), 

\begin{equation}\label{decomposition-2}
{\mathcal H}^{\g} \, (S_{\mathbf C}) = 0
\end{equation}

\noindent
for each $\g >0$. On the other hand, 
since ${\mathbf C}$ is a cone and $z \in {\rm sing} \, {\mathbf C} \cap {\mathbf S}^{n}$, we have 
that $\{ t\z \, : \, t > 0\} \subset {\rm sing} \, {\mathbf C}.$ In fact, we must have that 

\begin{equation}\label{decomposition-3}
\{ t\z \, : \, t > 0\} \cap B_{\s}^{n+1}(0) \subset S_{\mathbf C}.
\end{equation}

\noindent
For if not, $\z \in B_{\mathbf C}$ in which case 
${\mathbf C}$ would have a (unique) multiplicity 2 tangent plane at $\z$, and since 
$V_{j} \to {\mathbf C}$, by Theorem~\ref{maintheorem}, it follows that for all sufficiently large $j$,
${\rm spt} \, \|V_{j}\|$ is a 2-valued $C^{1, \alpha}$ graph in some neighborhood 
of $\z.$ But this contradicts the fact that $\z_{j} \in S_{V_{j}}.$ Hence we must have
(\ref{decomposition-3}), but this contradicts the dimension estimate (\ref{decomposition-2}).
This concludes the proof of the lemma.
\end{proof}

\bigskip

\section{Some further corollaries}\label{corollaries}
\begin{theorem}\label{curvest}
Let $\d \in (0, 1)$. There exist positive numbers $\G$ and $\s$ depending only on $\d$ such that if $2 \leq n \leq 6$ 
and $M$ is a 
an immersed, stable minimal hypersurface of $B_{2}^{n+1}(0)$ satisfying ${\mathcal H}^{n-2}({\rm sing} \, M) = 0$ and 
$\frac{{\mathcal H}^{n}(M)}{\omega_{n}2^{n}} \leq 3 - \d$, then ${\rm sing} \, M \cap B_{\s}^{n+1}(0) = \emptyset$ and

$$\mbox{\rm sup}_{M \cap B_{\s}^{n+1}(0)} \,|A| \leq \G$$ 

\noindent
where $A$ denotes the second fundamental form of $M$ and $|A|$ the length of $A$.\\ 
\end{theorem}

\medskip

\noindent
{\bf Remark:} If $M$ is assumed to be embedded, this result holds with mass bound arbitrary, and 
is due to R.Schoen and L. Simon~\cite{SS}.
In dimensions $2 \leq n \leq 5$, provided we {\em assume} ${\rm sing} \, M = \emptyset,$ the
result (for $M$ immersed) holds with mass bound arbitrary, and is due to 
R. Schoen, L. Simon and S. T. Yau~\cite{SSY}.\\

\begin{proof}
Set $\s_{1} = \min \,\{ \s(1, \d), \ldots, \s(6, \d)\}$ and $\s = \s_{1}/4,$ where 
$\s(n, \d)$ is as in Theorem~\ref{compactnessthm}. Then it follows directly 
by taking $M_{k} = M$ in Theorem~\ref{compactnessthm} that ${\rm sing} \, M \cap B_{\s_{1}}^{n+1}(0) = \emptyset$, 
so we only need to prove the curvature
estimate.\\ 

If there is no such $\G,$  then for some $n$ with $2 \leq n \leq 6$ and some 
$\d \in (0, 1)$, 
there exists a sequence $\{M_{k}\}$ of stable minimal hypersurfaces immersed in $B_{2}^{n+1}(0)$
with $0 \in M_{k},$ satisfying ${\mathcal H}^{n-2}({\rm sing} \, M_{k}) = 0$ (or we may assume 
${\rm sing} \, M_{k} \cap B_{1}^{n+1}(0) = \emptyset$ if we wish, in view of the preceding 
paragraph) and 
$\frac{{\mathcal H}^{n} \, (M_{k})}{\omega_{n} 2^{n}} \leq 3 - \d$ for each $k$, and yet
there exists a point $Z_{k} \in M_{k} \cap B_{\s}^{n+1}(0)$ for each $k$ with 

\begin{equation}\label{curvest-1}
|A_{k}|(Z_{k}) \to \infty,
\end{equation}

\noindent
where $A_{k}$ denotes the second fundamental form of $M_{k}$ and 
$|A_{k}|$ its length. By Theorem~\ref{compactnessthm}, there exists a stationary varifold $V$
of $B_{2}^{n+1}(0)$  such that after passing to a subsequence, which we continue to denote $\{M_{k}\}$, 
we have that $M_{k} \to V$ as varifolds, and that ${\rm spt} \, \|V\| \cap B_{\s_{1}}^{n+1}(0) = M$ where 
$M$ is a smooth (i.e. having ${\rm sing} \, M = \emptyset$) stable minimal hypersurface of $B_{\s_{1}}^{n+1}(0)$;
since varifold convergence implies convergence (of the supports of the weight measures) in 
Hausdorff distance, we also have that $Z_{k} \to Z$ for some $Z \in M \cap {\overline B_{\s}^{n+1}(0)}.$
But since $M$ is a regular immersed hypersurface, and the density of $M$ at $X$ is $\leq 3 - \d$ 
for every $X \in M \cap B_{\s_{1}/2}^{n+1}(0)$, the tangent cone to $M$ at $Z$ is either a multiplicity 1 plane, 
or a multiplicity 2 plane or a transversely intersecting pair of hyperplanes with $Z$ belonging to its
axis. Applying respectively Allard's regularity theorem, Theorem~\ref{maincorollary} or 
Theorem 1 of \cite{WN2}, we conclude that there exists a fixed radius $\r > 0$ independent of 
$k$ such that in each of these cases, for all sufficiently large $k$, we have that 

$${\rm sup} \, _{M_{k} \cap B_{\r}^{n+1}(Z)} \, |A_{k}| \leq \frac{C}{\r}$$

\noindent
for some fixed constant $C = C(n)$ independent of $k.$ But this contradicts (\ref{curvest-1}). The 
theorem is thus proved.  
\end{proof}

\medskip

\begin{theorem}[A Bernstein type theorem] Let $\d \in (0, 1)$. Suppose $2 \leq n \leq 6,$ $M$ is a 
complete, non-compact stable minimal hypersurfaces of 
${\mathbf R}^{n+1}$ satisfying $\frac{{\mathcal H}^{n} (M \cap B_{R}^{n+1}(0))}{\omega_{n}R^{n}} \leq 3 -\d$ for 
all $R> 0.$ Then $M$ must be a union of affine hyperplanes.
\end{theorem}

\medskip

\noindent
{\bf Remark}: This is a slight generalization of the Bernstein type theorem in \cite{WN3}, which asserts the existence 
of a number $\e \in (0,1)$ such that the conclusion of the theorem is true 
whenever $2 \leq n \leq 6$ and $\frac{{\mathcal H}^{n} M \cap B_{R}^{n+1}(0)}{\omega_{n}R^{n}} \leq 2 + \e$ for all $R > 0.$\\

\begin{proof}
By Theorem~\ref{curvest}, $\sup_{B_{\s R}^{n+1}(0)} \, |A| \leq \frac{\G}{R}$ for all 
$R > 0,$ where $\s>0$ and $\G$ are independent of $R$. Let $R \to \infty.$
\end{proof}

\medskip

The following result is an improvement of Lemma 1 of \cite{SS}. Note that our proof of it below uses the regularity theory; 
Lemma 1 of \cite{SS} 
on the other hand was used in {\em proving} the regularity theorem of \cite{SS}, and it would be interesting to see if 
the result below has a proof independent of regularity theory.\\
 
\begin{theorem}
Let $p \in (0, 4 + \sqrt{8/n})$, $\Lambda >0$ and $\th \in (0, 1).$ There exists a constant $C = C(n, p, \Lambda, \th)$ such that if 
$M$ is an embedded, stable minimal hypersurface of 
$B_{1}^{n+1}(0)$ with ${\mathcal H}^{n-2}(\mbox{\rm sing} \, M) < \infty$ and ${\mathcal H}^{n}(M) \leq \Lambda$, then 

$$\int_{M \cap B_{\th}^{n+1}(0)} |A|^{p} \leq C \left(\int_{M \cap B_{1}^{n+1}(0)} 1 - (\nu\cdot\nu_{0})^{2}\right)^{p/2}$$

\noindent
for any unit vector $\nu_{0} \in {\mathbf R}^{n+1}.$ Here 
$A$ denotes the second fundamental form of $M$ and $\nu$ the unit normal vector to $M.$\\

The estimate continues to hold if $M$ is immersed  
provided $\Lambda = \omega_{n}(3 -\d)$ for some $\d \in (0, 1)$ and ${\mathcal H}^{n-2} \, ({\rm sing} \, M) = 0.$
\end{theorem}

\begin{proof}
We argue by contradiction. If the estimate were not true for some $\Lambda,$ 
$p \in [4, 4 + \sqrt{8/n})$ and $\th \in (0, 1)$, then there exists a sequence of stable minimal hypersurfaces $M_{k}$ 
of $B_{1}^{n+1}(0)$ with ${\mathcal H}^{n}(M_{k}) \leq \Lambda$ and 

\begin{equation}\label{schoenest-1}
\int_{M_{k} \cap B_{\th}^{n+1}(0)} |A_{k}|^{p} \geq k \left(\int_{M_{k} \cap B_{1}^{n+1}(0)} 1 - (\nu_{k}\cdot\nu_{0}^{k})^{2}
\right)^{p/2}
\end{equation}

\noindent
where $\nu_{0}^{k}$ are unit vectors in ${\mathbf R}^{n+1}.$ Note that under the assumptions of 
the theorem, ${\mathcal H}^{n-7+\g} \, ({\rm sing} \, M_{k}) = 0$ for each $\g > 0$, which follows from
Theorem 3 of \cite{SS} in the embedded case, and from Theorem~\ref{compactnessthm} above in the immersed case. 
By Schoen-Simon-Yau integral curvature estimate (which was originally
proved for smooth, stable hypersurfaces but continues 
to hold for stable hypersurfaces $M$ with singularities provided
${\mathcal H}^{n-p}({\rm sing} \, M) < \infty$, as can be seen using an easy cut-off function argument), we have that 

\begin{equation}\label{schoenest-2}
\int_{M_{k} \cap B_{\th}^{n+1}(0)} |A_{k}|^{p} \leq C
\end{equation}

\noindent
where $C$ is a constant that depends only on $n$, $p,$ $\Lambda$ and $\th.$  
From (\ref{schoenest-1}) and (\ref{schoenest-2}), it follows that 

\begin{equation}\label{schoenest-3}
\int_{M_{k} \cap B_{1}^{n+1}(0)} 1 - (\nu_{k}\cdot\nu_{0}^{k})^{2} \to 0.
\end{equation}

\noindent
Since mass of $M_{k}$ is uniformly bounded, Allard's compactness theorem says that 
after passing to a subsequence, $M_{k} \to V$ for some stationary varifold $V$ of $B_{1}^{n+1}(0),$ and 
(\ref{schoenest-3}) says that $V$ must be a hyperplane (with some positive integer multiplicity.) Let us assume 
without loss of generality that this hyperplane is ${\mathbf R}^{n} \times \{0\}.$  
Now in case $M_{k}$ are embedded, by Schoen-Simon regularity theorem, this means that for 
all sufficiently large $k$, $M_{k} \cap (B_{\frac{1 + \th}{2}}(0) \times {\mathbf R})$ decomposes 
as the (disjoint) union of graphs of $m_{k}$ functions  $u^{k}_{i} \, : \, B_{\frac{1 + \th}{2}}(0) \to {\mathbf R},$ 
$1 \leq i \leq m_{k},$ (with $m_{k}$ bounded independently of $k$ by a number 
depending on $\Lambda$), 
each solving the minimal surface equation. In the immersed case under the stronger mass bound, 
by Theorem~\ref{maincorollary}, the same conclusion holds with $m_{k} \leq 2.$\\ 

Now let $L_{k}$ be the hyperplane determined by the unit vector $\nu_{0}^{k}$, and 
$l_{k} \, : \, {\mathbf R}^{n} \times \{0\} \to {\mathbf R}$ be the linear function whose graph is $L_{k}.$ (Note that 
$\nu_{0}^{k} \cdot e^{n+1} \to 1.$) 
Then $u^{k}_{i} - l_{k}$ solves an elliptic equation over $B_{\frac{1 + \th}{2}}(0),$ and 
so by elliptic estimates, we have a constant $C = C(n, \th)$ such that 
$\mbox{sup}_{B_{\frac{1 + 3\th}{4}}(0)} \, |D^{2}u^{k}_{i}| \leq C\|Du^{k}_{i} - Dl_{k}\|_{L^{2}(B_{\frac{1 + \th}{2}}(0))}$ 
and $\mbox{sup}_{B_{5/8}(0)} \, |Du^{k}_{i} - Dl_{k}| \leq C\|Du^{k}_{i} - Dl_{k}\|_{L^{2}(B_{3/4}(0))}$ for each $i$. 
But this means that $\mbox{sup}_{M_{k} \cap B_{\th}^{n+1}(0)} |A_{k}| \leq C\left(\int_{M_{k} \cap B_{1}^{n+1}(0)} 
1 - (\nu_{k}\cdot\nu_{0}^{k})^{2} \right)^{1/2}$ where $C = C(n, \Lambda),$ which contradicts (\ref{schoenest-1}).
\end{proof}

\begin{flushright}
{\sc Massachusetts Institute of Technology\\
77 Massachusetts Avenue, Cambridge, MA 02139.\\
Current address:\\
University of California, San Diego\\
9500, Gilman Drive\\
La Jolla, CA 92093-0112.}\\
\end{flushright}   

\end{document}